\documentclass[phd]{msthesis}
\usepackage{amsfonts, amssymb, amsmath, eucal, verbatim, amsthm, amscd, enumerate, graphics, lscape, tabls}

\newtheorem{thm}{Theorem}[chapter]
\newtheorem{prop}{Proposition}[chapter]
\newtheorem{lem}{Lemma}[chapter]
\newtheorem{defn}{Definition}[chapter]

\newtheorem{cor}{Corollary}[chapter]
\newtheorem{rem}{Remark}[chapter]
\newtheorem{ex}{Example}[chapter]
\newtheorem{claim}{Claim}[chapter]
\newtheorem{bst}{Bisymmetric Triple}[chapter]
\newtheorem{symp}{Symmetric Pair}[chapter]

\newenvironment{pr}{\textit{Proof: }}

\newcommand{\sg}{\sigma}

\newcommand{\al}{\alpha}
\newcommand{\be}{\beta}
\newcommand{\la}{\lambda}

\newcommand{\ga}{\gamma}

\newcommand{\de}{\delta}

\newcommand{\iso}{\cong}

\newcommand{\li}{\medskip}

\newcommand{\reals}{\mathbb{R}}

\newcommand{\complex}{\mathbb{C}}

\newcommand{\rationals}{\mathbb{Q}}

\newcommand{\R}{\mathcal{R}}

\newcommand{\mfe}{\mathfrak{e}}
\newcommand{\mff}{\mathfrak{f}}
\newcommand{\mfg}{\mathfrak{g}}
\newcommand{\mfh}{\mathfrak{h}}
\newcommand{\mfk}{\mathfrak{k}}
\newcommand{\mfl}{\mathfrak{l}}
\newcommand{\mfm}{\mathfrak{m}}
\newcommand{\mfn}{\mathfrak{n}}

\newcommand{\mfp}{\mathfrak{p}}
\newcommand{\mfs}{\mathfrak{s}}
\newcommand{\mfu}{\mathfrak{u}}

\newcommand{\mfsu}{\mathfrak{s}\mathfrak{u}}
\newcommand{\mfso}{\mathfrak{s}\mathfrak{o}}
\newcommand{\mfsp}{\mathfrak{s}\mathfrak{p}}

\def\bar{\begin{array}}
\def\ear{\end{array}}
\def\sbar{\begin{subarray}}
\def\sear{\end{subarray}}
\def\beq{\begin{equation}}
\def\eeq{\end{equation} }
\def\beqar{\begin{eqnarray}}
\def\eeqar{\end{eqnarray}}
\def\bal{\begin{align}}
\def\eal{\end{align}}
\def\bfig{\begin{figure}}
\def\efig{\end{figure}}
\def\bc{\begin{center}}
\def\ec{\end{center}}
\def\btab{\begin{table}}
\def\etab{\end{table}}
\def\bland{\begin{landscape}}
\def\eland{\end{landscape}}

\def\bproof{\begin{pr}}
\def\eproof{\end{pr}}
\def\bprop{\begin{prop}}
\def\eprop{\end{prop}}
\def\bthm{\begin{thm}}
\def\ethm{\end{thm}}
\def\blem{\begin{lem}}
\def\elem{\end{lem}}
\def\brem{\begin{rem}}
\def\erem{\end{rem}}
\def\bcor{\begin{cor}}
\def\ecor{\end{cor}}
\def\bex{\begin{ex}}
\def\eex{\end{ex}}
\def\bdfn{\begin{defn}}
\def\edfn{\end{defn}}
\def\bclaim{\begin{claim}}
\def\eclaim{\end{claim}}
\def\bbst{\begin{bst}}
\def\ebst{\end{bst}}
\def\bsymp{\begin{symp}}
\def\esymp{\end{symp}}

\newlength{\myVSpace}
\newcommand\xstrut{\raisebox{-.9\myVSpace}{\rule{0pt}{\myVSpace}}}

\makeatletter
\def\table{\@ifnextchar[{\table@i}{\table@i[\fps@table]}}
\def\table@i[#1]{\@float{table}[#1]\footnotesize}
\makeatother

\title{Einstein Homogeneous Riemannian Fibrations}
\author{F\'{a}tima Ara\'{u}jo}
\submityear{2008}

\begin{document}
\maketitle

\dedication {To my godparents Amadeu and Concei\c{c}\~{a}o}

\begin{acknowledgements}I owe my deepest gratitute to my
supervisor Prof. Dmitri Alekseevskii, who has invested
considerable time and energy into guiding me through my PhD. His
professional guidance and unfaltering support have been most
invaluable. I would like to thank Prof. Elmer Rees for his
supervision in a first stage of my PhD and Dr. Roger Bielawski for
his help and ideas during a critical
period of my research. \\ \\

I would like to thank the Funda\c{c}\~{a}o Portuguesa para a Ci\^{e}ncia e Tecnologia for sponsoring this work  and the Centro de Algebra da Universidade de Lisboa.\\ \\

Last but not least, I want to thank my family and close friends
for their unwavering emotional support and firm belief in me
throughout my research and life.\\ \\

\end{acknowledgements}

\standarddeclaration \noindent

\begin{abstract}
This thesis is dedicated to the study of the existence of homogeneous Einstein metrics on the total space of homogeneous fibrations such
that the fibers are totally geodesic manifolds. We obtain the Ricci curvature of an invariant metric with totally geodesic fibers and some necessary conditions for the existence of Einstein metrics with totally geodesic fibers in terms of Casimir operators. Some particular cases are studied, for instance, for normal base or fiber, symmetric fiber, Einstein base or fiber, for which the Einstein equations are manageable. We investigate the existence of such
Einstein metrics for invariant bisymmetric fibrations of maximal
rank, i.e., when both the base and the fiber are symmetric spaces
and the base is an isotropy irreducible space of maximal rank. We find this way new Einstein metrics. For
such spaces we describe explicitly the isotropy representation in
terms subsets of roots and compute the eigenvalues of the Casimir
operators of the fiber along the horizontal direction. Results for compact simply connected $4$-symmetric spaces of maximal rank follow from this. Also, new invariant Einstein metrics are found on Kowalski $n$-symmetric spaces.
\end{abstract}

\setcounter{page}{1} \tableofcontents

\begin{introduction}
The principal topic of study in this thesis is the existence of Einstein invariant  metrics on the total space of homogeneous Riemannian fibrations such that the fibers are totally geodesic submanifolds. In chapter 1 we introduce some main definitions and notation and deduce some essential formulas for the Ricci curvature of an invariant metric. Then we consider a fibration $G/L\rightarrow G/K$
such that $G$ is a compact connected semisimple Lie group $G$, and
$L \varsubsetneq K\varsubsetneq G$ are connected closed non-trivial subgroups of
$G$.  We assume that
$G/L$ has simple spectrum. On the total space $G/L$, we consider a $G$-invariant Riemannian metric such that the natural projection $G/L\rightarrow G/K$ is a Riemannian submersion and the fibers are totally geodesic submanifolds. We shall briefly call such a
metric an \textit{Einstein adapted metric}. We describe the Ricci curvature of any adapted metric in
terms of Casimir operators and obtain two necessary conditions for
existence of an Einstein adapted metric expressed only in terms of
the Casimir operators of the horizontal and vertical directions.
These will provide a tool to show that in many cases such a metric
cannot exist without further study, i.e.,  only by studying eigenvalues of
certain Casimir operators.

\li

In chapter 2 we restrict the object
of study to some special cases, where the Einstein equations are simpler. We consider
the cases where the metric on the fiber or on the base is a multiple
of the Killing form of $G$, in which cases are included those with
isotropy irreducible fiber or base. The case when both these two
conditions are satisfied gives rise to the study of the, throughout
called, Einstein binormal metrics. The existence of Einstein binormal metrics
translates into very simple algebraic conditions which shall allow
us to find out new Einstein metrics. Also we obtain necessary conditions for existence of an Einstein adapted metric such that the metric on the base space or the metric on the fiber are also Einstein. Finally, we apply the results obtained so far to the case when the fiber is a symmetric space and $N$ is isotropy irreducible.

\li

Chapter 3 is devoted to bisymmetric fibrations of maximal rank,
i.e., we consider a fibration $G/L\rightarrow G/K$, as in chapter
1, such that $L$ is a subgroup of maximal rank, $K/L$ is a
symmetric space and $G/L$ is an isotropy irreducible symmetric
space. We introduce the notion of a bisymmetric triple
$(\mfg,\mfk,\mfl)$ associated to a bisymmetric fibration. We
obtain all the bisymmetric triples $(\mfg,\mfk,\mfl)$ in the case
when $\mfg$ is a simple Lie algebra and classify them into two
different types, I and II. We classify all the Einstein adapted
metrics when $\mfg$ is an exceptional Lie algebra, for both Type I
and II. When $\mfg$ is a classical Lie algebra, we classify all
the Einstein adapted metrics for Type I. For Type II in the
classical case, we classify all Einstein binormal metrics and all
Einstein adapted metrics whose restriction to the fiber is also
Einstein; moreover, if one of these metrics exists we obtain all
the other Einstein adapted metrics. Finally, we apply the results
obtained to compact simply-connected irreducible $4$-symmetric
spaces. In appendix A we obtain all the necessary eigenvalues for
the Einstein equations for each bisymmetric triple considered in
this chapter.

\li

In chapter 4 we study the existence of Einstein adapted
metrics on the Kowalski $n$-symmetric spaces, i.e., we consider a fibration
$$\frac{\triangle^p G_0\times \triangle^q G_0}{\triangle ^nG_0}\rightarrow\frac{G_0^n}{\triangle^nG_0}\rightarrow
\frac{G_0^p}{\triangle^pG_0}\times \frac{G_0^q}{\triangle^qG_0},$$

where $G_0$ is compact connected simple Lie group and $\triangle^mG_0$ is the diagonal subgroup in $G_0^m$, for $m=p,q,n$. It is known that $\frac{G_0^n}{\triangle^nG_0}$ is a standard Einstein manifold and we prove that, for $n>4$, there exists another Einstein adapted metric, whereas, for $n=4$, the standard metric is the only Einstein adapted metric.
\end{introduction}

\chapter{}

In Section 1 we introduce some essential definitions and notation. We deduce a formula for the Ricci curvature of an
invariant metric on a reductive homogeneous space. In Section 2 we
obtain the Ricci curvature of an invariant metric with totally
geodesic fibers on the total space of a homogeneous fibration and
some necessary conditions for that metric to be Einstein.

\section{The Ricci Curvature of a Riemannian  Homogeneous Space}

A Riemannian metric $g$ is said to be Einstein if its Ricci curvature satisfies an equation of the form $Ric=Eg$, for some constant $E$, the Einstein constant of $g$ (\cite{Be}). This equation, commonly called the Einstein equation, is in general a very complicated system of partial differential equations of second order. Although so far no fully general results are known for existence of Einstein metrics, many results of existence and classification are known for many classes of spaces. Two examples of this are the  K\"{a}hler-Einstein metrics (\cite{Yau}, \cite{AP}, \cite{Sa}, \cite{Ti}) and the Sasakian-Einstein metrics (\cite{BG}). Many results are known on homogeneous Einstein metrics. For Riemannian homogeneous spaces the Einstein equation translates into a system of algebraic equations, which is an easier problem than its general version. However, even for this class of spaces we are far from knowing full answers.  Einstein normal homogeneous manifolds were classified by Wang and Ziller (\cite{WZ}). Nowadays, it is known that every compact simply connected homogeneous manifold with dimension less or equal to $11$ admits a homogeneous Einstein metric: any such manifold with dimension $2$ or $3$ has constant sectional curvature \cite{Be}; in dimension $4$, the result was shown by Jensen (\cite{Je1}), and by Alekseevsky, Dotti and Ferraris in dimension 5 (\cite{ADF}); in dimension $6$, the result is due to Nikonorov and Rodionov (\cite{NR2}), and in dimension $7$ it is due to Castellani, Romans and Warner (\cite{CRW}). All the $7$-dimensional homogeneous Einstein manifolds (\cite{Ni}) were obtained by Nikonorov. These results were extended to dimension up to $11$ by B\"{o}hm and Kerr (\cite{BK}). Many examples of homogeneous Einstein manifolds with dimension arbitrary big are known. Spheres and projective spaces are examples of this, where all homogeneous Einstein metrics were classified by Ziller (\cite{Zi2}). Also isotropy irreducible spaces (\cite{Wo}), symmetric spaces (\cite{He},\cite{Ke1}), flag manifolds, among many others, provide examples of Einstein manifolds with arbitrary big dimension. Einstein homogeneous fibrations have also been the object of study. We recall the work of Jensen on principal  fibers bundles (\cite{Je2}), where new invariant Einstein metrics are found on the total space of certain homogeneous fibrations, and the work of Wang and Ziller on principal torus bundles (\cite{WZ3}). Einstein homogeneous fibrations are the main focus of this thesis.

\li

Let $G$ be a Lie group and $L$ a closed subgroup. We denote by $\mfg$ and $\mfl$ the Lie algebras of $G$ and $L$, respectively. The homogeneous space $M=G/L$ is the space of all cosets $\{aL:a\in G\}$ endowed with the unique differentiable structure such that the canonical projection

\beq\bar{rrcl}\pi:& G & \rightarrow & M\\
&  a & \mapsto & aL \ear\eeq

is a submersion, i.e., $d\pi_a$ is onto for all $a\in G$, and with the natural transitive left action of $G$,

\beq\bar{rrcl}\tau:& G\times M & \rightarrow & M\\
&  (b,aL) & \mapsto & (ba)L \ear\eeq

Let  $X\in\mfg$ and let $exp\,tX$ be the one-parameter subgroup
generated by $X$. For every $a\in G$,

\beq d\pi_a(X)=\frac{d}{dt}(exp\,tX)aL\mid_{t=0}\eeq

and, in particular, for $o=\pi(e)=L$, this map yields an isomorphism

\beq\label{isomts}d\pi_e:\mfg/\mfl\iso T_oM.\eeq

For every $X\in\mfg$ we define a $G$-invariant vector field on $M$ by

\beq X^*_{aL}=d\pi_a(X)=\frac{d}{dt}(exp\,tX)gL\mid_{t=0}.\eeq

The homogeneous space $M$ is called \textit{reductive} if there exists a direct complement $\mfm$ of $\mfl$ in $\mfg$ which is $Ad\,L$-invariant, i.e.,

\beq\mfg=\mfl\oplus\mfm \textrm{ and }Ad\,L(\mfm)\subset\mfm.\eeq

The inclusion $Ad\,L(\mfm)\subset\mfm$ implies

\beq\label{invsubs2}[\mfl,\mfm]\subset\mfm\eeq

and the converse holds if $L$ is connected. The homogeneous space $M$ is reductive if $L$ is compact. Throughout, we suppose that $L$ is compact and we denote by $\mfm$ a reductive complement of $\mfl$ on $\mfg$. If $M$ is reductive we have an isomorphism

$$\mfm\iso T_oM$$

and the tangent space $T_oM$ is identified with $\mfm$ and consequently the vector field $X^*$ on $M$ is identified with $X\in\mfm$. Under this identification we shall simply write $X$ for $X^*_o$. Furthermore, the isotropy representation of $M$, $$Ad^{M}: L\rightarrow  GL(T_oM),$$ is equivalent to the adjoint representation of $L$ on $\mfm$. Consequently, there is a one-to-one correspondence between $G$-invariant objects on $M$ and $Ad\,L$-invariant objects on $\mfm$. In  particular, $G$-invariant metrics on $M$ correspond to $Ad\,L$-invariant scalar products on $\mfm$. More precisely, a metric  $g$ on $M$ is said to be $G$-invariant if, for every $a\in G$,

$$\tau_a^*g=g$$

and the correspondence between $G$-invariant metrics on $M$ and $Ad\,L$-scalar products on $\mfm$ is given

\beq\label{indmet}g_a(X^*_a,Y^*_a)=<X,Y>, \textrm{ for all }a\in G.\eeq                                                                                                                                                                                                                                                                                           \li

Let $Kill$ be the Killing form of $\mfg$. We recall that $Kill$ is
the bilinear form on $\mfg$ defined by \beq
Kill(X,Y)=tr(ad_Xad_Y), \,X,Y\in\mfg\eeq

where, for each $X\in \mfg$, $ad_X$ denotes the adjoint map

\beq\bar{rcl} \mfg & \rightarrow & \mfg\\
Y & \mapsto & [X,Y]\ear.\eeq

The Killing form of $\mfg$ is an $Ad\,G$-invariant bilinear form and it is non-degenerate if $\mfg$ is semisimple. Moreover, if $G$ is a compact connected semisimple Lie group, $Kill$ is negative definite. In this case, by (\ref{indmet}), the negative of the Killing form induces a $G$-invariant metric on $M$, the \textit{standard Riemannian metric} on $M$.

\li

With respect to the decomposition $\mfg=\mfl\oplus\mfm$, we write

\beq \label{adX} ad_X=\left(\bar{cc} 0 & C_X\\ B_X &
P_X\ear\right),\textrm{ for every  } X\in\mfm.\eeq

Hence, for $Y\in \mfm$,

\beq (ad_XY)_{\mfm}=P_XY \textrm{ and } (ad^2_X
Y)_{\mfm}=(B_XC_X+P_X^2)Y,\eeq

where the subscript $\mfm$ denotes the projection onto $\mfm$.

\li

Let $g_M$ be a $G$-invariant metric on $M$. As was explained above, there is a one-to-one
correspondence between $G$-invariant metrics on $M$ and
$Ad\,L$-invariant scalar products on $\mfm$. So let $<,>$ be the
$Ad\,L$-invariant scalar product on $\mfm$ corresponding to $g_M$.

\li

For every $X\in\mfm$, let $T_X$ be the endomorphism of $\mfm$
defined by

\beq\label{opT}<T_XY,Z>=<X,P_YZ>,\textrm{ for every } Y,Z\in\mfm.\eeq

\li

The \textit{Nomizu operator} of the scalar product $<,>$ (cf.
\cite{No}, \cite{NK}) is $$L_X\in End(\mfm),\,X\in\mfg,$$ defined by

\beq L_XY=-\nabla_{Y}X^*,\,Y\in\mfm\eeq

where $\nabla$ is the Riemannian connection of $g_M$. We have

\beq \label{no1}L_XY=\frac{1}{2}P_XY+U(X,Y),\eeq

where $U:\mfm\times\mfm\rightarrow \mfm$ is the operator

\beq\label{opU}U(X,Y)=-\frac{1}{2}(T_XY+T_YX),\,X,Y\in\mfm.\eeq

\li

The metric $g_M$ is called \textit{naturally reductive} if $U=0$. The curvature tensor, the sectional curvature and the Ricci curvature of $g_M$ are $G$-invariant tensors and thus they are determined by the following identities (\cite{No}), which represent their values at the point $o$. By using $L$ the curvature tensor of $g_M$ at $o$ can be written as

\beq
\label{curv}R(X,Y)=[L_X,L_Y]_{\mfm}-L_{[X,Y]_{\mfm}}-ad_{[X,Y]_{\mfl}},\eeq

for every $X,\,Y\in \mfm$. The sectional curvature $K$ of $g_M$ is defined by

\beq K(Z,X)=<R(Z,X)X,Z>,\eeq

for every $X,Z\in\mfm$ orthonormal with respect to $<,>$. The
Ricci curvature of $g_M$ is determined by

\beq Ric(X,X)=\sum_iK(Z_i,X),\,X\in\mfm\eeq

where $(Z_i)_i$ is an orthonormal basis of $\mfm$ with respect to $<,>$. The metric $g_M$ is said to be an \textit{Einstein metric} if

\beq Ric=Eg_M ,\eeq

for some constant $E$ called the \emph{Einstein constant} of
$g_M$. Equipped with a $G$-invariant Einstein metric, $M$ is
called an Einstein homogeneous manifold.

\li

Below we show some elementary properties of the Nomizu operator:

\blem (i) $L_X$ is skew-symmetric with respect to $<,>$, i.e.,

\beq\label{nomss}<L_XY,Z>+<Y,L_XZ>=0, \,X,Y,\in\mfm;\eeq

(ii) for every $X,Y\in\mfm$,

\beq \label{tf}L_XY-L_YX=[X,Y]_{\mfm}=P_XY.\eeq
\elem

\bproof From identities (\ref{opT}), (\ref{no1}) and  (\ref{opU}) we deduce

$$\bar{rl} 2<L_XY,Z> = & <P_XY,Z>-<T_XY,Z>-<T_YX,Z>\\ \\

= & <T_ZX,Y>+<T_XZ,Y>-<P_XZ,Y>=\\ \\

= & -2<Y,L_XZ>\ear$$

and this shows the skew-symmetry of $L_X$. To show (ii) we just
observe that $U(X,Y)=U(Y,X)$ and $P_XY=-P_YX$. Equivalently, this
assertion just means that the Levi-Civita connection is torsion
free since $\nabla_{X^*}Y^*=L_XY$.

$\Box$\eproof

For every $X,Y\in\mfm$, we define the operator

\beq R_XY=L_YX\eeq

and the vector

\beq V_X=L_XX.\eeq

\blem \label{R_X}For every $X\in\mfm$,

$$R_X=-\frac{1}{2}(P_X+P_X^*+T_X)\textrm{ and }
R_X^*=-\frac{1}{2}(P_X+P_X^*-T_X).$$\elem

\bproof Let $X$, $Y$, $Z\in\mfm$.

$$\bar{rl}2<R_XY,Z>= & 2<L_YX,Z>\\ \\

= & <[Y,X]_{\mfm},Z>+<Y,[Z,X]_{\mfm},>+<X,[Y,Z]_{\mfm}> \\ \\

= & -<P_XY,Z>-<P_X^*Y,Z>-<T_XY,Z>.\ear$$

Thus we obtain the required expression for $R_X$. The formula for
$R_X^*$ follows from the fact that $P_X+P_X^*$ is symmetric and
$T_X$ is skew symmetric with respect to $<,>$.

$\Box$

\eproof

\blem \label{seccurv}The sectional curvature of $g_M$ is

$$K(Z,X)=<\big(R_X^*R_X-P_X^2-P_X^*P_X-B_XC_X+P_{V_X}\big)Z,Z>,$$

for every $Z,\,X\in \mfm$ orthonormal with respect to $<,>$. \elem

\bproof Let $Z,\,X\in \mfm$. The proof is a straightforward
calculation by using the identities (\ref{curv}), (\ref{nomss})
and (\ref{tf}):

$$\bar{rl} K(Z,X) = &  <R(Z,X)X,Z> \\ \\

= &  <[L_Z,L_X]_{\mfm}X,Z> -<L_{[Z,X]_{\mfm}}X,Z> -<ad_{[Z,X]_{\mfl}}X,Z>  \\ \\

= & <L_ZL_XX,Z>-<L_XL_ZX,Z>-<L_X[Z,X]_{\mfm},Z>-\\ \\

 & -<[[Z,X]_{\mfm},X]_{\mfm},Z>-<[[Z,X]_{\mfl}X]_{\mfm},Z>\\ \\

= &
-<L_XX,L_ZZ,>+<L_ZX,L_XZ>+<[Z,X]_{\mfm},L_XZ>-\\ \\

& -<[[Z,X],X]_{\mfm},Z>\\ \\

 = &-<L_XX,L_ZZ,>+<L_ZX,L_ZX>+<L_ZX,[X,Z]_{\mfm}>+

\\ \\& +<[Z,X]_{\mfm},L_ZX>+<[Z,X]_{\mfm},[X,Z]_{\mfm}>-\\ \\

& -<[X,[X,Z]]_{\mfm},Z>\\ \\

=
&-<L_XX,L_ZZ,>+<L_ZX,L_ZX>-<[Z,X]_{\mfm},[Z,X]_{\mfm}>-\\ \\

&-<(P_X^2Z)_{\mfm},Z>\ear$$

$$\bar{rl}= & <L_{Z}X,L_{Z}X>-<[X,Z]_{\mfm},[X,Z]_{\mfm}>-<(P^2_X
Z)_{\mfm},Z>\\ \\

 & -<L_{Z}V_X,Z>\\ \\

= & <R_XZ,R_XZ>-<P_XZ,P_XZ>-<(P^2_X+B_XC_X)Z,Z>\\ \\

& -<L_{V_X}Z+[Z,V_X]_{\mfm},Z>\\
\\

= &
<(R_X^*R_X-P_X^*P_X-B_XC_X-P_X^2)Z,Z>+<(P_{V_X}-L_{V_X})Z,Z>.\ear$$

\li

Since $L_X$ is skew-symmetric with respect to $<,>$, we have
$<L_{V_X}Z,Z>=0$, for every $Z\in\mfm$, and we obtain

$$K(Z,X)=<\big(R_X^*R_X-P_X^*P_X-B_XC_X-P_X^2+P_{V_X}\big)Z,Z>,$$

as required.

$\Box$ \eproof

\li

\bthm\label{ricci} Let $X,Y\in\mfm$. Then

$$Ric(X,Y)=-\dfrac{1}{2}tr(2R_XR_Y+B_XC_Y+B_YC_X-2P_{U(X,Y)}).$$

In particular,

$$Ric(X,X)=-tr(R_X^2+B_X C_X-P_{V_X}).$$\ethm

\bproof We first compute $Ric(X,X)$ and then obtain $Ric(X,Y)$ by
polarization.

\li

Since $T_X$ is a skew-symmetric operator and $P_X+P_X^*$ is
symmetric,
$$tr\big((P_X+P_X^*)T_X\big)=-tr\big(T_X(P_X+P_X^*)\big)=tr\big((P_X+P_X^*)T_X)$$

and thus all the terms vanish. Therefore, by using Lemma \ref{R_X}
we obtain

$$tr(R_X^*R_X)=\frac{1}{4}\big((P_X+P_X^*)^2-T_X^2\big)=\frac{1}{4}tr(2P_X^2+2P_X^*P_X-T_X^2)$$

and

$$tr(R_X^2)=\frac{1}{4}\big((P_X+P_X^*)^2+T_X^2\big)=\frac{1}{4}tr(2P_X^2+2P_X^*P_X+T_X^2).$$

\li

Hence, $tr(R_X^*R_X-P_X^2-P_X^*P_X)=-tr(R_X^2)$.

\li

Let us suppose that $<X,X>=1$ and let $\{e_i\}_i$ be an
orthonormal basis  for $\mfm$ with respect to $<,>$ such that
$X=e_1$. Hence, we can apply Lemma \ref{seccurv} to obtain the
following:

$$\bar{rl}Ric(X,X)= & \sum_iK(e_i,X)\\ \\

= &
<\big(R_X^*R_X-P_X^*P_X-B_XC_X-P_X^2+P_{V_X}\big)e_i,e_i>\\
\\

= & tr\big(R_X^*R_X-P_X^2-P_X^*P_X-B_XC_X+P_{V_X}\big)\\ \\

= & -tr(R_X^2+B_XC_X-P_{V_X}).\ear$$

\li

Since both $Ric$ and the map $X\mapsto tr(R_X^2+B_XC_X-P_{V_X})$
are bilinear maps, the identity above holds even if $X$ is not
unit. Hence, for every $X\in\mfm$,

$$Ric(X,X)=-tr(R_X^2+B_XC_X-P_{V_X}).$$

Now we compute $Ric(X,Y)$. Since $Ric$ is a symmetric bilinear
operator we have $$2Ric(X,Y)= Ric(X+Y,X+Y)-Ric(X,X)-Ric(Y,Y).$$

By using the expression above for $Ric(X,X)$ we get

$$\bar{rl} 2Ric(X,Y)= & tr(-R^2_{X+Y}+R^2_X+R^2_Y)+\\ \\

+ & tr(-B_{X+Y}C_{X+Y}+B_XC_X+B_YC_Y)+\\ \\

+ & tr(P_{V_{X+Y}}-P_{V_X}-P_{V_Y}). \ear$$

\li

By bilinearity of $L$ and property (\ref{tf}) we obtain

$$\bar{rl}V_{X+Y}= & L_{X+Y}(X+Y)\\ \\

=& L_XX+L_YY+L_XY+L_YX\\ \\

=& V_X+V_Y+2L_XY-[X,Y]_{\mfm}\\ \\

=& V_X+V_Y+2L_XY+P_XY.\ear$$

\li

Therefore, $P_{V_{X+Y}}-P_{V_X}-P_{V_Y}=P_{2L_XY-P_XY}$.

\li

Also the identity $B_{X+Y}C_{X+Y}=B_XC_X+B_YC_Y+B_XC_Y+B_YC_X$
implies that $$-B_{X+Y}C_{X+Y}+B_XC_X+B_YC_Y=-B_XC_Y-B_YC_X.$$

\li

Moreover, $R^2_{X+Y}=R^2_X+R^2_Y+R_XR_Y+R_YR_X$ and thus

$$tr(-R^2_{X+Y}+R^2_X+R^2_Y)=-tr(R_XR_Y+R_YR_X)=-tr(2R_XR_Y).$$

\li

Therefore, we obtain
$$\bar{rl}2Ric(X,Y)= & -tr(2R_XR_Y+B_XC_Y+B_YC_X)+tr(P_{2L_XY-P_XY})\\ \\

= & -tr(2R_XR_Y+B_XC_Y+B_YC_X)+tr(P_{2U(X,Y)}).\ear$$

$\Box$\eproof

\bcor\label{ric2}Let $X,Y\in\mfm$.

$$Ric(X,Y)=-\frac{1}{4}tr(2P_X^*P_Y+T_XT_Y)-\frac{1}{2}Kill(X,Y)+tr(P_{U(X,Y)}).$$

\ecor

\bproof By using Theorem \ref{ricci} and Lemma \ref{R_X}, we write
$Ric$ as follows:

\beq\label{riccihere}Ric(X,Y)=-\frac{1}{2}\left(\frac{1}{2}(P_X+P_X^*+T_X)(P_Y+P_Y^*+T_Y)+B_XC_Y+B_YC_X-2P_{U(X,Y)}\right).\eeq

Since $P_X+P_X^*$ and $P_Y+P_Y^*$ are symmetric linear maps and
$T_X$ and $T_Y$ are skew-symmetric, we have

\beq\label{trace0}tr((P_X+P_X^*)T_Y)=tr(T_X(P_Y+P_Y^*))=0.\eeq

Moreover,

\beq\label{tracequal}tr(P_XP_Y)=tr(P_X^*P_Y^*))\textrm{ and }
tr(P_X^*P_Y)=tr(P_XP_Y^*)).\eeq

We can use (\ref{adX}) to write $Kill$ as follows:

\beq
\label{killequal}Kill(X,Y)=tr(ad_Xad_Y)=tr(C_XB_Y+B_XC_Y+P_XP_Y).\eeq

Finally, by using (\ref{trace0}), (\ref{tracequal}) and
(\ref{killequal}) we simplify (\ref{riccihere}) to obtain the
expression stated for $Ric$.

$\Box$\eproof

\bdfn A symmetric bilinear map $\be$ on $\mfm$ is said to be
associative if $\be([u,v]_{\mfm},w)=\be(u,[v,w]_{\mfm})$, for
every $u,v,w\in\mfm$. \edfn

\brem \label{ric2rem}If there exists on $\mfm$ an associative
symmetric bilinear form $\be$ such that $\be$ is non-degenerate,
then $tr(P_{U(X,Y)})=0$, for all $X,Y\in\mfm$. Indeed, if such a
bilinear form exists, $tr\,P_a=0$, for every $a\in\mfm$. Let
$\{w_i\}_i$ and $\{w_i'\}_i$ be bases of $\mfm$ dual with respect
to $\be$, i.e., $\be(w_i,w_j')=\de_{ij}$. Then, for every
$a\in\mfm$,

$$\bar{rl} \be(P_aw_i,w_i')= & \be([a,w_i]_{\mfm},w_i')\\ \\

= & -\be(w_i,[a,w_i']_{\mfm})\\ \\

= & -\be(P_aw_i',w_i).\ear$$

Hence, $tr(P_a)=0$. Also, if the metric $g_M$ on $M$ is naturally reductive, then $P_{U(X,Y)}=0$, for all $X,Y\in\mfm$, since, in this case, $U$ is identically zero.

$\diamond$

\erem

\bdfn \label{tracesdef}Let $U$, $V$ be $Ad\,L$-invariant vector
subspaces of $\mfg$. We define a bilinear map $Q^{UV}:\mfm\times\mfm
\rightarrow \reals$ by

$$Q^{UV}_{XY}=tr([X,[Y,\cdot]_V]_U), \,X,Y\in\mfm,$$

where the subscripts $U$ and $V$ denote the projections onto $U$
and $V$, respectively. \edfn

\bdfn \label{Casimirdef}Let $U$ be an $Ad\,L$-invariant vector
subspace of $\mfg$ such that the restriction of Kill to $U$ is
non-degenerate. The Casimir operator of $U$ with respect to the
Killing form of $\mfg$ is the operator

$$C_U=\sum_iad_{u_i}ad_{u_i'},$$

where $\{u_i\}_i$ and $\{u_i'\}_i$ are bases of $U$ which are dual
with respect to Kill, i.e., $Kill(u_i,u_j')=\de_{ij}$.\edfn

The Casimir operator $C_{U}$ is an $Ad\,L$-invariant linear map
and thus it is scalar on any irreducible $Ad\,L$-module. In
particular, if $\mfg$ is simple, then $C_{\mfg}=Id_{\mfg}$.

\blem \label{traces1}Suppose that $\mfg$ is semisimple. Let $U$,
$V$ be $Ad\,L$-invariant vector subspaces of $\mfg$ such that the
restrictions of the Killing form to $U$ and $V$ are both non
degenerate.

(i) $Q^{UV}$ is an $Ad\,L$-invariant symmetric bilinear map.
Hence, if $W\subset\mfg$ is any irreducible $Ad\,L$-submodule,
then $Q^{UV}\mid_{W\times W}$ is a multiple of
 $Kill\mid_{W\times W}$.

(ii) $Q^{UV}=Q^{VU}$. \elem

\bproof Since $\mfg$ is a semisimple Lie algebra its Killing form
is non-degenerate. As in addition $Kill\mid_{U\times U}$ and
$Kill\mid_{V\times V}$ are non-degenerate, we may consider the
orthogonal complements $U^{\perp}$ and $V^{\perp}$ of $U$ and $V$,
respectively, in $\mfg$ with respect to $Kill$. It follows that

$$Kill\mid_{U\times\mfg}=Kill(\cdot,p_U\cdot)\mid_{U\times\mfg}$$

and

$$Kill\mid_{V\times\mfg}=Kill(\cdot,p_V\cdot)\mid_{V\times\mfg},$$

where $p_U$ and $p_V$ are the projections onto $U$ and $V$,
respectively.

Also, we may consider bases $\{w_i\}_i$ and $\{w_i'\}_i$ of $U$
which are dual with respect to $Kill$. By using these facts we
have the following:

\li

(i) Let $X,Y\in\mfm$ and $g\in L$.

$$\bar{rl}Kill([X,[Y,w_i]_V]_U,w_i')= & Kill([X,[Y,w_i]_V],w_i')\\ \\

= & -Kill([Y,w_i]_V,[X,w_i'])\\ \\

= & -Kill([Y,w_i],[X,w_i']_V)\\ \\

= & Kill(w_i,[Y,[X,w_i']_V])\\ \\

= & Kill(w_i,[Y,[X,w_i']_V]_U). \ear$$

Therefore, $tr([X,[Y,\cdot]_V]_U)=tr([Y,[X,\cdot]_V]_U)$ and thus
$Q^{UV}_{XY}=Q^{UV}_{YX}$. So $Q^{UV}$ is symmetric.

To show the $Ad\,L$-invariance of $Q^{UV}$ we note that since $V$
and $V^{\perp}$ are $Ad\,L$-invariant subspaces and $\mfg=V\oplus
V^{\perp}$, the projections on $V$ and $V^{\perp}$ are also
$Ad\,L$-invariant linear maps.

$$\bar{rl}Kill([Ad_gX,[Ad_gY,w_i]_V]_U,w_i')= & Kill([Ad_gX,[Ad_gY,w_i]_V],w_i')\\ \\

= & Kill(Ad_{g^{-1}}[Ad_gX,[Ad_gY,w_i]_V],Ad_{g^{-1}}w_i')\\ \\

= & Kill([X,Ad_{g^{-1}}[Ad_gY,w_i]_V],Ad_{g^{-1}}w_i')\\ \\

= & Kill([X,[Y,Ad_{g^{-1}}w_i]_V],Ad_{g^{-1}}w_i')\\ \\

= & Kill([X,[Y,Ad_{g^{-1}}w_i]_V]_U,Ad_{g^{-1}}w_i').\ear$$

If $\{w_i\}_i$ and $\{w'_i\}_i$ are dual bases of $U$ with respect
to $Kill$, then $\{Ad_{g^{-1}}w_i\}_i$ and $\{Ad_{g^{-1}}w'_i\}_i$
are still dual bases as well since the Killing form is invariant
under inner automorphisms. So by the above we conclude that
$$tr([Ad_gX,[Ad_gY,\cdot]_V]_U)=tr([X,[Y,\cdot]_V]_U)$$ and thus
$Q^{UV}$ is $Ad\,L$-invariant.

\li

(ii) For $Z\in\mfm$, let $A_Z=\left(ad_Z\mid_U\right)_V$ and
$B_Z=\left(ad_Z\mid_V\right)_U$. We have

$$Q_{XY}^{VU}=tr(A_XB_Y)=tr(B_YA_X)=Q_{YX}^{UV}.$$

Hence, by symmetry of $Q^{UV}$, we conclude that
$Q_{XY}^{VU}=Q_{YX}^{UV}=Q_{XY}^{UV}$, for every $X,Y\in\mfm$.
Therefore, $Q^{UV}=Q^{VU}$.

$\Box$\eproof

\li

\blem \label{traces2}Suppose that $\mfg$ is semisimple. Let $U$,
$V$ be $Ad\,L$-invariant vector subspaces of $\mfg$ such that the
restrictions of the Killing form to $U$ and $V$ are both
non-degenerate. For every $X,Y\in\mfm$,

(i) if $ad_XU\subset V$ or $ad_YU\subset V$, then
$Q^{UV}_{XY}=Kill(C_UX,Y)=Kill(X,C_UY)$;

(ii) if $ad_XV\perp U$ or $ad_YV\perp U$, then $Q^{UV}_{XY}=0$;

(iv) if $ad_Xad_YU\perp U$ or $ad_Yad_XU\perp U$, then
$Q^{UV}_{XY}=0$.

\elem

\bproof (i) Let $C_U=\sum_iad_{w_i}ad_{w'_i}$ be the Casimir
operator of $U$. Since $Q^{UV}=Q^{VU}$ it suffices to suppose that
$ad_YU\subset V$. If $ad_YU\subset V$, then
$$Q^{UV}_{XY}=tr([X,[Y,\cdot]]_U)=tr(ad_Xad_Y\mid_U).$$

Since

$$\bar{rl}Kill([X,[Y,w_i]]_U,w'_i) = & Kill([X,[Y,w_i]],w'_i)\\ \\

= & -Kill([Y,w_i],[X,w'_i])\\ \\

= & Kill(Y,[w_i,[w'_i,X]]),\ear$$

we have $Q^{UV}_{XY}=\sum_iKill(Y,[w_i,[w'_i,X]])=Kill(Y,C_UX)$. By
symmetry of $Q^{UV}$ we also get $Q^{UV}_{XY}=Kill(X,C_UY)$.

\li

(ii) If $ad_XV\perp U$, then, for every $w,w'\in U$,
$Kill([X,[Y,w]_V],w')=0$ and thus $Q_{XY}^{UV}=0$, for every $Y\in
\mfm$. By symmetry, the same conclusion holds if $ad_YV\perp U$.

\li

(iii) If $ad_Xad_YU\perp U$, then, for every $w,w'\in U$,
$Kill([X,[Y,w]],w')=0$ and thus $Kill([X,[Y,w]_V]_U,w')=0$. Hence
$Q_{XY}^{UV}=0$. If $ad_Yad_XU\perp U$, then $Q_{XY}^{UV}=0$ by
symmetry.

$\Box$ \eproof

\brem \label{traces3}In Lemmas \ref{traces1} and \ref{traces2} the
condition that $\mfg$ is semisimple may be replaced by requiring
that there is on $\mfg$ a non-degenerate $Ad\,L$-invariant
symmetric bilinear form $\be$, since in the proofs above the
Killing form may be replaced by any such form $\be$. In this case,
the orthogonality conditions in Lemma \ref{traces2} should be
understood as conditions with respect to $\be$.\erem

\bthm \label{ric1}Let $\be$ be an associative $Ad\,L$-invariant
non-degenerate symmetric bilinear form on $\mfm$. Let
$\mfm=\mfm_1\oplus\ldots\oplus\mfm_m$ be a decomposition of $\mfm$
into $Ad\,L$-invariant subspaces such that $\be\mid_{\mfm_j\times
\mfm_i}=0$, if $i\neq j$. Let $g_M$ be the $G$-invariant
pseudo-Riemannian metric on $G/L$ induced by the scalar product of
the form

\beq\label{mdef1}<,>=\oplus_{j=1}^m\nu_j
\be\mid_{\mfm_j\times\mfm_j},\eeq

for $\nu_j>0$, for every $j=1,\ldots,m$. For every $X\in\mfm_a$ and
$Y\in\mfm_b$, the Ricci curvature of $g_M$ is given as follows:

$$Ric(X,Y)=\frac{1}{2}\sum_{j,k=1}^m\left(\frac{\nu_k}{\nu_j}-\frac{\nu_a\nu_b}{2\nu_k\nu_j}\right)Q_{XY}^{\mfm_j\mfm_k}-\frac{1}{2}Kill(X,Y).$$

\ethm

\bproof First we observe that the non-degeneracy of $\be$ and the
condition of pairwise orthogonality of the $\mfm_j$'s imply that
$\be\mid_{\mfm_j\times\mfm_j}$ is in fact non-degenerate. Let
$X\in\mfm_a$ and $Y\in\mfm_b$. By Corollary \ref{ric2} we have

$$Ric(X,Y)=-\frac{1}{4}tr(2P_X^*P_Y+T_XT_Y)-\frac{1}{2}Kill(X,Y)+tr(P_{U(X,Y)}).$$

According to Remark \ref{ric2rem}, we have $tr(P_{U(X,Y)})=0$.

\li

Let $j=1,\ldots,m$ and let $\{w_i\}_i$ and $\{w'_i\}_i$ be dual
bases for $\mfm_j$ with respect to $\be$. We note that such bases
exist as $\be\mid_{\mfm_j\times \mfm_j}$ is non-degenerate.

$$\bar{rl}<T_XT_Yw_i,w_i'>= & <X,[T_Yw_i,w'_i]_{\mfm}>\\ \\

= & \nu_a \be(X,[T_Yw_i,w'_i])\\ \\

= & -\nu_a \be(T_Yw_i,[X,w'_i])\\ \\

= & -\nu_a \sum_{k=1}^m \be(T_Yw_i,[X,w'_i]_{\mfm_k})\\ \\

= & -\nu_a \sum_{k=1}^m\dfrac{1}{\nu_k}<T_Yw_i,[X,w'_i]_{\mfm_k}>\\ \\

= & -\nu_a \sum_{k=1}^m\dfrac{1}{\nu_k}<Y,[w_i,[X,w'_i]_{\mfm_k}]_{\mfm}>\ear$$

$$\bar{rl}
= & -\nu_a\nu_b \sum_{k=1}^m\dfrac{1}{\nu_k}\be(Y,[w_i,[X,w'_i]_{\mfm_k}])\\
\\

= & -\nu_a\nu_b \sum_{k=1}^m\dfrac{1}{\nu_k}\be([Y,w_i],[X,w'_i]_{\mfm_k})\\
\\

= & -\nu_a\nu_b \sum_{k=1}^m\dfrac{1}{\nu_k}\be([Y,w_i]_{\mfm_k},[X,w'_i])\\
\\

= & \nu_a\nu_b \sum_{k=1}^m\dfrac{1}{\nu_k}\be(w'_i,[X,[Y,w_i]_{\mfm_k}])\\
\\

= & \nu_a\nu_b \sum_{k=1}^m\dfrac{1}{\nu_k}\be(w'_i,[X,[Y,w_i]_{\mfm_k}]_{\mfm_j})\\
\\

= & \nu_a\nu_b
\sum_{k=1}^m\dfrac{1}{\nu_k\nu_j}<w'_i,[X,[Y,w_i]_{\mfm_k}]_{\mfm_j}>.\ear$$

\li

Hence,

$$tr(T_XT_Y\mid_{\mfm_j})=\nu_a\nu_b
\sum_{k=1}^m\dfrac{1}{\nu_k\nu_j}tr([X,[Y,\cdot]_{\mfm_k}]_{\mfm_j})=\nu_a\nu_b
\sum_{k=1}^m\dfrac{1}{\nu_k\nu_j}Q_{XY}^{\mfm_j\mfm_k}$$

and thus

$$tr(T_XT_Y)=\nu_a\nu_b
\sum_{j,k=1}^m\frac{1}{\nu_k\nu_j}Q_{XY}^{\mfm_j\mfm_k}.$$

\li

$$\bar{rl}<P_X^*P_Yw_i,w'_i>= & <P_Yw_i,P_Xw'_i>\\ \\

= & \sum_{k=1}^m<[Y,w_i]_{\mfm_k},[X,w'_i]_{\mfm_k}>\\ \\

= & \sum_{k=1}^m\nu_k\be([Y,w_i]_{\mfm_k},[X,w'_i])\\ \\

= & -\sum_{k=1}^m\nu_k\be(w'_i,[X,[Y,w_i]_{\mfm_k}])\\ \\

= & -\sum_{k=1}^m\nu_k\be(w'_i,[X,[Y,w_i]_{\mfm_k}]_{\mfm_j})\\ \\

= &
-\sum_{k=1}^m\dfrac{\nu_k}{\nu_j}<w'_i,[X,[Y,w_i]_{\mfm_k}]_{\mfm_j}>.\ear$$

\li

Then
$$tr(P_X^*P_Y\mid_{\mfm_j})=-\sum_{k=1}^m\dfrac{\nu_k}{\nu_j}tr([X,[Y,\cdot]_{\mfm_k}]_{\mfm_j})=-
\sum_{k=1}^m\dfrac{\nu_k}{\nu_j}Q_{XY}^{\mfm_j\mfm_k}$$

and thus we get

$$tr(P_X^*P_Y)=-
\sum_{j,k=1}^m\frac{\nu_k}{\nu_j}Q_{XY}^{\mfm_j\mfm_k}.$$

By using Corollary \ref{ric2} we finally obtain the required
expression for $Ric(X,Y)$.

$\Box$ \eproof

We recall that a metric $g_M$ is said to be \textit{normal} if it is a
multiple of an associative $Ad\,L$-invariant non-degenerate
symmetric bilinear form on $\mfm$.

\bcor \label{ric1cor1}If $g_M$ is a normal metric, then
$Ric(\mfm_i,\mfm_j)=0$, for all $i\neq j$. For every $X\in\mfm_j$,

$$Ric(X,X)=-\frac{1}{4}Kill(X,X)-\frac{1}{2}Kill(C_{\mfl}X,X),$$

where $C_{\mfl}$ is the Casimir operator of $\mfl$ with respect to
the Killing form. Furthermore, if $\mfm_j$ is irreducible, then

$$Ric\mid_{\mfm_j\times\mfm_j}=-\frac{1}{2}\left(\frac{1}{2}+c_{\mfl,j}\right)Kill\mid_{\mfm_j\times\mfm_j},$$

where $c_{\mfl,j}$ is the eigenvalue of $C_{\mfl}$ on $\mfm_j$.\ecor

\bproof Let $X\in\mfm_a$ and $Y\in\mfm_b$. If $g_M$ is a normal
metric, then there exists an $Ad\,L$-invariant non-degenerate
symmetric bilinear form $\be$ on $\mfm$ which induces $g_M$.
Hence, in Theorem \ref{ric1} we can take $\nu_1=\ldots=\nu_m=1$
and obtain the following:

$$\bar{rl}Ric(X,Y)= &
\frac{1}{4}\sum_{j,k=1}^mQ_{XY}^{\mfm_j\mfm_k}-\frac{1}{2}Kill(X,Y)\\ \\

= & \frac{1}{4}Q_{XY}^{\mfm\mfm}-\frac{1}{2}Kill(X,Y)\\ \\

= & \frac{1}{4}Q_{XY}^{\mfm\mfg}-\frac{1}{4}Q_{XY}^{\mfm\mfl}-\frac{1}{2}Kill(X,Y)\\ \\

= & \frac{1}{4}Q_{XY}^{\mfm\mfg}-\frac{1}{4}Q_{XY}^{\mfl\mfm}-\frac{1}{2}Kill(X,Y)\\ \\

= & \frac{1}{4}Kill(C_{\mfm}X,Y)-\frac{1}{4}Kill(C_{\mfl}X,Y)-\frac{1}{2}Kill(X,Y)\\ \\

= & -\frac{1}{4}Kill(X,Y)-\frac{1}{2}Kill(C_{\mfl}X,Y).\ear$$

Since $C_{\mfl}(\mfm_a)\subset\mfm_a$, it is clear that
$Ric(X,Y)=0$ if $a\neq b$ and $Ric$ is well determined by elements
$Ric(X,X)$ with $X\in\mfm_a$. If $\mfm_j$ is irreducible, then
$C_{\mfl}$ is scalar on $\mfm_j$ and we obtain the identity given
for $Ric$.

$\Box$\eproof

The formula above for the Ricci curvature of a normal metric was first found by M.Y. Wang and W. Ziller in \cite{WZ}. From Corollary \ref{ric1cor1}, it is clear that a necessary and sufficient condition for a normal metric to be Einstein is that the Casimir operator of $\mfl$ is scalar on the isotropy space $\mfm$. For instance, this condition holds if $\mfm$ is irreducible.  Simply connected non-strongly isotropy irreducible homogeneous spaces which admit a normal Einstein metric were classified by M.Y. Wang and W. Ziller in \cite{WZ}, when $G$ is a compact connected simple group. Also, more generally, simply connected compact standard homogeneous manifolds were studied by E.D. Rodionov in \cite{Ro4}.

\li

We obtain a similar formula to that of Corollary \ref{ric1cor1}, in the case when the submodules $\mfm_1,\ldots,\mfm_m$ pairwise commute.

\bcor \label{ric1cor2}If $\mfm_1,\ldots,\mfm_m$ pairwise commute,
i.e., $[\mfm_a,\mfm_b]=0$, for all $a\neq b$, then
$Ric(\mfm_a,\mfm_b)=0$, for all $a\neq b$. For every $X\in\mfm_a$,

$$Ric(X,X)=-\frac{1}{4}Kill(X,X)-\frac{1}{2}Kill(C_{\mfl}X,X),$$

where $C_{\mfl}$ is the Casimir operator of $\mfl$ with respect to
the Killing form. Furthermore, if $\mfm_j$ is irreducible, then

$$Ric\mid_{\mfm_j\times\mfm_j}=-\frac{1}{2}\left(\frac{1}{2}+c_{\mfl,j}\right)Kill\mid_{\mfm_j\times\mfm_j},$$

where $c_{\mfl,j}$ is the eigenvalue of $C_{\mfl}$ on $\mfm_j$.\ecor

\bproof Let $X\in\mfm_a$ and $Y\in\mfm_b$. If
$\mfm_1,\ldots,\mfm_m$ pairwise commute, then, by Lemma
\ref{traces2}, we have $Q_{XY}^{\mfm_j\mfm_k}=0$, for every
$j,k\neq a,b$. In particular, all these bilinear maps vanish if
$a\neq b$. Hence, if $a\neq b$,
$Ric(X,Y)=-\frac{1}{2}Kill(X,Y)=0$. Therefore, the Ricci curvature
is well determined by elements of the form $Ric(X,X)$, with
$X\in\mfm_a$, and

$$Ric(X,X)= \frac{1}{4}Q_{XX}^{\mfm_a\mfm_a}-\frac{1}{2}Kill(X,X).$$

Since $Q_{XX}^{\mfm_a\mfm_a}=Q_{XX}^{\mfm\mfm}$, the rest of the
proof is similar to the proof of Corollary \ref{ric1cor1}.

$\Box$ \eproof

\newpage

\section{Homogeneous Riemannian Fibrations}\label{sectionRF}

\subsection{Notation and Hypothesis}\label{notation1}

In this Section we obtain the Ricci curvature of an invariant
metric with totally geodesic fibers on the total space of a
homogeneous fibration. We start by settling once and for all the
notation used throughout.

\li

Let $G$ be a compact connected semisimple Lie group and
$L\varsubsetneq K\varsubsetneq G$ connected closed non-trivial
subgroups of $G$. We denote $M=G/L$, $N=G/K$ and $F=K/L$. We
consider the natural fibration

\beq\bar{rrcl}\pi: & M & \rightarrow & N\\
& aL & \mapsto & aK  \ear\eeq

with fiber $F$ and structural group $K$. We denote by $\mfg$,
$\mfk$ and $\mfl$ the Lie algebras of $G$, $K$ and $L$,
respectively. By $Kill$ we denote the Killing form of $G$ and we
set $B=-Kill$. Also, we denote the Killing forms of $K$ and $L$ by
$Kill_{\mfk}$ and $Kill_{\mfl}$, respectively. As $G$ is compact
and semisimple, the Killing form of $G$ is negative definite and
thus $B$ is positive definite. We consider an orthogonal
decomposition of $\mfg$ with respect to $B$ given by

\beq\mfg=\mfl\oplus\underbrace{\mfp\oplus\mfn}_{\mfm},\eeq

where $\mfk=\mfl\oplus\mfp$. Clearly, $\mfg=\mfl\oplus\mfm$,
$\mfg=\mfk\oplus\mfn$ and $\mfk=\mfl\oplus\mfp$ are reductive
decompositions for $M$, $N$ and $F$, respectively. Hence, we have
the following inclusions

\beq [\mfk,\mfn]\subset\mfn,\,[\mfl,\mfn]\subset\mfn,\,
[\mfl,\mfp]\subset\mfp\textrm{ and } [\mfp,\mfn]\subset\mfn.\eeq

\li

An $Ad\,K$-invariant scalar product on $\mfn$ induces a
$G$-invariant Riemannian metric $g_N$ on $N$ and an
$Ad\,L$-invariant scalar product on $\mfp$ induces a $G$-invariant
Riemannian metric $g_F$ on $F$. The orthogonal direct sum of these
scalar products on $\mfm$ defines a $G$-invariant Riemannian
metric $g_M$ on $M$ which projects onto a $G$-invariant metric on
$N$. Moreover, if $\mfp$ and $\mfn$ do not contain any equivalent
$Ad\,L$-submodules, then any $G$-invariant metric which projects
onto a $G$-invariant metric on $N$ is constructed in this fashion.
We recall the following result due to L.B\'{e}rard-Bergery
(\cite{BB}, \cite{Be} 9 \S H):

\bthm The natural projection $\pi: M
 \rightarrow N$ is a Riemannian submersion from $(M,g_M)$ to $(N,g_N)$ with totally geodesic
 fibers.\ethm

Throughout this thesis we shall refer to such a metric $g_M$ as an
adapted metric:

\bdfn An \textbf{adapted} metric on $M$ is a $G$-invariant
Riemannian metric $g_M$  such that the natural projection $\pi: M
 \rightarrow N$ is a Riemannian submersion and consequently the fibers are totally geodesic submanifolds.   The fibration $M\rightarrow N$ equipped
with an adapted metric $g$ is then called a \textbf{Riemannian
fibration}. \edfn

An adapted metric on $M$ shall be denoted by $g_M$ and, as already
introduced above, $g_F$ and $g_N$ shall denote the projection of
$g_M$ onto the base space $N$ and $g_F$ its restriction to the
fiber $F$.

\li

We consider a decomposition $\mfp=\mfp_1\oplus\ldots\oplus\mfp_s$
into irreducible $Ad\,L$-submodules pairwise orthogonal with
respect to B. Also let $\mfn=\mfn_1\oplus\ldots\oplus\mfn_n$ be an
orthogonal decomposition into irreducible $Ad\,K$-submodules.
Throughout we assume the following hypothesis:

$$\bar{l}\textrm{(i) $\mfp_1,\ldots,\mfp_s$ are pairwise
inequivalent irreducible $Ad\,L$-submodules;}\\
\textrm{(ii) $\mfn_1,\ldots,\mfn_n$ are
pairwise inequivalent irreducible $Ad\,K$-submodules;}\\
\textrm{(iii) $\mfp$ and $\mfn$ do not contain equivalent
$Ad\,L$-submodules.} \ear$$

We shall refer to this hypothesis by saying that $M$ has
\textbf{simple spectrum}. Under this hypothesis, according to
Schur's Lemma, any $Ad\,L$-invariant scalar product on
$\mfm=\mfp\oplus\mfn$ which restricts to an $Ad\,K$-invariant
scalar product on $\mfn$ is of the form

\beq\label{mdef}<,>=\left(\oplus_{a=1}^s\la_aB\mid_{\mfp_a\times\mfp_a}\right)\oplus\left(\oplus_{k=1}^n\mu_kB\mid_{\mfn_k\times\mfn_k}\right),\eeq

for some $\la_1,\ldots,\la_s,\mu_1,\ldots,\mu_n>0$. Since an
adapted metric $g_M$ on $M$ projects onto a $G$-invariant
Riemannian metric on $N$, $g_M$ is necessarily induced by a scalar
product of the form (\ref{mdef}). To denote that $g_M$ is induced
by a scalar product as in (\ref{mdef}) we shall write

\beq\label{nuplemet}
g_M=g_M(\la_1,\ldots,\la_s;\,\mu_1,\ldots,\mu_n).\eeq

Similarly, we write

\beq g_F=g_F(\la_1,\ldots,\la_s) \textrm{ and
}g_N=g_N(\,\mu_1,\ldots,\mu_n).\eeq

\li

By $Ric$ we mean the Ricci curvature of $g_M$ and by $Ric^F$ and
$Ric^N$ the Ricci curvature of $g_N$ and $g_F$, respectively.

\li

In the following Sections we compute the Ricci curvature for $g_M$
and find some necessary conditions so that $g_M$ is an Einstein
metric. We recall that in Theorem \ref{ric1} we have shown that
the Ricci curvature of any metric on $M$ can be described using
the bilinear maps $Q_{XY}$ defined in \ref{tracesdef}.

\subsection{The Ricci Curvature in the Direction of the Fiber}

In this Section we obtain the Ricci curvature of the adapted metric

$$g_M=g_M(\la_1,\ldots,\la_s;\,\mu_1,\ldots,\mu_n)$$

in the vertical direction $\mfp$. We recall that $\mfp$ decomposes
into the direct sum of the pairwise inequivalent irreducible
$Ad\,L$-submodules $\mfp_1,\ldots,\mfp_s$ and, as explained above,
$g_M$ is induced by the scalar product (\ref{mdef})

$$\left(\oplus_{a=1}^s\la_aB\mid_{\mfp_a\times\mfp_a}\right)\oplus\left(\oplus_{k=1}^n\mu_kB\mid_{\mfn_k\times\mfn_k}\right),$$

while $g_F$ is the restriction of $g_M$ to the fiber, i.e.,

$$g_F=g_F(\la_1,\ldots,\la_s).$$

\blem \label{tracesp}Let $X\in\mfp$ and $Y\in\mfm$.

(i) $Q_{XY}^{\mfn_j\mfp_a}=Q_{XY}^{\mfp_a\mfn_j}=0$,
$a=1,\ldots,s$, $j=1,\ldots,n$;

(ii) $Q_{XY}^{\mfn_i\mfn_j}=0$, if $i\neq j$, $i,j=1,\ldots,n$;

(iii) $Q_{XY}^{\mfn_j\mfn_j}=Kill(C_{\mfn_j}X,Y)$, $j=1,\ldots,n$.
\elem

\bproof Since $ad_X\mfp\subset\mfk\perp\mfn$ we have
$Q_{XY}^{\mfn_j\mfp_a}=0$, from Lemma \ref{traces2}. From Lemma
\ref{traces1} (ii),
$Q_{XY}^{\mfp_a\mfn_j}=Q_{XY}^{\mfn_j\mfp_a}=0$.

As $ad_X\mfn_j\subset\mfn_j$ ,we have
$Q_{XY}^{\mfn_j\mfn_j}=Kill(C_{\mfn_j}X,Y)$. Moreover, since
$\mfn_j\perp\mfn_i$, for every $i\neq j$, we also conclude that
$Q_{XY}^{\mfn_i\mfn_j}=0$, if $i\neq j$.

$\Box$ \eproof

\blem\label{ricF}The Ricci curvature of
$g_F=g_F(\la_1,\ldots,\la_s)$ is of the form
$Ric^F=\oplus_{a=1}^sq_aB\mid_{\mfp_a\times\mfp_a}$, with

$$q_a=\frac{1}{2}\sum_{b,c=1}^s\left(\frac{\la_a^2}{2\la_c\la_b}-\frac{\la_c}{\la_b}\right)q^{cb}_a+\frac{\ga_a}{2}$$

where the constants $q^{cb}_a$ and $\ga_a$ are such that
$$Kill_{\mfk}\mid_{\mfp_a\times\mfp_a}=\ga_aKill\mid_{\mfp_a\times\mfp_a}$$
and
$$Q^{\mfp_b\mfp_c}\mid_{\mfp_a\times\mfp_a}=q^{cb}_aKill\mid_{\mfp_a\times\mfp_a}.$$

In particular, $Ric^F(\mfp_a,\mfp_b)=0$, if $a\neq b$. \elem

\bproof  Since $\mfp_1,\ldots,\mfp_s$ are pairwise inequivalent
irreducible $Ad\,L$-submodules and $Ric^F$ is an $Ad\,L$-invariant
symmetric bilinear form, we may write
$Ric^F=\oplus_{a=1}^sq_aB\mid_{\mfp_a\times\mfp_a}$, for some
constants $q_1,\ldots,q_s$. In particular, we have
$Ric^F(\mfp_a,\mfp_b)=0$, for every $a,b=1,\ldots,s$ such that
$a\neq b$. By Theorem \ref{ric1}, the Ricci curvature of $g_F$ is

$$Ric^F(X,X)=\dfrac{1}{2}\sum_{b,c=1}^s\left(\frac{\la_b}{\la_c}-\dfrac{\la_a^2}{2\la_c\la_b}\right)Q_{XX}^{\mfp_c\mfp_b}-\frac{1}{2}Kill_{\mfk}(X,X).$$

By Lemma \ref{traces1} the maps $Q^{\mfp_c\mfp_b}$ are
$Ad\,L$-invariant symmetric bilinear maps. Since $\mfp_a$ is
$Ad\,L$-irreducible, there are constants $q^{cb}_a$ such that

$$Q^{\mfp_c\mfp_b}\mid_{\mfp_a\times\mfp_a}=q^{cb}_aKill\mid_{\mfp_a\times\mfp_a}.$$

Similarly, the Killing form of $\mfk$, $Kill_{\mfk}$, is an
$Ad\,L$-invariant symmetric bilinear form on $\mfp_a$. So, by
irreducibility of $\mfp_a$, there is a constant $\ga_a$ such that

$$Kill(C_{\mfk}\cdot,\cdot)\mid_{\mfp_a\times\mfp_a}=\ga_aKill\mid_{\mfp_a\times\mfp_a}.$$

By the expression above for $Ric^F$, we must have

$$q_a=\frac{1}{2}\sum_{b,c=1}^s\left(\frac{\la_a^2}{2\la_c\la_b}-\frac{\la_c}{\la_b}\right)q^{cb}_a+\frac{\ga_a}{2}$$

and the result follows from this.

$\Box$ \eproof

\li

\bprop \label{riccip}Let
$g_M=g_M(\la_1,\ldots,\la_s;\,\mu_1,\ldots,\mu_n)$ be an adapted
metric on $M$.

For every $a,b=1,\ldots,s$ such that $a\neq b$,
$Ric(\mfp_a,\mfp_b)=0$. For every $X\in\mfp_a$, $a=1,\ldots,s$,

$$Ric(X,X)=\left(q_a+\dfrac{\la_a^2}{4}\sum_{j=1}^n\dfrac{c_{\mfn_j,a}}{\mu_j^2}\right)B(X,X),$$

where, for $j=1,\ldots,n$, the constants $c_{\mfn_j,a}$ are such
that

$$Kill(C_{\mfn_j}\cdot,\cdot)\mid_{\mfp_a\times\mfp_a}=c_{\mfn_j,a}Kill\mid_{\mfp_a\times\mfp_a}$$

and $C_{\mfn_j}$ is the Casimir operator of $\mfn_j$ with respect
to $Kill$. The constant $q_a$ is such that

$$q_a=\frac{1}{2}\sum_{b,c=1}^s\left(\frac{\la_a^2}{2\la_c\la_b}-\frac{\la_c}{\la_b}\right)q^{cb}_a+\frac{\ga_a}{2},$$

with $q^{cb}_a$ and $\ga_a$ defined by

$$\bar{l}Kill_{\mfk}\mid_{\mfp_a\times\mfp_a}=\ga_aKill\mid_{\mfp_a\times\mfp_a}\\
\\
Q^{\mfp_b\mfp_c}\mid_{\mfp_a\times\mfp_a}=q^{cb}_aKill\mid_{\mfp_a\times\mfp_a}.\ear$$\eprop

\bproof Since $\mfp_1,\ldots,\mfp_s$ are pairwise inequivalent
irreducible $Ad\,L$-submodules and $Ric\mid_{\mfp\times\mfp}$ is
an $Ad\,L$-invariant symmetric bilinear form, we have that
$Ric\mid_{\mfp\times\mfp}$ is diagonal, i.e.,

$$Ric\mid_{\mfp\times\mfp}=a_1B\mid_{\mfp_1\times\mfp_1}\oplus\ldots\oplus a_sB\mid_{\mfp_s\times\mfp_s},$$

for some constants $a_1,\ldots,a_s$. In particular, we have
$Ric(\mfp_a,\mfp_b)=0$, for every $a,b=1,\ldots,s$ such that $a\neq
b$. Hence, $Ric\mid_{\mfp\times\mfp}$ is determined by elements
$Ric(X,X)$ with $X\in\mfp_a$, $a=1,\ldots,s$.

\li

By Lemma \ref{tracesp} we obtain that only
$Q_{XX}^{\mfn_j\mfn_j}=Kill(C_{\mfn_j}X,X)$ and
$Q_{XX}^{\mfp_b\mfp_c}$ are non-zero. Therefore, by Theorem
\ref{ric1} we obtain that

$$Ric(X,X)=$$

$$\frac{1}{2}\sum_{j,k=1}^s\left(\frac{\la_k}{\la_j}-\frac{\la_a^2}{2\la_j\la_k}\right)Q_{XX}^{\mfp_j\mfp_k}+\frac{1}{2}\sum_{j=1}^n\left(1-\frac{\la_a^2}{2\mu_j^2}\right)Q_{XX}^{\mfn_j\mfn_j}-\frac{1}{2}Kill(X,X).$$

\li

We have
$$\bar{rl}\sum_{j=1}^mQ_{XX}^{\mfn_j\mfn_j}= & \sum_{j=1}^mKill(C_{\mfn_j}X,Y)\\ \\

= & Kill(C_{\mfn}X,X)\\
\\

= & Kill(X,X)-Kill(C_{\mfk}X,X)\\ \\

= & Kill(X,X)-Kill_{\mfk}(X,X).\ear$$

Hence we can rewrite $Ric(X,X)$ as follows:

$$\underbrace{\frac{1}{2}\sum_{a,b=1}^s\left(\frac{\la_b}{\la_c}-\frac{\la_a^2}{2\la_c\la_b}\right)Q_{XX}^{\mfp_c\mfp_b}-\frac{1}{2}Kill_{\mfk}(X,X)}_{(1)}-\frac{1}{2}\sum_{j=1}^n\frac{\la_a^2}{2\mu_j^2}Kill(C_{\mfn_j}X,X).$$

As we saw in the proof of Lemma \ref{ricF}, the summand (1) is
just $Ric^F(X,X)=q_aB(X,X)$.

Furthermore, since $Kill(C_{\mfn_j}\cdot,\cdot)=Q^{\mfn_j\mfn_j}$
are $Ad\,L$-invariant symmetric bilinear maps (Lemma \ref{traces1})
and $\mfp_a$ is $Ad\,L$-irreducible, there are constants
$c_{\mfn_j,a}$ such that

$$Kill(C_{\mfn_j}\cdot,\cdot)\mid_{\mfp_a\times\mfp_a}=c_{\mfn_j,a}Kill\mid_{\mfp_a\times\mfp_a}.$$

Therefore,

$$Ric(X,X)=\left(q_a+\frac{1}{4}\sum_{j=1}^n\frac{\la_a^2}{\mu_j^2}c_{\mfn_j,a}\right)B(X,X).$$

$\Box$\eproof

\subsection{The Ricci Curvature in the Horizontal Direction}

We obtain the Ricci curvature of an adapted metric
$g_M=g_M(\la_1,\ldots,\la_s;\,\mu_1,\ldots,\mu_n)$ in the
horizontal direction $\mfn$. We recall that $\mfn$ decomposes into
the direct sum of the pairwise inequivalent irreducible
$Ad\,K$-submodules $\mfn_1,\ldots,\mfn_n$ and, as explained above,
$g_M$ is induced by the scalar product (\ref{mdef})

$$\left(\oplus_{a=1}^s\la_aB\mid_{\mfp_a\times\mfp_a}\right)\oplus\left(\oplus_{k=1}^n\mu_kB\mid_{\mfn_k\times\mfn_k}\right)$$

and $g_N$ is the projection of $g_M$ onto the base space, i.e.,

$$g_N=g_N(\mu_1,\ldots,\mu_n).$$

\blem \label{tracesn}Let $X\in\mfn_k$ and $Y\in\mfm$.

(i) $Q_{XY}^{\mfn_j\mfp_a}=Q_{XY}^{\mfp_a\mfn_j}=0$, for $j\neq
k$;

(ii)
$Q_{XY}^{\mfp_a\mfn_k}=Q_{XY}^{\mfn_k\mfp_a}=Kill(C_{\mfp_a}X,Y)$;

(iii) $Q_{XY}^{\mfp_j\mfp_a}=0$, for $i,j=1,\ldots,s$. \elem

\bproof We have $ad_X\mfp_a\subset\mfn_k\perp\mfp,\mfn_j$, for
$j\neq k$. Thus, $Q_{XY}^{\mfn_j\mfp_a}=0$, for $j\neq k$ and
$Q_{XY}^{\mfp_j\mfp_a}=0$, for $i,j=1,\ldots,s$, from Lemma
\ref{traces2}. Also from $ad_X\mfp_a\subset\mfn_k$ we deduce that
$Q_{XY}^{\mfp_a\mfn_k}=Kill(C_{\mfp_a}X,Y)$. From Lemma
\ref{traces1} (ii), we also obtain
$Q_{XY}^{\mfp_a\mfn_j}=Q_{XY}^{\mfn_j\mfp_a}=0$, for $j\neq k$ and
$Q_{XY}^{\mfn_k\mfp_a}=Q_{XY}^{\mfp_a\mfn_k}=Kill(C_{\mfp_a}X,Y)$,
for $j=k$.

$\Box$ \eproof

\blem\label{ricN}The Ricci curvature of
$g_N=g_N(\mu_1,\ldots,\mu_n)$ is of the form
$Ric^N=\oplus_{k=1}^nr_kB\mid_{\mfn_k\times\mfn_k}$, where

$$r_k=\frac{1}{2}\sum_{j,i=1}^n\left(\frac{\mu_a^2}{2\mu_i\mu_j}-\frac{\mu_i}{\mu_j}\right)r^{ji}_k+\frac{1}{2}$$

and the constants $r^{ji}_k$ are such that
$$Q^{\mfn_j\mfn_i}\mid_{\mfn_k\times\mfn_k}=r^{ji}_kKill\mid_{\mfn_k\times\mfn_k}.$$

In particular, $Ric^N(\mfn_k,\mfn_j)=0$, if $k\neq j$. \elem

\bproof The metric $g_N$ is induced by an $Ad\,K$-invariant scalar
product on $\mfn$. Hence, $Ric^N$ is an $Ad\,K$-invariant symmetric
bilinear form on $\mfn$. Since the subspaces $\mfn_j$,
$j=1,\ldots,n$, are irreducible pairwise inequivalent
$Ad\,K$-submodules, we may write

$$Ric^N=\oplus_{k=1}^mr_kB\mid_{\mfn_k\times\mfn_k},$$

for some constants $r_1,\ldots,r_n$. It is clear that
$Ric^N(\mfn_k,\mfn_j)=0$, if $k\neq j$.

It follows from Theorem \ref{ric1} that, for every $X\in\mfn_k$,

$$Ric^N(X,X)=\frac{1}{2}\sum_{j,i=1}^n\left(\frac{\mu_i}{\mu_j}-\frac{\mu_k^2}{2\mu_i\mu_j}\right)Q_{XX}^{\mfn_j\mfn_i}-\frac{1}{2}Kill(X,X).$$

Since each subspace $\mfn_j$ is $Ad\,K$-invariant, the bilinear
maps $Q_{XY}^{\mfn_j\mfn_i}$ are $Ad\,K$-invariant symmetric
bilinear maps (Lemma \ref{traces1}). Hence, by irreducibility of
the $\mfn_k$'s it follows that
$Q^{\mfn_j\mfn_i}\mid_{\mfn_k\times\mfn_k}=r^{ji}_k
Kill\mid_{\mfn_k\times\mfn_k},$ for some constants $r^{ji}_k$,
$j,i,k=1,\ldots,n$.

\li

By definition of the $r_k$'s it must be

$$r_k=\frac{1}{2}\sum_{j,i=1}^n\left(\frac{\mu_k^2}{2\mu_i\mu_j}-\frac{\mu_i}{\mu_j}\right)r^{ji}_a+\frac{1}{2}.$$

$\Box$\eproof

\bprop \label{riccin}Let
$g_M=g_M(\la_1,\ldots,\la_s;\,\mu_1,\ldots,\mu_n)$ be an adapted
metric on $M$. For every $k,j=1,\ldots,n$, such that $j\neq k$,
$Ric(\mfn_k,\mfn_j)=0$. For every $X\in\mfn_k$,

$$Ric(X,X)=-\dfrac{1}{2\mu_k}\sum_{a=1}^s\la_aB(C_{\mfp_a}X,X)+r_kB(X,X),$$

where, for every $a=1,\ldots,s$, $C_{\mfp_a}$ is the Casimir
operator of $\mfp_a$ with respect to $Kill$,

$$r_k=\frac{1}{2}\sum_{j,i=1}^n\left(\frac{\mu_a^2}{2\mu_i\mu_j}-\frac{\mu_i}{\mu_j}\right)r^{ji}_k+\frac{1}{2}$$

and the constants $r^{ji}_k$ are such that

$$Q^{\mfn_j\mfn_i}\mid_{\mfn_k\times\mfn_k}=r^{ji}_kKill\mid_{\mfn_k\times\mfn_k}.$$\eprop

\bproof Let $X\in\mfn_k$ and $Y\in\mfn_{k'}$. By Lemma
\ref{tracesn} we have $Q_{XY}^{\mfp_a\mfp_b}=0$, for every
$a,b=1,\ldots,s$. Also,
$Q_{XY}^{\mfn_j\mfp_a}=Q_{XY}^{\mfp_a\mfn_j}=0$, if $j\neq k,k'$.
Therefore, it follows from Theorem \ref{ric1} and Lemma \ref{ricN}
that, if $k\neq k'$, then

$$Ric(X,Y)=\frac{1}{2}\sum_{j,i=1}^n\left(\frac{\mu_i}{\mu_j}-\frac{\mu_k\mu_{k'}}{2\mu_i\mu_j}\right)Q_{XY}^{\mfn_j\mfn_i}-\frac{1}{2}Kill(X,Y)=Ric^N(X,Y)=0.$$

Hence, $Ric\mid_{\mfn\times\mfn}$ is determined by elements
$Ric(X,X)$ with $X\in\mfn_k$, $k=1,\ldots,n$. For $X\in\mfn_k$, by
Theorem \ref{ric1}, we get

$$Ric(X,X)=$$

$$\frac{1}{2}\sum_{k=1}^s\left(\frac{\mu_k}{\la_a}-\frac{\mu_k^2}{2\mu_
k\la_a}\right)Q_{XX}^{\mfp_a\mfn_k}+\frac{1}{2}\sum_{a=1}^s\left(\frac{\la_a}{\mu_k}-\frac{\mu_k^2}{2\mu_
k\la_a}\right)Q_{XX}^{\mfn_k\mfp_a}+Ric^N(X,X).$$

From Lemma \ref{tracesn}, we know that
$Q_{XX}^{\mfn_k\mfp_a}=Q_{XX}^{\mfp_a\mfn_k}=Kill(C_{\mfp_a}X,X)$.
Hence, we simplify the expression above obtaining

$$Ric(X,X)=\dfrac{1}{2}\sum_{a=1}^s\dfrac{\la_a}{\mu_k}Kill(C_{\mfp_a}X,X)+Ric^N(X,X).$$
Finally, since $Ric^N(X,X)=r_kB(X,X)$, using the notation of Lemma
\ref{ricN}, we conclude that

$$Ric(X,X)=-\dfrac{1}{2}\sum_{a=1}^s\dfrac{\la_a}{\mu_k}B(C_{\mfp_a}X,X)+r_kB(X,X).$$

$\Box$\eproof

\subsection{The Ricci Curvature in the Mixed Direction}

In the previous two sections we determined the Ricci curvature of

$$g_M=g_M(\la_1,\ldots,\la_s;\,\mu_1,\ldots,\mu_n)$$

on the directions of $\mfp$ and $\mfn$. Here we obtain the Ricci
curvature in the direction of $\mfp\times\mfn$.

\bprop \label{riccipn}Let
$g_M=g_M(\la_1,\ldots,\la_s;\,\mu_1,\ldots,\mu_n)$ be an adapted
metric on $M$. For every $X\in\mfp_a$ and $Y\in\mfn_k$,
$$Ric(X,Y)=\dfrac{\la_a\mu_k}{4}\sum_{j=1}^n\dfrac{B(C_{\mfn_j}X,Y)}{\mu_j^2},$$

where, for every $j=1,\ldots,n$, $C_{\mfn_j}$ is the Casimir
operator of $\mfn_j$ with respect to $Kill$. \eprop

\bproof For $X\in\mfp$ we know from Lemma \ref{tracesp} that
$Q_{XY}^{\mfn_j\mfp_a}=Q_{XY}^{\mfp_a\mfn_j}=0$, for every
$a=1,\ldots,s,\,j=1,\ldots,n$, and $Q_{XY}^{\mfn_i\mfn_j}=0$, if
$i\neq j$, $i,j=1,\ldots,n$, whereas
$Q_{XY}^{\mfn_j\mfn_j}=Kill(C_{\mfn_j}X,Y)$, $j=1,\ldots,n$.
Moreover, for $Y\in\mfn_k$, since
$ad_Xad_Y\mfp\subset\mfn_k\perp\mfp$, from Lemma \ref{tracesp} we
also obtain that $Q_{XY}^{\mfp_b\mfp_c}=0$, for every
$b,c=1,\ldots,s$. Therefore, only
$Q_{XY}^{\mfn_j\mfn_j}=Kill(C_{\mfn_j}X,Y)$, $j=1,\ldots,n$, may
not be zero. Furthermore, $Kill(X,Y)=0$. Hence, from Theorem
\ref{ric1} we get

$$Ric(X,Y)=\frac{1}{2}\sum_{j=1}^n\left(1-\frac{\la_a\mu_k}{\mu_j^2}\right)Q_{XY}^{\mfn_j\mfn_j}=\frac{1}{2}\sum_{j=1}^n\left(1-\frac{\la_a\mu_k}{\mu_j^2}\right)Kill(C_{\mfn_j}X,Y).$$

On the other hand,
$$\sum_{j=1}^nKill(C_{\mfn_j}X,Y)=Kill(C_{\mfn}X,Y)=Kill(X,Y)-Kill(C_{\mfk}X,Y)=0,$$
since $C_{\mfk}\mfp\subset\mfk\perp\mfn$. Therefore,

$$Ric(X,Y)=-\frac{\la_a\mu_k}{4}\sum_{j=1}^n\frac{Kill(C_{\mfn_j}X,Y)}{\mu_j^2}.$$

$\Box$\eproof

\subsection{Necessary Conditions for the Existence of an Adapted Einstein Metric}

From the expressions obtained previously for the Ricci curvature
in the horizontal direction and in the direction of
$\mfp\times\mfn$ we obtain two necessary conditions for the
existence of an adapted Einstein metric on $M$. These are
restrictions on Casimir operators and shall be extremely useful in
the chapters ahead.

\bcor \label{cond1}Let
$g_M=g_M(\la_1,\ldots,\la_s;\,\mu_1,\ldots,\mu_n)$ be an adapted
metric on $M$. If $g_M$ is Einstein, then the operator
$\sum_{a=1}^s\la_aC_{\mfp_a}\mid_{\mfn_k}$ is scalar. \ecor

\bproof Let $g_M$ be an adapted metric as defined in (\ref{mdef}).
If $g_M$ is Einstein with Einstein constant $E$, then,
$Ric\mid_{\mfn\times\mfn}=E<,>\mid_{\mfn\times\mfn}$ and thus
$Ric\mid_{\mfn\times\mfn}$ is $Ad\,K$-invariant. Therefore, by
Proposition \ref{riccin}, we conclude that
$\sum_{a=1}^s\la_aC_{\mfp_a}\mid_{\mfn}$ has to be
$Ad\,K$-invariant. Hence, $\sum_{a=1}^s\la_aC_{\mfp_a}\mid_{\mfn_k}$
is scalar, by irreducibility of $\mfn_k$.

$\Box$ \eproof

\bcor \label{cond2}Let
$g_M=g_M(\la_1,\ldots,\la_s;\,\mu_1,\ldots,\mu_n)$  be an adapted
metric on $M$. The orthogonality condition $Ric(\mfp,\mfn)=0$
holds if and only if

\beq\label{ort1}\sum_{j=1}^n\frac{1}{\mu_j^2}C_{\mfn_j}(\mfp)\subset
\mfk.\eeq

Moreover, if $g_M$ is Einstein, then (\ref{ort1}) holds. \ecor

\bproof Let $g_M$ be an adapted metric of the form
$g_M(\la_1,\ldots,\la_s;\,\mu_1,\ldots,\mu_n)$. From Proposition
\ref{riccipn}, we obtain that $Ric(\mfp,\mfn)=0$ if and only if,
for every $X\in\mfp_a$ and $Y\in\mfn_b$,

$$Kill\left(\sum_{j=1}^n\frac{C_{\mfn_j}}{\mu_j^2}X,Y\right)=0.$$

This holds if only if
$\sum_{j=1}^n\frac{C_{\mfn_j}}{\mu_j^2}X\subset \mfk$, for every
$X\in\mfp$.

If $g_M$ is Einstein with Einstein constant $E$, then
$Ric(\mfp,\mfn)=E<\mfp,\mfn>=0$. This shows the last assertion of
the Corollary.

$\Box$ \eproof

The two previous Corollaries may be restated in the following way,
which emphasizes the fact that the two necessary conditions
obtained for existence of an Einstein adapted metric are just
algebraic conditions on the Casimir operators of the submodules
$\mfp_a$ and $\mfn_k$.

\bcor \label{cond3}If there exists on $M$ an Einstein adapted
metric, then there are positive constants $\la_1,\ldots,\la_s$
such that the operator $\sum_{a=1}^s\la_aC_{\mfp_a}\mid_{\mfn_k}$
is scalar. Furthermore, if $g_N$ is not the standard metric, then
there are positive constants $\nu_1,\ldots,\nu_n$, not all equal,
such that

$$\sum_{j=1}^n\nu_jC_{\mfn_j}(\mfp)\subset \mfk.$$\ecor

\bproof The assertions follow from Corollaries \ref{cond1}
 and \ref{cond2}.  In \ref{cond2} we set $\nu_k=1/\mu_k^2$. Hence, $\nu_1=\ldots=\nu_n$ occurs when $g_N$ is the standard metric. Moreover, if $\nu_1=\ldots=\nu_n$, the inclusion in \ref{cond2} is equivalent to $C_{\mfn}(\mfp)\subset\mfk$, which always holds since $C_{\mfn}=Id-C_{\mfk}$ and $C_{\mfk}$ maps $\mfp$ into $\mfk$. So we obtain a condition on the $C_{\mfn_j}$'s only when $g_N$ is not standard.

 $\Box$\eproof

\newpage

\chapter{}

\addtolength{\myVSpace}{0.1cm}

As in Chapter 1, we consider a homogeneous fibration $F\rightarrow
M \rightarrow N$, for $M=G/L$, $N=G/K$ and $F=K/L$, where $G$ is a
compact connected semisimple Lie group and $L\varsubsetneq
K\varsubsetneq G$ are connected closed non-trivial subgroups, and
an adapted metric $g_M$ on $M$. We consider some particular cases
by imposing restrictions on the metric $g_M$, which shall lead to
simpler expressions of the Ricci curvature and thus allow us to
determine further conditions for the existence of an Einstein
adapted metric. Unless otherwise stated we follow the notation
used in Chapter 1. Thus, as before, $\mfp_1,\ldots,\mfp_s$ are the
irreducible pairwise inequivalent $Ad\,L$-submodules of $\mfp$,
the tangent space to the fiber, and $\mfn_1,\ldots,\mfn_n$ are the
irreducible pairwise inequivalent $Ad\,K$-submodules of $\mfn$,
the tangent space to the base. An adapted metric $g_M$ on $M$ is
written as

$$g_M=g_M(\la_1,\ldots,\la_s;\,\mu_1,\ldots,\mu_n)$$

meaning that $g_M$ is induced by the scalar product

$$\left(\oplus_{a=1}^s\la_aB\mid_{\mfp_a\times\mfp_a}\right)\oplus\left(\oplus_{k=1}^n\mu_kB\mid_{\mfn_k\times\mfn_k}\right),$$

on the tangent space $\mfm=\mfp\oplus\mfn$ of $M$. We assume that
$M$ has simple spectrum as in Section \ref{notation1}.

\section{Riemannian Fibrations with Normal
Fiber}\label{normalfiber}

In this section we consider an adapted metric $g_M$ whose
restriction to the fiber, $g_F$, is a multiple of the Killing form
of $\mfg$. Hence, we have

\beq\label{mdefnf}g_M=g_M(\underbrace{\la,\ldots,\la}_s;\,\mu_1,\ldots,\mu_n)
\eeq

and \beq g_F=g_F(\underbrace{\la,\ldots,\la}_s)\eeq

by setting $\la_1=\ldots=\la_s=\la$ in (\ref{mdef}) and
(\ref{nuplemet}). Clearly, when equipped with $g_F$, $F$ becomes a
normal Riemannian manifold. In particular, if the Killing form of
$\mfk$ is a multiple of the Killing form of $\mfg$, then $F$ is a
standard Riemannian manifold. This shall be the case when, for
instance, $\mfp$ is irreducible, but it will not be the case in
general.

\bprop \label{riccinf}Let $g_M$ be an adapted metric on $M$ of the
form
$$g_M(\underbrace{\la,\ldots,\la}_s;\,\mu_1,\ldots,\mu_n).$$ The
Ricci curvature of $g_M$ is as follows:

\li

(i) For every $X\in\mfp_a$,

$$Ric(X,X)=\left(q_a+\frac{\la^2}{4}\sum_{j=1}^n\dfrac{c_{\mfn_j,a}}{\mu_j^2}\right)B(X,X),$$

with $$q_a=\frac{1}{2}\left(c_{\mfl,a}+\frac{\ga_a}{2}\right),$$

where $c_{\mfl,a}$ is the eigenvalue of the Casimir operator of
$\mfl$ on $\mfp_a$, $\ga_a$ is defined by

$$Kill_{\mfk}\mid_{\mfp_a\times\mfp_a}=\ga_aKill\mid_{\mfp_a\times\mfp_a}$$

and $c_{\mfn_j,a}$ is given by

$$Kill(C_{\mfn_j}\cdot,\cdot)\mid_{\mfp_a\times\mfp_a}=c_{\mfn_j,a}Kill\mid_{\mfp_a\times\mfp_a},$$

where $C_{\mfn_j}$ is the Casimir operator of $\mfn_j$ with
respect to $Kill$.

(ii) For every $X\in\mfn_k$,

$$Ric(X,X)=-\frac{\la}{2\mu_k}B(C_{\mfp}X,X)+r_kB(X,X),$$

where $r_k$ is as defined in Lemma \ref{ricN};

(iii) For every $X\in\mfp_a$ and $Y\in\mfn_k$,
$$Ric(X,Y)=\dfrac{\la\mu_k}{4}\sum_{j=1}^n\dfrac{B(C_{\mfn_j}X,Y)}{\mu_j^2};$$

(iv) $Ric(\mfp_a,\mfp_b)=0$, for every $a,b=1,\ldots,s$ such that
$a\neq b$, and $Ric(\mfn_i,\mfn_j)=0$, for every $i,j=1,\ldots,n$
such that $i\neq j$. \eprop

\bproof (i) By Corollary \ref{ric1cor1}, if $g_F$ is a multiple of
$B$, then we obtain that

$$Ric^F(X,X)=-\frac{1}{2}\left(\frac{1}{2}+c_{\mfl,a}'\right)Kill_{\mfk}(X,X),$$

for $X\in\mfp_a$, where $c_{\mfl,a}'$ is the eigenvalue of the
Casimir operator of $\mfl$ with respect to the Killing form of
$\mfk$ on $\mfp_a$. Clearly,
$$\bar{rl}c_{\mfl,a}'Kill_{\mfk}(X,X)= & Kill_{\mfk}(C'_{\mfl}X,X)\\ \\
= & tr(ad_X^2\mid_{\mfp_a})\\ \\

= & Kill(C_{\mfl}X,X)\\ \\

= & c_{\mfl,a}Kill(X,X).\ear$$

By recalling that

$Kill_{\mfk}\mid_{\mfp_a\times\mfp_a}=\ga_aKill\mid_{\mfp_a\times\mfp_a}$,
we write
$Ric^F=\frac{1}{2}\left(\frac{\ga_a}{2}+c_{\mfl,a}\right)B$ and by
following the notation in Lemma \ref{ricF}, we have

$$q_a=\frac{1}{2}\left(\frac{\ga_a}{2}+c_{\mfl,a}\right).$$

The result then follows from this and Proposition \ref{riccip}.

\li

(ii) follows directly from Proposition \ref{riccin}, by observing
that $\sum_{a=1}^sC_{\mfp_a}=C_{\mfp}$.

(iii) The Ricci curvature in the direction $\mfp\times\mfn$
essentially remains unchanged; the expression given is just that of
Proposition \ref{riccipn} after replacing $\la_1,\ldots,\la_s$ by
$\la$.

(iv) These orthogonality conditions are satisfied by any adapted
metric on $M$ and were shown to hold in Propositions \ref{riccip}
and \ref{riccin}.

$\Box$ \eproof

\bprop \label{riccip2} Let
$g_M(\underbrace{\la,\ldots,\la}_s;\,\mu_1,\ldots,\mu_n)$ be any
adapted metric on $M$ and suppose that $\mfp_1,\ldots,\mfp_s$
pairwise commute, i.e., $[\mfp_a,\mfp_b]=0$ if $a\neq b$. For every
$X\in\mfp_a$,

$$Ric(X,X)=\left(q_a+\frac{\la_a^2}{4}\sum_{j=1}^n\dfrac{c_{\mfn_j,a}}{\mu_j^2}\right)B(X,X),$$

where all the constants are as in Proposition \ref{riccinf}.\eprop

\bproof This follows from the fact that, in this case, Corollary
\ref{ric1cor2} gives also
$$Ric^F(X,X)=-\frac{1}{2}\left(\frac{1}{2}+c_{\mfl,a}'\right)Kill_{\mfk}(X,X),$$
as in the previous proof. Hence, the expression for $Ric(X,X)$ is
exactly the same obtained above for (i) in Proposition
\ref{riccinf}.

$\Box$\eproof

\bcor \label{condnf}If there exists on $M$ an Einstein adapted
metric of the form
$g_M(\underbrace{\la,\ldots,\la}_s;\,\mu_1,\ldots,\mu_n)$, then
$C_{\mfp}$ and $C_{\mfl}$ are scalar on $\mfn_k$,
$k=1,\ldots,n$.\ecor

\bproof Since $\sum_{a=1}^sC_{\mfp_a}=C_{\mfp}$, the necessary
condition for $g_M$ to be Einstein given in Corollary \ref{cond1},
translates into the condition that $C_{\mfp}$ is scalar on
$\mfn_j$, if $\la_1=\ldots=\la_s=\la$. We have that
$C_{\mfk}=C_{\mfp}+C_{\mfl}$ is scalar on $\mfn_j$, since $\mfn_j$
is irreducible as a $K$-module. Then $C_{\mfp}$ is scalar on
$\mfn_j$ if and only if $C_{\mfl}$ is.

$\Box$ \eproof

\section{Riemannian Fibrations with Standard Base}\label{normalbase}

In this section we consider an adapted metric $g_M$ whose projection
onto the base space, $g_N$, is a multiple of the Killing form of
$\mfg$. Hence, we have

\beq\label{mdefnf}g_M=g_M(\la_1,\ldots,\la_s;\,\underbrace{\mu,\ldots,\mu}_n)
\eeq

and \beq g_N=g_N(\underbrace{\mu,\ldots,\mu}_n),\eeq

by setting $\mu_1=\ldots=\mu_n=\mu$ in (\ref{mdef}) and
(\ref{nuplemet}). In this particular case, when equipped with
$g_N$, $N$ is a standard Riemannian manifold.

\bprop\label{riccinb}Let $g_M$ be an adapted metric on $M$ of the
form $$g_M(\la_1,\ldots,\la_s;\,\underbrace{\mu,\ldots,\mu}_n).$$
The Ricci curvature of $g_M$ is as follows:

(i) For every $X\in\mfp_a$,

$$Ric(X,X)=\left(q_a+\frac{\la_a^2}{4\mu^2}(1-\ga_a)\right)B(X,X),$$

where $q_a$ and $\ga_a$ are as defined in Lemma \ref{ricF}, i.e.,
they are defined by the identities

$$Kill_{\mfk}\mid_{\mfp_a\times\mfp_a}=\ga_aKill\mid_{\mfp_a\times\mfp_a}\,and\,Ric^F\mid_{\mfp_a\times\mfp_a}=q_aB\mid_{\mfp_a\times\mfp_a};$$

(ii) For every  $X\in\mfn_k$,
$$Ric(X,X)=-\dfrac{1}{2}\sum_{a=1}^s\dfrac{\la_a}{\mu}B(C_{\mfp_a}X,X)+r_kB(X,X),$$

with $$r_k=\frac{1}{2}\left(\frac{1}{2}+c_{\mfk,k}\right),$$

where $c_{\mfk,k}$ is the eigenvalue of the Casimir operator
$C_{\mfk}$ on $\mfn_k$;

(iii) $Ric(\mfp,\mfn)=0$;

(iv) $Ric(\mfp_a,\mfp_b)=0$, for every $a,b=1,\ldots,s$ such that
$a\neq b$, and $Ric(\mfn_i,\mfn_j)=0$, for every $i,j=1,\ldots,n$
such that $i\neq j$. \eprop

\bproof (i) From the fact that $\ga_a+\sum_{j=1}^n c_{\mfn_j,a}=1$,
we obtain

$$\sum_{j=1}^n\frac{\la_a^2}{\mu^2}c_{\mfn_j,a}=\frac{\la_a^2}{\mu^2}(1-\ga_a).$$

The required expression follows immediately from Proposition
\ref{riccip}.

\li

(ii) From Corollary \ref{ric1cor1} we obtain that
$r_k=\frac{1}{2}\left(\frac{1}{2}+c_{\mfk,k}\right)$, where $r_k$ is
as defined in Lemma \ref{ricN}. The expression then follows from
Proposition \ref{riccin}.

\li

(iii) By using the fact that $C_{\mfn}=\sum_{j=1}^nC_{\mfn_j}$, from
Proposition \ref{riccipn} it follows that

$$Ric(X,Y)=\frac{\la_a}{4\mu}B(C_{\mfn}X,Y),$$

for every $X\in\mfp_a$ and $Y\in\mfn_k$. Moreover, since
$C_{\mfn}=C_{\mfg}-C_{\mfk}=Id-C_{\mfk}$ and
$C_{\mfk}(\mfp)\subset\mfk$, we have that $C_{\mfn}(X)\in\mfk$ is
orthogonal to $Y\in\mfn$ with respect to $B$. Hence, $Ric(X,Y)=0$.

\li

(iv) these orthogonality conditions are simply those in
Propositions \ref{riccip} and \ref{riccin}.

$\Box$\eproof

\bcor \label{riccinb2}Let $g_M$ be any adapted metric on $M$ and
suppose that $\mfn_1,\ldots,\mfn_n$ pairwise commute, i.e.,
$[\mfn_j,\mfn_k]=0$, for $k\neq j$. Then, for every $X\in\mfn_k$,
$$Ric(X,X)=-\dfrac{1}{2}\sum_{a=1}^s\dfrac{\la_a}{\mu_k}B(C_{\mfp_a}X,X)+r_kB(X,X),$$

where $r_k=\frac{1}{2}\left(\frac{1}{2}+c_{\mfk,k}\right)$ and
$c_{\mfk,k}$ is the eigenvalue of the Casimir operator $C_{\mfk}$ on
$\mfn_k$. \ecor

\bproof The proof is immediate by using Corollary \ref{ric1cor2}
and Proposition \ref{riccin}.

$\Box$\eproof

\section{Binormal Riemannian Fibrations}

A $G$-invariant metric $g_M$ on $M$ of the form

\beq\label{mdefbi}g_M=(\underbrace{\la,\ldots,\la}_s;\,\underbrace{\mu,\ldots,\mu}_n)
\eeq
 is called
\textbf{binormal}. That is, a binormal metric is induced by the
scalar product

$$\la B\mid_{\mfp\times\mfp}\oplus\mu B\mid_{\mfn\times\mfn}$$

on $\mfm$. The fibration $F\rightarrow M\rightarrow N$ is then
called a \textbf{binormal Riemannian fibration}. Clearly, a binormal
metric projects onto an invariant metric on the base space $N$ and
thus it is an adapted metric. For a binormal metric $g_M$, both
$g_N$ and $g_F$ are multiples of the Killing form of $\mfg$

\beq g_F=(\underbrace{\la,\ldots,\la}_s)\textrm{ and
}g_N=(\underbrace{\mu,\ldots,\mu}_n) \eeq

and thus $F$ is a normal Riemannian manifold, which is standard if
$Kill_{\mfk}$ is a multiple of $Kill$, and $N$ is a standard
Riemannian manifold.

\li

In this Section we obtain the Ricci curvature of a binormal metric
$g_M$ on $M$ and conditions for such a metric to be Einstein. As
we shall see, the conditions for the existence of an Einstein
binormal metric translate in very simple conditions on the Casimir
operators of $\mfk$, $\mfl$ and $\mfp_a$, $a=1,\ldots,s$. The
results that we found in Sections \ref{normalfiber} and
\ref{normalbase} yield the following description of the Ricci
curvature:

\bcor \label{riccitypeI} Let
$g_M=g_M(\underbrace{\la,\ldots,\la}_s;\,\underbrace{\mu,\ldots,\mu}_n)$
be a binormal metric on $M$.

(i) For every $X\in\mfp_a$,

$$Ric(X,X)=\left(q_a+\frac{\la^2}{4\mu^2}(1-\ga_a)\right)B(X,X),$$

where $q_a=\frac{1}{2}\left(\frac{\ga_a}{2}+c_{\mfl,a}\right)$,
$c_{\mfl,a}$ is the eigenvalue of $C_{\mfl}$ on $\mfp_a$ and $\ga_a$
is determined by

$$Kill_{\mfk}\mid_{\mfp_a\times\mfp_a}=\ga_aKill\mid_{\mfp_a\times\mfp_a};$$

(ii) For every $Y\in\mfn_j$,
$$Ric(Y,Y)=-\frac{\la}{2\mu}B(C_{\mfp}Y,Y)+r_jB(Y,Y),$$

where $r_j=\frac{1}{2}\left(\frac{1}{2}+c_{\mfk,j}\right)$ and
$c_{\mfk,j}$ is the eigenvalue of $C_{\mfk}$ on $\mfn_j$;

(iii) Moreover, $Ric(\mfp,\mfn)=0$;

(iv) $Ric(\mfp_a,\mfp_b)=0$, for every $a,b=1,\ldots,s$ such that
$a\neq b$, and $Ric(\mfn_i,\mfn_j)=0$, for every $i,j=1,\ldots,n$
such that $i\neq j$.

\ecor

\bproof For a binormal metric on $M$, $g_F$ and $g_N$ are multiples
of $Kill$, so we use Propositions \ref{riccinf} and \ref{riccinb}.

$\Box$\eproof

\bdfn \label{deltas1}For $i,j=1,\ldots,n$ and $a,b=1,\ldots,s$, we
set

(i) $\de_{ij}^{\mfk}=c_{\mfk,i}-c_{\mfk,j}$ and
$\de_{ij}^{\mfl}=c_{\mfl,i}-c_{\mfl,j}$;

(ii) $\de_{ab}^{\mfk}=\ga_a-\ga_b$ and
$\de_{ab}^{\mfl}=c_{\mfl,a}-c_{\mfl,b}$.

\edfn

\bthm \label{binormal1}(i) If $C_{\mfp}$ is not scalar on each
$\mfn_j$, then there are no binormal Einstein metrics on $M$;

(ii) Suppose that $C_{\mfp}$ is scalar on each $\mfn_j$ and write
$C_{\mfp}\mid_{\mfn_j}=b^jId_{\mfn_j}$, for $j=1,\ldots,n$. Then
there is a one-to-one correspondence, up to homothety, between
binormal Einstein metrics on $M$ and positive solutions of the
following set of equations on the unknown $X$:

\beq \label{einI1}\de_{ij}^{\mfk}(1-X)=\de_{ij}^{\mfl}, \textrm{ if
}n>1, \eeq

\beq
\label{einI2}(2\de_{ab}^{\mfl}+\de_{ab}^{\mfk})X^2=\de_{ab}^{\mfk},
\textrm{ if $s>1$}, \eeq

\beq\label{einI3}\left(\ga_a+2c_{\mfl,a}\right)X^2-\left(1+2c_{\mfk,j}\right)X+(1-\ga_a+2b^j)=0.\eeq

for every $a,b=1,\ldots,s$ and $i,j=1,\ldots,n$, where $c_{\mfl,a}$
is the eigenvalue of $C_{\mfl}$ on $\mfp_a$, $\ga_a$ is determined
by

$$Kill_{\mfk}\mid_{\mfp_a\times\mfp_a}=\ga_aKill\mid_{\mfp_a\times\mfp_a},$$

$c_{\mfk,j}$ is the eigenvalue of $C_{\mfk}$ on $\mfn_j$ and the
$\de$'s are as in Definition \ref{deltas1}. If such a positive
solution $X$ exists, then binormal Einstein metrics are, up to
homothety, given by

$$<,>=B\mid_{\mfp\times\mfp}\oplus XB\mid_{\mfn\times\mfn}.$$
\ethm

\bproof Let
$g_M(\underbrace{\la,\ldots,\la}_s;\,\underbrace{\mu,\ldots,\mu}_n)$
be a binormal metric on $M$ and $$X=\frac{\mu}{\la}.$$

By Lemma \ref{condnf}, we have that, if $g_M$ is Einstein, then
$C_{\mfp}$ and $C_{\mfl}$ are scalar on $\mfn_j$, for every
$j=1,\ldots,n$. Say

$$C_{\mfp}\mid_{\mfn_j}=b^jId \textrm{ and } C_{\mfl}\mid_{\mfn_j}=c_{\mfl,j}Id.$$

Suppose that $g$ is Einstein with constant $E$. From Corollary
\ref{riccitypeI}, we obtain the Einstein equations

\beq \label{e1}-\dfrac{\la}{2\mu}b^j+r_j=\mu E,\,j=1,\ldots,n\eeq

\beq\label{e2}\dfrac{1}{2}\left(\frac{\ga_a}{2}+c_{\mfl,a}+\frac{\la^2}{2\mu^2}(1-\ga_a)\right)=\la
E,\,a=1,\ldots,s.\eeq

\li

If $n>1$, from Equation (\ref{e1}) we obtain the following:

\beq\label{difij}
\dfrac{\la}{2\mu}(b^i-b^j)=r_i-r_j,\,i,j=1,\ldots,n.\eeq

By using Lemma \ref{ricN} we have
$$r_i-r_j=\dfrac{1}{2}\left(\dfrac{1}{2}+c_{\mfk,i}\right)-\dfrac{1}{2}\left(\dfrac{1}{2}+c_{\mfk,j}\right)=\dfrac{1}{2}(c_{\mfk,i}-c_{\mfk,j}),$$

whereas

$$b^i-b^j=(c_{\mfk,i}-c_{\mfk,j})-(c_{\mfl,i}-c_{\mfl,j}).$$

Therefore, Equation (\ref{difij}) becomes
$$-\frac{\la}{\mu}(\underbrace{c_{\mfl,i}-c_{\mfl,j}}_{\de_{ij}^{\mfl}})=\left(1-\frac{\la}{\mu}\right)(\underbrace{c_{\mfk,i}-c_{\mfk,j}}_{\de_{ij}^{\mfk}}).$$

By using the variable $X$, we rewrite the equation above as
$-\frac{1}{X}\de_{ij}^{\mfl}=\left(1-\frac{1}{X}\right)\de_{ij}^{\mfk}$,
and this yields $\de_{ij}^{\mfl}=(1-X)\de_{ij}^{\mfk}$.

\li

Equation (\ref{e2}) may be rewritten as

\beq\label{e3}\dfrac{1}{2}\left(\frac{\ga_a}{2}+c_{\mfl,a}\right)X+(1-\ga_a)\frac{1}{4X}=\mu
E.\eeq

Hence, if $s>1$, for $a,b=1,\ldots,s$, we get

$$\dfrac{1}{2}\left(\frac{\ga_a}{2}+c_{\mfl,a}\right)X+(1-\ga_a)\frac{1}{4X}=\dfrac{1}{2}\left(\frac{\ga_b}{2}+c_{\mfl,b}\right)X+(1-\ga_b)\frac{1}{4X},$$

which yields

$$\underbrace{c_{\mfl,a}-c_{\mfl,b}}_{\de_{ab}^{\mfl}}=\frac{1}{2}\left(\frac{1}{X^2}-1\right)(\underbrace{\ga_a-\ga_b}_{\de_{ab}^{\mfk}}).$$

By solving this equation we obtain

$$(2\de_{ab}^{\mfl}+\de_{ab}^{\mfk})X^2=\de_{ab}^{\mfk}.$$

\li

Finally, by using Equations (\ref{e1}) and (\ref{e3}) we obtain the
equality

$$\dfrac{1}{2}\left(\frac{\ga_a}{2}+c_{\mfl,a}\right)X+(1-\ga_a)\frac{1}{4X}=-\dfrac{b^j}{2X}+\frac{1}{2}\left(\frac{1}{2}+c_{\mfk,j}\right),$$

which rearranged gives

$$\left(\frac{\ga_a}{2}+c_{\mfl,a}\right)X^2-\left(\frac{1}{2}+c_{\mfk,j}\right)X+\frac{1}{2}(1-\ga_a+2b^j)=0.$$

$\Box$\eproof

An immediate Corollary is the following:

\bcor\label{bincor1} Suppose that $F$ and $N$ are isotropy
irreducible spaces such that $dim\,F>1$. There exists on $M$ an
Einstein adapted metric if and only if $C_{\mfp}$ is scalar on
$\mfn$ and $\triangle \geq 0$, where
$$\triangle=(1+2c_{\mfk,\mfn})^2-4(\ga+2c_{\mfl,\mfp})(1-\ga+2b),$$

$C_{\mfk,\mfn}$ is the eigenvalue of $C_{\mfk}$ on $\mfn$,
$C_{\mfl,\mfp}$ is the eigenvalue of $C_{\mfl}$ on $\mfp$, $b$ is
the eigenvalue of $C_{\mfp}$ on $\mfn$ and $\ga$ is such that
$Kill_{\mfk}\mid_{\mfp\times\mfp}=\ga Kill\mid_{\mfp\times\mfp}$.

If all these conditions are satisfied, then Einstein adapted
metrics are, up to homothety, given by

$$g_M=B\mid_{\mfp\times\mfp}\oplus XB\mid_{\mfn\times\mfn}, \textrm{
where
}X=\frac{1+2c_{\mfk,\mfn}\pm\sqrt{\triangle}}{2(\ga+2c_{\mfl,\mfp})}.$$\ecor

\bproof Since $\mfp$ is an irreducible $Ad\,L$-module and $\mfn$
is an irreducible $Ad\,K$-module, then any adapted metric on $M$
is binormal. Hence, we use Theorem \ref{binormal1}. By the
irreducibility of $\mfp$ and $\mfn$, we have $s=1$ and $n=1$ and
thus Einstein binormal metrics are given by positive solutions of
(\ref{einI3}), if $C_{\mfp}$ is scalar on $\mfn$. Hence, from
Theorem \ref{binormal1} we conclude that there exists on $M$ an
Einstein binormal metric if and only if $C_{\mfp}$ is scalar on
$\mfn$ and $\triangle\geq 0$, where
$$\triangle=(1+2c_{\mfk,\mfn})^2-4(\ga+2c_{\mfl,\mfp})(1-\ga+2b).$$

\li

Since $F$ is isotropy irreducible and $dim\,F>1$, we have
$\ga+2c_{\mfl,p}\neq 0$ and the polynomial in (\ref{einI3}) has
exactly degree two. In fact, if $\ga+2c_{\mfl,p}=0$, then
$\ga=c_{\mfl,p}=0$ and thus, in particular, $\mfp$ lies in the
center of $\mfk$. But the hypothesis that $\mfp$ is irreducible
and abelian implies that $\mfp$ is $1$-dimensional which
contradicts the hypothesis that $dim\,F>1$. Therefore,
$\ga+2c_{\mfl,p}\neq 0$. In this case, the solutions of
(\ref{einI3}) are

$$X=\frac{1+2c_{\mfk,\mfn}\pm\sqrt{\triangle}}{2(\ga+2c_{\mfl,\mfp})}.$$

$\Box$\eproof

In the case when $F$ is 1-dimensional, the fibration $M\rightarrow
N$ is a principal circle bundle, since $F$ is an abelian compact connected
1-dimensional group. We recall that Einstein metrics on principal fiber bundles have been widely studied (\cite{Je2},\cite{WZ3}) and, in particular, homogeneous Einstein metrics on circle bundles were classified McKenzie Y. Wang and Wolfgang Ziller in \cite{WZ3}. We revisit metrics on circle bundles by stating the following:

\bcor\label{bincor2} Suppose that $N$ is isotropy irreducible and
$F$ is isomorphic to the circle group. There exists on $M$ exactly
one $G$-invariant Einstein metric, up to homothety, given by

$$g_M=B\mid_{\mfp\times\mfp}\oplus XB\mid_{\mfn\times\mfn},\,with\,X=\frac{2+m}{m(1+2c_{\mfk,\mfn})},$$

where $c_{\mfk,\mfn}$ is the eigenvalue of $C_{\mfk}$ on $\mfn$
and $m=dim\,G/K$. \ecor

\bproof The fact that $\mfp$ is 1-dimensional implies that $\mfp$
lies in the center of $\mfk$. Hence, in the notation of Corollary
\ref{bincor1}, $\ga=c_{\mfl,\mfp}=0$. On the other hand, if $\mfn$
is $Ad\,K$-irreducible then, the semisimple part of $K$ acts
transitively on $\mfn$. Moreover, since $\mfp$ lies in the center
of $\mfk$, then the semisimple part of $\mfl$ coincides with the
semisimple part of $\mfk$. Hence, $L$ also acts transitively on
$\mfn$ and $\mfn$ is an irreducible $Ad\,L$-module as well.
Consequently, any $G$-invariant metric on $M$ is adapted and
moreover is binormal, by the irreducibility of $\mfp$ and $\mfn$.
Furthermore, $C_{\mfp}$ must be scalar on $\mfn$, since $C_{\mfk}$
and $C_{\mfl}$ are scalar on $\mfn$. Therefore, $G$-invariant
Einstein metrics are given by positive solutions of (\ref{einI3})
in Theorem \ref{binormal1}. Since $\ga=c_{\mfl,\mfp}=0$,
(\ref{einI3}) is just a degree-one equation whose solution is

\beq\label{solcircle1} X=\frac{1+2b}{1+2c_{\mfk,\mfn}},\eeq

where $b$ is the eigenvalue of $C_{\mfp}$ on $\mfn$. Now we
compute $b$, which is the eigenvalue of $C_{\mfp}$ on $\mfn$.
Since $\mfg$ is simple we have $tr(C_{\mfp})=dim\,\mfp=1$. Since
$\mfp$ lies in the center of $\mfk$, $C_{\mfp}$ vanishes on $\mfk$
and thus
$tr(C_{\mfp})=tr(C_{\mfp}\mid_{\mfn})=b\,\textrm{dim}\,\mfn=bm$.
Hence,
$$b=\frac{1}{m}.$$ By replacing $b$ on (\ref{solcircle1}) we obtain
the desired expression for $X$.

$\Box$\eproof

\bex  \textbf{Circle Bundles over Compact Irreducible Hermitian
Symmetric Spaces. } An application of Corollary \ref{bincor2} occurs when the base space is an irreducible symmetric
space. So let us consider a fibration $F\rightarrow M \rightarrow
N$ where $F$ is isomorphic to the circle group and $N$ is an
isotropy irreducible symmetric space. Since $F$ is the circle
group, $\mfp$ lies in the center of $\mfk$. Hence, $K$ has
one-dimensional center, since for a compact irreducible symmetric
space the center of $K$ has at most dimension 1. Moreover, in this
case $N$ is a compact irreducible Hermitian symmetric space. In
particular, $L$ must coincide with the semisimple part of $K$. Compact irreducible Hermitian symmetric spaces $G/K$ are classified (see
e.g. \cite{He}). All the possible $G$, $K$ and $L$ are listed in Table \ref{sf2}, together with the coefficient $X$ of
the, unique, Einstein adapted metric on $G/L$, as in Corollary \ref{bincor2}.

\li

\btab[h!]\caption{Circle Bundles Over irreducible Hermitian
Symmetric Spaces.}\label{sf2}
$$\bar{|lll|c|}\hline\xstrut  G & K & L & X \\ \hline

\xstrut SU(n) & S(U(p)\times U(n-p)) & SU(p)\times SU(n-p) &
\frac{p(n-p)+1}{2p(n-p)}\\

\xstrut SO(2n) &  U(n) &   SU(n) & \frac{n(n-1)+2}{2n(n-1)}\\

\xstrut SO(n) & SO(2)\times SO(n-2) & SO(n-2) & \frac{n-1}{n-2}\\

\xstrut Sp(n) & U(n) & SU(n) & \frac{n(n+1)+2}{2n(n+1)}\\

\xstrut E_6 & SO(10)\times U(1) & SO(10) & \frac{17}{32}\\

\xstrut E_7 & E_6\times U(1) & E_6 & \frac{14}{27}\\ \hline\ear$$
\etab

\eex

Finally, if $F$ is not isotropy irreducible, under some hypothesis
we can show the following:

\bcor \label{binormal3}Suppose $F$ is not isotropy irreducible and
that there exists a constant $\al$ such that

$$Kill_{\mfl}\mid_{\mfp\times\mfp}=\al
Kill_{\mfk}\mid_{\mfp\times\mfp}.$$

For $a=1,\ldots,s$, let $\ga_a$ be the constant determined by
$$Kill_{\mfk}\mid_{\mfp_a\times\mfp_a}=\ga_aKill\mid_{\mfp_a\times\mfp_a}.$$

If for some $a,b=1,\ldots,s$, $\ga_a\neq \ga_b$, then there exists a
binormal Einstein metric on $M$ if and only if, for every
$j=1,\ldots,n$,

\beq\label{caskjlj}c_{\mfl,j}=\left(1-\frac{1}{\sqrt{2\al+1}}\right)\left(c_{\mfk,j}+\frac{1}{2}
\right)\eeq

and $C_{\mfp}$ is scalar on each $\mfn_j$, where $c_{\mfl,j}$ and
$c_{\mfk,j}$ are the eigenvalues of $C_{\mfl}$ and $C_{\mfk}$,
respectively, on $\mfn_j$. In this case, there is a unique binormal
Einstein metric, up to homothety, given by
$$B\mid_{\mfp\times\mfp}\oplus
\frac{1}{\sqrt{2\al+1}}B\mid_{\mfn\times\mfn}.$$

\ecor

\bproof If $Kill_{\mfl}\mid_{\mfp\times\mfp}=\al
Kill_{\mfk}\mid_{\mfp\times\mfp}$, then
\beq\label{eqcas1}c_{\mfl,a}=\al \ga_a, \textrm{ for every
}a=1,\ldots,s.\eeq

Therefore, for any $a,b=1,\ldots,s$, if $s>1$,
$2\de_{ab}^{\mfl}+\de_{ab}^{\mfk}=(2\al+1)\de_{ab}^{\mfk}$ and,
thus, Equation (\ref{einI2}) in Theorem \ref{binormal1} becomes

\beq\label{eqcas2}(2\al+1)\de_{ab}^{\mfk}X^2=\de_{ab}^{\mfk}.\eeq

In particular, (\ref{eqcas1}) implies that $c_{\mfl,a}=0$ if and
only if $\ga_a=0$ (thus if $\mfp$ has submodules where $L$ acts
trivially, then $Kill_{\mfk}$ vanish on those submodules and then
they lie in the center of $\mfk$. If $K$ is semisimple, then the
isotropy representation of $K/L$ is faithful). The fact that the
isotropy representation of $K/L$ is not irreducible implies that
$\mfp$ decomposes as a direct sum $\mfp_1\oplus\ldots\oplus\mfp_s$
with $s>1$. For the indices for which
$\ga_a\neq\ga_b$\footnote{This condition implies that the
representation of $L$ on at least one of the $\mfp_a$'s is
faithful. }, we have $\de_{ab}^{\mfk}\neq 0$ and (\ref{eqcas2})
implies that
$$X=\frac{1}{\sqrt{2\al+1}}.$$

Hence, $X=\frac{1}{\sqrt{2\al+1}}$ must be a root of the polynomial
in (\ref{einI3}). By using the fact that $c_{\mfl,a}=\al \ga_a$ and
$b_j=c_{\mfk,j}-c_{\mfl,j}$, simple calculations show that

\beq\label{casj1}c_{\mfl,j}=\left(1-\frac{1}{\sqrt{2\al+1}}\right)\left(c_{\mfk,j}+\frac{1}{2}
\right).\eeq

We observe that this condition implies (\ref{einI1}) in Theorem
\ref{binormal1}, as we can see by the equalities below:

$$\de_{ij}^{\mfl}=c_{\mfl,i}-c_{\mfl,j}=\left(1-\frac{1}{\sqrt{2\al+1}}\right)(c_{\mfk,i}-c_{\mfk,j})=(1-X)\de_{ij}^{\mfk}.$$

Hence, there is a binormal Einstein metric if and only if
(\ref{casj1}) is satisfied and the operator $C_{\mfp}$ is scalar
on $\mfn_j$, for every $j=1,\ldots,n$. In this case, according
also to Theorem \ref{binormal1} such metric is, up to homothety,
given by $B\mid_{\mfp\times\mfp}\oplus
\frac{1}{\sqrt{2\al+1}}B\mid_{\mfn\times\mfn}$.

$\Box$\eproof

\bcor \label{binormal4}Suppose $F$ is not isotropy irreducible and
that there exists a constant $\al$ such that

$$Kill_{\mfl}\mid_{\mfp\times\mfp}=\al
Kill_{\mfk}\mid_{\mfp\times\mfp}.$$

For $a=1,\ldots,s$, let $\ga_a$ be the constant determined by
$$Kill_{\mfk}\mid_{\mfp_a\times\mfp_a}=\ga_aKill\mid_{\mfp_a\times\mfp_a}.$$

If for some $a,b=1,\ldots,s$, $\ga_a\neq \ga_b$ and there exists on
$M$ an Einstein binormal metric, then the number $\sqrt{2\al+1}$ is
a rational. \ecor

\bproof This follows from the fact that the eigenvalues of
$C_{\mfk}$ and $C_{\mfl}$ on $\mfn_j$ are rational numbers. Since
$\mfk$ is a compact algebra, the eigenvalue of its Casimir operator
on the complex representation on $\mfn_j^{\complex}$ is given by

$$\frac{<\la_j,\la_j+2\de>}{2h^*(\mfg)}\in \rationals,$$

where $\la_j$ is the highest weight for  $\mfn_j^{\complex}$,
$2\de$ is the sum of all positive roots of $\mfk$ and $h^*(\mfg)$
is the dual Coxeter number of $\mfg$ (\cite{Hu}, \cite{Pa}). A
similar formula holds for $C_{\mfl,j}$ and we conclude that
$C_{\mfl,j}$ and $C_{\mfk,j}$ are rational numbers. If there
exists a binormal Einstein metric on $M$, then $C_{\mfl,j}$ and
$C_{\mfk,j}$ are related by formula (\ref{caskjlj}) in Corollary
\ref{binormal3}. This implies that $\sqrt{2\al+1}$ is a rational
number.

$\Box$\eproof

\section{Riemannian Fibrations with Einstein Fiber and Einstein Base}

In this section we investigate conditions for the existence of an
Einstein adapted metric
$g_M=(\la_1,\ldots,\la_s;\,\mu_1,\ldots,\mu_n)$ on $M$ such that
$g_F$ or $g_N$ are also Einstein.

\bthm \label{gnein} Let
$g_M=g_M(\la_1,\ldots,\la_s;\,\mu_1,\ldots,\mu_n)$ be an adapted
metric on $M$. If $g_M$ and $g_N$ are both Einstein and $n>1$,
then

$$\frac{\mu_j}{\mu_k}=\frac{r_j}{r_k}=\left(\frac{b^j}{b^k}\right)^{\frac{1}{2}}, \textrm{ for
}j,k=1,\ldots,n,$$

where $b^j$ is the eigenvalue of the operator
$\sum_{a=1}^s\la_aC_{\mfp_a}$ on $\mfn_j$, for $j=1,\ldots,n$, and
the $r_j$'s are determined by
$Ric^N=\oplus_{k=1}^nr_kB\mid_{\mfn_k\times\mfn_k}$ as in Lemma
\ref{ricN}. Up to homothety, there exists at most one Einstein
metric $g_N$ on $N$ such that the corresponding $g_M$ on $M$ is
Einstein.\ethm

\bproof Let $g_M(\la_1,\ldots,\la_s;\,\mu_1,\ldots,\mu_n)$ be an
adapted metric on $M$. From Corollary \ref{cond1} we know that if
$g_M$ is Einstein, then there are constants $b^j$ such that

$$\sum_{a=1}^s\la_aC_{\mfp_a}\mid_{\mfn_j}=b^jId_{\mfn_j}.$$

We recall from Lemma \ref{ricN} that
$Ric^N=\oplus_{k=1}^nr_kB\mid_{\mfn_k\times\mfn_k}$. Hence, if $g_N$
is Einstein, then

\beq\label{rmus1}\frac{r_1}{\mu_1}=\ldots=\frac{r_n}{\mu_n}.\eeq

From this equalities we obtain that

\beq\label{rmus2}\frac{\mu_j}{\mu_k}=\frac{r_j}{r_k}, \textrm{ for
}j,k=1,\ldots,n.\eeq

From Proposition \ref{riccin}, for $X\in\mfn_k$, the Ricci curvature
of $g_M$ is

$$Ric(X,X)=-\dfrac{1}{2\mu_k}\sum_{a=1}^s\la_aB(C_{\mfp_a}X,X)+r_kB(X,X)=\big(-\dfrac{b^k}{2\mu_k}+r_k\big)B(X,X).$$

If $g_M$ is Einstein, then from the expression above we obtain the
following Equations

\beq\label{rmus3}
-\dfrac{b^k}{2\mu_k^2}+\frac{r_k}{\mu_k}=-\dfrac{b^j}{2\mu_j^2}+\frac{r_j}{\mu_j}.\eeq

The identities (\ref{rmus1}) and (\ref{rmus3}) imply that
$$\dfrac{b^k}{\mu_k^2}=\dfrac{b^j}{\mu_j^2}$$

and consequently, by using (\ref{rmus2}),

$$\left(\frac{r_j}{r_k}\right)^2=\left(\frac{\mu_j}{\mu_k}\right)^2=\frac{b^j}{b^k}.$$

Finally, we observe that although the fact that $g_N$ is Einstein
implies the equalities $\dfrac{\mu_j}{\mu_k}=\dfrac{r_j}{r_k}$,
there might be more then one solution for the $n$-tuples
$(\mu_1,\ldots,\mu_n)$, up to scalar multiplication, since the
$r_j$'s in general depend on the $\mu_i$'s. This is obvious since
clearly there might be many distinct Einstein metrics on $N$, up to
homothety. This is explicit in the formula given in Lemma
\ref{ricN}. However, as the eigenvalues $b^j$ are independent of the
constants $\mu_1,\ldots,\mu_n$, the ratios
$\frac{\mu_j}{\mu_k}=\left(\frac{b^j}{b^k}\right)^{\frac{1}{2}}$
imply that there is at most one possible choice for $g_N$, up to
scalar multiplication.

$\Box$\eproof

\bthm \label{gfein} Let
$g_M(\la_1,\ldots,\la_s;\,\mu_1,\ldots,\mu_n)$ be an adapted
metric on $M$. If $g_M$ and $g_F$ are both Einstein and $s>1$,
then

$$\frac{\la_a}{\la_b}=\frac{q_a}{q_b}= \sum_{j=1}^n\frac{C_{\mfn_j,b}}{\mu_j^2}\Big/\sum_{j=1}^n\frac{C_{\mfn_j,a}}{\mu_j^2} ,\textrm{ for
}a,b=1,\ldots,s,$$

where $c_{\mfn_j,a}$ is such that
$Kill(C_{\mfn_j}\cdot,\cdot)\mid_{\mfp_a\times\mfp_a}=c_{\mfn_j,a}Kill\mid_{\mfp_a\times\mfp_a}$,
for $a=1,\ldots,s$ and the $q_a$'s are determined by
$Ric^F=\oplus_{a=1}^sq_aB\mid_{\mfp_a\times\mfp_a}$ as in Lemma
\ref{ricF}. Up to scalar multiplication, there exists at most one
Einstein metric $g_F$ on $F$ such that the corresponding metric
$g_M$ on $M$ is Einstein.\ethm

\bproof The proof is similar to that of Theorem \ref{gnein}, by
using Lemma \ref{ricF} and Proposition \ref{riccip}.

$\Box$\eproof

\bcor  Let $g_M=g_M(\la_1,\ldots,\la_s;\,\mu_1,\ldots,\mu_n)$ be
an adapted metric on $M$. If $g_M$, $g_N$ and $g_F$ are Einstein,
then

$$\frac{r_j}{r_k}=\left(\frac{b^j}{b^k}\right)^{\frac{1}{2}}, \textrm{ for
}j,k=1,\ldots,n,$$

and

$$\frac{q_a}{q_b}= \sum_{j=1}^n\frac{C_{\mfn_j,b}}{b_j}\Big/\sum_{j=1}^n\frac{C_{\mfn_j,a}}{b_j} ,\textrm{ for
}a,b=1,\ldots,s,$$

where all the constants are as in Theorems \ref{gnein} and
\ref{gfein}.\ecor

\bproof Using Theorem \ref{gnein}, we write
$\mu_j^2=\frac{b^j}{b^1}\mu_1^2$. The second formula follows
immediately from this and Theorem \ref{gfein}.

$\Box$\eproof

\bthm Let $g_M=g_M(\la_1,\ldots,\la_s;\,\mu_1,\ldots,\mu_n)$ be an
adapted metric on $M$. Suppose that $g_M$, $g_N$ and $g_F$ are
Einstein and let $E$, $E_F$ and $E_N$ be the corresponding
Einstein constants. If $E\neq E_N$, then

$$\mu_j=\left(\frac{b^j}{2(E_N-E)}\right)^{\frac{1}{2}},$$

$$\la_a=2\frac{E-E_F}{E_N-E}\left(\frac{C_{\mfn_j,a}}{b^j}\right)^{-1}.$$

where $b^j$ is the eigenvalue of the operator
$\sum_{a=1}^s\la_aC_{\mfp_a}$ on $\mfn_j$ and $c_{\mfn_j,a}$ is
such that
$Kill(C_{\mfn_j}\cdot,\cdot)\mid_{\mfp_a\times\mfp_a}=c_{\mfn_j,a}Kill\mid_{\mfp_a\times\mfp_a}$.\ethm

\bproof Let $g_M=(\la_1,\ldots,\la_s;\,\mu_1,\ldots,\mu_n)$ be an
adapted metric on $M$. If $g_M$, $g_N$ and $g_F$ are all Einstein,
from Propositions \ref{riccip} and \ref{riccin} we get

$$-\frac{1}{2\mu_j}b^j+\mu_j E_N=\mu_j E ,$$

from which, if $E_N\neq E$, we deduce

\beq\label{eqmuj1}\mu_j^2=\frac{b^j}{2(E_N-E)}\eeq

and

$$\la_a E_F+\frac{\la_a^2}{4}\frac{C_{\mfn_j,a}}{\mu_j^2}=\la_a E.$$

From this we get

\beq \la_a\frac{C_{\mfn_j,a}}{\mu_j^2}=4(E-E_F).\eeq

We obtain the required formula by replacing (\ref{eqmuj1}) in the
equation above.

$\Box$\eproof

\section{Riemannian Fibrations with Symmetric
Fiber}\label{sectionsymF}

In this section we consider a fibration $F\rightarrow M\rightarrow
N$ such that $F$ is a symmetric space and $N$ is isotropy
irreducible. We specify the Ricci curvature of an adapted metric
on $M$ and obtain the Einstein equations in some particular cases.

\li

If $F=K/L$ is a symmetric space, then we consider its DeRham
decomposition

\beq\label{deRham} K/L=K_0/L_0\times K_1/L_1\times\ldots\times
K_s/L_s,\eeq

where $K_0$ is the center of $K$ and, for $a=1,\ldots,s$, $K_a$ is
simple. By $\mfk_a$ and $\mfl_a$ we denote the Lie algebras of
$K_a$ and $L_a$, respectively. In particular, for $a=1,\ldots,s$,
$K_a/L_a$ is an irreducible symmetric space. Thus $\mfp_a$ may be
chosen as a symmetric reductive complement of $\mfl_a$ in
$\mfk_a$. Since $\mfk_a$ is simple, the Casimir operator of $\mfk$
is scalar on $\mfk_a$. Hence, in the equality
$Kill_{\mfk}\mid_{\mfp_a\times\mfp_a}=\ga_a
Kill\mid_{\mfp_a\times\mfp_a}$, the constant $\ga_a$ is simply the
eigenvalue of the Casimir operator of $\mfk$ on $\mfk_a$, because
$Kill_{\mfk}=Kill(C_{\mfk}\cdot,\cdot)$. For $a=0$ this is still
true with $\ga_0=0$.

\bprop \label{ricciFsym}Suppose that $F$ is a symmetric space and
$N$ is isotropy irreducible and let
$g_M=(\la_0,\ldots,\la_s;\,\mu)$ be an adapted metric on $M$.  The
Ricci curvature of $g_M$ is as follows:

(i) For every $X\in\mfp_a$, $a=0,\ldots,s$,

$$Ric(X,X)=\left(\frac{\ga_a}{2}+\dfrac{\la_a^2}{4\mu^2}(1-\ga_a)\right)B(X,X),$$

where $\ga_a$ is the eigenvalue of $C_{\mfk}$ on $\mfk_a$;

(ii) For every $X\in\mfn$,

$$Ric(X,X)=-\frac{1}{2}\sum_{a=1}^s\la_aB(C_{\mfp_a}X,X)+rB(X,X),$$

where $r=\dfrac{1}{2}\big(\dfrac{1}{2}+c_{\mfk,\mfn}\big)$ and
$c_{\mfk,\mfn}$ is the eigenvalue of $C_{\mfk}$ on $\mfn$;

(iii) $Ric(\mfp,\mfn)=0$;

(iv) For every $a,b=0,\ldots,s$ such that $a\neq b$,
$Ric(\mfp_a,\mfp_b)=0$ and for every $i,j=1,\ldots,n$ such that
$i\neq j$, $Ric(\mfn_i,\mfn_j)=0$.

\eprop

\bproof Since $N$ is isotropy irreducible the expressions for the
Ricci curvature of $g_M$ are given by Proposition \ref{riccinb}. In
particular, for $X\in\mfp_a$,

$$Ric(X,X)=\left(q_a+\frac{\la_a^2}{4\mu^2}(1-\ga_a)\right)B(X,X).$$

If we consider the DeRham decomposition of $F$ as in (\ref{deRham}),
$[\mfp_a,\mfp_b]=0$, for every $a\neq b$. Hence, from Proposition
\ref{riccip2}, we have
$q_a=\frac{1}{2}\left(\frac{\ga_a}{2}+c_{\mfl,a}\right)$. Since $F$
is a symmetric space, then $Kill_{\mfl}\mid_{\mfp\times\mfp}=
\frac{1}{2}Kill_{\mfk}\mid_{\mfp\times\mfp}$ and thus
$C_{\mfl}\mid_{\mfp\times\mfp}=\frac{1}{2}C_{\mfk}\mid_{\mfp\times\mfp}$.
Hence, $c_{\mfl,a}=\frac{\ga_a}{2}$. Therefore
$q_a=\dfrac{\ga_a}{2}$.

\li

Since $N$ is irreducible, by using Proposition \ref{riccinb}, for
$X\in\mfn$, we write

$$Ric(X,X)=-\dfrac{1}{2}\sum_{a=1}^s\dfrac{\la_a}{\mu}B(C_{\mfp_a}X,X)+rB(X,X),$$

where $$r=\frac{1}{2}\left(\frac{1}{2}+c_{\mfk,\mfn}\right)$$ and
$c_{\mfk,\mfn}$ is the eigenvalue of the Casimir operator $C_{\mfk}$
on $\mfn$.

(iii) and (iv) follow directly from Proposition \ref{riccinb} as
well.

$\Box$ \eproof

\bthm \label{eqFsym}Suppose that $F$ is a symmetric space and $N$
is an isotropy irreducible space. Moreover, suppose that
$C_{\mfp_a}\mid_{\mfn}=b_aId_{\mfn}$, for some constants $b_a$,
for every $a=1,\ldots,s$. There exists on $M$ an Einstein adapted
metric if and only if there are positive solutions of the
following system of $s$ algebraic equations in the unknowns
$X_1,\ldots,X_s$:\footnote{$\widehat{X_a}$ means that $X_a$ does
not occur in the product. }

\beqar 2\ga_1X_1^2X_a+(1-\ga_1)X_a-2\ga_aX_1X_a^2-(1-\ga_a)X_1=0,
\,a=2,\ldots,s\nonumber\\
2\sum_{a=1}^sb_aX_1\ldots\widehat{X_a} \ldots X_s-4rX_1\ldots
X_s+2\ga_1X_1^2X_2\ldots X_s+(1-\ga_1)X_2\ldots X_s=0,\nonumber
\eeqar

where $\ga_a$ is the eigenvalue of $C_{\mfk}$ on $\mfp_a$,
$r=\frac{1}{2}\left(\frac{1}{2}+c_{\mfk,\mfn}\right)$ and
$c_{\mfk,\mfn}$ is the eigenvalue of $C_{\mfk}$ on $\mfn$. To each
$s$-tuple $(X_1,\ldots,X_s)$ corresponds a family of Einstein
adapted metrics on $M$ given, up to homothety, by

$$g_M=\oplus_{a=1}^s\frac{1}{X_a}B\mid_{\mfp_a\times\mfp_a}\oplus
B\mid_{\mfn\times\mfn}.$$ \ethm

\bproof Let $g_M=(\la_1,\ldots,\la_s;\,\mu)$ be an adapted metric
on $M$. First we observe that the hypothesis
$C_{\mfp_a}\mid_{\mfn}=b_aId_{\mfn}$, for every $a=1,\ldots,s$,
implies that $\sum_{a=1}^s\la_aC_{\mfp_a}$ is scalar for any
choice of $\la_a$'s. Hence, the necessary condition for the
existence of an Einstein adapted metric on $M$ given by Corollary
\ref{cond1} is satisfied. Moreover, (iii) and (iv) of Proposition
\ref{ricciFsym} imply that for $g_M$ to be Einstein, it suffices
to analyze the equations

\beqar
\label{eqFsym1}Ric\mid_{\mfp_a\times\mfp_a}=\la_aEB\mid_{\mfp_a\times\mfp_a},\,a=1,\ldots,s\\
\label{eqFsym2}Ric\mid_{\mfn\times\mfn}=\mu
EB\mid_{\mfn\times\mfn},\eeqar

where $E$ is the Einstein constant of $g_M$.

We introduce the unknowns

$$X_a=\frac{\mu}{\la_a},\,a=1,\ldots,s.$$

By using $C_{\mfp_a}\mid_{\mfn}=b_aId_{\mfn}$ and the $X_a$'s,
Equation \ref{eqFsym2} may be rewritten as

\beq\label{eqFsym3}-\sum_{a=1}^s\frac{b_a}{2X_a}+r=\mu E.\eeq

Also, by using Proposition \ref{ricciFsym} and the $X_a$'s, Equation
\ref{eqFsym1} may be rewritten as

\beq\label{eqFsym4}\frac{\ga_a}{2}+\dfrac{1-\ga_a}{4X_a^2}=\la_a
E.\eeq

By multiplying (\ref{eqFsym4}) by $X_a$ we get

\beq\label{eqFsym5}\frac{2\ga_aX_a^2+1-\ga_a}{4X_a}=\mu E.\eeq

Therefore, the Einstein Equations are just

\beqar\label{eqFsym6}\frac{2\ga_aX_a^2+1-\ga_a}{4X_a}=\frac{2\ga_1X_1^2+1-\ga_1}{4X_1},\,a=1,\ldots,s\\
\label{eqFsym7}-\sum_{a=1}^s\frac{b_a}{2X_a}+r=\frac{2\ga_1X_1^2+1-\ga_1}{4X_1}.\eeqar

We obtain the equations stated in the theorem simply by
rearranging (\ref{eqFsym6}) and (\ref{eqFsym7}). We recall that
since $N$ is irreducible, we have
$r=\frac{1}{2}\left(\frac{1}{2}+c_{\mfk,\mfn}\right)$, as in
Proposition \ref{ricciFsym}, and thus $r$ does not depend on
$\mu$. So $X_1,\ldots,X_s$ are actually the only unknowns of the
system above.

$\Box$\eproof

\bcor\label{binSym1} Suppose that $F$ and $N$ are irreducible
symmetric spaces and $dim\,F>1$. There exists on $M$ an Einstein
adapted metric if and only if $C_{\mfp}$ is scalar on $\mfn$ and
$\triangle' \geq 0$, where
$$\triangle'=1-2\ga(1-\ga+2b),$$

$\ga$ is the eigenvalue of $C_{\mfk}$ on $\mfp$ and $b$ is the
eigenvalue of $C_{\mfp}$ on $\mfn$. If these two conditions are
satisfied, then Einstein adapted metrics are homothetic to
$g_M=B\mid_{\mfp\times\mfp}\oplus XB\mid_{\mfn\times\mfn}$, where

$$X=\frac{1\pm\sqrt{\triangle'}}{2\ga}.$$\ecor

\bproof It follows from Corollary \ref{bincor1} and from the fact
that, since $F$ and $N$ are irreducible symmetric spaces,
$c_{\mfk,\mfn}=\dfrac{1}{2}$ and $c_{\mfl,\mfp}=\dfrac{\ga}{2}$.

$\Box$\eproof

\bcor \label{binormal5}Suppose that $F$ is a symmetric space and
$N$ is isotropy irreducible.

(i) If $C_{\mfp}$ is not scalar on $\mfn$ or $C_{\mfk}$ is not
scalar on $\mfp$, then there is no binormal Einstein metric on
$M$.

(ii) Suppose that $C_{\mfp}$ is scalar on $\mfn$ and $C_{\mfk}$ is
scalar on $\mfp$, and write $C_{\mfp}\mid_{\mfn}=bId_{\mfn}$ and
$C_{\mfk}\mid_{\mfp}=\ga Id_{\mfp}$. There is an one-to-one
correspondence between binormal Einstein metrics on $M$ and
positive roots of the polynomial

\beq\label{polbin5}2\ga
X^2-\left(1+2c_{\mfk,\mfn}\right)X+(1-\ga+2b)=0.\eeq

for every $a,b=1,\ldots,s$ and $i,j=1,\ldots,n$, where
$c_{\mfk,\mfn}$ is the eigenvalue of $C_{\mfk}$ on $\mfn$. If such
a positive solution $X$ exists, then binormal Einstein metrics
are, up to homothety, given by

$$<,>=B\mid_{\mfp\times\mfp}\oplus XB\mid_{\mfn\times\mfn}.$$
\ecor

\bproof If $F$ is a symmetric space, then
$Kill_{\mfl}\mid_{\mfp\times\mfp}=\al
Kill_{\mfk}\mid_{\mfp\times\mfp}$, with $\al=\frac{1}{2}$. The
number $\sqrt{2\al+1}=\sqrt{2}$ is not a rational. Hence,
Corollary \ref{binormal4} implies that if there exists a binormal
Einstein metric on $M$, then $\ga_1=\ldots=\ga_s=\ga$ for some
constant $\ga$, i.e., the Casimir operator of $\mfk$ is scalar on
$\mfp$.

Moreover, we know from Theorem \ref{binormal1} that the condition
that $C_{\mfp}$ is scalar on $\mfn$ is also a necessary condition
for the existence of a binormal Einstein metric.

The polynomial (\ref{polbin5}) is just (\ref{einI3}) from Theorem
\ref{binormal1}, for $c_{\mfl,\mfp}=\frac{\ga}{2}$ and $n=1$.
Also, the condition (\ref{einI2}) from Theorem \ref{binormal1} is
satisfied since for $\ga_1=\ldots=\ga_s$, we have
$\de_{ab}^{\mfk}=\de_{ab}^{\mfl}=0$.

$\Box$ \eproof

\bcor \label{binormal6}Suppose that $F$ is a symmetric space. If
there exists on $M$ a binormal Einstein metric $g_M$, then $g_F$
is Einstein. The converse holds if $C_{\mfk}$ is scalar on
$\mfp$.\ecor

\bproof If $F$ is irreducible, then any metric on $F$ is Einstein.
So let us suppose that $F$ is a reducible symmetric space. Then, by
Corollary \ref{binormal5}, the existence of a binormal Einstein
metric $g_M$ on $M$ implies that $\ga_1=\ldots=\ga_s=\ga$ for some
$\ga$. From Proposition \ref{riccinf}, we have
$Ric^F\mid_{\mfp_a\times\mfp_a}=q_aB(X,X)$, where
$q_a=\frac{1}{2}\left(\frac{\ga_a}{2}+c_{\mfl,a}\right)=\frac{\ga_a}{2}=\frac{\ga}{2}$,
for every $a=1,\ldots,s$. Therefore, $g_F$ is Einstein with Einstein
constant $E_F=\frac{\ga}{2\la}.$

Conversely, let $g_M=(\la_1,\ldots,\la_s,\mu)$ be any Einstein
adapted metric on $M$. If $C_{\mfk}$ is scalar on $\mfp$, then
$\ga_1=\ldots=\ga_s$. Hence, if $g_F$ is Einstein, we have from
Theorem \ref{gfein} that

$$\frac{\la_a}{\la_b}=\frac{q_a}{q_b}=\frac{\ga_a/2}{\ga_b/2}=1$$

and $g_M$ is binormal.

$\Box$\eproof

Hence, binormal Einstein metrics are such that the restriction to
the fiber is Einstein. As the next two results show, there might
exists other Einstein adapted metrics satisfying this property.

\bcor\label{gfeincor}Suppose that $F$ is a symmetric space and $N$
is isotropy irreducible. Let $\ga_a$ be the eigenvalue of $C_{\mfk}$
on $\mfp_a$, $a=1,\ldots,s$, and $g_M$ an adapted metric on $M$. If
$g_M$ and $g_F$ are both Einstein, then

$$\ga_a=\ga_b\textrm{ or }\ga_a=1-\ga_b,\,a,b=1,\ldots,s.$$
\ecor

\bproof If $F$ is a symmetric space, we have $q_a=\frac{\ga_a}{2}$,
for every $a=1,\ldots,s$. On the other hand, if $N$ is isotropy
irreducible, then $\mfn$ is an irreducible $Ad\,K$-module, and since
$C_{\mfn,a}=1-\ga_a$, the identity in Theorem \ref{gfein} becomes

$$\frac{\ga_a/2}{\ga_b/2}=\frac{1-\ga_b}{\mu^2}\Big/\frac{1-\ga_a}{\mu^2}.$$

Hence, we obtain the equation

$$\ga_a(1-\ga_a)=\ga_b(1-\ga_b),\textrm{ for
}a,b=1,\ldots,s,$$

whose solutions are $\ga_a=\ga_b$ or $\ga_a=1-\ga_b$.

$\Box$\eproof

\bcor \label{gFeincor2}Suppose that $F$ is a symmetric space such
that $\mfp=\mfp_1\oplus\mfp_2$, where $\mfp_1$ and $\mfp_2$ are
non-abelian, and $N$ is isotropy irreducible. Suppose that
$C_{\mfp_a}\mid_{\mfn}=b_aId_{\mfn}$, for some constants $b_a$,
$a=1,2$.  Let $\ga_a$ be the eigenvalue of $C_{\mfk}$ on $\mfp_a$,
$a=1,2$, and $c_{\mfk,\mfn}$ be the eigenvalue of $C_{\mfk}$ on
$\mfn$.

If there exists on $M$ an Einstein adapted metric $g_M$ such that
$g_F$ is also Einstein, then one of the following cases holds:

(i) $\ga_2=\ga_1$ and $\triangle\geq 0$, where

$$\triangle=(1+c_{\mfk,\mfn})^2-8\ga_1(1-\ga_1+2b).$$

If these two conditions are satisfied the metric, then $g_M$ is
the binormal metric given, up to homothety, by

$$g_M=B\mid_{\mfp\times\mfp}\oplus
XB\mid_{\mfn\times\mfn},\,
where\,X=\frac{1+c_{\mfk,\mfn}\pm\sqrt{\triangle}}{2\ga_1}.$$

(ii) $\ga_2=1-\ga_1$ and $D(\ga_1)\geq 0$, where

$$D(\ga_1)=4r^2-4b_1\ga_1-4b_2(1-\ga_1)-2\ga_1(1-\ga_1)$$

and $r=\frac{1}{2}\left(\frac{1}{2}+c_{\mfk,\mfn}\right)$. If
these two conditions are satisfied the metric $g_M$ is given, up
to homothety, by

$$g_M=\frac{1}{X_1}B\mid_{\mfp_1\times\mfp_1}\oplus\frac{1}{X_2}B\mid_{\mfp_2\times\mfp_2}\oplus
B\mid_{\mfn\times\mfn},$$

where $$X_2=\frac{\ga_1 X_1}{1-\ga_1}\textrm{  and
}X_1=\frac{2r\pm\sqrt{D(\ga_1)}}{2\ga_1}.$$

\ecor

\bproof First we observe that the hypothesis that $\mfp_1$ and
$\mfp_2$ are non-abelian implies that $\ga_1,\,\ga_2\neq 0$. Let
$g_M=g_M(\la_1,\ldots,\la_s,\mu)$ be an Einstein adapted metric on
$M$ such that $g_F$ is also Einstein. Corollary \ref{gfeincor},
implies that either $\ga_2=\ga_1$ or $\ga_2=1-\ga_1$.

In the case $\ga_2=\ga_1$, the statement follows from Corollaries
\ref{binormal5} and \ref{binormal6}. In the case $\ga_2=1-\ga_1$,
we obtain from Theorem \ref{gfein}, that

\beq\label{rellambdas}\frac{\la_1}{\la_2}=\frac{\ga_1/2}{\ga_2/2}=\frac{\ga_1}{1-\ga_1}.\eeq

On the other hand, according to Theorem \ref{eqFsym}, an adapted
Einstein metric on $M$ corresponds to positive solutions of the
equations

\beqar
\label{eqII11}2\ga_1X_1^2X_2+(1-\ga_1)X_2-2\ga_2X_1X_2^2-(1-\ga_2)X_1=0
\\
\label{eqII12}2b_1X_2+2b_2X_1-4rX_1X_2+2\ga_1X_1^2X_2+(1-\ga_1)X_2=0,
\eeqar

where $X_a=\frac{\mu}{\la_a}$. By using the identity
(\ref{rellambdas}), we solve the system of equations
(\ref{eqII11}) and (\ref{eqII12}) for $X_2=
\frac{\ga_1}{1-\ga_1}X_1$ and $\ga_2=1-\ga_1$, in order to obtain
the solutions stated in (ii).

$\Box$ \eproof

\li

The following two Corollaries classify all the Einstein adapted
metrics in the cases when $\ga_2=\ga_1$ or $\ga_2=1-\ga_1$. These
results follow immediately from Corollary \ref{gFeincor2} and from
solving the equations (\ref{eqII11}) and (\ref{eqII12}) for
$\ga_2=\ga_1$ or $\ga_2=1-\ga_1$, respectively.

\bcor \label{contcor1}Suppose that $F$ is a symmetric space such
that $\mfp=\mfp_1\oplus\mfp_2$, where $\mfp_1$ and $\mfp_2$ are
non-abelian, and $N$ is isotropy irreducible. Suppose that
$C_{\mfp_a}\mid_{\mfn}=b_aId_{\mfn}$, for some constants $b_a$,
$a=1,2$, and let $\ga_a$ be the eigenvalue of $C_{\mfk}$ on
$\mfp_a$, $a=1,2$, and $c_{\mfk,\mfn}$ the eigenvalue of
$C_{\mfk}$ on $\mfn$.

Suppose that $\ga_2=\ga_1$, i.e., $C_{\mfk}$ is scalar on $\mfp$.
If there exists on $M$ an Einstein adapted metric $g_M$, then one
of the following two cases holds:

(i) $g_F$ is also Einstein and $g_M$ is a binormal metric given by
Corollary \ref{gFeincor2} (i).

(ii) $D(\ga_1)\geq 0$, where

$$D(\ga_1)=4r^2(1-\ga_1)-2\ga_1(2b_2+1-\ga_1)(2b_1+1-\ga_1)$$

and $r=\frac{1}{2}\left(\frac{1}{2}+c_{\mfk,\mfn}\right)$. The
metric $g_M$ is given, up to homothety, by

$$g_M=\frac{1}{X_1}B\mid_{\mfp_1\times\mfp_1}\oplus\frac{1}{X_2}B\mid_{\mfp_2\times\mfp_2}\oplus
B\mid_{\mfn\times\mfn},$$

where $$X_2=\frac{1-\ga_1}{2\ga_1 X_1}\textrm{  and
}X_1=\frac{2r(1-\ga_1)\pm\sqrt{(1-\ga_1)D(\ga_1)}}{2\ga_1(2b_2+1-\ga_1)}.$$
In this second case, $g_F$ is not Einstein and $g_M$ is not
binormal.\ecor

\bcor\label{contcor2}Suppose that $F$ is a symmetric space such
that $\mfp=\mfp_1\oplus\mfp_2$, where $\mfp_1$ and $\mfp_2$ are
non-abelian, and $N$ is isotropy irreducible. Suppose that
$C_{\mfp_a}\mid_{\mfn}=b_aId_{\mfn}$, for some constants $b_a$,
$a=1,2$, and let $\ga_a$ be the eigenvalue of $C_{\mfk}$ on
$\mfp_a$, $a=1,2$, and $c_{\mfk,\mfn}$ the eigenvalue of
$C_{\mfk}$ on $\mfn$.

Suppose that $\ga_2=1-\ga_1$. If there exists on $M$ an Einstein
adapted metric $g_M$, then one of the following two cases holds:

(i) $g_F$ is also Einstein and $g_M$ is the metric given by
Corollary \ref{gFeincor2} (ii).

(ii) $D(\ga_1)\geq 0$, where

$$D(\ga_1)=4r^2-2(2b_2+\ga_1)(2b_1+1-\ga_1)$$

and $r=\frac{1}{2}\left(\frac{1}{2}+c_{\mfk,\mfn}\right)$. The
metric $g_M$ is given, up to homothety, by

$$g_M=\frac{1}{X_1}B\mid_{\mfp_1\times\mfp_1}\oplus\frac{1}{X_2}B\mid_{\mfp_2\times\mfp_2}\oplus
B\mid_{\mfn\times\mfn},$$

where
$$X_2=\frac{1}{2X_1}\textrm{  and
}X_1=\frac{2r\pm\sqrt{D(\ga_1)}}{2(2b_2+\ga_1)}.$$

$g_M$ is never binormal and in the second case $g_F$ is not
Einstein.\ecor

\newpage

\chapter{}\label{chpaterbs}

\addtolength{\myVSpace}{0.1cm}

As in the previous chapters, we consider a homogeneous fibration
$F\rightarrow M \rightarrow N$, for $M=G/L$, $N=G/K$ and $F=K/L$,
where $G$ is a compact connected semisimple Lie group and
$L\varsubsetneq K\varsubsetneq G$ connected closed non-trivial
subgroups. In this chapter we suppose that both the fiber $F$ and
the base space $N$ are symmetric spaces of maximal rank and,
moreover, $N$ is isotropy irreducible.  The triple formed by the Lie
algebras of $G$, $K$, $L$, denoted by $(\mfg,\mfk,\mfl)$, shall be
called a bisymmetric triple of maximal rank. We classify all the
bisymmetric triple of maximal rank when $\mfg$ is simple and obtain
formulas to compute the eigenvalues which are necessary to decide
about the existence of Einstein adapted metrics. For each triple, we
present the eigenvalues of the Casimir operators of the irreducible
$L$-invariant subspaces of the fiber on the horizontal direction and the
eigenvalues of the Casimir operator of $\mfk$ on the vertical
direction. The computations for these eigenvalues are in Appendix
\ref{cpproofs} as well as a description of the isotropy
representation in terms of subset of roots for each triple. Finally,
we study the existence of adapted Einstein metrics by using the
results in previous chapters. Tables are presented in the end of this chapter. We use the notation
used in previous chapters unless stated otherwise.

\section{Introduction}

For all the definitions and properties concerning the roots system
of a Lie algebra please see \cite{He} or \cite{OV}. Let $G$ be a
compact connected semisimple Lie group and $L\varsubsetneq K
\varsubsetneq G$ connected closed non-trivial subgroups such that
$N=G/K$ is isotropy irreducible. As in section \ref{sectionRF} of
Chapter 1, $\mfn$ and $\mfp$ denote the reductive complements of
$\mfk$ in $\mfg$ and of $\mfl$ in $\mfk$, respectively. The
subspace $\mfn$ is irreducible as an $Ad\,K$-module and $\mfp$ may
decompose into the direct sum
$\mfp=\mfp_1\oplus\ldots\oplus\mfp_s$ of irreducible
$Ad\,L$-modules. We suppose that $M$ has simple spectrum, i.e.,
$\mfn$ do not contain any $Ad\,L$-submodule equivalent to any of
the  $\mfp_1,\ldots,\mfp_s$ and $\mfp_1,\ldots,\mfp_s$ are
pairwise inequivalent.

\li

Initially, we only suppose that $L$ is a subgroup of maximal rank
in $G$. We choose a Cartan subalgebra $\mfh$ of $\mfg^{\complex}$
such that $\mfh\subset\mfl^{\complex}$. Let $\R$ be a system of
nonzero roots for $\mfg^{\complex}$ with respect to $\mfh$. As
usual we have a decomposition of $\mfg^{\complex}$ into root
subspaces

$$\mfg^{\complex}=\mfh\oplus(\oplus_{\al\in \R}\mfg^{\al}),$$

where $\mfg^{\al}=\{u\in\mfg^{\complex}:ad\,h(u)=\al u, \forall
h\in\mfh\}$.

We have $Kill(\mfg^{\al},\mfg^{\be})\neq 0$ if and only if
$\al+\be=0$ and thus, for every $\al\in \R$, we can take
$E_{\al}\in \mfg^{\al}$ such that $Kill(E_{\al},E_{-\al})=1$.
Since $[\mfg^{\al},\mfg^{-\al} ]\subset\mfh$, each pair
$E_{\al},E_{-\al}$ determines an element $H_{\al}$ in the Cartan
subalgebra $\mfh$ given by $H_{\al}=[E_{\al},E_{-\al}]$. The
vectors $H_{\al}$ are such that $Kill(H_{\al},h)=\al(h)$, for
every $h\in\mfh$. In particular, the length $|\al|$ of a root
$\al\in\R$ in $\mfg$ is defined by

\beq |\al|^2=\al(H_{\al})=Kill(H_{\al},H_{\al}).\eeq

For every $\al,\be\in\R$ such that $\al+\be\in\R$, since
$[\mfg^{\al},\mfg^{\be} ]\subset\mfg^{\al+\be}$, we define numbers
$N_{\al,\be}\in \complex$ by

\beq\label{sc1}[E_{\al},E_{\be} ]=N_{\al,\be}E_{\al+\be},\eeq

called the \emph{structure constants}. The $N_{\al,\be}$'s satisfy
the following properties:

\beqar \label{scprop1}N_{\al,\be}=-N_{\be,\al}\\
\label{scprop2}N_{-\al,\be+\al}=N_{-\be,-\al}=N_{\al,\be},\eeqar

for every $\al,\be\in\R$ such that $\al+\be\in\R$.

A basis $\{E_{\al}\}_{\al\in \R}$ of $\oplus_{\al\in R}\mfg^{\al}$
formed by elements chosen as above is called a \emph{standard
normalized basis} and that is what we will use throughout. By
using such a basis we construct the elements

\beq\label{basiscf}
X_{\al}=\frac{E_{\al}-E_{-\al}}{\sqrt{2}}\textrm{ and
}Y_{\al}=\frac{i(E_{\al}+E_{-\al})}{\sqrt{2}}.\eeq

The vectors $X_{\al}$ and $Y_{\al}$ are unit vectors with respect
to $B$. Together with the maximal toral subalgebra $i\mfh_{\reals}$, $X_{\al}$ and
$Y_{\al}$ generate a compact real form for $\mfg^{\complex}$ which
we identify with $\mfg$ (see e.g. \cite{He}, ch.III).

\li

Since $L\subset K$, $\mfh$ is also a Cartan subalgebra for
$\mfk^{\complex}$. We define the following subsets of roots

\beqar \R_{\mfl}=\{\al\in\R:E_{\al}\in \mfl^{\complex}\}\\
\R_{\mfk}=\{\al\in\R:E_{\al}\in \mfk^{\complex}\}\\
\R_{\mfn}=\R-\R_{\mfk}=\{\al\in\R:E_{\al}\in \mfn^{\complex}\}\\
\R_{\mfp}=\R_{\mfk}-\R_{\mfl}=\{\al\in\R:E_{\al}\in
\mfp^{\complex}\}\eeqar

We may also consider the subsets of roots
\beq\R_{\mfp_a}=\{\al\in\R:E_{\al}\in\mfp_a^{\complex}\},\,a=1,\ldots,s,\eeq

and, since $\mfl$ has maximal rank and $\mfk=\mfl\oplus\mfp$,
$\mfp_a=<X_{\al},Y_{\al} :\al\in\R_{\mfp_a}^+>$. Since
$Kill(E_{\al},E_{-\al})=1$, the bases
$\{E_{\al}\}_{\al\in\R_{\mfp_a}}$ and
$\{E_{-\al}\}_{\al\in\R_{\mfp_a}}$ of $\mfp_a^{\complex}$ are dual
with respect to $Kill$. Moreover,
$\{X_{\al},Y_{\al}\}_{\al\in\R_{\mfp_a}^+}$ is an orthonormal
basis for $\mfp_a$ with respect to $B$. Consequently, the Casimir
operators of $\mfp_a^{\complex}$ and of $\mfp_a$ are

\beqar
\label{cpacomplex}C_{\mfp_a^{\complex}}=\sum_{\al\in\R_{\mfp_a}}ad_{E_{\al}}ad_{E_{-\al}}\\
\label{cpareal}C_{\mfp_a}=-\sum_{\al\in\R_{\mfp_a}^+}\big(ad_{X_{\al}}^2+ad_{Y_{\al}}^2\big).
\eeqar

\li

Since $\mfk$ has maximal rank $\mfg=\mfk\oplus\mfn$, we have
$\mfn^{\complex}=<E_{\al}:\al\in\R_{\mfn}>$ and
$\mfn=<X_{\al},Y_{\al}:\al\in\R_{\mfn}^+>$. The subspace $\mfn$ is
by hypothesis irreducible as an $Ad\,K$-submodule. If
$\mfn=\oplus_j\mfn^j$ is a decomposition of $\mfn$ into
irreducible $Ad\,L$-modules, we write

\beq \R_{\mfn^j}=\{\phi\in\R: E_{\phi}\in(\mfn^j)^{\complex}
\}\textrm{ and }\mfn^j=\{X_{\phi},
Y_{\phi}:\phi\in\R_{\mfn^j}^+\}.\eeq

\li

We recall that one of the conditions for existence of an Einstein
adapted metric  on $M$, given in Corollary \ref{cond1}, is that
there are $\la_1,\ldots,\la_s>0$ such that the operator
$\sum_{a=1}^s\la_aC_{\mfp_a}$ is scalar on $\mfn$. The Casimir
operator $C_{\mfp_a}$ is necessarily scalar on the irreducible
$Ad\,L$-submodules $\mfn^j$. Since $B(X_{\phi},X_{\phi})=1$, the
eigenvalue of $C_{\mfp_a}$ on $\mfn^j$ is given by
$B(C_{\mfp_a}X_{\phi},X_{\phi})$, for any $\phi\in\R_{\mfn^j}^+$. We
shall write $b_a^{\phi}$ for this eigenvalue, i.e.,

\beqar b_a^{\phi}=B(C_{\mfp_a}X_{\phi},X_{-\phi})\\
C_{\mfp_a}\mid_{\mfn^j}=b_a^{\phi} Id_{\mfn^j},\,
\forall\,\phi\in\R_{\mfn^j} .\eeqar

Furthermore, the eigenvalue of $C_{\mfp_a}$ on $\mfn^j$,
$b_a^{\phi}$, must coincide with the eigenvalue of
$C_{\mfp_a^{\complex}}$ on $(\mfn^j)^{\complex}$. Hence, we also
have

\beq \label{cpa2}
b_a^{\phi}=Kill(C_{\mfp_a^{\complex}}E_{\phi},E_{-\phi}).\eeq

\brem We observe that in previous Sections the notation $\mfn_j$
was used to denote $Ad\,K$-irreducible submodules of $\mfn$, while
in here we use the similar notation $\mfn^j$ to denote
$Ad\,L$-irreducible submodules, whereas $\mfn$ is
$Ad\,K$-irreducible. Similarly, $b_a^j$ was used before to denote
the eigenvalue of $C_{\mfp_a}$, if this operator was scalar, on
$\mfn_j$, while in here $b_a^{\phi}$ is the eigenvalue of this
same operator on $\mfn^j$. No confusion should arise from this
since we shall use this second notation only when $\mfn$ is an
irreducible $Ad\,K$-module.

\erem

The necessary condition for existence of an adapted Einstein metric
on $M$ given in Corollary \ref{cond1} can now be rewritten as
follows:

\bcor\label{cond1mr} If there exists on $M$ an Einstein adapted
metric, then there are positive constants $\la_1,\ldots,\la_s$
such that
$$\sum_{a=1}^s\la_a(b_a^{\phi_1}-b_a^{\phi_2})=0,$$

for every $\phi_1,\phi_2\in\R_{\mfn}$.\ecor

The condition in Corollary \ref{cond1mr} shall play a fundamental
role as a preliminary test for existence of Einstein adapted
metrics. It is a very restrictive condition which is not satisfied
by many of the spaces under study in this Chapter.

\li

For any roots $\phi$ and $\al$ let $\phi+n\al$, $p_{\al\phi}\leq
n\leq q_{\al\phi}$, be the $\al$-series containing $\phi$. By
definition, the $\al$-series containing $\phi$ is the set of all
roots of the form $\phi+n\al$ where $n$ is an integer. It is known
that $\phi+n\al$ is an interrupted series (\cite{He}, Chap.III, \S
4). For roots $\al$ and $\phi$ the square of the structure
constant $N_{\al\phi}$ is given by

\beq\label{sc1}N_{\al,\phi}^2=\frac{q_{\al\phi}(1-p_{\al\phi})}{2}\al(H_{\al}).\eeq

\li

\bprop \label{bejs}Suppose that $rank\,L=rank\,G$. For every
$\phi\in \R_{\mfn}$ and $a=1,\ldots,s$,

$$b_a^{\phi}=\dfrac{1}{2}\sum _{\al\in
\R_{\mfp_a}^+}d_{\al\phi}|\al|^2,$$

where $d_{\al\phi}=q_{\al\phi}-p_{\al\phi}-2p_{\al\phi}q_{\al\phi}$
and $\phi+n\al$, $p_{\al\phi}\leq n\leq q_{\al\phi}$ is the
$\al$-series containing $\phi$.

\eprop

\bproof By using (\ref{cpa2}) and (\ref{cpacomplex}) we obtain the
following:

$$\bar{rl}b_a^{\phi}= & Kill(C_{\mfp_a^{\complex}}E_{\phi},E_{-\phi})\\ \\

= &
\sum_{\al\in\R_{\mfp_a}}Kill([E_{-\al},[E_{\al},E_{\phi}]],E_{-\phi})\\
\\

= &
\sum_{\al\in\R_{\mfp_a}}N_{\al,\phi}Kill([E_{-\al},E_{\phi+\al}],E_{-\phi})\\
\\
= &
\sum_{\al\in\R_{\mfp_a}}N_{\al,\phi}N_{-\al,\phi+\al}Kill(E_{\phi},E_{-\phi})\\
\\

= & \sum_{\al\in\R_{\mfp_a}}N_{\al,\phi}N_{-\al,\phi+\al}\ear$$

From (\ref{scprop2}) we have $N_{-\al,\phi+\al}=N_{\al,\phi}$ and
we get

$$b_a^{\phi}=\sum_{\al\in\R_{\mfp_a}}N_{\al,\phi}^2=\sum_{\al\in\R_{\mfp_a}^+}(N_{\al,\phi}^2+N_{-\al,\phi}^2).$$

Now let $\phi+n\al$, $p_{\al\phi}\leq n\leq q_{\al\phi}$, be the
$\al$-series containing $\phi$. It is known that

$$N_{\al,\phi}^2=\frac{q_{\al\phi}(1-p_{\al\phi})}{2}\al(H_{\al}),$$

as mentioned in (\ref{sc1}).

On the other hand, to compute $N_{-\al,\phi}^2$ we need the
$(-\al)$-series containing $\phi$. Clearly, this series is
$\phi-n'\al$, where $-q_{\al\phi}\leq n'\leq -p_{\al\phi}$. Hence,
we obtain the following:

$$N_{-\al,\phi}^2=\frac{-p_{\al\phi}(1-(-q_{\al\phi}))}{2}(-\al)(H_{-\al})=\frac{-p_{\al\phi}(1+q_{\al\phi})}{2}\al(H_{\al}).$$

Hence,

$$b_a^{\phi}=\sum_{\al\in\R_{\mfp_a}^+}\left(\frac{q_{\al\phi}(1-p_{\al\phi})}{2}-\frac{-p_{\al\phi}(1+q_{\al\phi})}{2}\right)\al(H_{\al}), $$

which yields the required formula.

$\Box$\eproof

\li

Let us consider a decomposition of $\mfk$ into its center $\mfk_0$
and simple ideals $\mfk_a$, for $a=1,\ldots, t$,

\beq\label{dec1}\mfk=\mfk_0\oplus\mfk_1\oplus\ldots\oplus\mfk_t,\eeq

and let $\ga_a$ denote the eigenvalue of the Casimir operator of
$\mfk$ on $\mfk_a$. We present a formula to compute the eigenvalues
$\ga_a$'s by making use of dual Coxeter numbers. We start by
recalling some facts about roots. On this topic we refer to
(\cite{BD}, V.5) and (\cite{Hu}, 10.4).

There are at most two different lengths in a given irreducible root
system, and the corresponding roots are designated by \textbf{long}
and \textbf{short} roots. If there is only one length it is
conventional to say that all the roots are long. If $\al$ is a long
root and $\be$ is short, then

\beq\label{rl1}\bar{ll}\dfrac{|\al|^2}{|\be|^2}=3, &\textrm{ in the case of }G_2 \textrm{ and}\\
\dfrac{|\al|^2}{|\be|^2}=2, &\textrm{ in the case of
}B_n,C_n\textrm{ and } F_4. \ear\eeq

In the remaining cases, $A_n$, $D_n$, $E_6$, $E_7$ and $E_8$, there
is only one length. These facts can be read off from the
corresponding Dynkin diagrams.

We also recall that a length of a root $\al$ is given by
$$|\al|^2=\al(H_{\al})=Kill(H_{\al},H_{\al}).$$

The \textbf{dual Coxeter number} of a simple Lie algebra $\mfg$ is
the number given by

$$h^*(\mfg)=\dfrac{1}{|\al|^2},$$

where $\al$ is a long root (see e.g. \cite{Pa}). The dual Coxeter
numbers of each irreducible root system are given in Table
\ref{tabcoxeter}.

We may suppose that $\mfh_a=\mfh\cap\mfk_a$ is a Cartan subalgebra
of $\mfk_a$ and thus a root of $\mfk_a$ can be viewed as a root for
$\mfg$. Hence we can compare lengths of roots of $\mfg$ with lengths
of roots of $\mfk_a$. So let $\de_a$ be the ratio of the square
length of a long root for $\mfg$ to that of $\mfk_a$, i.e.,

$$\de_a=\dfrac{|\al|^2_{\mfg}}{|\be|^2_{\mfg}}=\dfrac{Kill(H_{\al},H_{\al})}{Kill(H_{\be},H_{\be})},$$

where $\al$ is a long root of $\mfg$ and $\be$ is a long root of
$\mfk_a$. Clearly, $\de_a=1$ if there exists only one length for
$\mfg$ or if both $\mfg$ and $\mfk_a$ have two lengths. If
$\de_a\neq 1$, then, according to (\ref{rl1}), $\de_a$ is equal to
either 2 or 3. We recall the following result by D. Panyushev:

\bprop \label{gajs}\cite{Pa} Suppose that $\mfg$ is simple. Then

$$\ga_a=\dfrac{h^*(\mfk_a)}{\de_a . h^*(\mfg)},\,a=1,\ldots,s,$$

where $h^*(\mfk_a)$ and $h^*(\mfg)$ are the dual Coxeter numbers of
$\mfk_a$ and $h^*(\mfg)$, respectively. \eprop

We observe that so far we used only the fact that $L$ has maximal
rank. Now if, in addition, $F=K/L$ is a symmetric space, then we
consider its DeRham decomposition

\beq\label{deRham2} K/L= K_1/L_1\times\ldots\times K_s/L_s,\eeq

where, for $a=1,\ldots,s$, $K_a$ is simple, as in Section
\ref{sectionsymF}. We observe that since $L$ has maximal rank in the
deRham decomposition of $K/L$ the factor $K_0/L_0$ must be trivial.
The Lie algebras of $K_1,\ldots, K_s$, denoted by
$\mfk_1,\ldots,\mfk_s$, are some of the ideals in the decomposition
(\ref{dec1}). As explained at the beginning of Section
\ref{sectionsymF}, the irreducible $Ad\,L$-submodules $\mfp_1,
\ldots,\mfp_s$ are chosen as the symmetric reductive complements of
$\mfl_a$ in $\mfk_a$, for $a=1,\ldots,s$. Hence, the constant
$\ga_a$ defined throughout by the equality
$Kill_{\mfk}\mid_{\mfp_a\times\mfp_a}=\ga_a
Kill\mid_{\mfp_a\times\mfp_a}$ in previous sections, is now just
 the eigenvalue of the Casimir operator of
$\mfk$ on $\mfk_a$, and thus can be determined by the formula in
Proposition \ref{gajs}.

In Appendix A we compute the eigenvalues $b_a^{\phi}$'s and
$\ga_a$'s using Propositions \ref{bejs} and \ref{gajs}. Their
values are indicated in Tables \ref{eigIexc}, \ref{eigIclass},
\ref{eigIIexc} and \ref{eigIIclass} in Section \ref{alltabs} for
each bisymmetric triple.

\section{Bisymmetric Triples of Maximal Rank - Classification}

Let us consider a homogeneous fibration $F\rightarrow M \rightarrow
N$, for $M=G/L$, $N=G/K$ and $F=K/L$, where $G$ is a compact
connected semisimple Lie group and $L\varsubsetneq K\varsubsetneq G$
connected closed non-trivial subgroups such that $F$ and $N$ are
symmetric spaces. We shall call such a fibration a \textbf{bisymmetric fibration}. In particular, the pairs $(G,K)$ and $(K,L)$ are symmetric pairs of compact type \cite{He}. With a slight abuse of terminology, we shall also say that the pairs of Lie algebras $(\mfg,\mfk)$ and $(\mfk,\mfl)$ are symmetric pairs of compact type whenever the corresponding pairs $(G,K)$ and $(K,L)$ are.

\bdfn A \textbf{bisymmetric triple} is a triple $(\mfg,\mfk,\mfl)$ where $\mfg$, $\mfk$ and $\mfl$ are Lie algebras satisfying the following conditions:

(i) $\mfl \varsubsetneq \mfk \varsubsetneq \mfg$;

(ii) $(\mfg,\mfk)$ and $(\mfk,\mfl)$ are symmetric pairs of compact type.

A bisymmetric triple is said to be irreducible if $\mfg$ is a simple Lie algebra and said to be of maximal rank if $\mfl$ has maximal rank in $\mfg$, i.e., it contains a maximal toral subalgebra of $\mfg$.

\edfn

Clearly, there is a one-to-one correspondence between bisymmetric fibrations, up to cover, and bisymmetric triples. \textbf{All the bisymmetric triples $\mathbf{(\mfg,\mfk,\mfl)}$ considered in this chapter are irreducible and of maximal rank}, even when this is not explicitly stated. Consequently, \textbf{any bisymmetric fibration $\mathbf{F\rightarrow M \rightarrow
N}$ here considered is such that $\mathbf{L}$ is a subgroup of maximal rank in $\mathbf{G}$ and  $\mathbf{N=G/K}$ is an irreducible symmetric space}.

\li

\bdfn A bisymmetric triple $(\mfg,\mfk,\mfl)$ is said to be of

(i) \textbf{Type I} if $F$ is an isotropy irreducible symmetric
space; equivalently, if $\mfp$ is an irreducible $Ad\,L$-module;

(ii) \textbf{Type II} if $F$ is the direct product of two isotropy
irreducible symmetric spaces; equivalently, if
$\mfp=\mfp_1\oplus\mfp_2$, where $\mfp_1$ and $\mfp_2$ are
nontrivial irreducible $Ad\,L$-modules.

\li

A bisymmetric fibration $F\rightarrow M\rightarrow N$ of is said to
be of Type I or II if the corresponding bisymmetric triple
$(\mfg,\mfk,\mfl)$ is either of Type I or II, respectively.\edfn

As we shall see any irreducible bisymmetric triple of maximal rank with $\mfg$
simple is of Type I or II.

\li

Isotropy irreducible symmetric spaces have been classified and a
classification can be found in \cite{He}. By using this we obtain a
list of all possible triples $(\mfg,\mfk,\mfl)$ such that $\mfl$ and
$\mfk$ are subalgebras of maximal rank of $\mfg$ and $(\mfg,\mfk)$
and $(\mfk,\mfl)$ are symmetric pairs of compact type. By inspection
of the classification of symmetric pairs $(\mfg,\mfk)$ of compact
type in \cite{He} we obtain that those of maximal rank are the pairs
in Tables \ref{spexc} and \ref{spclass}.

We observe that the cases when $\mfk$ is the centralizer of a torus
are only the cases

$(\mfe_6,\mfso_{10}\oplus \reals )$,
$(\mfe_7,\mfe_6\oplus\reals)$, $(\mfso_{2n},\mfu_n)$,
$(\mfso_n,\reals\oplus\mfso_{n-2})$, $(\mfsp_{n},\mfu_n)$ and
$(\mfsu_n,\mfsu_p\oplus\mfsu_{n-p}\oplus\reals)$. This follows
from the fact that these are the only subalgebras $\mfk$
corresponding to painted Dynkin diagrams of the Dynkin diagram of
$\mfg$ or as they are the only ones such that $\mfk$ is not
centerless. In all the other cases $\mfk$ is semisimple. If $\mfk$
is simple then $(\mfk,\mfl)$ shall be an irreducible symmetric
pair, i.e., $\mfp$ is an irreducible $L$-invariant subspace. Thus,
$(\mfg,\mfk,\mfl)$ is of type I. In the cases where
$\mfk=\mfk_1\oplus\reals$ with $\mfk_1$ a simple ideal of $\mfk$,
since we require $\mfl$ to be of maximal rank, we have
$\mfl=\mfl_1\oplus\reals$, where $\mfl_1$ is a subalgebra of
$\mfk_1$ with maximal rank and $(\mfk,\mfl)\iso(\mfk_1,\mfl_1)$ is
an irreducible symmetric pair. Thus, in this case,
$\mfp\iso\mfp_1$ is also an irreducible $L$-invariant subspace and
$(\mfg,\mfk,\mfl)$ is of type I. In the cases where
$\mfk=\mfk_1\oplus\mfk_2$, with $\mfk_1$ and $\mfk_2$ simple
ideals of $\mfk$, we have $\mfl=\mfl_1\oplus\mfl_2$, where, for
$i=1,2$, $\mfl_i$ is a subalgebra of $\mfk_i$ of maximal rank.
Clearly, one of the $\mfl_i$'s must be proper as we require that
$\mfl$ is a proper subalgebra of $\mfk$. If both $\mfl_1$ ad
$\mfl_2$ are proper, then $\mfp=\mfp_1\oplus\mfp_2$, where
$\mfp_1$ and $\mfp_2$ are nonzero irreducible $\mfl$-invariant
subspaces. Hence, in this case, $(\mfg,\mfk,\mfl)$ is of type II.
If exactly one of $\mfl_i$'s coincides with $\mfk_i$, then
$(\mfk,\mfl)\iso(\mfk_j,\mfl_j)$ and $\mfp=\mfp_j$, for that $j$
satisfying $\mfl_j\neq\mfk_j$, and once again $(\mfg,\mfk,\mfl)$
is of type I. Finally, we have the case of the spaces
$(\mfsu_n,\mfsu_p\oplus\mfsu_{n-p}\oplus\reals)$,
$p=1,\ldots,n-1$. Clearly, $\mfl$ must be of the form
$\mfl=\mfl_1\oplus\mfl_2\oplus\reals$, where $\mfl_1$ and $\mfl_2$
are maximal rank subalgebras of $\mfsu_p$ and $\mfsu_{n-p}$,
respectively. We obtain a triple of type $I$ if exactly one of the
$\mfl_i$'s is proper and a triple of type II if both $\mfl_1$ and
$\mfl_2$ are proper. This proves the following:

\blem An irreducible bisymmetric triple of maximal rank
$(\mfg,\mfk,\mfl)$ such that $\mfg$ is simple is either of Type I
or II. Moreover, all such bisymmetric triples $(\mfg,\mfk,\mfl)$
are those in Tables and \ref{eigIexc}, \ref{eigIclass},
\ref{eigIIexc} and \ref{eigIIclass}.\elem

\blem \label{Ilemma1} For an irreducible bisymmetric triple of maximal rank
$(\mfg,\mfk,\mfl)$ of Type I, let $\ga$ be the eigenvalue of the
Casimir operator of $\mfk$ on $\mfp$ and $b^{\phi}$'s the
eigenvalues of the Casimir operator of $\mfp$ on $\mfn$. For each
bisymmetric triple of maximal rank these eigenvalues have the values
listed in Tables \ref{eigIexc} and \ref{eigIclass}. \elem

\blem\label{IIlemma1} For an irreducible bisymmetric triple of maximal rank
$(\mfg,\mfk,\mfl)$ of Type II, let $\ga_a$ be the eigenvalue of
the Casimir operator of $\mfk$ on $\mfp_a$ and $b_a^{\phi}$'s the
eigenvalues of the Casimir operator of $\mfp_a$ on $\mfn$,
$a=1,2$. For each bisymmetric triple of maximal rank these
eigenvalues have the values listed in Tables \ref{eigIIexc} and
\ref{eigIIclass}. \elem

\section{Einstein Adapted Metrics for Type I}

In this Section we determine all the bisymmetric Riemannian
fibrations $F\rightarrow M \rightarrow N$ of maximal rank of Type
I which admit an Einstein adapted metric. We recall that for Type
I, $\mfp$ is an irreducible $Ad\,L$-submodule. Moreover, since
$\mfn$ is an irreducible $Ad\,K$-module, any adapted metric is
binormal. As in addition $F$ and $N$ are symmetric spaces, we may
apply Corollary \ref{binSym1}. We recall that according to
Corollary \ref{binSym1}, there exists on $M$ an Einstein
(binormal) adapted metric if and only if

\beq\bar{l}(i)\textrm{ the Casimir operator of $\mfp$ is scalar on
$\mfn$
and}\\ \\
(ii)\label{tI1} \triangle=1-2\ga(1-\ga+2b)\geq 0.\ear\eeq

If these two conditions are satisfied, Einstein binormal metrics are, up to homothety,
given by

\beq \label{metTypeI}g_M=B\mid_{\mfp\times\mfp}\oplus
XB\mid_{\mfn\times\mfn},\textrm{ where
}X=\frac{1\pm\sqrt{\triangle}}{2\ga}.\eeq

We also recall that $b$ is the eigenvalue of $C_{\mfp}$ on $\mfn$,
in the case when this operator is scalar, and $\ga$ is the
eigenvalue of the Casimir operator of $\mfk$ on $\mfp$. These
constants are computed in Appendix \ref{cpproofs} and their values
are indicated in Tables \ref{eigIexc}  and \ref{eigIclass}, as has
been stated in Lemma \ref{Ilemma1}. We recall that condition (i)
translates into $b^{\phi}=b$, for every $\phi\in\R_{\mfn}$, for
Type I, according to Corollary \ref{cond1mr}. Hence, the first
test for existence of an adapted Einstein metric shall be to
observe if there exists only one eigenvalue $b^{\phi}$ in the
corresponding columns of Tables \ref{eigIexc}  and
\ref{eigIclass}.

\bthm \label{mtypeI}The bisymmetric fibrations $F\rightarrow M\rightarrow N$ of Type I such that there exists on
$M$ an Einstein adapted metric are those whose bisymmetric triples are listed in Tables
\ref{mIexc} and \ref{mIclass}. For each case there are exactly two
Einstein adapted metrics. Furthermore, these Einstein metrics
are, up to homothety, given by

$$g_M=B\mid_{\mfp\times\mfp}\oplus XB\mid_{\mfn\times\mfn},$$ where
$X$ is indicated in the Tables mentioned above.  \ethm

\bproof As explained in the discussion above, by Corollary
\ref{binSym1}, the existence of an adapted Einstein metric implies
that the Casimir operator of $\mfp$ is scalar on $\mfn$. By
inspection we conclude from Tables \ref{eigIclass}  and
\ref{eigIexc} that the only spaces satisfying this condition are
those corresponding to the labels

$$\ref{cpbn3},\,\ref{cpdn5},\,\ref{cpcn5},\,\ref{cpf41},\,\ref{cpf42}\,\ref{cpg21},\,\ref{cpg22},\,\ref{cpe81},\,\ref{cpe83},\,\ref{cpe88},\,\ref{cpe71},\,\ref{cpe74},\,\ref{cpe62},\,\ref{cpe63}$$

and

\ref{cpan1} for $l=\frac{p}{2}$ with $p$ even; in this case,
$b=\frac{p}{4n}$;

\ref{cpbn2} for $s=\frac{n-p}{2}$ with $n-p$ even; in this case,
$b=\frac{n-p}{2(2n-1)}$;

\ref{cpdn2} for $l=\frac{p}{2}$ with $p$ even; in this case,
$b=\frac{p}{4(n-1)}$;

\ref{cpcn2} for $l=\frac{p}{2}$ with $p$ even; in this case,
$b=\frac{p}{8(n+1)}$;

\ref{cpe78} for $p=1$; in this case $b=\frac{1}{9}$;

\ref{cpe64} for $p=1$; in this case $b=\frac{1}{8}$.

\li

We compute $\triangle$ given in formula \ref{tI1}, and the values
obtained are as follows

{\footnotesize{
$$\bar{ll} & \triangle \\ \hline

\xstrut\ref{cpan1} & \left(\frac{n-p}{n}\right)^2>0\\

\xstrut\ref{cpbn2} & \frac{4p^2+8p-4n+5}{(2n-1)^2}>0,\,\forall
p=\lfloor \frac{\sqrt{4n-1}}{2}\rfloor,\ldots,n-1 \\

\xstrut\ref{cpbn3} & \left(\frac{2p+1}{2n-1}\right)^2>0\\

\xstrut\ref{cpdn2} & \frac{p^2-(2n+1)p+n^2+1}{(n-1)^2}>0 \\

\xstrut\ref{cpdn5} & \left(\frac{p-n}{n-1}\right)^2>0
\\

\xstrut\ref{cpcn2} &
\frac{3p^2+(3-4n)p+2(n^2+1)}{2(n+1)^2}>0\\

\xstrut\ref{cpcn5} & \left(\frac{n-p}{n+1}\right)^2>0\\

\xstrut\ref{cpf41} &  \frac{106-63p+7p^2}{162}> 0,\,\textrm{ iff
}p=1,\,7\\

\xstrut\ref{cpf42} & \frac{49}{81}>0\\

\xstrut\ref{cpg21} & \frac{1}{4}>0 \\

\xstrut\ref{cpg22} & \frac{11}{18}>0\\

\xstrut\ref{cpe81} & \frac{7p^2-56p+113}{225}>0\\

\xstrut\ref{cpe83} & \frac{196}{225}>0\\

\xstrut\ref{cpe88} & -\frac{2}{25}<0 \\

\xstrut\ref{cpe71} & \frac{64}{81}>0 \\

\xstrut\ref{cpe74} & \frac{164-60p+5p^2}{324}>0,\,p=2,4 \\

\xstrut\ref{cpe78} & \frac{25}{81}>0\\

\xstrut\ref{cpe62} & \frac{1}{9}>0,\, p=2; -\frac{1}{9},\,p=4\\

\xstrut\ref{cpe63} & \frac{25}{36}>0\\

\xstrut\ref{cpe64} & \frac{1}{4}>0
 \ear$$

 }}

For $\ref{cpf41}$, since $p=1,3,5,7$, then $\triangle>0$ for $p=1,7$
and $\triangle<0$ otherwise.

For $\ref{cpe74}$, we have $p=2,4,6$. Then, $\triangle>0$ for
$p=2,4$ and $\triangle<0$ for $p=6$.

For $\ref{cpe81}$, we have $p=1,\ldots,4$ and thus $\triangle>0$ for
every $p$.

For $\ref{cpdn2}$, $p=1,\ldots,\lfloor \frac{n}{2}\rfloor$. We have
that $\triangle\geq 0$ if and only if $$p\in
\left(\big(-\infty,\frac{2n+1-\sqrt{4n-3}}{2}\big)\cup
\big(\frac{2n+1+\sqrt{4n-3}}{2},+\infty\big)\right)\cap\{1,\ldots,\lfloor
\frac{n}{2}\rfloor\}$$

We can show that $\frac{2n+1-\sqrt{4n-3}}{2}\geq \frac{n}{2}$ and
thus $\triangle> 0$, for every $p=1,\ldots,\lfloor
\frac{n}{2}\rfloor$.

For $\ref{cpbn2}$, $p=1,\ldots,n-1$. We show that $\triangle\geq 0$,
if and only if
$$p\in \left(-1+\frac{\sqrt{4n-1}}{2},+\infty\right)\cap\{1,\ldots,n-1\}$$

Since, for every $n$, $-1+\frac{\sqrt{4n-1}}{2}<n-1$ we conclude
that $\triangle> 0$, if and only if $p=\lfloor
-1+\frac{\sqrt{4n-1}}{2}\rfloor+1,\ldots,n-1=\lfloor
\frac{\sqrt{4n-1}}{2}\rfloor,\ldots,n-1$.

\li

Finally, we compute $X=\frac{1\pm\sqrt{\triangle}}{2\ga}$ for those
cases when $\triangle>0$. The values of $X$ are indicated in Table
\ref{mIexc}.

$\Box$\eproof

\brem For the bisymmetric triples of Type I such that $C_{\mfp}$ is
scalar on $\mfn$ but $\triangle<0$, there is still a complex non
real solution $X=\frac{1\pm\sqrt{\triangle}}{2\ga}$. So even in
these cases we can conclude that there exists a complex Einstein
adapted metric on $M$ as in (\ref{metTypeI}). These are just the
cases \ref{cpbn2} for $p\geq \lfloor
-1+\frac{\sqrt{4n-1}}{2}\rfloor$, \ref{cpf41} for $p=3,5$ and
\ref{cpe88}.

\erem

\li

\section{Einstein Adapted Metrics for Type II}

In this Section we study the existence of Einstein adapted metrics
on bisymmetric fibrations of Type II. Whereas for Type I any
adapted metric was binormal this is clearly not true for Type II,
since $\mfp$ is not an irreducible $Ad\,L$-module. We shall
classify all the bisymmetric triples which admit an Einstein
binormal metric. Since for bisymmetric fibrations, in particular,
$F$ is a symmetric space, we know from Corollary \ref{binormal6}
that for any Einstein binormal metric $g_M$, $g_F$ is also
Einstein. We shall also classify all the bisymmetric triples which
admit an Einstein non-binormal adapted metric $g_M$ whose
restriction $g_F$ is also Einstein.

\li

Since for type II $\mfp=\mfp_1+\mfp_2$, where $\mfp_i$, $i=1,2$, are
irreducible $L$-modules and $\mfn$ is an irreducible $K$-module, the
existence of an Einstein adapted metric implies that there exist
positive constants $\la_1$, $\la_2$, such that the operator
$\la_1C_{\mfp_1}+\la_2C_{\mfp_2}$ is scalar on $\mfn$, according to
Corollary \ref{cond1}. As rephrased  in Corollary \ref{cond1mr}, the
condition above translates into

\beq\label{condIIb}
\la_1(b_1^{\phi_1}-b_1^{\phi_2})+\la_2(b_2^{\phi_1}-b_2^{\phi_2})=0,\,\textrm{
for every }\phi_1,\phi_2\in\R_{\mfn}\eeq

By using Tables \ref{eigIIexc} and \ref{eigIIclass} we conclude the
following:

\blem \label{1testeII}The only bisymmetric triples satisfying
condition (\ref{condIIb}) are the cases \ref{cpdn7}, \ref{cpcn7},
\ref{cpg23}, \ref{cpe89}, \ref{cpe75} and

$$\bar{l}\xstrut\ref{cpan3}\, for\, $p=2l$,\, $n-p=2s$\\

\xstrut\ref{cpdn4}\, for\, $p=2l$,\, $n-p=2s$\\

\xstrut\ref{cpdn8}\, for\, $p=2l$\\

\xstrut\ref{cpcn4}\, for\, $p=2l$,\, $n-p=2s$\\

\xstrut\ref{cpcn8}\, for\, $p=2l$\\

\xstrut\ref{cpe65}\, for\, $p=1$\ear$$

For all other bisymmetric triples of Type II we can conclude that
there exists no Einstein adapted metric on $M$. \elem

Furthermore, for all the triples listed in Lemma \ref{1testeII}, we
observe that $C_{\mfp_1}$ and $C_{\mfp_2}$ are scalar on $\mfn$ and
thus $C_{\mfp}$ is also scalar on $\mfn$. We shall write

\beqar \label{expcp1}C_{\mfp_i}\mid_{\mfn}=b_iId_{\mfn},\, i=1,2\\
\label{expcp2}C_{\mfp}\mid_{\mfn}=bId_{\mfn},\textrm{ for
}b=b_1+b_2,\eeqar

following the notation used in previous chapters.

\bthm \label{biII}The bisymmetric fibrations $F\rightarrow M\rightarrow N$ of Type II such that there exists on
$M$ an Einstein binormal metric are those whose bisymmetric triples are listed
in Table \ref{bimII}. Furthermore, the binormal Einstein metrics
are, up to homothety, given by

$$g_M=B\mid_{\mfp\times\mfp}\oplus XB\mid_{\mfn\times\mfn},$$ where
$X$ is indicated in Table \ref{bimII}. In all the cases, $g_N$ and
$g_F$ are also Einstein.\ethm

\bproof Binormal Einstein metrics are in this case given by
Corollary \ref{binormal5}. First we observe that in order to exist
an Einstein binormal metric on $M$, $C_{\mfp}$ must be scalar on
$\mfn$ and $C_{\mfk}$ must be scalar on $\mfp$. The triples which
satisfy the first condition are those listed in Lemma
\ref{1testeII}. Furthermore, the second condition implies that
$\ga_2=\ga_1$. From the cases in Lemma \ref{1testeII} we conclude
from Tables \ref{eigIIexc} and \ref{eigIIclass}  that the spaces
which satisfy the condition $\ga_2=\ga_1$ are those listed below:

\li

\ref{cpan3} for $s=l=2p$, $n=4l$; in this case,
$\ga_1=\ga_2=\frac{1}{2}$ and $b_1=b_2=\frac{1}{8}$;

\li

\ref{cpdn4} for $s=l=2p$, $n=4l$, $l\geq 2$; in this case,
$\ga_1=\ga_2=\frac{2l-1}{4l-1}$ and $b_1=b_2=\frac{l}{2(4l-1)}$;

\li

\ref{cpdn7} for $n=2p$, $p\geq 2$; in this case,
$\ga_1=\ga_2=\frac{p-1}{2p-1}$ and $b_1=b_2=\frac{p-1}{4(2p-1)}$;

\li

\ref{cpdn8} for $p=2l$, $n=4l$; in this case,
$\ga_1=\ga_2=\frac{2l-1}{4l-1}$ and $b_1=\frac{l}{2(4l-1)}$,
$b_2=\frac{2l-1}{4(4l-1)}$;

\li

\ref{cpcn4} for $s=l=2p$, $n=4l$; in this case,
$\ga_1=\ga_2=\frac{2l+1}{4l+1}$ and $b_1=b_2=\frac{l}{4(4l+1)}$;

\li

\ref{cpcn7} for $n=2p$; in this case, $\ga_1=\ga_2=\frac{p+1}{2p+1}$
and  $b_1=b_2=\frac{p+1}{4(2p+1)}$;

\li

\ref{cpcn8} for $p=2l$, $n=4l$; in this case,
$\ga_1=\ga_2=\frac{2l+1}{4l+1}$ and  $b_1=\frac{l}{4(4l+1)}$,
$b_2=\frac{2l+1}{4(4l+1)}$.

Now Einstein binormal metrics are given by positive solutions of
(\ref{polbin5}). Since $c_{\mfk,\mfn}=\frac{1}{2}$, we have that
there exists an Einstein binormal metric if and only if

\beq\label{tII1}\triangle=1-2\ga(1-\ga+2b)\geq 0,\eeq

where $\ga=\ga_1=\ga_2$ and $b$ is the eigenvalue of $C_{\mfp}$ on
$\mfn$, i.e., $b=b_1+b_2$. In such a case these Einstein metrics are
given by homotheties of

$$g_M=B\mid_{\mfp\times\mfp}\oplus XB\mid_{\mfn\times\mfn},\textrm{
where }X=\frac{1\pm\sqrt{\triangle}}{2\ga}.$$

We compute $\triangle$:

$$\bar{ll} & \triangle \\ \hline

\xstrut\ref{cpan3} & 0\\

\xstrut\ref{cpdn4} & \frac{1}{(4l-1)^2}>0\\

\xstrut\ref{cpdn7} & \frac{1}{2p-1}>0\\

\xstrut\ref{cpdn8} & \frac{2l}{(4l-1)^2}>0\\

\xstrut\ref{cpcn4} &\frac{4l^2+2l+1}{(4l+1)^2}>0\\

\xstrut\ref{cpcn7} & -\frac{1}{2p+1}<0\\

\xstrut\ref{cpcn8} & \frac{l(2l-1)}{(4l+1)^2}>0\\

\ear$$

Except in the case \ref{cpcn7}, there exists an Einstein adapted
metric. The values for $X$ are indicated in Table \ref{bimII}.

$g_N$ is Einstein because $N$ is irreducible and $g_F$ is Einstein
due to Corollary \ref{binormal6}.

$\Box$\eproof

\brem In the case \ref{cpcn7}, where $C_{\mfp}$ is scalar on $\mfn$
but $\triangle=-\frac{1}{2p+1}<0$, for every $p$, we can still
consider the non-real complex solution
$X=\frac{1\pm\sqrt{\triangle}}{2\ga}$ which gives rise to a non-real
complex Einstein binormal metric on $M$. \erem

\bthm\label{gFeinII}The bisymmetric fibrations $F\rightarrow
M\rightarrow N$ of Type II such that there exists on $M$ an
Einstein adapted metric such that $g_F$ is Einstein are those with
an Einstein binormal metric, as in Theorem \ref{biII} and Table
\ref{bimII}, and the fibration corresponding to the bisymmetric
triple

$$(\mfsu_{2(l+s)},\mfsu_{2l}\oplus\mfsu_{2s}\oplus\reals,\mfsu_l\oplus\mfsu_l\oplus\mfsu_s\oplus\mfsu_s\oplus\reals^3),$$

whose Einstein adapted metric is, up to homothety, given by

$$g_M=\frac{2l}{l+s}B\mid_{\mfp_1\times\mfp_1}\oplus\frac{2s}{l+s}B\mid_{\mfp_2\times\mfp_2}\oplus B\mid_{\mfn\times\mfn}.$$

This metric is binormal if and only if $l=s$. \ethm

\bproof The only cases which may admit an Einstein adapted metric
are those listed in Lemma \ref{1testeII}, since they are the only
that satisfy the necessary condition (\ref{condIIb}) for existence
of such a metric. Furthermore, if we require that $g_F$ is also
Einstein, then we must have one of the two conditions in Corollary
\ref{gfeincor2}. In the first case, $\ga=1=\ga_2$, $g_M$ is a
binormal metric according to Corollary \ref{gfeincor2}. These are
the cases given by Theorem \ref{biII} and listed in Table
\ref{bimII}. In the second case, we have $\ga_2=1-\ga_1$. By
inspection of the triples listed in Lemma \ref{1testeII}, we
conclude that this is possible only for the triple \ref{cpan3}
when $p=2l$ and $n-p=2s$. In this case $\ga_1=\frac{l}{l+s}$,
$\ga_2=\frac{s}{l+s}=1-\ga_1$, $b_1=\frac{\ga_1}{4}$ and
$b_2=\frac{\ga_2}{4}=\frac{1-\ga_1}{4}$. Using Corollary
\ref{gFeincor2}, we obtain $D(\ga_1)=0$ and thus
$X_1=\frac{l+s}{2l}$, $X_2=\frac{l+s}{2s}$.

$\Box$ \eproof

\li

For bisymmetric fibrations of Type II the classifications of all Einstein adapted metrics is a difficult problem due to the high complexity of the Einstein equations. The classification of Einstein binormal metrics and Einstein adapted metrics such that $g_F$ is Einstein as well, as done above, is a way to restrict the problem in a way that the Einstein equations are manageable. It can be seen by Theorems \ref{gFeinII} and \ref{biII} that no bisymmetric fibration of Type II in the exceptional admits neither an
Einstein binormal metric nor an Einstein adapted metric with $g_F$ Einstein. However, for this spaces of exceptional type is is possible to classify all Einstein adapted metrics with the help of Maple. Once again, since $C_{\mfp_1}$ and $C_{\mfp_2}$ must be scalar on $\mfn$, from Lemma \ref{1testeII} we know that the only cases of Type II and exceptional are the fibrations \ref{cpe65} for $p=1$, \ref{cpg23}, \ref{cpe75} and \ref{cpe89}. The results obtained for these cases are synthetized in Theorem \ref{genII} and  Table \ref{tabgenII}.

\bthm \label{genII}The only bisymmetric fibrations $F\rightarrow M\rightarrow N$ of Type II, such that $G$ is an exceptional Lie group, which admit an Einstein adapted metric are those whose bisymmetric triples are listed in Table \ref{tabgenII}. None of these metrics is neither binormal nor such that $g_F$ is Einstein.\ethm

\bproof  As mentioned above, it follows from Lemma \ref{1testeII} that the only cases of Type II with $\mfg$ exceptional wich may admit an Einstein adapted metric are the cases \ref{cpe65} for $p=1$, \ref{cpg23}, \ref{cpe75} and \ref{cpe89}. We recall that for each of these spaces any Einstein adapted metric is of the form

$$g_M=\frac{1}{X_1}B\mid_{\mfp_1\times\mfp_1}\oplus\frac{1}{X_2}B\mid_{\mfp_2\times\mfp_2}\oplus B\mid_{\mfn\times\mfn},$$ where
$X_1$ and $X_2$ are positive solutions of the system of equations given in Theorem \ref{eqFsym} which are as follows:

\beqar
\label{eqII112}2\ga_1X_1^2X_2+(1-\ga_1)X_2-2\ga_2X_1X_2^2-(1-\ga_2)X_1=0,
\\
\label{eqII122}2b_1X_2+2b_2X_1-2X_1X_2+2\ga_1X_1^2X_2+(1-\ga_1)X_2=0.
\eeqar

Also we recall that the eigenvalues $b_i$ and $\ga_i$, for $i=1,2$, can be found in Table \ref{eigIIexc}. To show the result we use Maple and so many details are ommited.

\textbf{\ref{cpg23})} For the bisymmetric triple \ref{cpg23} the
non-zero solutions of the equations (\ref{eqII112}) and
(\ref{eqII122}) are of the form

$$X_1=\al_i,\,X_2=-\frac{7}{4}\al_i^3+12\al_i^2-\frac{5899}{36}\al_i+19,$$

where $\al_i$ is a root of the polynomial $$t(z)=63z^4-432z^3+1088z^2-1224z+513.$$

Since $X_1=\al_i$ we are interested only in positive real roots of $t$. This polynomial has exactly two positive roots which are

$$\bar{c}\al_1=\frac{12\xi\be}{7}+\frac{\be^3}{126}-\frac{\sqrt{6}}{126}\frac{(320\xi^2\be^2+7\xi^4\be^2-25781\be^2+43416\xi^3)^{\frac{1}{2}}}{\xi\be}  \\

\al_2=\frac{12\xi\be}{7}+\frac{\be^3}{126}+\frac{\sqrt{6}}{126}\frac{(320\xi^2\be^2+7\xi^4\be^2-25781\be^2+43416\xi^3)^{\frac{1}{2}}}{\xi\be} \\

\ear$$

where $\xi=(17756+81\sqrt{7662443})^{\frac{1}{6}}$ and $\be=(960\xi^2-42\xi^4+154686)^{\frac{1}{4}}$.

Simple calculations show that $\al_i$, for $i=1,2$, yields positive real values for $X_2$ as well. Approximations for the corresponding values of $X_1$ and $X_2$ are given below:

$$\bar{ccc}i & X_1 & X_2 \\

     1 & 0.5526 & 3.6958\\

     2 & 0.7432 & 4.7185 \ear$$

\li

Hence, there are on $M$ exactly two Einstein adapted metrics.

\textbf{\ref{cpe65}, for $p=1$)} In this case the non-zero
solutions of the equations (\ref{eqII112}) and (\ref{eqII122}) are
given by

$$X_1=\al_i,\,X_2=-\frac{156}{7}\al_i^3+\frac{552}{7}\al_i^2-\frac{571}{7}\al_i+\frac{176}{7},$$

where $\al_i$ is a root of the polynomial $$t(z)=234z^4-828z^3+993z^2-474z+77.$$ The polynomial $t$ has exactly four positive roots which are indicated below:

$$\bar{c}\al_1=\frac{23}{26}-\frac{\sqrt{2}\be}{156}-\frac{1}{156}\left(\frac{-3664\xi\be+26\xi^2\be+71786\be+26460\sqrt{2}\xi}{\xi\be}\right)^{\frac{1}{2}}\\

\al_2=\frac{23}{26}-\frac{\sqrt{2}\be}{156}+\frac{1}{156}\left(\frac{-3664\xi\be+26\xi^2\be+71786\be+26460\sqrt{2}\xi}{\xi\be}\right)^{\frac{1}{2}}\\

\al_3=\frac{23}{26}+\frac{\sqrt{2}\be}{156}-\frac{1}{156}\left(\frac{-3664\xi\be+26\xi^2\be+71786\be+26460\sqrt{2}\xi}{\xi\be}\right)^{\frac{1}{2}}\\

\al_4=\frac{23}{26}+\frac{\sqrt{2}\be}{156}+\frac{1}{156}\left(\frac{-3664\xi\be+26\xi^2\be+71786\be+26460\sqrt{2}\xi}{\xi\be}\right)^{\frac{1}{2}}\\
\ear$$

where $\xi=(136819+36i\sqrt{1796295})^{\frac{1}{3}}$ and $\be=\left(\frac{13\xi^2+916\xi+35893}{\xi}\right)^{\frac{1}{2}}$.

Simple calculations show that $\al_i$, for $i=1,\ldots,4$, yields positive real values for $X_1$ and $X_2$ whose approximations are given below:

$$\bar{ccc}i & X_1 & X_2 \\

     1 & 0.3702 & 4.6215\\

     2 & 0.5345 & 0.6682\\

     3 & 1.0499 & 0.6338\\

     4 & 1.5838 & 5.2195 \ear$$

\li

Hence, there exist exactly four Einstein adapted metrics on $M$.

\li

\textbf{\ref{cpe75}, for $p=2$)}  For the bisymmetric triple
\ref{cpe75}, in the case $p=2$, the non-zero solutions of the
equations (\ref{eqII112}) and (\ref{eqII122}) are of the form

$$X_1=\frac{1}{2}\al_i,\,X_2=-\frac{140}{3}\al_i^3+148\al_i^2-\frac{681}{5}\al_i+\frac{184}{5},$$

where $\al_i$ is a root of the polynomial $$t(z)=350z^4-1110z^3+1179z^2-492z+69.$$ The polynomial $t$ has exactly four positive real roots which are

$$\bar{c}\al_1=\frac{111}{140}-\frac{\be}{140}-\frac{1}{140}\left(\frac{2634\xi\be-14\xi^2\be-64526\be-35250\xi}{\xi\be}\right)^{\frac{1}{2}}\\

\al_2=\frac{111}{140}-\frac{\be}{140}+\frac{1}{140}\left(\frac{2634\xi\be-14\xi^2\be-64526\be-35250\xi}{\xi\be}\right)^{\frac{1}{2}}\\

\al_3=\frac{111}{140}+\frac{\be}{140}-\frac{1}{140}\left(\frac{2634\xi\be-14\xi^2\be-64526\be-35250\xi}{\xi\be}\right)^{\frac{1}{2}}\\

\al_4=\frac{111}{140}+\frac{\be}{140}+\frac{1}{140}\left(\frac{2634\xi\be-14\xi^2\be-64526\be-35250\xi}{\xi\be}\right)^{\frac{1}{2}}\\
\ear$$

where $\xi=(290727+500i\sqrt{53545})^{\frac{1}{3}}$ and $\be=\left(\frac{14\xi^2+1317\xi+64526}{\xi}\right)^{\frac{1}{2}}$.

Simple calculations show that $\al_i$, $i=1,\ldots,4$, yields positive real values of $X_1$ and $X_2$ whose approximations are given below:

$$\bar{ccc}i & X_1 & X_2 \\

     1 & 0.3086 & 7.4890\\

     2 & 0.4686 & 0.6737\\

     3 & 0.9326 & 0.6496\\

     4 & 1.4616 & 8.1878 \ear$$

\li

Hence, there are on $M$ exactly four Einstein adapted metrics.

\li

\textbf{\ref{cpe75}, for $p=4$)} For the bisymmetric triple
\ref{cpe75}, in the case $p=2$, the non-zero solutions of the
equations (\ref{eqII112}) and (\ref{eqII122}) are given by

$$X_1=\al_i,\,X_2=-\frac{100}{3}\al_i^3+100\al_i^2-\frac{262}{3}\al_i+26,$$

where $\al_i$ is a root of the polynomial $$t(z)=200z^4-600z^3+614z^2-264z+39.$$

Since $X_1=\al_i$, only positive roots of $t$ yield positive solutions of the equations above. The polynomial $t$ has exactly two positive roots which are indicated below:

$$\bar{c}\al_1=\frac{3\xi\be}{4}+\frac{\be^3}{60}-\frac{\sqrt{3}}{60}\frac{122\xi^2\be^2+\xi^4\be^2-1151\be^2+1620\xi^3}{\xi\be}\\

\al_2=\frac{3\xi\be}{4}+\frac{\be^3}{60}+\frac{\sqrt{3}}{60}\frac{122\xi^2\be^2+\xi^4\be^2-1151\be^2+1620\xi^3}{\xi\be}\\
\ear$$

where $\xi=(109457+180\sqrt{416842})^{\frac{1}{6}}$ and $\be=(14\xi^2+1317\xi+64526)^{\frac{1}{4}}$.

Simple calculations show that $\al_i$, $i=1,\ldots,4$, yields positive real values of $X_1$ and $X_2$ whose approximations are given below:

$$\bar{ccc}i & X_1 & X_2 \\

     1 & 0.3143 & 7.3931\\

     2 & 1.4375 & 8.0839\ear$$
\li

Therefore, there are on $M$ exactly two Einstein adapted metrics.

\li

\textbf{\ref{cpe75}, for $p=6$)} For the bisymmetric triple
\ref{cpe75}, in the case $p=6$, the non-zero solutions of the
equations (\ref{eqII112}) and (\ref{eqII122}) are of the form

$$X_1=3\al_i,\,X_2=-\frac{2500}{3}z^3+820z^2-235z+24,$$

where $\al_i$ is a root of the polynomial $$t(z)=1250z^4-1230z^3+415z^2-60z+3.$$

This polynomial has exactly two positive roots which are

$$\bar{c}\al_1=\frac{41\times3^{\frac{3}{4}}\xi\be}{1500}+\frac{5^{\frac{2}{3}}3^{\frac{3}{4}}\be^3}{22500}-\frac{5^{\frac{5}{6}}3^{\frac{3}{4}}\sqrt{6}}{22500}\frac{(3887\sqrt{3}5^{\frac{1}{3}}\xi^2\be^2+125\sqrt{3}\xi^4\be^2-20875\sqrt{3}5^{\frac{2}{3}}\be^2+527553\times 5^{\frac{2}{3}}\xi^{3})^{\frac{1}{2}}}{\xi\be} \\

\al_2=\frac{41\times3^{\frac{3}{4}}\xi\be}{1500}+\frac{5^{\frac{2}{3}}3^{\frac{3}{4}}\be^3}{22500}+\frac{5^{\frac{5}{6}}3^{\frac{3}{4}}\sqrt{6}}{22500}\frac{(3887\sqrt{3}5^{\frac{1}{3}}\xi^2\be^2+125\sqrt{3}\xi^4\be^2-20875\sqrt{3}5^{\frac{2}{3}}\be^2+527553\times 5^{\frac{2}{3}}\xi^{3})^{\frac{1}{2}}}{\xi\be} \\
\ear$$

where $\xi=(14027+18\sqrt{2}\sqrt{483323})^{\frac{1}{6}}$ and $\be=(3887\times 5^{\frac{2}{3}}\xi^2-250\times 5^{\frac{1}{3}}\xi^4+208750)^{\frac{1}{4}}$.

Simple calculations show that $\al_i$, $i=1,\ldots,2$, yields positive real values of $X_1$ and $X_2$ whose approximations are given below:

$$\bar{ccc}i & X_1 & X_2 \\

     1 & 0.3163 & 7.3606\\

     2 & 1.4292 & 8.0485\ear$$
\li

Hence, there are on $M$ exactly two Einstein adapted metrics.

\li

\textbf{\ref{cpe89})} In this case there is no Einstein adapted metric on $M$. We get that, in particular, $X_2$ would be a root of the polynomial

$$t(z)=9z^4-195z^3+1198z^2-1395z+464,$$

but $t$ does not admit any positive root.

$\Box$\eproof

In the classical case, similar methods as those briefly exposed in
the proof above may be attempted to obtain solutions for the
Einstein equations. However, as it can be read from Table
\ref{eigIIclass} the eigenvalues depend on parameters in general
which would retrieve very complicated equations. Though for the
bisymmetric triples which satisfy one of the conditions
$\ga_1=\ga_2$ or $\ga_2=1-\ga_1$ it is possible to classify all
the Einstein adapted metrics by using Corollaries \ref{contcor1}
and \ref{contcor2}, respectively. For these spaces, the Einstein
adapted metrics $g_M$ whose restriction $g_F$ is also Einstein
were classified in \ref{gFeinII}. So it remains to obtain all the
other possible Einstein adapted metrics.

\bthm Let $F\rightarrow M\rightarrow N$ be a bisymmetric fibration
of Type II such that $\ga_2=\ga_1$ or $\ga_2=1-\ga_1$. If there
exists on $M$ an Einstein adapted metric such that $g_F$ is not
Einstein, then the corresponding bisymmetric triple
$(\mfg,\mfk,\mfl)$ is one of the triples in Table
\ref{nonbimII}.\ethm

\bproof In the only case when $\ga_2=1-\ga_1$, \ref{cpan3} for
$p=2l$ and $n=2(l+s)$, we have $\ga_1=\frac{l}{l+s}$,
$\ga_2=\frac{s}{l+s}$, $b_1=\frac{l}{4(l+s)}=\frac{\ga_1}{4}$ and
$b_2=\frac{s}{4(l+s)}=\frac{1-\ga_1}{4}$. The required metric
should be given by Corollary \ref{contcor2} (ii). Simple
calculation show that $D(\ga_1)=-\frac{1}{2}\ga_1(1-\ga_1)<0$, for
every $l,s$, since $0<\ga_1<1$. Hence in this case there are no
other Einstein adapted metrics besides those found previously.

\li

The cases such that $\ga_2=\ga_1$ are those listed in the proof of
Theorem \ref{biII}. In this case, there exists an Einstein adapted
metric such that $g_F$ is not Einstein if and only if
$D(\ga_1)\geq 0$, where

$$D(\ga_1)=4r^2(1-\ga_1)-2\ga_1(2b_2+1-\ga_1)(2b_1+1-\ga_1),$$

according to Corollary \ref{contcor1} (ii).

\li

In the cases \ref{cpan3} for $s=l=2p$, $n=4l$, \ref{cpdn7} for
$n=2p$, \ref{cpcn7}, for $n=2p$, as in the proof of \ref{bimII} we
have $b_1=b_2=\frac{\ga}{4}$. Hence, we simplify the expression
for $D(\ga_1)$ given above as

\beq D(\ga_1)=\frac{1}{2}(-\ga_1^3+4\ga_1^2-6\ga_1+2)\eeq

For \ref{cpan3}, $s=l=2p,\,n=4l$, we have $\ga_1=\frac{1}{2}$ and
$D(\frac{1}{2})<0$; for \ref{cpdn7}, $n=2p,\, p\geq 2$, we have
$\ga_1=\frac{p-1}{2p-1}$ and $D(\ga_1)\geq 0$ only for
$p=2,\ldots,6$; for \ref{cpcn7}, $n=2p$, $\ga=\frac{p+1}{2p+1}$
and $D<0$, for every $p\geq 1$.

\li

For the case \ref{cpdn4}, with $s=l=2p$, $n=4l$, $l\geq 2$, as in
the proof of \ref{bimII}, we have $\ga_1=\frac{2l-1}{4l-1}$ and
$b_1=b_2=\frac{l}{2(4l-1)}=\frac{1-\ga_1}{4}$. For this case, we
have

$$D(\ga_1)=-\frac{1}{2}(\ga_1-1)(3\ga_1-1)(3\ga_1-2).$$

Since $\ga\in (\frac{3}{7},\frac{1}{2})$,  we have $D(\ga_1)<0$,
for every $\ga$ and thus there is no positive solution.

\li

For \ref{cpdn8}, where $p=2l$, $n=4l$, we have
$\ga_1=\frac{2l-1}{4l-1}$ and
$b_1=\frac{l}{2(4l-1)}=\frac{1-\ga_1}{4}$,
$b_2=\frac{2l-1}{4(4l-1)}=\frac{\ga_1}{4}$. We rewrite $D(\ga_1)$
as follows:

$$ D(\ga_1)=\frac{1}{2}(1-\ga_1)(3\ga_1^2-6\ga_1+2).$$

Since $\ga_1=\frac{2l-1}{4l-1}\in(\frac{1}{3},\frac{1}{2})$,
simple calculations show that $3\ga_1^2-6\ga_1+2\geq 0$ only for
$l=1$. So only for $l=1$ exists a metric with the desired
properties.

\li

For \ref{cpcn4}, with $s=l=2p$, $n=4l$, we have
$\ga_1=\frac{2l+1}{4l+1}$ and
$b_1=b_2=\frac{l}{4(4l+1)}=\frac{1-\ga_1}{8}$. Hence

$$ D(\ga_1)=\frac{1}{8}(1-\ga_1)(25\ga_1^2-25\ga_1+8).$$

As $25\ga_1^2-25\ga_1+8>0$, for every $\ga_1$, there exists a
metric with the required properties for every $l\geq 1$.

\li

In the case \ref{cpcn8}, for $p=2l$, $n=4l$, we have
$\ga_1=\frac{2l+1}{4l+1}$ and
$b_1=\frac{l}{4(4l+1)}=\frac{1-\ga_1}{8}$,
$b_2=\frac{2l+1}{4(4l+1)}=\frac{\ga_1}{4}$. Thus

$$D(\ga_1)=\frac{1}{4}(1-\ga_1)(5\ga_1^2-10\ga_1+4).$$

Since $\ga_1\in(\frac{1}{2},\frac{3}{5})$, we conclude that
$5\ga_1^2-10\ga_1+4\geq 0$, for every $l\geq 3$. Hence, there
exists an adapted Einstein metric for every $l\geq 3$.

$\Box$\eproof

\brem For the triples in Lemma \ref{1testeII}, where $C_{\mfp_i}$
is scalar on $\mfn$, such that $D(\ga_1)<0$, we can still consider
the non-real complex solutions $X_1,X_2$ which gives rise to a non
binormal Einstein adapted metric on $M$ with non-real complex
coefficients. The spaces in these conditions are the cases
\ref{cpan3} for solutions of the form $X_2=\frac{1}{2X_1}$, and
for solutions of the form $X_2=\frac{\ga X_1}{1-\ga}$ we have the
cases \ref{cpan3}, $s=p=2l,\,n=4l$, for every $l$, \ref{cpdn7},
$n=2p$, for $p\geq 5$, \ref{cpcn7}, $n=2p$, for every $p$,
\ref{cpdn4}, with $s=p=2l$, $n =4l$, for $l\geq 2$, \ref{cpdn8},
$p=2l$, $n=4l$, for $l\geq 2$, \ref{cpcn8}, $p=2l$, $n=4l$, for
$l=1,2$.

\erem

\section{Application to $4$-symmetric Spaces}

A homogeneous space $G/L$ is said to be a $4$-symmetric space if
there exists $\sg\in Aut(G)$ such that

$$(G_{\sg})_0\subset L\subset G_{\sg}$$

and $\sg$ has order $4$. Compact simply connected irreducible
$4$-symmetric spaces have been classified by J.A.Jimenez in
\cite{Ji} following the previous work of V.Ka\v{c} (see e.g.
\cite{He}, Chap.X), J.A Wolf and A.Gray \cite{WG}. It is shown in
\cite{Ji} that any compact simply connected irreducible $4$-symmetric space is the
total space of a fiber bundle whose fiber and base space are
symmetric spaces and the base is an isotropy irreducible space of
maximal rank. These spaces are fully described in Tables III, IV and
V in \cite{Ji}. Hence, for each compact simply connected irreducible
$4$-symmetric space $M$ there is a bisymmetric fibration
$F\rightarrow M \rightarrow N$ of maximal rank whose base space $N$
is isotropy irreducible. The bisymmetric triples $(\mfg,\mfk,\mfl)$
corresponding to $4$-symmetric spaces of maximal rank must be some
of Tables \ref{eigIexc}, \ref{eigIclass}, \ref{eigIIexc} and
\ref{eigIIclass}. Hence, a simple comparison between the Tables in
this chapter and the classification in \cite{Ji} allow us to easily
conclude about the existence of Einstein metrics on $4$-symmetric
spaces.

From Theorem \ref{mtypeI} and Tables \ref{mIexc} and \ref{mIclass} we conclude the following:

\bcor Let $M=G/L$ be a compact simply connected irreducible $4$-symmetric spaces of
Type I and $(\mfg,\mfk,\mfl)$ a bisymmetric triple corresponding to $M$ such that $L\subset K$. If $M$ admits an Einstein adapted metric, then $(\mfg,\mfk,\mfl)$ is one of the triples listed below.

(i) In the exceptional case:

$$\bar{l}(\mff_4,\mfso_9,\mfso_7\oplus\reals),\\

(\mff_4, \mfsp_3\oplus\mfsu_2,\mfsp_3\oplus\reals),\\

(\mfg_2,\mfsu_2\oplus\mfsu_2, \mfsu_2\oplus\reals), (\mfg_2,\mfsu_2\oplus\mfsu_2, \reals\oplus\mfsu_2),\\

(\mfe_8,\mfso_{16},\mfso_{2p}\oplus\mfso_{16-2p}),\,p=1,3,\\

(\mfe_7,\mfso_{12}\oplus\mfsu_2, \mfso_{12}\oplus\reals),\\

(\mfe_6,\mfso_{10}\oplus\reals, \mfso_8\oplus\reals\oplus\reals),\\

(\mfe_6,\mfsu_6\oplus\mfsu_2, \mfsu_6\oplus\reals),\\

(\mfe_6, \mfsu_6\oplus\mfsu_2, \mfsu_5\oplus\reals\oplus\mfsu_2),\ear$$

(ii) in the classical case:

$$\bar{l}(\mfso_{2n+1},\mfso_{2p+1}\oplus\mfso_{2(n-p)}, \mfso_{2p+1}\oplus\mfu_{n-p}),\\

(\mfso_{2n},\mfso_{2p}\oplus \mfso_{2(n-p)},\mfu_p\oplus\mfso_{2(n-p)}),\\

(\mfsp_{n},\mfsp_{p}\oplus \mfsp_{n-p},\mfu_p\oplus\mfsp_{n-p}).\ear$$\ecor

\li

From Theorem \ref{gFeinII} and Tables \ref{bimII} and \ref{nonbimII} we obtain the following results for Type II.

\bcor Let $M=G/L$ be a compact simply connected irreducible $4$-symmetric spaces of
Type II and $(\mfg,\mfk,\mfl)$ a bisymmetric triple corresponding to $M$ such that $L\subset K$. If $M$ admits an Einstein adapted metric $g_M$ such that $g_F$ is also Einstein, then  $(\mfg,\mfk,\mfl)$ is either

$$(\mfso_{8l},\mfso_{4l}\oplus\mfso_{4l}, \mfso_{2l}\oplus\mfso_{2l}\oplus\mfso_{2l}\oplus\mfso_{2l})$$

or

$$(\mfsu_{2(l+s)},\mfsu_{2l}\oplus\mfsu_{2s}\oplus\reals, \mfsu_{l}\oplus\mfsu_{l}\oplus\mfsu_{s}\oplus\mfsu_{s}\oplus\reals\oplus\reals\oplus\reals).$$

The Einstein metric is binormal in the first
case and in the second for $l=s$.\ecor

Finally, the result below follows from Theorem \ref{genII} and Table \ref{tabgenII}.

\bcor Let $M=G/L$ be a compact simply connected irreducible $4$-symmetric spaces of
Type II, such that $G$ is an exceptional Lie group, and $(\mfg,\mfk,\mfl)$ a bisymmetric triple corresponding to $M$ such that $L\subset K$. If $M$ admits an Einstein adapted metric, then $(\mfg,\mfk,\mfl)$ is one of the following three cases:

$$\bar{l}(\mfe_6,\mfsu_6\oplus\mfsu_2, \mfsu_5\oplus\reals\oplus\reals)\\

(\mfe_7,\mfso_{12}\oplus\mfsu_2,\mfso_{10}\oplus\reals\oplus\reals)\\

(\mfe_7,\mfso_{12}\oplus\mfsu_2,\mfso_{6}\oplus\mfso_6\oplus\reals).\ear$$

There are 4 Einstein adapted metrics for each of the first two
cases and 2 for the last case. None of these metrics is binormal
or such that the restriction to the fiber is Einstein.\ecor

\brem We observe that for the $4$-symmetric spaces corresponding to the first two cases the base space considered by Jimenez is different from the one indicated above. In this two cases, the bisymmetric triples considered in \cite{Ji} are

$$\bar{l}(\mfe_6,\mfso_{10}\oplus\reals,\mfsu_5\oplus\reals\oplus\reals)\\

(\mfe_7,\mfe_6\oplus\reals,\mfso_{10}\oplus\reals\oplus\reals).\ear$$

For these triples there is no Einstein adapted metric since the Casimir operator of $\mfp$ is not scalar.\erem

\newpage

\section{Tables}\label{alltabs}

\btab[h]\caption{Dual Coxeter Numbers}\label{tabcoxeter}

$$\bar{cc}

\textrm{Coxeter group} & \textrm{Dual Coxeter number}\\\hline

A_n & n+1\\

B_n & 2n-1\\

C_n & n+1\\

D_n & 2n-2\\

E_6 & 12\\

E_7 &  18\\

E_8 & 30\\

F_4 & 9\\

G_2 & 4

\ear$$ \etab

\btab[h]\caption{Symmetric pairs of compact type of maximal rank -
Exceptional Spaces}\label{spexc}
$$\bar{|l|l||l|l|} \hline \mfg & \mfk & \mfg & \mfk
\\\hline

\xstrut\mff_4 & \mfsp_3\oplus\mfsu_2 & \mfe_7 & \mfsu_8\\

\xstrut\mff_4 & \mfso_9 & \mfe_7 & \mfe_6\oplus\reals\\

\xstrut\mfg_2 & \mfsu_2\oplus \mfsu_2 & \mfe_7 & \mfso_{12}\oplus \mfsu_2\\

\xstrut\mfe_6 & \mfso_{10}\oplus \reals & \mfe_8 & \mfso_{16}\\

\xstrut\mfe_6 & \mfsu_6\oplus\mfsu_2 & \mfe_8 & \mfe_7\oplus
\mfsu_2\\ \hline\ear$$\etab

\btab[h]\caption{Symmetric pairs of compact type of maximal rank -
Classical Spaces}\label{spclass}
$$\bar{|l|l|} \hline \mfg & \mfk
\\\hline

\xstrut\mfso_{2n} & \mfu_n\\

\xstrut\mfso_n & \mfso_{2p}\oplus\mfso_{n-2p}\\

\xstrut\mfsp_{n} & \mfu_n\\

\xstrut\mfsp_{n} & \mfsp_p\oplus\mfsp_{n-p}\\

\xstrut\mfsu_n &
\mfsu_p\oplus\mfsu_{n-p}\oplus\reals\\\hline\ear$$\etab

\btab[c]\caption{Bisymmetric triples of type I and their eigenvalues
- Exceptional spaces}\label{eigIexc}
$$\bar{|c|lll|c|c|}\hline \xstrut Appendix & \mfg & \mfk  & \mfl & \ga & b^{\phi}\\ \hline

\xstrut\ref{cpf41} & \mff_4 & \mfso_9 &
\mfso_p\oplus\mfso_{9-p},\,p=1,3,5,7 & \frac{7}{9} &
\frac{p(9-p)}{72}\\ \hline

\xstrut\ref{cpf42} &\mff_4 & \mfsp_3\oplus\mfsu_2 &
\mfsp_3\oplus\reals &
\frac{2}{9} & \frac{1}{18}\\

\xstrut\ref{cpf43} & & & \mfu_3\oplus\mfsu_2 & \frac{4}{9} &
\frac{1}{4},\,\frac{2}{9}\\

\xstrut\ref{cpf44} & & & \mfsp_2\oplus\mfsu_2\oplus \mfsu_2 &
\frac{4}{9} & \frac{1}{9},\,\frac{1}{18}\\\hline

\xstrut\ref{cpg21} & \mfg_2 & \mfsu_2\oplus \mfsu_2 &
\reals\oplus\mfsu_2 &
\frac{1}{2} & \frac{1}{8}\\

\xstrut\ref{cpg22} & &  & \mfsu_2\oplus\reals  & \frac{1}{6} &
\frac{1}{6}\\\hline

\xstrut\ref{cpe81} & \mfe_8 & \mfso_{16} &
\mfso_{2p}\oplus\mfso_{16-2p},\,p=1,\ldots,4 & \frac{1}{5} &
\frac{p(8-p)}{60}\\

\xstrut\ref{cpe82} & & & \mfu_8 & \frac{1}{5} &
\frac{4}{15},\,\frac{3}{15},\,\frac{7}{15} \\\hline

\xstrut\ref{cpe83} & \mfe_8 & \mfe_7\oplus \mfsu_2 & \mfe_7\oplus \reals & \frac{1}{15} & \frac{1}{60}\\

\xstrut\ref{cpe84} & & & \mfe_6\oplus\reals\oplus \mfsu_2 &
\frac{3}{5} &
\frac{11}{60},\frac{9}{20}\\

\xstrut\ref{cpe86} & & & \mfso_{12}\oplus\mfsu_2\oplus\mfsu_2 &
\frac{3}{5}
& \frac{4}{15},\frac{1}{5}\\

\xstrut\ref{cpe88} & & & \mfsu_8\oplus \mfsu_2 & \frac{3}{5} &
\frac{1}{4}\\\hline

\xstrut\ref{cpe71} & \mfe_7 & \mfso_{12}\oplus \mfsu_2 &
\mfso_{12}\oplus\reals & \frac{1}{9} & \frac{1}{36}\\

\xstrut\ref{cpe72} & & & \mfu_6\oplus \mfsu_2 & \frac{5}{9} &
\frac{1}{6},\,\frac{5}{18}\\

\xstrut\ref{cpe74} & & & \mfso_p\oplus\mfso_{12-p}\oplus \mfsu_2,
\,p=2,4,6 & \frac{5}{9} &  \frac{p(12-p)}{144}\\\hline

\xstrut\ref{cpe76} & \mfe_7 & \mfe_6\oplus\reals &
\mfso_{10}\oplus\reals\oplus\reals & \frac{2}{3}
&\frac{2}{9},\frac{1}{6},\,\frac{4}{9}\\

\xstrut\ref{cpe77} &  &  & \mfsu_6\oplus\mfsu_2\oplus\reals &
\frac{2}{3} & \frac{5}{18},\frac{2}{9}\\\hline

\xstrut\ref{cpe78} & \mfe_7 & \mfsu_8 &
\mfsu_p\oplus\mfsu_{8-p}\oplus\reals,\,1\leq p\leq 4 & \frac{4}{9} &
\bar{cc} \xstrut\frac{1}{9}, & p=1\\ \xstrut\frac{2}{9},\,\frac{1}{6}, & p=2\\
\xstrut\frac{2}{9},\,\frac{1}{3}, & p=3\\ \xstrut\frac{2}{9},
\frac{4}{9},\,\frac{11}{36}, & p=4 \ear
\\\hline

\xstrut\ref{cpe61} & \mfe_6 & \mfso_{10}\oplus \reals &
\mfu_5\oplus\reals
& \frac{2}{3} & \frac{5}{12},\,\frac{1}{6},\,\frac{1}{4}\\

\xstrut\ref{cpe62} & & & \mfso_p\oplus\mfso_{10-p}\oplus\reals,
\,p=2,4 & \frac{2}{3} & \frac{p(10-p)}{96}
\\\hline

\xstrut\ref{cpe63} & \mfe_6 & \mfsu_6\oplus\mfsu_2 &
\mfsu_6\oplus\reals &
\frac{1}{6} & \frac{1}{24}\\

\xstrut\ref{cpe64} & & &
\mfsu_p\oplus\mfsu_{6-p}\oplus\reals\oplus\mfsu_2 & \frac{1}{2} &
\frac{p+2}{24},\,\frac{p}{8}\\ \hline \ear$$\etab

\li

\bland\btab[c]\caption{Bisymmetric triples of type I and their
eigenvalues - Classical spaces}\label{eigIclass}
$$\bar{|c|lll|c|c|}\hline Appendix & \mfg & \mfk  & \mfl & \ga & b^{\phi}\\ \hline

\xstrut\ref{cpan1} & \mfsu_n & \mfsu_p\oplus\mfsu_{n-p}\oplus\reals
& \mfsu_l\oplus\mfsu_{p-l}\oplus\reals\oplus\mfsu_{n-p}\oplus\reals
& \frac{p}{n} & \frac{p-l}{2n},\,\frac{l}{2n}\\ \hline

\xstrut\ref{cpbn1} & \mfso_{2n+1} &
\mfso_{2p+1}\oplus\mfso_{2(n-p)},\,p=0,\ldots,n-1 &
\mfso_{2l+1}\oplus\mfso_{2(p-l)}\oplus\mfso_{2(n-p)} &
\frac{2p-1}{2n-1} & \frac{p-l}{2n-1},\,\frac{4l+1}{4(2n-1)}\\

\xstrut\ref{cpbn2} & & &
\mfso_{2p+1}\oplus\mfso_{2s}\oplus\mfso_{2(n-p-s)}
& \frac{2(n-p-1)}{2n-1} & \frac{n-p-s}{2n-1},\,\frac{s}{2n-1}\\

\xstrut\ref{cpbn3} & & & \mfso_{2p+1}\oplus\mfu_{n-p} &
\frac{2(n-p-1)}{2n-1} & \frac{n-p-1}{2(2n-1)}\\ \hline

\xstrut\ref{cpdn1} & \mfso_{2n} & \mfu_n &
\mfu_p\oplus\mfu_{n-p},\,p=1,\ldots,n-1 & \frac{n}{2(n-1)} &
\frac{n-p}{2(n-1)},\,\frac{p}{2(n-1)},\,\frac{n-2}{4(n-1)}\\ \hline

\xstrut\ref{cpdn2} & \mfso_{2n} &
\mfso_{2p}\oplus\mfso_{2(n-p)},\,p=1,\ldots,\lfloor
\frac{n}{2}\rfloor  &
\mfso_{2l}\oplus\mfso_{2(p-l)}\oplus\mfso_{2(n-p)}& \frac{p-1}{n-1}
& \frac{p-l}{2(n-1)},\,\frac{l}{2(n-1)}\\

\xstrut\ref{cpdn5} & & & \mfu_p\oplus\mfso_{2(n-p)} &
\frac{p-1}{n-1} & \frac{p-1}{4(n-1)}\\ \hline

\xstrut\ref{cpcn1} & \mfsp_{n} & \mfu_n &
\mfu_p\oplus\mfu_{n-p},\,p=1,\ldots,n-1 & \frac{n}{2(n+1)} &
\frac{n-p}{2(n+1)},\,\frac{n-p}{n+1},\,\frac{p}{2(n+1)},\,\frac{p}{n+1},\,\frac{n+2}{2(n+1)}
 \\ \hline

\xstrut\ref{cpcn2} & \mfsp_n & \mfsp_p\oplus\mfsp_{n-p} &
\mfsp_l\oplus\mfsp_{p-l}\oplus\mfsp_{n-p} & \frac{p+1}{n+1} &
\frac{p-l}{4(n+1)},\,\frac{l}{4(n+1)}\\

\xstrut\ref{cpcn5} &  &  & \mfu_p\oplus\mfsp_{n-p} & \frac{p+1}{n+1}
& \frac{p+1}{4(n+1)}\\ \hline \ear$$\etab \eland

\bland\btab[c]\caption{Bisymmetric triples of type II and their
eigenvalues - Exceptional
spaces}\label{eigIIexc}$$\bar{|c|lll|c|c|c|c|}\hline Appendix &\mfg &
\mfk  & \mfl & \ga_1 & \ga_2 & b_1^{\phi} & b_2^{\phi} \\\hline

\xstrut\ref{cpf45} & \mff_4 & \mfsp_3\oplus\mfsu_2 & \mfu_3\oplus
\reals & \frac{4}{9} & \frac{2}{9} &
\big(\frac{1}{4},\frac{2}{9}\big) &
\frac{1}{18}\\

\xstrut\ref{cpf46} & & & \mfsu_2\oplus\mfsp_2\oplus\reals &
\frac{4}{9} &
\frac{2}{9} & \big(\frac{1}{9},\frac{1}{18}\big) & \frac{1}{18}\ \\
\hline

\xstrut\ref{cpg23} & \mfg_2 & \mfsu_2\oplus \mfsu_2  &
\reals\oplus\reals & \frac{1}{2} & \frac{1}{6} & \frac{1}{8} &
\frac{1}{6}\\\hline

\xstrut\ref{cpe85} & \mfe_8 & \mfe_7\oplus \mfsu_2 &
\mfe_6\oplus\reals\oplus\reals & \frac{3}{5} & \frac{1}{15} &
\big(\frac{11}{60},\frac{11}{60},\frac{9}{20}\big) & \frac{1}{60}\\

\xstrut\ref{cpe87} & & & \mfso_{12}\oplus\mfsu_2\oplus\reals &
\frac{3}{5} & \frac{1}{15} & \big(\frac{4}{15},\frac{1}{5}\big)  &
\frac{1}{60}
\\

\xstrut\ref{cpe89} & & & \mfsu_8\oplus\reals & \frac{3}{5} &
\frac{1}{15} & \frac{1}{4}& \frac{1}{60} \\ \hline

\xstrut\ref{cpe73} & \mfe_7 & \mfso_{12}\oplus \mfsu_2 &
\mfu_6\oplus \reals & \frac{5}{9} & \frac{1}{9} &
\big(\frac{5}{18},\frac{1}{6},\frac{5}{18}\big) & \frac{1}{36}\\

\xstrut\ref{cpe75} & & & \mfso_p\oplus\mfso_{12-p}\oplus \reals,
\,p\textrm { even} & \frac{5}{9} & \frac{1}{9} & \frac{1}{36}&
\frac{p(12-p)}{144}\\ \hline

\xstrut\ref{cpe65} & \mfe_6 & \mfsu_6\oplus\mfsu_2 &
\mfsu_p\oplus\mfsu_{6-p}\oplus \reals\oplus\reals & \frac{1}{2} &
\frac{1}{6}
 & \frac{1}{24} & \frac{p+2}{24},\, \frac{p}{8}\\ \hline
\ear$$\etab\eland

\bland \btab[c]\caption{Bisymmetric triples of type II and their
eigenvalues -
Classical spaces}\label{eigIIclass}$$\bar{|c|lll|c|c|c|c|}\hline Appendix & \mfg & \mfk  & \mfl & \ga_1 & \ga_2 & b_1^{\phi} & b_2^{\phi} \\
\hline
\xstrut\ref{cpan3} & \mfsu_n & \mfsu_p\oplus\mfsu_{n-p}\oplus\reals
&
\mfsu_l\oplus\mfsu_{p-l}\oplus\mfsu_s\oplus\mfsu_{n-p-s}\oplus\reals\oplus\reals\oplus\reals
& \frac{p}{n} & \frac{n-p}{n} &
\big(\frac{p-l}{2n},\frac{l}{2n}\big) &
\big(\frac{n-p-s}{2n},\frac{s}{2n}\big) \\ \hline

\xstrut\ref{cpbn5} & \mfso_{2n+1} & \mfso_{2p+1}\oplus
\mfso_{2(n-p)} & \mfso_{2l+1}\oplus
\mfso_{2(p-l)}\oplus\mfso_{2s}\oplus\mfso_{2(n-p-s)} &
\frac{2p-1}{2n-1} & \frac{2(n-p-1)}{2n-1} &
\big(\frac{p-l}{2n-1},\frac{4l+1}{4(2n-1)}\big) &
\big(\frac{n-p-s}{2n-1},\frac{s}{2n-1}\big)\\

\xstrut\ref{cpbn4} &  & & \mfso_{2l+1}\oplus
\mfso_{2(p-l)}\oplus\mfu_{n-p} & \frac{2p-1}{2n-1} &
\frac{2(n-p-1)}{2n-1}  &
\big(\frac{p-l}{2n-1},\frac{4l+1}{4(2n-1)}\big) &
\frac{n-p-1}{2(2n-1)}\\ \hline

\xstrut\ref{cpdn4} & \mfso_{2n} & \mfso_{2p}\oplus \mfso_{2(n-p)} &
\mfso_{2l}\oplus
\mfso_{2(p-l)}\oplus\mfso_{2s}\oplus\mfso_{2(n-p-s)} &
\frac{p-1}{n-1}& \frac{n-p-1}{n-1} &
\big(\frac{p-l}{2(n-1)},\frac{l}{2(n-1)}\big) &
\big(\frac{n-p-s}{2(n-1)},\frac{s}{2(n-1)}\big)\\

\xstrut\ref{cpdn7} &  & & \mfu_p\oplus \mfu_{n-p} & \frac{p-1}{n-1}&
\frac{n-p-1}{n-1} & \frac{p-1}{4(n-1)} & \frac{n-p-1}{4(n-1)}\\

\xstrut\ref{cpdn8} &  & & \mfso_{2l}\oplus
\mfso_{2(p-l)}\oplus\mfu_{n-p} & \frac{p-1}{n-1}& \frac{n-p-1}{n-1}
& \big(\frac{p-l}{2(n-1)},\frac{l}{2(n-1)}\big) &
\frac{n-p-1}{4(n-1)}\\ \hline

\xstrut\ref{cpcn4} & \mfsp_n & \mfsp_p\oplus \mfsp_{n-p} &
\mfsp_l\oplus \mfsp_{p-l}\oplus\mfsp_{s}\oplus\mfsp_{n-p-s} &
\frac{p+1}{n+1}& \frac{n-p+1}{n+1} &
\big(\frac{p-l}{4(n+1)},\frac{l}{4(n+1)}\big) &
 \big(\frac{n-p-s}{4(n+1)},\frac{s}{4(n+1)}\big)\\

\xstrut\ref{cpcn7} &  &  & \mfu_p\oplus \mfu_{n-p} &
\frac{p+1}{n+1}&
\frac{n-p+1}{n+1} & \frac{p+1}{4(n+1)}& \frac{n-p+1}{4(n+1)}\\

\xstrut\ref{cpcn8} &  &  & \mfsp_l\oplus \mfsp_{p-l}\oplus\mfu_{n-p}
& \frac{p+1}{n+1}& \frac{n-p+1}{n+1} &
\big(\frac{p-l}{4(n+1)},\frac{l}{4(n+1)}\big) &
\frac{n-p+1}{4(n+1)}\\ \hline

 \ear$$\etab \eland

\btab[c]\caption{Einstein Bisymmetric triples of type I -
Exceptional spaces}\label{mIexc}
$$\bar{|lll|c|}\hline\mfg & \mfk  & \mfl  & X\\ \hline

\xstrut\mff_4 & \mfsp_3\oplus\mfsu_2 & \mfsp_3\oplus\reals &
\frac{1}{2},\,4\\\hline

\xstrut\mff_4 & \mfso_9 & \mfso_8 & 1,\,\frac{2}{7}\\

\xstrut& & \mfso_7\oplus\reals & \frac{9\pm\sqrt{8}}{14}\\\hline

\xstrut\mfg_2 & \mfsu_2\oplus \mfsu_2 & \reals\oplus\mfsu_2 &
\frac{1}{2},\,\frac{3}{2}\\

\xstrut&  & \mfsu_2\oplus\reals & \frac{6\pm\sqrt{22}}{2}\\\hline

\xstrut\mfe_6 & \mfso_{10}\oplus \reals &
\mfso_{8}\oplus\reals\oplus\reals &  1,\,\frac{1}{2}
\\\hline

\xstrut\mfe_6 & \mfsu_6\oplus\mfsu_2 & \reals\oplus\mfsu_5\oplus\mfsu_2 & \frac{1}{2},\,\frac{3}{2}\\

\xstrut & & \mfsu_6\oplus\reals &
\frac{1}{2},\,\frac{11}{2}\\\hline

\xstrut\mfe_7 & \mfso_{12}\oplus \mfsu_2 & \mfso_{12}\oplus\reals  & \frac{17}{2},\,\frac{1}{2}\\

\xstrut & & \reals\oplus\mfso_{10}\oplus \mfsu_2 &
\frac{1}{2},\,\frac{13}{10}\\

\xstrut & & \mfso_4\oplus\mfso_{8}\oplus \mfsu_2  & 1,\, \frac{4}{5} \\\hline

\xstrut\mfe_7 & \mfsu_8 & \mfsu_7\oplus\reals &
\frac{1}{2},\,\frac{7}{4}\\\hline

\xstrut\mfe_8 & \mfe_7\oplus \mfsu_2 & \mfe_7\oplus \reals &
\frac{1}{2},\,\frac{29}{2}\\\hline

\xstrut\mfe_8 & \mfso_{16} & \mfso_{2p}\oplus\mfso_{16-2p} &
 \frac{15\pm\sqrt{7p^2-56p+113}}{14}  \\ \hline
 \ear$$\etab

\bland\btab[c]\caption{Einstein Bisymmetric triples of type I -
Classical spaces}\label{mIclass}
$$\bar{|lll|c|}\hline\mfg & \mfk  & \mfl & X\\ \hline

\xstrut\mfso_{2n} & \mfso_{2p}\oplus\mfso_{2(n-p)} &
\mfso_{p}\oplus\mfso_{p}\oplus\mfso_{2(n-p)},\,p\,even &
\frac{n-1\pm\sqrt{p^2-(2n+1)p+n^2+1}}{2(p-1)}\\

\xstrut& & \mfu_p\oplus\mfso_{2(n-p)} &
\frac{1}{2},\,\frac{n}{p-1}-\frac{1}{2}\\\hline

\xstrut\mfso_{2n+1} & \mfso_{2p+1}\oplus\mfso_{2(n-p)} &
\mfso_{2p+1}\oplus\mfso_{n-p}\oplus\mfso_{n-p},\,n-p\,even &
\frac{2n-1\pm\sqrt{4p^2+8p-4n+5}}{4(n-p-1)}\\

\xstrut& & \mfso_{2p+1}\oplus\mfu_{n-p} &
\frac{1}{2},\,\frac{n+p}{2(n-p-1)}\\\hline

\xstrut\mfsp_n & \mfsp_{2l}\oplus\mfsp_{n-2l} &
\mfsp_l\oplus\mfsp_l\oplus\mfsp_{n-2l} &
\frac{n+1\pm\sqrt{6l^2+(3-4n)l+n^2+1}}{2(2l+1)}\\

\xstrut &  & \mfu_p\oplus\mfsp_{n-p} & \frac{1}{2},\,\frac{1}{2}+\frac{n-p}{p+1}\\\hline

\xstrut\mfsu_n & \mfsu_{2l}\oplus\mfsu_{n-2l}\oplus\reals &
\mfsu_l\oplus\mfsu_l\oplus\reals\oplus\mfsu_{n-2l}\oplus\reals
 & \frac{1}{2},\,\frac{n}{2l}-\frac{1}{2}\\ \hline\ear$$\etab\eland

\btab[c]\caption{Bisymmetric triples of type II with Einstein
metric such that $g_F$ is also Einstein}\label{bimII}

\bc{$g_M$ binormal}\ec

$$\bar{|lll|c|}\hline\mfg & \mfk  & \mfl & X\\ \hline

\xstrut\mfsu_{4l} & \mfsu_{2l}\oplus\mfsu_{2l}\oplus\reals &
\mfsu_l\oplus\mfsu_l\oplus\oplus\mfsu_l\oplus\mfsu_l\oplus\reals^3 & 1\\

\xstrut\mfso_{8l} & \mfso_{4l}\oplus\mfso_{4l} & \mfso_{2l}
\oplus\mfso_{2l} \oplus\mfso_{2l} \oplus\mfso_{2l} &
1,\,\frac{2l}{2l-1}\\

\xstrut\mfso_{8l} & \mfso_{4l}\oplus\mfso_{4l} & \mfso_{2l}
\oplus\mfso_{2l} \oplus\mfu_{2l} &
\frac{4l-1\pm\sqrt{2l}}{2(2l-1)}\\

\xstrut\mfso_{4l} & \mfso_{2l}\oplus\mfso_{2l} & \mfu_{l}
\oplus\mfu_{l},\,l\geq 2 & \frac{2l-1\pm\sqrt{2l-1}}{2(l-1)}\\

\xstrut\mfsp_{4l} & \mfsp_{2l}\oplus\mfsp_{2l} &
\mfsp_l\oplus\mfsp_l\oplus\mfsp_l\oplus\mfsp_l &
\frac{4l+1\pm\sqrt{4l^2+2l+1}}{2(2l+1)}\\

\xstrut\mfsp_{4l} & \mfsp_{2l}\oplus\mfsp_{2l} &
\mfsp_l\oplus\mfsp_l\oplus\mfu_{2l} &
\frac{4l+1\pm\sqrt{l(2l-1)}}{2(2l+1)}
\\ \hline\ear$$

\bc{$g_M$ non-binormal}\ec

$$\bar{|lll|c|c|}\hline\mfg & \mfk  & \mfl & X_1 & X_2\\ \hline

\xstrut\mfsu_{2(l+s)} & \mfsu_{2l}\oplus\mfsu_{2s}\oplus\reals &
\mfsu_l\oplus\mfsu_l\oplus\mfsu_s\oplus\mfsu_s\oplus\reals^3 &
\frac{l+s}{2l} & \frac{l+s}{2s}\\\hline\ear$$\etab

\bland\btab[c]\caption{All other Einstein adapted metrics for the
bisymmetric triples of Type II which admit an EAM $g_M$ such that
$g_F$ is Einstein}\label{nonbimII}
$$\bar{|lll|c|c|}\hline\mfg & \mfk  & \mfl & X_1 & X_2\\ \hline

\xstrut\mfso_{4l}& \mfso_{2l}\oplus\mfso_{2l} &
\mfu_l\oplus\mfu_l,
\,l=2,\ldots,6 & \frac{2l(l-1)\pm\sqrt{(-l^4+7l^3-5l^2+l)/2}}{2(l-1)(3l-1)} & \frac{l}{2(l-1)}.\frac{1}{X_1}\\

\xstrut\mfso_{8} & \mfso_4\oplus\mfso_4 &
\reals\oplus\reals\oplus\mfu_2 & \frac{4\pm\sqrt{6}}{5}& \frac{1}{X_1}\\

\xstrut\mfsp_{4l} & \mfsp_{2l}\oplus\mfsp_{2l} &
\mfsp_l\oplus\mfsp_l\oplus\mfsp_l\oplus\mfsp_l,\,l\geq 1 & \frac{4l+1\pm\sqrt{14l^2+7l+4}}{5(2l+1)}& \frac{l}{2l+1}.\frac{1}{X_1} \\

\xstrut\mfsp_{4l} & \mfsp_{2l}\oplus\mfsp_{2l} &
\mfsp_l\oplus\mfsp_l\oplus\mfu_{2l},\,l\geq 3 & \frac{2(4l+1)\pm\sqrt{4l^2-8l-1}}{5(2l+1)}& \frac{l}{2l+1}.\frac{1}{X_1}\\

\hline\ear$$\etab\eland

\btab[c]\label{tabgenII}\caption{Einstein bisymmetric fibrations of Type II with $\mfg$ exceptional}
$$\bar{|lll|c|} \hline\mfg & \mfk & \mfl & Number\, of\, Einstein\, adapted\, metrics\\

 \hline\xstrut\mfg_2 & \mfsu_2\oplus\mfsu_2 & \reals\oplus\reals & 2\\

 \xstrut \mfe_6 & \mfsu_6\oplus\mfsu_2 & \mfsu_5\oplus\reals\oplus\reals & 4\\

 \xstrut \mfe_7 & \mfso_{12}\oplus\mfsu_2 & \reals\oplus\mfso_{10}\oplus\reals & 4\\

 \xstrut \mfe_7 & \mfso_{12}\oplus\mfsu_2& \mfso_{4}\oplus\mfso_{8}\oplus\reals & 2\\

 \xstrut \mfe_7 & \mfso_{12}\oplus\mfsu_2& \mfso_{6}\oplus\mfso_{6}\oplus\reals & 2\\ \hline\ear$$\etab

 \newpage

\newpage

\chapter{}

In this Chapter we consider a fibration
$$\frac{\triangle^p G_0\times \triangle^q G_0}{\triangle^n G_0}\rightarrow\frac{G_0^n}{\triangle^nG_0}\rightarrow
\frac{G_0^p}{\triangle^pG_0}\times \frac{G_0^q}{\triangle^qG_0},$$

where $G_0$ is compact connected simple Lie group and $\triangle^mG_0$ is the diagonal subgroup in $G_0^m$, for $m=p,q,n$. The spaces
$\frac{G_0^n}{\triangle^nG_0}$ are $n$-symmetric spaces and it has
been proved by Kowalski that under some conditions they are not
$k$-symmetric for $k<n$ (see e.g. \cite{Ko}). Hence, throughout this
Chapter we shall designate such a space by a Kowaslki $n$-symmetric
space. McKenzie Y. Wang and Wolfgang Ziller have shown that these spaces are
standard Einstein manifolds (\cite{WZ2}, \cite{Ro1}). We obtain new Einstein metrics with totally geodesic fibers. In section 1 we describe the isotropy subspaces and compute the necessary eigenvalues to obtain the Ricci curvature of an adapted metric on $M$. In Section 2 we show that, for $n>4$, there exists at least one non-standard Einstein adapted metric on $M$, which is binormal or such that the metric on the base space is also Einstein if and only if  $p=q$. We prove that for $n=4$ the standard metric is the only Einstein adapted metric.  We remark that for $n=4$ the fibration above is a
bisymmetric fibration of non-maximal rank and whose base space is
isotropy reducible in opposition to the cases studied in Chapter 3.

\section{Kowalski N-Symmetric Spaces - The Isotropy Representation and the Casimir Operators}

Let $G_0$ be a compact connected simple Lie group and $\mfg_0$ its
Lie algebra. For any positive integer $m$ we denote by $G_0^m$ (or
$\mfg_0^m$) the direct product of $G_0$ ($\mfg_0$, resp.) by itself
$m$ times. By $\triangle ^mG_0$ (or $\triangle^m \mfg_0$)  we denote
the diagonal in $G_0^m$ (in $\mfg_0^m$, resp.). Clearly, the Lie
algebras of $G_0^m$  and $\triangle ^mG_0$ are $\mfg_0^m$ and
$\triangle^m \mfg_0$, respectively.

\li

Let $n,p,q$ be positive integers such that $p+q=n$ and
$2\leq p\leq q\leq n-2$. Set $G=G_0^n$ and consider the following closed
subgroups of $G$:

$$K=\triangle^p G_0\times \triangle^q G_0,$$

$$L=\triangle ^nG_o\subset K.$$

\li

The Lie algebras of $G$, $K$ and $L$ are, respectively,

$$\mfg=\mfg_0^n,$$

$$\mfk=\triangle^p \mfg_0\times \triangle^q \mfg_0,$$

$$\mfl=\triangle^n \mfg_0.$$

\li

Following the notation of previous chapters, we write $M=G/L$, $N=G/K$ and $F=K/L$.
We consider the fibration $F \rightarrow M \rightarrow N$. We note that
$N=\frac{G_0^p}{\triangle^pG_0}\times
\frac{G_0^q}{\triangle^qG_0}$ and thus this fibration is

$$\frac{G_0^n}{\triangle^nG_0}\rightarrow \frac{G_0^p}{\triangle^pG_0}\times
\frac{G_0^q}{\triangle^qG_0}$$

with fiber $$F=\frac{\triangle^p G_0\times \triangle^q
G_0}{\triangle^nG_0}.$$

\li

Let $\Phi_0$ be the Killing form of $\mfg_0$. Then, the Killing form
of $\mfg$ is

$$\Phi=\underbrace{\Phi_0+\ldots+\Phi_0}_n.$$

Following the notation of chapters 2 and 3, let $\mfn$ be the orthogonal complement of $\mfk$ in $\mfg$ and
$\mfp$ be an orthogonal complement of $\mfl$ in $\mfk$, with respect
to $\Phi$. Then,
$$\mfg=\underbrace{\mfl\oplus\mfp}_{\mfk}\oplus\mfn$$

and we write $\mfm=\mfp\oplus\mfn$.

\li

\blem (i)
$\mfp=\{(\underbrace{qX,\ldots,qX}_p,\underbrace{-pX,\ldots-pX}_q):X\in\mfg_0
\}$ and $\mfp$ is an irreducible $Ad\,L$-module;

\li

(ii) $\mfn=\mfn_1\oplus\mfn_2$, where

$$\mfn_1=\{(X_1,\ldots,X_p,0,\ldots,0):X_j\in\mfg_0,\sum_{j=1}^pX_j=0
\}\subset \mfg_0^p\times 0_q$$

$$\mfn_2=\{(0,\ldots,0,X_1,\ldots,X_q):X_j\in\mfg_0,\sum_{j=1}^qX_j=0
\}\subset 0_p\times\mfg_0^q;$$

\elem

\bproof Let $(\underbrace{Z,\ldots,Z}_n)\in\mfl$ and
$(\underbrace{X,\ldots,X}_p,\underbrace{Y,\ldots,Y}_q)\in\mfk$,
where $X$, $Y$ and $Z$ are arbitrary elements in $\mfg_0$.

$$\bar{rl}0 = & \Phi((X,\ldots,X,Y,\ldots,Y) ,(Z,\ldots,Z))\\ \\

= & p\Phi_0(X,Z)+q\Phi_0(Y,Z)\\ \\

= & \Phi_0(pX+qY,Z)\ear$$

Since $\Phi_0$ is non-degenerate on $\mfg_0$ and $Z$ is arbitrary,
the identity above is possible if and only if $pX+qY=0$. Hence, we
conclude that $\mfp$ is formed by elements of the form

$$(\underbrace{qX,\ldots,qX}_p,\underbrace{-pX,\ldots-pX}_q),$$

where $X\in\mfg_0$. Moreover, it is clear that $Ad\,L$-submodules of
$\mfp$ correspond to ideals of $\mfg_0$. Since $\mfg_0$ is simple we
conclude that $\mfp$ is irreducible.

\li

Clearly, we may write $\mfn=\mfn_1\oplus\mfn_2$, where $\mfn_1$ is
an orthogonal complement of $\triangle^p\mfg_0\times 0$ in
$\mfg_0^p\times 0$ and $\mfn_2$ is an orthogonal complement of
$0\times \triangle^q\mfg_0$ in $0\times\mfg_0^q$, with respect to
$\Phi$.

Let $(X_1,\ldots,X_p,0,\ldots,0)\in\mfg_0^p\times 0$ and
$(Z,\ldots,Z,0,\ldots,0)\in\triangle^p\mfg_0\times 0$, where $Z$ and
$X_j$ are arbitrary elements in $\mfg_0$.

$$\bar{rl}0 = & \Phi((X_1,\ldots,X_p,0,\ldots,0) ,(Z,\ldots,Z,0,\ldots,0))\\ \\

= & \Phi_0(X_1,Z)+\ldots+\Phi_0(X_p,Z)\\ \\

= & \Phi_0(\sum_{j=1}^pX_j,Z).\ear$$

Also by the nondegeneracy of $\Phi_0$, we conclude that the identity
above holds if and only if $\sum_{j=1}^pX_j=0$. This gives the
required expression for $\mfn_1$. To prove the expression for
$\mfn_2$ is similar.

\li

$\Box$ \eproof

\li

As usual we denote the Casimir operator of a subspace $V$
with respect to $\Phi$, the Killing form of $\mfg$, by $C_{V}$ (see Definition \ref{Casimirdef}). Also, we denote the identity map and the zero map of $V$, by $Id_V$
and $0_V$, respectively.

\bprop \label{ck1}(i) $C_{\mfg}=Id_{\mfg}$;

\li

(ii) $C_{\mfl}=\dfrac{1}{n}Id_{\mfg}$;

\li

(iii) $C_{\mfp}=\dfrac{q}{np}Id_{\mfg_0^p}\times
\dfrac{p}{nq}Id_{\mfg_0^q} $;

\li

(iv) $C_{\mfk}=\dfrac{1}{p}Id_{\mfg_0^p}\times
\dfrac{1}{q}Id_{\mfg_0^q} $;

\li

(v) $C_{\mfn_1}=\left(1-\dfrac{1}{p}\right)Id_{\mfg_0^p}\times
0_{\mfg_0^q}$ and  $C_{\mfn_2}=0_{\mfg_0^p}\times
\left(1-\dfrac{1}{q}\right)Id_{\mfg_0^q}$.
 \eprop

\bproof Let $(u_i)_i$ and $(u'_i)_i$ be bases of $\mfg_0$ dual with respect to
$\Phi_0$. Then the Casimir operator of $\mfg_0$ with respect to
$\Phi_0$ is $C_{\mfg_0}=\sum_iad_{u_i}ad_{u'_i}$. We observe that
since $\mfg_0$ is a simple Lie algebra, $C_{\mfg_0}=Id_{\mfg_0}$.

\li

(i)
$$\{(u_i,0,\ldots,0),(0,u_i,\ldots,0),\ldots,(0,0,\ldots,u_i)\}_i$$
and
$$\{(u'_i,0,\ldots,0),(0,u'_i,\ldots,0),\ldots,(0,0,\ldots,u'_i)\}_i$$
are bases for $\mfg$. We have

$$\Phi((0,\ldots,\underbrace{u_i}_{k^{th}},\ldots,0),
(0,\ldots,\underbrace{u'_j}_{l^{th}},\ldots,0))=
\de_{kl}\Phi_0(u_i,u'_j)=\de_{kl}\de_{ij}.$$

Hence, the bases above are dual for $\Phi$ and thus we may write

$$\bar{rl}C_{\mfg}= & \sum_iad_{(u_i,0,\ldots,0)}ad_{(u'_i,0,\ldots,0)}+\ldots +\sum_iad_{(0,\ldots,0,u_i)}ad_{(0,\ldots,0,u'_i)}\\ \\

= & (\sum_iad_{u_i}ad_{u'_i},0,\ldots,0)+\ldots+(0,0,\ldots,\sum_iad_{u_i}ad_{u'_i})\\ \\

= & (C_{\mfg_0},\ldots,C_{\mfg_0})\\ \\

= & (Id_{\mfg_0},\ldots,Id_{\mfg_0})\\ \\

= & Id_{\mfg}.\ear$$

\li

(ii) $\{(\underbrace{u_i,\ldots,u_i}_n)\}_i$ and
$\{(\underbrace{u'_i,\ldots,u'_i}_n)\}_i$ are bases for $\mfl$ and
we have
$$\Phi((u_i,\ldots,u_i),(u'_j,\ldots,u'_j))=n\Phi_0(u_i,u'_j)=n\de_{ij}.$$

Hence, $\{\frac{1}{\sqrt{n}}(\underbrace{u_i,\ldots,u_i}_n)\}_i$
and $\{\frac{1}{\sqrt{n}}(\underbrace{u'_i,\ldots,u'_i}_n)\}_i$
are bases for $\mfl$ dual with respect to $\Phi$. So we have

$$\bar{rl}C_{\mfl}= &
\dfrac{1}{n}\sum_iad_{(u_i,\ldots,u_i)}ad_{(u'_i,\ldots,u'_i)}\\ \\

= &
\dfrac{1}{n}(\sum_iad_{u_i}ad_{u'_i},\ldots,\sum_iad_{u_i}ad_{u'_i}
)\\ \\

= & \dfrac{1}{n}(C_{\mfg_0},\ldots,C_{\mfg_0})\\ \\

= & \dfrac{1}{n}(Id_{\mfg_0},\ldots,Id_{\mfg_0})\\ \\

= & \dfrac{1}{n}Id_{\mfg}. \ear$$

\li

(iii)
$\{(\underbrace{qu_i,\ldots,qu_i}_p,\underbrace{-pu_i,\ldots,-pu_i}_q)\}_i$
and
$\{(\underbrace{qu'_i,\ldots,qu'_i}_p,\underbrace{-pu'_i,\ldots,-pu'_i}_q)\}_i$
are bases of $\mfp$.

$$\bar{rl} &  \Phi((qu_i,\ldots,qu_i,-pu_i,\ldots,-pu_i),(qu'_j,\ldots,qu'_j,-pu'_j,\ldots,-pu'_j))\\ \\

= & q^2p\Phi_0(u_i,u'_j)+p^2q\Phi_0(u_i,u'_j)\\ \\

= & (p+q)pq\de_{ij}\\ \\

= & npq \de_{ij}. \ear$$

Hence dual bases of $\mfp$ for $\Phi$ are as follows:

\bc{$\left\{\left(\left(\frac{q}{np}\right)^{\frac{1}{2}}u_i,\ldots,\left(\frac{q}{np}\right)^{\frac{1}{2}}u_i,-\left(\frac{p}{nq}\right)^{\frac{1}{2}}u_i,\ldots,-\left(\frac{p}{nq}\right)^{\frac{1}{2}}u_i\right)\right\}_i$}\ec
\bc{and
$\left\{\left(\left(\frac{q}{np}\right)^{\frac{1}{2}}u'_i,\ldots,\left(\frac{q}{np}\right)^{\frac{1}{2}}u'_i,-\left(\frac{p}{nq}\right)^{\frac{1}{2}}u'_i,\ldots,-\left(\frac{p}{nq}\right)^{\frac{1}{2}}u'_i\right)\right\}_i$}.\ec

Similar calculations as those done above show that
$$C_{\mfp}=\frac{q}{np}(C_{\mfg_0},\ldots, C_{\mfg_0})\times \frac{p}{nq}(C_{\mfg_0},\ldots,
C_{\mfg_0})=\dfrac{q}{np}Id_{\mfg_0^p}\times
\dfrac{p}{nq}Id_{\mfg_0^q}.$$

(iv) $$C_{\mfk}=
C_{\mfl}+C_{\mfp}=\left(\dfrac{1}{n}+\dfrac{q}{np}\right)Id_{\mfg_0^p}\times
\left(\dfrac{1}{n}+\dfrac{p}{nq}\right)Id_{\mfg_0^q} =
\dfrac{1}{p}Id_{\mfg_0^p}\times \dfrac{1}{q}Id_{\mfg_0^q}.$$

\li

(v) $$C_{\mfn_1}+C_{\mfn_2}=
C_{\mfn}=C_{\mfg}-C_{\mfk}=\left(1-\dfrac{1}{p}\right)Id_{\mfg_0^p}\times
\left(1-\dfrac{1}{q}\right)Id_{\mfg_0^q}$$

Since $[\mfn_1,\mfg]\subset \mfg_0^p\times 0$ and
$[\mfn_2,\mfg]\subset 0\times \mfg_0^q$, we conclude that

$$C_{\mfn_1}=\left(1-\dfrac{1}{p}\right)Id_{\mfg_0^p}\times
0_{\mfg_0^q} \textrm{ and }C_{\mfn_2}=0_{\mfg_0^p}\times
\left(1-\dfrac{1}{q}\right)Id_{\mfg_0^q}.$$

$\Box$\eproof

For the eigenvalues of the Casimir operators of $\mfl$, $\mfk$,
$\mfp$ and $\mfn$ we use notation similar to that used in previous
chapters. We recall that $c_{\mfl,\mfp}$ is the eigenvalue of
$C_{\mfl}$ on $\mfp$, $c_{\mfk,i}$ is the eigenvalue of $C_{\mfk}$
on $\mfn_i$. The Casimir operator of $\mfp$ is scalar on $\mfn_i$, as we can see from Corollary \ref{ck2}, and $b^i$ denotes the eigenvalue of $C_{\mfp}$ on $\mfn_i$,
for $i=1,2$. Also $c_{\mfn_i,\mfp}$ and $\ga$ are the constants
defined by

\beqar
\Phi(C_{\mfn_i}.,.)\mid_{\mfp\times\mfp}=c_{\mfn_i,\mfp}\Phi\mid_{\mfp\times\mfp},\,i=1,2\\
\Phi(C_{\mfk}.,.)\mid_{\mfp\times\mfp}=\ga
\Phi\mid_{\mfp\times\mfp}.\eeqar

\bcor \label{ck2}(i) $c_{\mfl,\mfp}=\dfrac{1}{n}$;

(ii) $C_{\mfp}$ is scalar on $\mfn_j$, $j=1,2$ and
$b^1=\dfrac{q}{np}$ and $b^2=\dfrac{p}{nq}$;

(iii) $c_{\mfk,1}=\dfrac{1}{p}$ and $c_{\mfk,2}=\dfrac{1}{q}$;

(iv) $\ga=\dfrac{q^2+p^2}{npq}$;

(v) $c_{\mfn_1,\mfp}=\dfrac{(p-1)q}{pn}$ and
$c_{\mfn_2,\mfp}=\dfrac{(q-1)p}{qn}$.

\ecor

\bproof The number $c_{\mfl,\mfp}$ is the eigenvalue of the Casimir
operator of $\mfl$ on $\mfp$. Thus, it follows from Proposition \ref{ck1} (ii),
that $c_{\mfl,\mfp}=\dfrac{1}{n}$. Since $\mfn_1\subset
\mfg_0^p\times 0$ and $\mfn_2\subset 0\times \mfg_0^q$, we conclude
from Proposition \ref{ck1} (iii) that $C_{\mfp}$ is scalar on $\mfn_1$ and on
$\mfn_2$; moreover its eigenvalues on these two spaces are
$b^1=\dfrac{q}{np}$ and $b^2=\dfrac{p}{nq}$, respectively.
Similarly, it follows from Proposition \ref{ck1} (iv) that the eigenvalues of
$C_{\mfk}$ on $\mfn_1$ and on $\mfn_2$ are $c_{\mfk,1}=\dfrac{1}{p}$
and $c_{\mfk,2}=\dfrac{1}{q}$, respectively. Now to show (iv) we
recall that $\ga$ is defined by the identity

$$\Phi(C_{\mfk}\cdot,\cdot)\mid_{\mfp\times\mfp}=\ga\Phi\mid_{\mfp\times\mfp}.$$

Let
$(\underbrace{qX,\ldots,qX}_p,\underbrace{-pX,\ldots-pX}_q)\in\mfp$.

$$\bar{rl}& \Phi((qX,\ldots,qX,-pX,\ldots-pX),(qX,\ldots,qX,-pX,\ldots-pX))\\
\\

= & \left(q^2.p+p^2.q\right)\Phi_0(X,X)\\ \\

= & npq\Phi_0(X,X)\ear$$

By using Proposition \ref{ck1} (iv) we get the following:

$$\bar{rl}& \Phi(C_{\mfk}(qX,\ldots,qX,-pX,\ldots-pX),(qX,\ldots,qX,-pX,\ldots-pX))\\
\\

= & \Phi((\frac{q}{p}X,\ldots,\frac{q}{p}X,-\frac{p}{q}X,\ldots-\frac{p}{q}X),(qX,\ldots,qX,-pX,\ldots-pX))\\
\\

= & \left(\frac{q^2}{p}.p+\frac{p^2}{q}.q\right)\Phi_0(X,X)\\ \\

= & (p^2+q^2)\Phi_0(X,X)\ear$$

Therefore, $\ga=\dfrac{q^2+p^2}{npq}$.

\li

Finally, for $j=1,2$, the numbers $c_{\mfn_j,\mfp}$ are defined by

$$\Phi(C_{\mfn_j}\cdot,\cdot)\mid_{\mfp\times\mfp}=c_{\mfn_j,\mfp}\Phi\mid_{\mfp\times\mfp}.$$

From Proposition \ref{ck1} (v), we obtain the following:

$$\bar{rl}& \Phi(C_{\mfn_1}(qX,\ldots,qX,-pX,\ldots-pX),(qX,\ldots,qX,-pX,\ldots-pX))\\
\\

= & \Phi\left(\left(\big(1-\frac{1}{p}\big)qX,\ldots,\big(1-\frac{1}{p}\big)qX,0,\ldots,0\right),(qX,\ldots,qX,-pX,\ldots-pX)\right)\\
\\

= & \dfrac{(p-1)q^2}{p}.p\Phi_0(X,X)\\ \\

= & (p-1)q^2\Phi_0(X,X)\ear$$

Therefore,
$c_{\mfn_1,\mfp}=\dfrac{(p-1)q^2}{npq}=\dfrac{(p-1)q}{np}$.
Similarly, we show that $c_{\mfn_2,\mfp}=\dfrac{(q-1)p}{qn}$.

$\Box$\eproof

\section{Existence of Einstein Adapted Metrics}

In this Section we investigate the existence of adapted metrics on
$M$ with respect to the fibration $$M=\frac{G_0^n}{\triangle^nG_0}\rightarrow
\frac{G_0^p}{\triangle^pG_0}\times \frac{G_0^q}{\triangle^qG_0},$$
as
in previous Section. We shall consider adapted metrics of the form

\beq\label{mdefknss}g_M=g_M(\la,\mu_1,\mu_2)\eeq

with  respect to the decomposition $\mfm=\mfp\oplus\mfn_1\oplus\mfn_2$, i.e., $g_M$ is induced by the scalar product

\beq\la B\mid_{\mfp\times\mfp}\oplus\mu_1B\mid_{\mfn_1\times\mfn_1}\oplus\mu_2B\mid_{\mfn_2\times\mfn_2}\eeq

on $\mfm$, where $B=-\Phi$, the negative of the Killing form. We observe that $\mfn_1$ and $\mfn_2$ are inequivalent
$Ad\,K$-modules, but they are not irreducible, for $n>4$. Hence, adapted
metrics on $M$ are not necessarily of the form (\ref{mdefknss}). However, throughout we shall focus only on adapted
metrics of the form $g_M=g_M(\la,\mu_1,\mu_2)$, unless it is explicitly stated
otherwise.

We recall that it has been proved by McKenzie Y. Wang and Wolfgang
Ziller (\cite{WZ2}, \cite{Ro1}) that the homogeneous space
$M=\frac{G_0^n}{\triangle^nG_0}$, $n\geq 2$, is a standard
Einstein manifold. Hence, there exists at least one Einstein
adapted metric of the form (\ref{mdefknss}). Indeed, by Corollary
\ref{ric1cor1} and Corollary \ref{ck2}, the Ricci curvature of the
standard metric is simply

\beq
Ric=\frac{1}{2}\left(\frac{1}{2}+c_{\mfl,\mfm}\right)B=\left(\frac{1}{4}+\frac{1}{2n}\right)B.\eeq

Consequently, the standard metric is Einstein. The standard metric is an example of a binormal metric on $M$ with
respect to the fibration $M\rightarrow N$. Below we shall classify
all the binormal Einstein metrics on $M$.

\li

Although the submodules $\mfn_1$ and $\mfn_2$ are not $Ad\,K$-irreducible for $n> 4$, the Casimir
operators of $\mfl$, $\mfk$ and $\mfp$ are always scalar on  $\mfn_1$ and
on $\mfn_2$. Hence, it is enough to consider one irreducible submodule in
$\mfn_1$ and one irreducible submodule in $\mfn_2$, for effects of
Ricci curvature. Therefore, according to Theorem \ref{binormal1} (2.19), there is an
one-to-one correspondence, up to homothety, between binormal adapted
Einstein metrics on $M$ and positive solutions of the following
set of equations:

\beq \label{eink1}\de_{12}^{\mfk}(1-X)=\de_{12}^{\mfl}\eeq

\beq\label{eink2}(\ga+2c_{\mfl,\mfp})X^2-\left(1+2c_{\mfk,j}\right)X+(1-\ga+2b^j)=0,\,j=1,2\eeq

Given a positive solution $X$, then binormal adapted Einstein
metrics are given, up to homothety, by

$$g_M=g_M(1,X),$$

i.e., are induced by scalar products of the form $<,>=B\mid_{\mfp\times\mfp}\oplus XB\mid_{\mfn\times\mfn}$ on $\mfm$.

\bthm\label{binormalknss} Let us consider the fibration
$$M=\frac{G_0^n}{\triangle^nG_0}\rightarrow
\frac{G_0^p}{\triangle^pG_0}\times \frac{G_0^q}{\triangle^qG_0}=N,$$

where $p+q=n$ and $2\leq p\leq q\leq n-2$.

If $p\neq q$ or $n=4$, then there exists on $M$ precisely one
binormal Einstein metric, up to homothety, which is the standard
metric. For $n>4$ and $p=q$, then there are on $M$ precisely two
binormal Einstein metrics, up to homothety, which are the standard
metric and the metric induced by the scalar product

$$B\mid_{\mfp\times\mfp}\oplus \frac{n}{4}B\mid_{\mfn\times\mfn}.$$
 \ethm

\bproof From Corollary \ref{ck2} we obtain that
$\de^{\mfk}_{12}=c_{\mfk_1}-c_{\mfk,2}=\frac{1}{p}-\frac{1}{q}$
whereas
$\de_{12}^{\mfl}=c_{\mfl_1}-c_{\mfl,2}=\frac{1}{n}-\frac{1}{n}=0$.
Hence, Equation (\ref{eink1}) implies that $X=1$ or $p=q$. So if
$p\neq q$, if there exists a binormal Einstein metric it must be the
standard metric. This we already know it is Einstein by \cite{Ro1}.
Therefore, if $p\neq q$, then there exists, up to homothety, exactly
one binormal Einstein metric on $G/L$ which is the standard one.

\li

By using Corollary \ref{ck2} Equation (\ref{eink2}) may be rewritten
as

\beq\label{eink3}nX^2-q(p+2)X+pq+q-p=0,  \textrm{ for }j=1 \eeq

and

\beq\label{eink4}nX^2-p(q+2)X+pq+p-q=0,  \textrm{ for }j=2 \eeq

It is clear that $X=1$ is actually a solution of both equations,
and thus the standard metric is in fact Einstein.

\li

Now suppose that $p=q$. As $n=p+q$, then $p=q=\frac{n}{2}$.
Therefore, (\ref{eink3}) and (\ref{eink4}) become equivalent to

\beq 4X^2-(n+4)X+n=0\eeq

The polynomial above has two positive roots, $1$ and $\frac{n}{4}$.
Therefore, for $p=q$ and $n>4$, there exist precisely two binormal
Einstein metrics.

$\Box$\eproof

\bthm \label{form2}Let $g_M$ be an Einstein adapted metric on $M$
of the form $g_M(\la;\mu_1,\mu_2)$. The projection $g_N$ onto the
base space is also Einstein if and only if $p=q$ and $g_M$ is
binormal. \ethm

\bproof By Theorem \ref{gnein} we know that if $g_M$ and $g_N$ are
Einstein then we must have the relation

\beq\label{rbsk}\frac{r_1}{r_2}=\left(\frac{b^1}{b^2}\right)^{\frac{1}{2}}.\eeq

From Lemma \ref{ck2} (ii) we obtain
$\left(\frac{b^1}{b^2}\right)^{\frac{1}{2}}=\frac{q}{p}$. Since
$[\mfn_1,\mfn_2]=0$, from Corollary \ref{riccinb2}, we get
$r_i=\frac{1}{2}\left(\frac{1}{2}+c_{\mfk,i}\right)$, $i=1,2$.
Hence, $\frac{r_1}{r_2}=\frac{(p+2)q}{(q+2)p}$, by using Lemma
\ref{ck2} (iii). Therefore, (\ref{rbsk}) is possible if and only if
$p=q$. Also from the proof of Theorem \ref{gnein}, if $g_N$ and
$g_M$ are Einstein, then
$\frac{\mu_1}{\mu_2}=\left(\frac{b^1}{b^2}\right)^{\frac{1}{2}}=\frac{p}{q}=1$
and $g_M$ is binormal. Conversely, if $g_M$ is binormal and $p=q=n/2$, then by the above we also get

$$\frac{r_1}{\mu_1}=\frac{r_2}{\mu_2}$$ and $g_N$ is Einstein.

$\Box$\eproof

We observe that since $\mfp$ is an irreducible $Ad\,L$-submodule, $g_F$ is always Einstein.

Theorems \ref{binormalknss} and \ref{form2} classify all Einstein binormal metrics on $M$ and all Einstein adapted metrics such that $g_N$ is also Einstein. It still remains to understand if there are other Einstein adapted metrics besides these. The Einstein equations in general for arbitrary $p$ and $q$ are extremely complicated. However with the help of Maple it is still possible to solve the problem in general. Next we shall classify all the Einstein adapted metrics on $M$ of the form $g_M(\la,\mu_1,\mu_2)$.

\blem \label{knsseeq}Consider the fibration
$$M=\frac{G_0^n}{\triangle^nG_0}\rightarrow
\frac{G_0^p}{\triangle^pG_0}\times \frac{G_0^q}{\triangle^qG_0}=N,$$

where $p+q=n$ and $2\leq p\leq q\leq n-2$. There is a one-to-one
correspondence between Einstein adapted metrics on $M$ of the form
$g_M=g_M(\la;\mu_1,\mu_2)$, up to homothety, and positive
solutions of the following system of Equations:

\beqar -2q^2X_1^2+nq(p+2)X_1+2p^2X_2^2-np(q+2)X_2=0 \\
n^2+q^2(p+1)X_1^2+p^2(q-1)X_2^2-nq(p+2)X_1=0. \eeqar

To a positive solution $(X_1,X_2)$ corresponds an Einstein metric
of the form $g_M=g_M\big(1;\frac{1}{X_1},\frac{1}{X_2}\big)$.\elem

\bproof Let $g_M$ be an adapted metric on $M$ of the form
$g_M(\la;\mu_1,\mu_2)$. We set

\beq\label{un1} X_i=\dfrac{\la}{\mu_i},\,i=1,2.\eeq

Since the fiber $F$ is irreducible we may use Proposition \ref{riccinf} to obtain the Ricci curvature. For $X\in\mfp$,

$$Ric(X,X)=\left(\frac{1}{2}\left(c_{\mfl,\mfp}+\frac{\ga}{2}\right)+\frac{\la^2}{4}\sum_{j=1}^n\dfrac{c_{\mfn_j,\mfp}}{\mu_j^2}\right)B(X,X).$$

Hence by using the eigenvalues in Corollary \ref{ck2} and the unknowns $X_1$ and $X_2$ defined in (\ref{un1}), we obtain

$$\frac{1}{2}\left(c_{\mfl,\mfp}+\frac{\ga}{2}\right)=\frac{1}{2n}+\frac{q^2+p^2}{4npq}=\frac{n}{4pq}$$

and

\beq Ric(X,X)=\left(\frac{n}{4pq}+\frac{(p-1)q}{4pn}X_1^2+\frac{(q-1)p}{4qn}X_2^2\right)B(X,X).\eeq

For $Y\in\mfn_k$, the ricci curvature is given by

$$-\frac{\la}{2\mu_k}B(C_{\mfp}Y,Y)+r_kB(Y,Y).$$

The Casimir operator of $\mfp$ is scalar on $\mfn_i$ with
eigenvalues $\frac{q}{np}$, for $i=1$, and $\frac{p}{nq}$, for
$i=2$, as we can see from Corollary \ref{ck2}. Since
$[\mfn_1,\mfn_2]=0$, we use Corollary \ref{riccinb2} to obtain
$r_k$:

$$r_k=\frac{1}{2}\left(\frac{1}{2}+c_{\mfk,k}\right).$$

Hence, from Corollary \ref{ck2}, we get

\beqar  Ric(Y,Y)=\left(-\frac{q}{2np}X_1+\frac{p+2}{4p}\right)B(Y,Y),\,Y\in\mfn_1\\
Ric(Y,Y)=\left(-\frac{p}{2nq}X_2+\frac{q+2}{4q}\right)B(Y,Y),\,Y\in\mfn_2.\eeqar

Finally, as $C_{\mfn_i}(\mfp)\subset \mfk$, for $i=1,2$, from Proposition \ref{riccinf}, we conclude that $Ric(\mfp,\mfn)=0$. Therefore, the Einstein equations for $g_M$ are just

\beqar\label{knsseeq1} \frac{n}{4pq}+\frac{(p-1)q}{4pn}X_1^2+\frac{(q-1)p}{4qn}X_2^2=\la E\\
\label{knsseeq2}-\frac{q}{2np}X_1+\frac{p+2}{4p}=\mu_1E\\
\label{knsseeq3}-\frac{p}{2nq}X_2+\frac{q+2}{4q}=\mu_2E\eeqar

where $E$ is the Einstein constant. We obtain the system of equations stated in this lemma by eliminating $E$ from the system above and rearranging the resulting equations.

$\Box$\eproof

\bthm\label{genknss} Let us consider the fibration
$$M=\frac{G_0^n}{\triangle^nG_0}\rightarrow
\frac{G_0^p}{\triangle^pG_0}\times \frac{G_0^q}{\triangle^qG_0}=N,$$

where $p+q=n$ and $2\leq p\leq q\leq n-2$. If $n>4$, there exist on $M$ exactly two Einstein adapted metrics of the form $g_M=g_M(\la,\mu_1,\mu_2)$ and one is the standard metric. For $n=4$ the only Einstein adapted metric is the standard one.

For $p=q$ the non-standard Einstein adapted metric is
binormal.\ethm

\bproof According to Lemma \ref{knsseeq}, the Einstein equations for $g_M$ are as follows:

\beqar -2q^2X_1^2+nq(p+2)X_1+2p^2X_2^2-np(q+2)X_2=0 \\
n^2+q^2(p+1)X_1^2+p^2(q-1)X_2^2-nq(p+2)X_1=0. \eeqar

By using Maple we obtain that the solutions of the system above are $X_1=X_2=1$ and

\beq \label{solgen2}X_1=\alpha,\,X_2=\left(\frac{-q^2(p+1)\alpha^2+nq(p+2)\alpha-n^2}{p^2(q-1)}\right)^{\frac{1}{2}},\eeq

where $\alpha$ is a root of the polynomial

$$t(Z)=4q^2Z^3-4q(n+pq+2)Z^2+n(q(q+2)(p+1)+n+8)Z-(q+3)n^2.$$

The solution $X_1=X_2=1$  corresponds to a standard metric and we recover the result that $M$ is an Einstein standard manifold. We are now interested in analysing the existence of other metrics. First we observe that from the expression for $X_2$ in (\ref{solgen2}) we conclude that,

$$X_2\in\reals\textrm{ if and only if }\al\in\left(\frac{n}{q(p+1)},\frac{n}{q}\right).$$

For this we compute the roots of the polynomial $-q^2(p+1)\alpha^2+nq(p+2)\alpha-n^2$, as in (\ref{solgen2}).

Simple calculations show that

\beqar t\left(\frac{n}{q}\right)=\frac{p(q-1)^2n^2}{q}>0\nonumber\\
t\left(\frac{n}{q(p+1)}\right)=-\frac{p(p+3)^2(q-1)n^2}{q(p+1)^3}<0\nonumber
\eeqar

and thus $t$ has at least one (positive) root in the interval $\left(\frac{n}{q(p+1)},\frac{n}{q}\right)$, by the Bolzano Theorem. By the explained above, to this root corresponds an Einstein adapted metric on $M$. Furthermore, we show that this root is unique and distinct from 1. From this we conclude that there exists a non-standard Einstein adapted metric on $M$. We observe that the derivative of $t$,
$\frac{dt}{dZ}$ has no real zeros. Simple calculations show that the zeros of $\frac{dt}{dZ}$ are

$$\frac{n+pq+2}{3q}\pm\frac{\sqrt{\de}}{6q},$$

where $$\de=(q+1)^2p^2-(q-1)(3q^2+4q-8)p-(q-1)(3q^2+8q+16).$$

We show that $\de<0$. For $p=q$, $\de=-2q^4-2q^3+8q^2-16q+16<0$, for every $q\geq 2$. So we suppose that $p<q$. In this case, since $p^2\leq (q-1)p$, we have

$$\bar{rl}\de\leq & (q-1)\big(-(2p+3)q^2-(2p+8)q+(9p-16)\big)\\

< & (q-1)(-2p^3-5p^2+p-16)\\

< & 0,\ear$$

for every $p\geq 2$.

With this we conclude that $\frac{dt}{dZ}$ has no real zeros and thus the root of $t$ found above is the unique real root of $t$. Moreover, we must guarantee that this root does not yield the solution $X_1=X_2=1$. If $X_1=X_2=1$, then $\al=1$ is a root of $t$. This may be possible since $1\in\left(\frac{n}{q(p+1)},\frac{n}{q}\right)$. Since

$$t(1)=p(q+2)(q-1)(n-4),$$

and, consequently, $\al=1$ is a root of $t$ if and only if $n=4$.
By using (\ref{solgen2}) we get that if $X_1=1$ when $n=4$, then
$X_2=1$ as well. Since non-standard Einstein adapted metrics are
given by pairs of the form (\ref{solgen2}), with $\al\neq 1$, we
conclude that there exists a unique non-standard Einstein adapted
metric of the form $g_M(\la,\mu_1,\mu_2)$ if and only if $n>4$; in
the case $n=4$, the standard metric is the unique Einstein adapted
metric of the form $g_M(\la,\mu_1,\mu_2)$. Finally, we observe
that, if $n=4$, the subspaces $\mfn^1$ and $\mfn^2$ are
irreducible $Ad\,L$-submodules. Hence, any adapted metric on $M$
is of the form $g_M(\la,\mu_1,\mu_2)$. Therefore, we conclude
that, for $n=4$, there exists a unique Einstein adapted metric on
$M$ and it is the standard one.

Since there is a unique non-standard Einstein adapted metric on
$M$, it follows from Theorem \ref{binormalknss} that this metric
is binormal if and only if $p=q$.

$\Box$\eproof

\li

As it has been observed previously the modules $\mfn_k$, $k=1,2$, are not irreducible $Ad\,K$-modules. From Lemma \ref{ck1} we deduce a possible decomposition for $\mfn_1$ and $\mfn_2$ into irreducible $Ad\,K$-submodules:

\blem \label{decmfn}$\mfn_1=\oplus_{j=1}^{p-1}\mfn_{1,j}$ and
$\mfn_2=\oplus_{j=1}^{q-1}\mfn_{2,j}$ , where

$$\mfn_{1,j}=\{(\underbrace{X,\ldots,X}_j,-jX,0,\ldots,0)\in\mfg_0^n:X\in\mfg_0\}, \textrm{ for } j=1,\ldots,p-1$$

$$\mfn_{2,j}=\{(\underbrace{0,\ldots,0}_p,\underbrace{X,\ldots,X}_j,-jX,0,\ldots,0)\in\mfg_0^n:X\in\mfg_0\}, \textrm{ for } j=1,\ldots,q-1.$$

Furthermore, the $\mfn_{i,j}$'s are irreducible
$Ad\,K$-submodules. \elem

The fact that the modules above are irreducible follows from the fact that $\mfg_0$ is simple. Also by similar calculations to those in Lemma \ref{ck2} we obtain the Casimir operators of the submodules $\mfn_{1,j}$ and $\mfn_{2,j}$:

\blem $C_{\mfn_{1,j}}=\Big(\underbrace{\frac{1}{j(j+1)},\ldots,
\frac{1}{j(j+1)}}_j,\dfrac{j}{j+1},0,\ldots,0\Big), \textrm{ for }
j=1,\ldots,p-1$, and

$C_{\mfn_{2,j}}=\Big(\underbrace{0,\ldots,0}_p,\underbrace{\frac{1}{j(j+1)},\ldots,
\frac{1}{j(j+1)}}_j,\dfrac{j}{j+1},0,\ldots,0\Big), \textrm{ for }
j=1,\ldots,q-1$.
\elem

\bproof $\{(\underbrace{u_i,\ldots,u_i}_j,-ju_i,0\ldots,0\}_i$ and
$\{(\underbrace{u'_i,\ldots,u'_i}_j,-ju'_i,0,\ldots,0)\}_i$ are
bases of $\mfn_{1,j}$.

$$\bar{rl} &  \Phi((u_i,\ldots,u_i,-ju_i,0,\ldots,0),(u'_j,\ldots,u'_j,-ju'_j,0,\ldots,0))\\ \\

= & k\Phi_0(u_i,u'_j)+k^2\Phi_0(u_i,u'_j)\\ \\

= & (k+1)k\de_{ij}\ear$$

Hence dual bases of $\mfn_{1,j}$ for $\Phi$ are as follows:

$\{\frac{1}{\sqrt{k(k+1)}}(u_i,\ldots,u_i,-ju_i,0\ldots,0\}_i$ and
$\{\frac{1}{\sqrt{k(k+1)}}(u'_i,\ldots,u'_i,-ju'_i,0,\ldots,0)\}_i.
$

By using these bases we obtain the required expression for the
Casimir operator of $\mfn_{1,j}$. For $C_{\mfn_{2,j}}$ is similar.

$\Box$\eproof

Hence, we may consider on $M$ an adapted metric of the form

\beq \label{mdefknss2} g_M=g_M(\la;\mu_{1,1},\ldots,\mu_{1,p-1},\mu_{2,1},\ldots,\mu_{2,q-1}).\eeq

Clearly, whereas the submodules $\mfn_{i,j}$ are $Ad\,K$-irreducible, the reader should note that they are not pairwise inequivalent. Hence, adapted metrics on $M$ are not necessarily of the form (\ref{mdefknss2}).

\blem \label{form1}Let $g_M$ be an adapted metric on $M$ of the form

$$g_M(\la;\mu_{1,1},\ldots,\mu_{1,p-1},\mu_{2,1},\ldots,\mu_{2,q-1}).$$

If $g_M$ is Einstein, then it is of the form

\beq\label{mdefknss2}g_M(\la,\mu_1,\mu_2).\eeq

\elem

\bproof By Corollary \ref{cond2}, if exists on $M$ an Einstein
adapted metric, then

\beq\label{cond2knss}\big(\sum_{j=1}^{p-1}\nu_{1,j}C_{\mfn_{1,j}}+\sum_{j=1}^{q-1}\nu_{2,j}C_{\mfn_{2,j}}\big)(\mfp)\subset
\mfk,\eeq

for some $\nu_{1,j},\,\nu_{2,j}>0$. The inclusion (\ref{cond2knss})
implies
\beq\label{cond3knss}\sum_j\nu_{i,j}C_{\mfn_{i,j}}(\mfp)\subset
\mfk,\eeq

for $i=1,2$, since $C_{\mfn_{1,j}}(\mfg)\subset \mfg_0^p\times 0_q$
and $C_{\mfn_{2,j}}(\mfg)\subset 0_p\times\mfg_0^q$. For
$k=1,\ldots,p$, let $P_k$ denote the $k$th component of
$P=\sum_j\nu_{1,j}C_{\mfn_{1,j}}$. By \ref{ck1} (vi), we have

$$P_k=\frac{k-1}{k}\nu_{1,k-1}+\frac{1}{k(k+1)}\sum_{j=k}^{p-1}\nu_{1,k},$$

for $k=2,\ldots,p-1$. We may then write

\beq\label{pk1}P_{k+1}=P_k+\frac{k-1}{k}(\nu_{1,k}-\nu_{1,k-1}).\eeq

On the other hand, the condition $P(\mfp)\subset\mfk$, implies that
$P_k=P_{k+1}$. Hence, from (\ref{pk1}) we obtain that
$\nu_{1,k}=\nu_{1,k-1}$. Similarly we show that
$\nu_{2,k}=\nu_{2,k-1}$. Now from the proof of Corollary
\ref{cond2}, we can see that for an adapted Einstein metric on $M$
as in (\ref{mdefknss}), then the constants $\nu_{1,k}=1/\mu_{1,k}^2$
and $\nu_{2,k}=1/\mu_{2,k}^2$ must satisfy (\ref{cond2knss}). This
concludes the proof.

$\Box$\eproof

Therefore, Theorem \ref{genknss} may be understood as a classification of Einstein metrics of the form (\ref{mdefknss2}).

\bcor \label{genknss} Let us consider the fibration
$$M=\frac{G_0^n}{\triangle^nG_0}\rightarrow
\frac{G_0^p}{\triangle^pG_0}\times \frac{G_0^q}{\triangle^qG_0}=N,$$

where $p+q=n$ and $2\leq p\leq q\leq n-2$. Let

\beq \label{mdeffim}g_M(\la;\mu_{1,1},\ldots,\mu_{1,p-1},\mu_{2,1},\ldots,\mu_{2,q-1})\eeq

be an adapted metric on $M$, corresponding to the decomposition of $\mfn$ given in Lemma \ref{decmfn}. If $n>4$, there exist on $M$ exactly two Einstein adapted metrics of the form (\ref{mdeffim}) and one is the standard metric. \ecor

\newpage

\section{Closing Remarks}

The main aim of this thesis was to establish conditions for
existence of homogeneous Einstein metrics with totally geodesic
fibers and to bring new existence and non-existence results for a
significant class of homogeneous spaces. Necessary and sufficient
conditions for existence of such metrics were obtained under some
hypothesis and, for some classes of homogeneous fibrations, new
Einstein metrics with totally geodesic fibers were found and in
other cases all such metrics were classified. Nevertheless, some
questions remain open for further research.

\li

For irreducible bisymmetric fibrations of maximal rank all the
Einstein adapted metrics were classified in the case when $G$ is
an exceptional Lie group. If $G$ is a classical group, all the
Einstein adapted metrics were classified for Type I, whereas for
Type II we classified only those whose restriction to the fiber is
still Einstein. For these bisymmetric fibrations which admit an
Einstein adapted metric which satisfies this condition we can
still classify all the other Einstein adapted metrics. However, in
this classical case, it remains to obtain a general classification
of  Einstein adapted metrics.

\li

Furthermore, the techniques and results for irreducible
bisymmetric fibrations of maximal rank suggest that a similar
research can be developed for homogeneous fibrations whose fiber
and base space are $p$-symmetric spaces of higher order. For
instance, the classification of compact simply-connected
$3$-symmetric spaces would allow us to consider fibrations whose
fiber is a $3$-symmetric space and the base is still an
irreducible symmetric space.

\li For Kowalski $n$-symmetric spaces we have shown that there
exists a non-standard Einstein adapted metric, if $n>4$. For $n=4$
we have proved that the standard metric is the only Einstein
metric with totally geodesic fibers. It remains to classify all
the Einstein adapted metrics, if $n>4$.

\appendix
\chapter{}\label{cpproofs}

A classification of bisymmetric triples of maximal rank was given
in Chapter 3 and they are listed in Tables \ref{eigIexc},
\ref{eigIclass}, \ref{eigIIexc} and \ref{eigIIclass}. In this
Appendix we determine the isotropy representation for each
bisymmetric triple and compute the eigenvalues $\ga_a$'s and
$b_a^{\phi}$ of the Casimir operators $C_{\mfk}$ and $C_{\mfp_a}$
along the subspaces $\mfp_a$ and $\mfn$, respectively. We use the
formulas presented in Chapter 3 in Propositions \ref{bejs} and
\ref{gajs}. We shall systematically use the roots systems and the
dual Coxeter numbers. The roots systems used can be found in
\cite{He} and \cite{OV} and the dual Coxeter numbers are given in
Table \ref{tabcoxeter}.

\li

As introduced in Chapter 3, if $\mfn=\oplus_j\mfn^j$ is a
decomposition of $\mfn$ into irreducible $Ad\,L$-modules, we write

\beq \R_{\mfn^j}=\{\phi\in\R: E_{\phi}\in(\mfn^j)^{\complex}
\}\textrm{ and }\mfn^j=<X_{\phi},
Y_{\phi}:\phi\in\R_{\mfn^j}^+>.\eeq

We recall that $b_a^{\phi}$, for $\phi\in\R_{\mfn^j}$, is the eigenvalue of $C_{\mfp_a}$ on $\mfn^j$:

$$C_{\mfp_a}\mid_{\mfn^j}=b_a^{\phi}Id_{\mfn^j}.$$

If $\mfp$ is $Ad\,L$-irreducible, we write simply $b^{\phi}$ for the eigenvalue of $C_{\mfp}$ on $\mfn^j$. These eigenvalues are given by the formula from Proposition
\ref{bejs}:

\beq \label{bejs2}b_a^{\phi}=\dfrac{1}{2}\sum _{\al\in
\R_{\mfp_a}^+}d_{\al\phi}|\al|^2,\eeq

where

\beq\label{dalphi}d_{\al\phi}=q_{\al\phi}-p_{\al\phi}-2p_{\al\phi}q_{\al\phi}\eeq

and $\phi+n\al$, $p_{\al\phi}\leq n\leq q_{\al\phi}$ is the
$\al$-series containing $\phi$.

\li

Also if $\mfk_a$ is a simple ideal of $\mfk$, $\ga_a$ is the eigenvalue of $C_{\mfk}$ on $\mfk_a$:

$$C_{\mfk}\mid_{\mfk_a}=\ga_aId_{\mfk_a}.$$

And finally to compute the
$\ga_a$'s, when $\mfk_a$ is a simple ideal of $\mfk$, we use Proposition
\ref{gajs}:

\beq\label{gajs2}\ga_a=\dfrac{h^*(\mfk_a)}{\de_a .
h^*(\mfg)},\,a=1,\ldots,s\eeq . where $h^*(\mfk_a)$ and
$h^*(\mfg)$ are the dual Coxeter numbers of $\mfk_a$ and $\mfg$,
respectively, and $\de_a=|\al|^2/|\be|^2$, for $\al$ a long root
of $\mfg$ and $\be$ a long root of $\mfk_a$. If there is only one
root length on $\mfg$ or both $\mfg$ and $\mfk_a$ have two root
lengths, $\de_a=1$. If $\de_a\neq 1$ then $\de_a$ it is equal to
$2$ or $3$. If $\mfp$ is $Ad\,L$-irreducible, then we write simply
$\ga$ for the eigenvalue of $C_{\mfk}$ on the corresponding simple
ideal of $\mfk$.

Since the computations of the eigenvalues $\ga_a$ and $b_a^{\phi}$
consist of a systematic use of formulas (\ref{bejs2}) and
(\ref{gajs2}) some details of these computations are omitted. For
each simple Lie algebra $\mfg$ we present the set of roots of
$\mfg$ and the corresponding length of the roots. We consider each
symmetric pair  of maximal rank $(\mfg,\mfk)$ and present the
subsets of roots of each simple ideal $\mfk_a$ of $\mfk$, the
corresponding value of $\ga_a$ and the subset of roots of the
isotropy space $\mfn$. We recall that $\R_{\mfn}=\R-\R_{\mfk}$.
Finally, for each bisymmetric triple $(\mfg,\mfk,\mfl)$, we
indicate the subset of roots $\R_{\mfp}$ for the symmetric
complement $\mfp$ of the symmetric pair $(\mfk,\mfl)$ and the
subset of roots of each irreducible $Ad\,L$-submodule $\mfn^i$ of
$\mfn$. For bisymmetric triples of Type I, we present the
essential information to compute the eigenvalues $b^{\phi}$ on
each subspace $\mfn^i$. Since $\mfn^i$ is $Ad\,L$-irreducible, it
suffices to choose any root $\phi$ in $\R_{\mfn^i}$. Thus, we
choose a root $\phi$ in each $\mfn^i$ and indicate all the roots
$\al\in\R_{\mfp}^+$ such that the string $\phi+n\al$ is not
singular, i.e., it contains other roots besides $\phi$; only for
these $\al$'s the coefficients $d_{\al\phi}$ in formula
(\ref{bejs2}) are non-zero. We indicate the elements in the
non-singular string $\phi+n\al$; only three cases occur: this
string is formed either by

(i) $\phi,\,\phi+\al$, in which case $p_{\al\phi}=0$ and $q_{\al\phi}=1$; thus, $d_{\al\phi}=1$;

(ii) $\phi,\,\phi-\al$, in which case $p_{\al\phi}=-1$ and $q_{\al\phi}=0$; thus, $d_{\al\phi}=1$;

(iii) $\phi,\,\phi\pm \al$, in which case $p_{\al\phi}=-1$ and $q_{\al\phi}=1$; thus, $d_{\al\phi}=4$.

We indicate the length $|\al|^2$ of each root $\al$ listed
previously. Once all this information is obtained, we apply
(\ref{bejs2}) to compute $b^{\phi}$.

The reader may notice throughout that, in some cases, the root
$\al$ indicated is not a positive root. We observe that since
$d_{\al\phi}=d_{-\al\phi}$, we may choose the $\al$'s in
$\R_{\mfp}$ independently of the sign, as long as only one of
$\pm\al$ is chosen. This allows  us to chose the necessary $\al$'s
without any considerations about the order of the roots in the Lie
algebra $\mfg$.

For bisymmetric triples of Type II, we obtain the eigenvalues
$b_a^{\phi}$ from the bisymmetric triples of Type I. For instance,
$b_1^{\phi}$ and $b_2^{\phi}$ for
$(\mfg_2,\mfsu_2\oplus\mfsu_2,\reals\oplus\reals)$ are given by
$b^{\phi}$ of $(\mfg_2,\mfsu_2\oplus\mfsu_2, \reals\oplus\mfsu_2)$
and by $b^{\phi}$ of $(\mfg_2,\mfsu_2\oplus\mfsu_2,
\mfsu_2\oplus\reals)$, respectively.

\li

In the cases of the classical Lie algebras, it is easy to
understand whether the sum of two roots is a root, due to the
simplicity of their root systems. However, in the exceptional
cases, this may be rather complicated, mainly in the cases of the
Lie algebras $\mfe_6$, $\mfe_7$ and $\mfe_8$. In these three
cases, auxiliary lemmas are provided, where conditions under which
the sum of two roots is a root are stated. These are the Lemmas
\ref{sumse8}, \ref{sumse7} and \ref{sumse6}.

\newpage

\section{$A_{n-1}$} In this Section we consider bisymmetric
triples of the form
$(\mfsu_n,\mfsu_p\oplus\mfsu_{n-p}\oplus\reals,\mfl)$, $1\leq
p\leq n-1$. We set $\mfk_1=\mfsu_p$, $\mfk_2=\mfsu_{n-p}$ and
$\mfk_0=\reals$.

For a root system of type $A_{n-1}$ for $\mfg$ we take

\beq\label{ran}\R=\{\pm(e_i-e_j):i\leq i<j\leq n\}.\eeq

In $\mfg$ there is only one root length and is

\beq\label{rlan}|\al|^2=\frac{1}{n}.\eeq

\li

{\footnotesize{

\bsymp $(\mfsu_n,\mfsu_p\oplus\mfsu_{n-p}\oplus\reals)$,
$p=1,\ldots,n-1$.

$$\bar{ccc} \mfk_i & \R_{\mfk_i} & \ga_i\\ \hline\hline

\mfsu_p & \{\pm(e_i-e_j):1\leq i<j\leq p\} & \frac{p}{n},\,p\geq
2\\
\mfsu_{n-p} & \{\pm(e_i-e_j):p+1\leq i<j\leq n\} &
\frac{n-p}{n},\,p\leq n-2\ear$$

\li

$$\R_{\mfn}=\{\pm(e_i-e_j):1\leq i\leq p,\,p+1\leq
j\leq n\}$$\esymp

\li

\bbst\label{cpan1}
$(\mfsu_n,\mfsu_p\oplus\mfsu_{n-p}\oplus\reals,\mfsu_l\oplus\mfsu_{p-l}\oplus\reals\oplus\mfsu_{n-p}\oplus\reals)$,
$1\leq p\leq n-1$, $1\leq l\leq p-1$. (Type I)

$$\bar{l}\R_{\mfp}=\{\pm(e_i-e_j):1\leq i\leq l,\,l+1\leq j\leq p\}\\

\mfn=\mfn^1\oplus\mfn^2\\

\R_{\mfn^1}=\{\pm(e_i-e_j):1\leq i\leq l,\,p+1\leq j\leq
n\}\\
\R_{\mfn^2}=\{\pm(e_i-e_j):l+1\leq i\leq p,\,p+1\leq j\leq
n\} \ear$$

$$\bar{cccccccc}\mfn^i & \phi\in\R_{\mfn^i} & \al\in\R_{\mfp}^+  & \phi+n\al & d_{\al\phi} & \textrm{No of } \al's & |\al|^2 & b^{\phi}\\ \hline\hline

\xstrut \mfn^1 & e_1-e_n & e_1-e_j,\,l+1\leq j\leq p &
\phi,\,\phi-\al & 1 & p-l & \frac{1}{n} & \frac{p-l}{2n}\\  \hline

\xstrut \mfn^2 & e_p-e_n & e_j-e_p,\,1\leq j\leq l &
\phi,\,\phi+\al & 1 & l & \frac{1}{n} & \frac{l}{2n}\ear$$\ebst

\li

 \bbst\label{cpan2}
$(\mfsu_n,\mfsu_p\oplus\mfsu_{n-p}\oplus\reals,\mfsu_p\oplus\mfsu_{s}\oplus\mfsu_{n-p-s}\oplus\reals\oplus\reals)$,\,$1\leq
p\leq n-1$, $1\leq s\leq n-p-1$. (Type I)

$$\bar{l}\R_{\mfp}=\{\pm(e_i-e_j):p+1\leq i\leq p+s,\,p+s+1\leq j\leq n\}\\

\mfn=\mfn^1\oplus\mfn^2\\

\R_{\mfn^1}=\{\pm(e_i-e_j):1\leq i\leq p,\,p+1\leq j\leq
p+s\}\\
\R_{\mfn^2}=\{\pm(e_i-e_j):1\leq i\leq p,\,p+s+1\leq j\leq
n\} \ear$$

$$\bar{cccccccc}\mfn^i & \phi\in\R_{\mfn^i} & \al\in\R_{\mfp}^+  & \phi+n\al & d_{\al\phi} & \textrm{No of } \al's & |\al|^2 & b^{\phi}\\ \hline\hline

\xstrut \mfn^1 & e_1-e_{p+s} & e_{p+s}-e_i,\,p+s+1\leq j\leq n &
\phi,\,\phi+\al & 1 & n-p-s & \frac{1}{n} & \frac{n-p-s}{2n}\\
\hline

\xstrut \mfn^2 & e_p-e_n & e_i-e_n,\,p+1\leq j\leq p+s &
\phi,\,\phi-\al & 1 & s & \frac{1}{n} & \frac{s}{2n}\ear$$\ebst

\li

\bbst\label{cpan3}
$(\mfsu_n,\mfsu_p\oplus\mfsu_{n-p}\oplus\reals,\mfsu_l\oplus\mfsu_{p-l}\oplus\mfsu_{s}\oplus\mfsu_{n-p-s}\oplus\reals\oplus\reals\oplus\reals)$,
$1\leq p\leq n-1$, $1\leq l\leq p-1$, $1\leq s\leq n-p-1$. (Type
II)

$$\bar{lll}\R_{\mfp_1}=\{\pm(e_i-e_j):1\leq i\leq l,\,l+1\leq j\leq
p\}\\

\R_{\mfp_2}=\{\pm(e_i-e_j):p+1\leq i\leq p+s,\,p+s+1\leq j\leq
n\}\\

\mfn=\mfn^1\oplus\mfn^2\oplus\mfn^3\oplus\mfn^4\ear$$

$$\bar{lcc}\R_{\mfn^i} & b_1^{\phi} & b_2^{\phi} \\ \hline
\{\pm(e_i-e_j):1\leq i\leq l,\,p+1\leq j\leq
p+s\} & \frac{p-l}{2n} & \frac{n-p-s}{2n}\\

\{\pm(e_i-e_j):1\leq i\leq 1,\,p+s+1\leq j\leq n\} & \frac{p-l}{2n} & \frac{s}{2n}\\

\{\pm(e_i-e_j):l+1\leq i\leq p,\,p+1\leq j\leq p+s\} & \frac{l}{2n} & \frac{ n-p-s}{2n}\\

\{\pm(e_i-e_j):l+1\leq i\leq p,\,p+s+1\leq j\leq
n\} & \frac{l}{2n} & \frac{s}{2n}
\ear$$
\ebst

}}

\newpage

\section{$B_n$}

In this Section we consider all the bisymmetric triples of the form
$(\mfso_{2n+1},\mfso_{2p+1}\oplus\mfso_{2(n-p)},\mfl)$, for $0\leq
p\leq n-1$.

A root system for $\mfg$ is

\beq\label{rlbn}\R=\{\pm e_i:1\leq i\leq n;\,\pm e_i\pm e_j:1\leq
i<j\leq n\}\eeq

and the length of a root is

\beq\label{rlbn}|\al|^2=\left\{\bar{ll}\frac{1}{2(2n-1)},& \al=\pm
e_i\\\frac{1}{2n-1}, & \al=\pm e_i \pm e_j\ear.\right.\eeq

\li

{\footnotesize{

\bsymp $(\mfso_{2n+1},\mfso_{2p+1}\oplus\mfso_{2(n-p)})$, $0\leq
p\leq n-1$.

$$\bar{ccc} \mfk_i & \R_{\mfk_i} & \ga_i\\ \hline\hline

\mfso_{2p+1} & \{\pm e_i:1\leq i\leq p;\,\pm e_i\pm e_j:1\leq
i<j\leq p\} & \frac{2p-1}{2n-1},\,p\geq 1\\
\mfso_{2(n-p)} & \{\pm e_i\pm e_j:p+1\leq i<j\leq n\} &
\frac{2(n-p-1)}{2n-1},\,p\leq n-2\ear$$

\li

$$\R_{\mfn}=\{\pm e_i\pm e_j:1\leq i\leq ,\,p+1\leq
j\leq n\}$$\esymp

\li

\bbst \label{cpbn1}
$(\mfso_{2n+1},\mfso_{2p+1}\oplus\mfso_{2(n-p)},\mfso_{2l+1}\oplus\mfso_{2(p-l)}\oplus\mfso_{2(n-p)})$,
$0\leq p\leq n-1$, $0\leq l\leq p-1$. (Type I)

$$\bar{l}\R_{\mfp}=\{\pm e_i: l+1\leq i\leq p,\,\pm e_i\pm e_j:1\leq i\leq l,\,l+1\leq j\leq p\}\\

\mfn=\mfn^1\oplus\mfn^2\\

\R_{\mfn^1}=\{\pm e_i:p+1\leq i\leq n;\,\pm e_i\pm e_j:1\leq i\leq l,\,p+1\leq j\leq
n\}\\

\R_{\mfn^2}=\{\pm e_i\pm e_j:l+1\leq i\leq p,\,p+1\leq j\leq
n\}\ear$$

$$\bar{cccccccc}\mfn^i & \phi\in\R_{\mfn^i} & \al\in\R_{\mfp}^+  & \phi+n\al & d_{\al\phi} & \textrm{No of } \al'ss & |\al|^2 & b^{\phi}\\ \hline\hline

\xstrut \mfn^1 & e_n & e_i,\,l+1\leq i\leq p & \phi,\,\phi\pm\al &
4 & p-l & \frac{1}{2(2n-1)} & \frac{p-l}{2n-1}\\  \hline

\xstrut \mfn^2 & e_p+e_n & \bar{c} e_i+e_p,\,1\leq i\leq l\\
e_p\ear & \phi,\,\phi-\al & 1 & \bar{c}2l\\1\ear & \bar{c}
\frac{1}{2n-1}\\\frac{1}{2(2n-1)}\ear &
\frac{4l+1}{4(2n-1)}\ear$$\ebst

\li

\bbst \label{cpbn2}
$(\mfso_{2n+1},\mfso_{2p+1}\oplus\mfso_{2(n-p)},\mfso_{2p+1}\oplus\mfso_{2s}\oplus\mfso_{2(n-p-s)})$,
$0\leq p\leq n-1$, $1\leq s\leq n-p-1$. (Type I)

$$\bar{l}\R_{\mfp}=\{\pm e_i\pm e_j:p+1\leq i\leq p+s,\,p+s+1\leq j\leq n\}\\

\mfn=\mfn^1\oplus\mfn^2\\

\R_{\mfn^1}=\{\pm e_i:p+1\leq i\leq p+s;\,\pm e_i\pm e_j:1\leq i\leq p,\,p+1\leq j\leq
p+s\}\\

\R_{\mfn^2}=\{\pm e_i:p+s+1\leq i\leq n;\,\pm e_i\pm e_j:1\leq
i\leq p,\,p+s+1\leq j\leq n\}\ear$$

$$\bar{cccccccc}\mfn^i & \phi\in\R_{\mfn^i} & \al\in\R_{\mfp}^+  & \phi+n\al & d_{\al\phi} & \textrm{No of } \al's & |\al|^2 & b^{\phi}\\ \hline\hline

\xstrut \mfn^1 & e_{p+1} & e_{p+1}\pm e_i,\,p+s+1\leq i\leq n &
\phi,\,\phi-\al & 1 & 2(n-p-s) & \frac{1}{2n-1} &
\frac{n-p-s}{2n-1}\\ \hline

\xstrut \mfn^2 & e_n & \pm e_i+e_n,\,p+1\leq i\leq p+s &
\phi,\,\phi+\al & 1 & 2s & \frac{1}{2n-1} &
\frac{s}{2n-1}\ear$$\ebst

\li

 \bbst\label{cpbn3}
$(\mfso_{2n+1},\mfso_{2p+1}\oplus\mfso_{2(n-p)},\mfso_{2p+1}\oplus\mfu_{n-p})$,
$0\leq p\leq n-1$. (Type I)

$$\bar{l}\R_{\mfp}=\{\pm (e_i+ e_j):p+1\leq i<j\leq n\}\\

\mfn \textrm{ irreducible }Ad\,L\textrm{-module}\ear$$

$$\bar{cccccccc} \phi\in\R_{\mfn} & \al\in\R_{\mfp}^+  & \phi+n\al & d_{\al\phi} & \textrm{No of } \al's & |\al|^2 & b^{\phi}\\ \hline\hline

\xstrut  e_n & e_i+e_n,\,p+1\leq i\leq n-1 & \phi,\,\phi-\al & 1 &
n-p-1 & \frac{1}{2n-1} & \frac{n-p-1}{2(2n-1)}\ear$$\ebst

\li

\bbst\label{cpbn4}
$(\mfso_{2n+1},\mfso_{2p+1}\oplus\mfso_{2(n-p)},\mfso_{2l+1}\oplus\mfso_{2(p-l)}\oplus\mfu_{n-p})$, $0\leq p\leq n-1$, $0\leq l\leq p-1$. (Type II)

$$\bar{lll}R_{\mfp_1}=\{\pm e_i: l+1\leq i\leq p,\,\pm e_i\pm e_j:1\leq i\leq l,\,l+1\leq j\leq
p\}\\

\R_{\mfp_2}=\{\pm (e_i+ e_j):p+1\leq i<j\leq n\}\\

\mfn=\mfn^1\oplus\mfn^2\ear$$

$$\bar{lcc}\R_{\mfn^i} & b_1^{\phi} & b_2^{\phi} \\ \hline
\{\pm e_i:p+1\leq i\leq n;\,\pm e_i\pm e_j:1\leq i\leq l,\,p+1\leq j\leq
n\} & \frac{p-l}{2n-1} & \frac{n-p-1}{2(2n-1)}\\

\{\pm e_i\pm e_j:l+1\leq i\leq p,\,p+1\leq j\leq
n\} & \frac{4l+1}{4(2n-1)} & \frac{n-p-1}{2(2n-1)}\ear$$\ebst

\li

\bbst\label{cpbn5}
$(\mfso_{2n+1},\mfso_{2p+1}\oplus\mfso_{2(n-p)},\mfso_{2l+1}\oplus\mfso_{2(p-l)}\oplus\mfso_{2s}\oplus\mfso_{2(n-p-s)})$,
$0\leq p\leq n-1$, $0\leq l\leq p-1$, $1\leq s\leq n-p-1$. (Type II)

$$\bar{lll}\R_{\mfp_1}=\{\pm e_i: l+1\leq i\leq p,\,\pm e_i\pm e_j:1\leq i\leq l,\,l+1\leq j\leq
p\},\\

\R_{\mfp_2}=\{\pm e_i\pm e_j:p+1\leq i\leq p+s,\,p+s+1\leq j\leq
n\}\\

\mfn=\mfn^1\oplus\mfn^2\oplus\mfn^3\oplus\mfn^4\ear$$

$$\bar{lcc}\R_{\mfn^i} & b_1^{\phi} & b_2^{\phi} \\ \hline
\{\pm e_i:p+1\leq i\leq p+s;\,\pm e_i\pm e_j:1\leq i\leq l,\,p+1\leq j\leq
p+s\} & \frac{p-l}{2n-1} & \frac{n-p-s}{2n-1}\\

\{\pm e_i:p+s+1\leq i\leq n;\,\pm e_i\pm e_j:1\leq i\leq
l,\,p+s+1\leq j\leq n\} & \frac{p-l}{2n-1} & \frac{s}{2n-1}\\

\{\pm e_i\pm e_j:l+1\leq i\leq p,\,p+1\leq j\leq
p+s\} & \frac{4l+1}{4(2n-1)} & \frac{n-p-s}{2n-1}\\

\{\pm e_i\pm e_j:l+1\leq i\leq p,\,p+s+1\leq j\leq
n\} & \frac{4l+1}{4(2n-1)} & \frac{s}{2n-1}
\ear$$
\ebst

}}

\newpage

\section{$D_n$}In this Section we consider all the bisymmetric
triples of the form $(\mfso_{2n},\mfu_n,\mfl)$ and
$(\mfso_{2n},\mfso_{2p}\oplus\mfso_{2n-2p},\mfl)$, $1\leq p\leq
n-1$.

A root system for $\mfg$ is

\beq\label{rldn}\R=\{\pm e_i\pm e_j:1\leq i<j\leq n\}\eeq

and the length of any root is

\beq\label{rldn}|\al|^2=\frac{1}{2(n-1)}.\eeq

\li

{\footnotesize{

\bsymp $(\mfso_{2n},\mfu_n)$.

$$\bar{ccc} \mfk & \R_{\mfk} & \ga\\ \hline\hline

\mfu_n & \{\pm (e_i- e_j):1\leq
i<j\leq n\} & \frac{n}{2(n-1)}\\
\ear$$

\li

$$\R_{\mfn}=\{\pm (e_i+e_j):1\leq
i<j\leq n\}$$\esymp

\li

\bsymp $(\mfso_{2n},\mfso_{2p}\oplus\mfso_{2n-2p})$,
$p=1,\ldots,n-1$.

$$\bar{ccc} \mfk_i & \R_{\mfk_i} & \ga_i\\ \hline\hline

\mfso_{2p} & \{\pm e_i\pm e_j:1\leq i<j\leq p\} & \frac{p-1}{n-1},\,p\geq 2\\

\mfso_{2n-2p} & \{\pm e_i\pm e_j:p+1\leq i<j\leq n\} &
\frac{n-p-1}{n-1},\,p\leq n-2\ear$$

\li

$$\R_{\mfn}=\{\pm e_i\pm e_j:1\leq i\leq p ,\,p+1\leq
j\leq n\}$$\esymp

\li

\bbst\label{cpdn1} $(\mfso_{2n},\mfu_n,\mfu_p\oplus\mfu_{n-p})$,
$0\leq p\leq n-1$. (Type I)

$$\bar{l}\R_{\mfp}=\{\pm (e_i- e_j):1\leq i\leq p,\,p+1\leq j\leq n\}\\

\mfn=\mfn^1\oplus\mfn^2\oplus\mfn^3\\

\R_{\mfn^1}=\{\pm(e_i+e_j):1\leq i<j\leq p\},\\

\R_{\mfn^2}=\{\pm(e_i+e_j):p+1\leq i<j\leq n,\},\\

\R_{\mfn^3}=\{\pm(e_i+e_j):1\leq i\leq p,\,p+1\leq j\leq n\}\ear$$

$$\bar{cccccccc}\mfn^i & \phi\in\R_{\mfn^i} & \al\in\R_{\mfp}^+  & \phi+n\al & d_{\al\phi} & \textrm{No of } \al'ss & |\al|^2 & b^{\phi}\\ \hline\hline

\xstrut \mfn^1 & e_1+e_p & \bar{c}e_1-e_i,\,p+1\leq i\leq
n\\e_p-e_i,\,p+1\leq i\leq n\ear & \phi,\,\phi-\al &
1 & \bar{c}n-p\\n-p\ear & \frac{1}{2(n-1)} & \frac{n-p}{2(n-1)}\\
\hline

\xstrut \mfn^2 & e_{p+1}+e_n & \bar{c} e_i-e_{p+1},\,1\leq i\leq
p\\ e_i-e_n,\,1\leq i\leq p\ear & \phi,\,\phi+\al
& 1 & \bar{c}p\\p\ear & \frac{1}{2(n-1)} & \frac{p}{2(n-1)}\\
\hline

\xstrut \mfn^2 & e_1+e_n & \bar{c} e_1-e_i,\,p+1\leq i\leq n-1\\
e_i-e_n,\,2\leq i\leq p\ear & \bar{c}\phi,\,\phi-\al\\
\phi,\,\phi+\al\ear & 1 & \bar{c}n-p-1\\p-1\ear & \frac{1}{2(n-1)}
& \frac{n-2}{4(n-1)}

\ear$$\ebst

\li

\bbst\label{cpdn2}
$(\mfso_{2n},\mfso_{2p}\oplus\mfso_{2(n-p)},\mfso_{2l}\oplus\mfso_{2(p-l)}\oplus\mfso_{2(n-p)})$, $1\leq p\leq n-1$, $1\leq l\leq p-1$. (Type I)

$$\bar{l}\R_{\mfp}=\{\pm e_i\pm e_j:1\leq i\leq l,\,l+1\leq j\leq p\}\\

\mfn=\mfn^1\oplus\mfn^2\\

\R_{\mfn^1}=\{\pm e_i\pm e_j:1\leq i\leq l,\,p+1\leq j\leq
n\}\\

\R_{\mfn^2}=\{\pm e_i\pm e_j:l+1\leq i\leq p,\,p+1\leq j\leq
n\} \ear$$

$$\bar{cccccccc}\mfn^i & \phi\in\R_{\mfn^i} & \al\in\R_{\mfp}^+  & \phi+n\al & d_{\al\phi} & \textrm{No of } \al's & |\al|^2 & b^{\phi}\\ \hline\hline

\xstrut \mfn^1 & e_1+e_n & e_1\pm e_i,\,l+1\leq i\leq p &
\phi,\,\phi-\al & 1 & 2(p-l) & \frac{1}{2(n-1)} &
\frac{p-l}{2(n-1)}\\ \hline

\xstrut \mfn^2 & e_p+e_n & \pm e_i+e_p,\,1\leq i\leq l &
\phi,\,\phi+\al & 1 & l & \frac{1}{2(n-1)} &
\frac{l}{2(n-1)}\ear$$\ebst

\li

\bbst\label{cpdn3}
$(\mfso_{2n},\mfs_{2p}\oplus\mfso_{2(n-p)},\mfso_{2p}\oplus\mfso_{2s}\oplus\mfso_{2(n-p-s)})$, $1\leq p\leq n-1$, $1\leq s\leq n-p-1$. (Type I)

$$\bar{l}\R_{\mfp}=\{\pm e_i\pm e_j:p+1\leq i\leq p+s,\,p+s+1\leq j\leq n\}\\

\mfn=\mfn^1\oplus\mfn^2\\

\R_{\mfn^1}=\{\pm e_i\pm e_j:1\leq i\leq p,\,p+1\leq j\leq
p+s\}\\

\R_{\mfn^2}=\{\pm e_i\pm e_j:1\leq i\leq p,\,p+s+1\leq j\leq
n\}\ear$$

$$\bar{cccccccc}\mfn^i & \phi\in\R_{\mfn^i} & \al\in\R_{\mfp}^+  & \phi+n\al & d_{\al\phi} & \textrm{No of } \al's & |\al|^2 & b^{\phi}\\ \hline\hline

\xstrut \mfn^1 & e_1+e_{p+1} & e_{p+1}\pm e_i,\,p+s+1\leq i\leq n
& \phi,\,\phi-\al & 1 & 2(n-p-s) & \frac{1}{2(n-1)} &
\frac{n-p-s}{2(n-1)}\\ \hline

\xstrut \mfn^2 & e_1+e_n & \pm e_i+e_n,\,p+1\leq i\leq p+s &
\phi,\,\phi-\al & 1 & 2s & \frac{1}{2(n-1)} &
\frac{s}{2(n-1)}\ear$$\ebst

\li

\bbst\label{cpdn4}
$(\mfso_{2n},\mfs_{2p}\oplus\mfso_{2(n-p)},\mfso_{2l}\oplus\mfso_{2(p-l)}\oplus\mfso_{2s}\oplus\mfso_{2(n-p-s)})$,
$1\leq p\leq n-1$, $1\leq l\leq p-1$, $1\leq s\leq n-p-1$ . (Type
II)

$$\bar{lll}\R_{\mfp_1}=\{\pm e_i\pm e_j:1\leq i\leq l,\,l+1\leq j\leq
p\}\\

\R_{\mfp_2}=\{\pm e_i\pm e_j:p+1\leq i\leq p+s,\,p+s+1\leq j\leq
n\}\\

\mfn=\mfn^1\oplus\mfn^2\oplus\mfn^3\oplus\mfn^4\ear$$

$$\bar{lcc}\R_{\mfn^i} & b_1^{\phi} & b_2^{\phi} \\ \hline
\{\pm e_i\pm e_j:1\leq i\leq l,\,p+1\leq j\leq p+s\} & \frac{p-l}{2(n-1)} & \frac{n-p-s}{2(n-1)}\\

\{\pm e_i\pm e_j:1\leq i\leq 1,\,p+s+1\leq j\leq n\} & \frac{p-l}{2(n-1)} & \frac{s}{2(n-1)}\\

\{\pm e_i\pm e_j:l+1\leq i\leq p,\,p+1\leq j\leq
p+s\} & \frac{l}{2(n-1)}& \frac{n-p-s}{2(n-1)}\\

\{\pm e_i\pm e_j:l+1\leq i\leq p,\,p+s+1\leq j\leq
n\}& \frac{l}{2(n-1)}& \frac{s}{2(n-1)}
\ear$$\ebst

\li

\bbst\label{cpdn5}
$(\mfso_{2n},\mfso_{2p}\oplus\mfso_{2(n-p)},\mfu_{p}\oplus\mfso_{2(n-p)})$,
$1\leq p\leq n-1$. (Type I)

$$\bar{l}\R_{\mfp}=\{\pm (e_i+ e_j):1\leq i<j\leq p\}\\

\mfn\textrm{ irreducible }Ad\,L\textrm{-module}\ear$$

$$\bar{cccccccc} \phi\in\R_{\mfn} & \al\in\R_{\mfp}^+  & \phi+n\al & d_{\al\phi} & \textrm{No of } \al's & |\al|^2 & b^{\phi}\\ \hline\hline

\xstrut e_1+e_n & e_1+ e_i,\,2\leq i\leq p & \phi,\,\phi-\al & 1 &
p-1 & \frac{1}{2(n-1)} & \frac{p-1}{4(n-1)}\\ \ear$$\ebst

\li

\bbst\label{cpdn6}
$(\mfso_{2n},\mfso_{2p}\oplus\mfso_{2(n-p)},\mfso_{2p}\oplus\mfu_{n-p})$, $1\leq p\leq n-1$. (Type I)

$$\bar{l}\R_{\mfp}=\{\pm (e_i+ e_j):p+1\leq i<j\leq n\}\\

\mfn\textrm{ irreducible }Ad\,L\textrm{-module}\ear$$

$$\bar{cccccccc} \phi\in\R_{\mfn} & \al\in\R_{\mfp}^+  & \phi+n\al & d_{\al\phi} & \textrm{No of } \al's & |\al|^2 & b^{\phi}\\ \hline\hline

\xstrut e_1+e_n & e_i+ e_n,\,p+1\leq i\leq n-1 & \phi,\,\phi-\al
& 1 & n-p-1 & \frac{1}{2(n-1)} & \frac{n-p-1}{4(n-1)}\\
\ear$$\ebst

\li

\bbst\label{cpdn7}
$(\mfso_{2n},\mfso_{2p}\oplus\mfso_{2(n-p)},\mfu_p\oplus\mfu_{n-p})$, $1\leq p\leq n-1$. (Type II)

$$\bar{lll}\R_{\mfp_1}=\{\pm (e_i+ e_j):1\leq i<j\leq p\}\\

\R_{\mfp_2}=\{\pm (e_i+ e_j):p+1\leq i<j\leq n\}\\

\mfn\textrm{ irreducible }Ad\,L\textrm{-module}\ear$$

$$\bar{cc} b_1^{\phi} & b_2^{\phi} \\ \hline
 \frac{p-1}{4(n-1)} & \frac{n-p-1}{4(n-1)}\ear$$\ebst

\li

\bbst\label{cpdn8}
$(\mfso_{2n},\mfso_{2p}\oplus\mfso_{2(n-p)},\mfso_{2l}\oplus\mfso_{2(p-l)}\oplus\mfu_{n-p})$,
$1\leq p\leq n-1$, $1\leq l\leq p-1$. (Type II)

$$\bar{lll}\R_{\mfp_1}=\{\pm e_i\pm e_j:1\leq i\leq l,\,l+1\leq j\leq p\}\\

\R_{\mfp_2}=\R_{\mfp_2}=\{\pm (e_i+ e_j):p+1\leq i<j\leq n\}\\

\mfn=\mfn^1\oplus\mfn^2\ear$$

$$\bar{lcc}\R_{\mfn^i} & b_1^{\phi} & b_2^{\phi} \\ \hline

\{\pm e_i\pm e_j:1\leq i\leq l,\,p+1\leq j\leq n\} & \frac{p-l}{2(n-1)} & \frac{n-p-1}{4(n-1)}\\

\{\pm e_i\pm e_j:l+1\leq i\leq p,\,p+1\leq j\leq n\} & \frac{l}{2(n-1)} & \frac{n-p-1}{4(n-1)}\\

\ear$$\ebst

}}

\newpage

\section{$C_n$} We consider the bisymmetric triples of the form
$(\mfsp_n,\mfu_n,\mfl)$ and
$(\mfsp_n,\mfsp_p\oplus\mfsp_{n-p},\mfl)$, for $1\leq p\leq n-2$.
The root system for $\mfg=\mfsp_n$ is

\beq\label{rcn}\R=\{\pm 2e_i:1\leq i\leq n;\,\pm e_i\pm e_j:1\leq
i<j\leq n\}.\eeq

In $\mfg$ there are two root lengths:

\beq\label{rlcn}|\al|^2=\left\{\bar{ll}\frac{1}{n+1}, & \al=\pm
2e_i\\ \frac{1}{2(n+1)}, & \al=\pm e_i\pm e_j\ear.\right.\eeq

{\footnotesize{

\bsymp $(\mfsp_n,\mfu_n)$.

$$\bar{ccc} \mfk & \R_{\mfk} & \ga\\ \hline\hline

\mfu_n & \{\pm (e_i- e_j):1\leq
i<j\leq n\} & \frac{n}{2(n+1)}\\
\ear$$

\li

$$\R_{\mfn}=\{\pm 2e_i,\,\pm (e_i+e_j):1\leq
i<j\leq n\}$$\esymp

\li

\bsymp$(\mfsp_n,\mfsp_p\oplus\mfsp_{n-p})$, $p=1,\ldots,n-2$.

$$\bar{ccc} \mfk_i & \R_{\mfk_i} & \ga_i\\ \hline\hline

\mfsp_p & \{\pm 2e_i:1\leq i\leq p;\,\pm e_i\pm e_j:1\leq i<j\leq
p\} & \frac{p+1}{n+1}\\

\mfsp_{n-p} & \{\pm 2e_i:p+1\leq i\leq n;\,\pm e_i\pm e_j:p+1\leq
i<j\leq n\} & \frac{n-p+1}{n+1}\ear$$

\li

$$\R_{\mfn}=\{\pm e_i\pm e_j:1\leq i\leq p,\,p+1\leq
j\leq n\}$$\esymp

\bbst\label{cpcn1}
$(\mfsp_n,\mfu_n,\mfu_p\oplus\mfu_{n-p})$, $1\leq p\leq n-1$. (Type I)

$$\bar{l}\R_{\mfp}=\{\pm (e_i- e_j):1\leq i\leq p,\,p+1\leq j\leq n\}\\

\mfn=\mfn^1\oplus\mfn^2\oplus\mfn^3\\

\R_{\mfn^1}=\{\pm 2e_i,\,\pm(e_i+e_j):1\leq i<j\leq p\},\\

\R_{\mfn^2}=\{\pm 2e_i,\,\pm(e_i+e_j):p+1\leq i<j\leq n,\},\\

\R_{\mfn^3}=\{\pm(e_i+e_j):1\leq i\leq p,\,p+1\leq j\leq n\}\ear$$

$$\bar{cccccccc} \mfn^i & \phi\in\R_{\mfn^i} & \al\in\R_{\mfp}^+  & \phi+n\al & d_{\al\phi} & \textrm{No of } \al's & |\al|^2 & b^{\phi}\\ \hline\hline

\xstrut\mfn^1 &  2e_1 & e_1-e_i,\,p+1\leq i\leq n &
\phi,\,\phi-\al & 1 & n-p & \frac{1}{2(n+1)} &
\frac{n-p}{4(n+1)}\\ \hline

 \xstrut\mfn^2 &  2e_n & e_i-e_n,\,1\leq i\leq p & \phi,\,\phi+\al & 1 & p & \frac{1}{2(n+1)} & \frac{p}{4(n+1)}\\ \hline

 \xstrut\mfn^3 &  e_1+e_n & \bar{c} e_i-e_n,\,2\leq i\leq p\\e_1-e_i,\,p+1\leq i\leq n-1\\e_1-e_n\ear & \bar{c}\phi,\,\phi+\al\\ \phi,\,\phi-\al\\ \phi,\,\phi\pm\al\ear & \bar{c}1\\1\\4\ear & \bar{c}p-1\\n-p-1\\1\ear & \frac{1}{2(n+1)} & \frac{n+2}{4(n+1)} \ear$$\ebst

\li

\bbst\label{cpcn2}
$(\mfsp_n,\mfsp_p\oplus\mfsp_{n-p},\mfsp_l\oplus\mfsp_{p-l}\oplus\mfsp_{n-p})$,
$1\leq p\leq n-1$, $1\leq l\leq p-1$. (Type I)

$$\bar{l}\R_{\mfp}=\{\pm e_i\pm e_j:1\leq i\leq l,\,l+1\leq j\leq p\}\\

\mfn=\mfn^1\oplus\mfn^2\\

\R_{\mfn^1}=\{\pm e_i\pm e_j:1\leq i\leq l,\,p+1\leq j\leq
n\}\\

\R_{\mfn^2}=\{\pm e_i\pm e_j:l+1\leq i\leq p,\,p+1\leq j\leq
n\}\ear$$

$$\bar{cccccccc}\mfn^i & \phi\in\R_{\mfn^i} & \al\in\R_{\mfp}^+  & \phi+n\al & d_{\al\phi} & \textrm{No of } \al's & |\al|^2 & b^{\phi}\\ \hline\hline

\xstrut \mfn^1 & e_1+e_n & e_1\pm e_i,\,l+1\leq i\leq p &
\phi,\,\phi-\al & 1 & p-l & \frac{1}{2(n+1)} &
\frac{p-l}{4(n+1)}\\ \hline

\xstrut \mfn^2 & e_p+e_n & \pm e_i+e_p,\,1\leq i\leq l &
\phi,\,\phi-\al & 1 & l & \frac{1}{2(n+1)} &
\frac{l}{4(n+1)}\ear$$\ebst

\li

\bbst\label{cpcn3}
$(\mfsp_n,\mfsp_p\oplus\mfsp_{n-p},\mfsp_p\oplus\mfsp_{s}\oplus\mfsp_{n-p-s})$,
$1\leq p\leq n-1$, $1\leq s\leq n-p-1$. (Type I)

$$\bar{l}\R_{\mfp}=\{\pm e_i\pm e_j:p+1\leq i\leq p+s,\,p+s+1\leq j\leq n\}\\

\mfn=\mfn^1\oplus\mfn^2\\

\R_{\mfn^1}=\{\pm e_i\pm e_j:1\leq i\leq p,\,p+1\leq j\leq
p+s\}\\

\R_{\mfn^2}=\{\pm e_i\pm e_j:l\leq i\leq p,\,p+s+1\leq j\leq
n\}\ear$$

$$\bar{cccccccc}\mfn^i & \phi\in\R_{\mfn^i} & \al\in\R_{\mfp}^+  & \phi+n\al & d_{\al\phi} & \textrm{No of } \al's & |\al|^2 & b^{\phi}\\ \hline\hline

\xstrut \mfn^1 & e_1+e_{p+1} & e_{p+1}\pm e_i,\,p+s+1\leq i\leq n
& \phi,\,\phi-\al & 1 & n-p-s & \frac{1}{2(n+1)} &
\frac{n-p-s}{4(n+1)}\\ \hline

\xstrut \mfn^2 & e_p+e_n & \pm e_i+e_n,\,p+1\leq i\leq p+s &
\phi,\,\phi-\al & 1 & s & \frac{1}{2(n+1)} &
\frac{s}{4(n+1)}\ear$$\ebst

\li

 \bbst\label{cpcn4}
$(\mfsp_n,\mfsp_p\oplus\mfsp_{n-p},\mfsp_l\oplus\mfsp_{p-l}\oplus\mfsp_s\oplus\mfsp_{n-p-s})$,
$1\leq p\leq n-1$, $1\leq l\leq p-1$, $1\leq s\leq n-p-1$. (Type
II)

$$\bar{lll}\R_{\mfp_1}=\{\pm e_i\pm e_j:1\leq i\leq l,\,l+1\leq j\leq
p\}\\ \\

\R_{\mfp_2}=\{\pm e_i\pm e_j:p+1\leq i\leq p+s,\,p+s+1\leq j\leq
n\}\\

\mfn=\mfn^1\oplus\mfn^2\oplus\mfn^3\oplus\mfn^4\ear$$

$$\bar{lcc}\R_{\mfn^i} & b_1^{\phi} & b_2^{\phi} \\ \hline
\{\pm e_i\pm e_j:1\leq i\leq l,\,p+1\leq j\leq
p+s\} & \frac{p-l}{4(n+1)} & \frac{n-p-s}{4(n+1)}\\

\{\pm e_i\pm e_j:1\leq i\leq l,\,p+s+1\leq j\leq n\}& \frac{p-l}{4(n+1)} & \frac{s}{4(n+1)}\\

\{\pm e_i\pm e_j:l+1\leq i\leq p,\,p+1\leq j\leq
p+s\} & \frac{l}{4(n+1)} & \frac{n-p-s}{4(n+1)}\\

\{\pm e_i\pm e_j:l+1\leq i\leq p,\,p+s+1\leq j\leq
n\} & \frac{l}{4(n+1)} & \frac{s}{4(n+1)}
\ear$$\ebst

\li

\bbst\label{cpcn5}
$(\mfsp_n,\mfsp_p\oplus\mfsp_{n-p},\mfu_p\oplus\mfsp_{n-p})$, $1\leq p\leq n-1$. (Type I)

$$\bar{l}\R_{\mfp}=\{\pm 2e_i:1\leq i\leq p;\,\pm (e_i+ e_j):1\leq i<j\leq
p\}\\

\mfn\textrm{ irreducible }Ad\,L\textrm{-module}\ear$$

$$\bar{cccccccc} \phi\in\R_{\mfn} & \al\in\R_{\mfp}^+  & \phi+n\al & d_{\al\phi} & \textrm{No of } \al's & |\al|^2 & b^{\phi}\\ \hline\hline

\xstrut e_1+e_n & \bar{c}2e_1\\e_1+e_i,\,2\leq i\leq p\ear &
\phi,\,\phi-\al & 1 & \bar{c}1\\p-1\ear &
\bar{c}\frac{1}{n+1}\\\frac{1}{2(n+1)}\ear & \frac{p+1}{4(n+1)}\\
\ear$$\ebst

\li

 \bbst\label{cpcn6}
$(\mfsp_n,\mfsp_p\oplus\mfsp_{n-p},\mfsp_p\oplus\mfu_{n-p})$, $1\leq
p\leq n-1$. (Type I)

$$\bar{l}\R_{\mfp}=\{\pm 2e_i:p+1\leq i\leq n;\,\pm (e_i+ e_j):p+1\leq i<j\leq
n\}\\

\mfn\textrm{ irreducible }Ad\,L\textrm{-module}\ear$$

$$\bar{cccccccc} \phi\in\R_{\mfn} & \al\in\R_{\mfp}^+  & \phi+n\al & d_{\al\phi} & \textrm{No of } \al's & |\al|^2 & b^{\phi}\\ \hline\hline

\xstrut e_1+e_n & \bar{c}2e_n\\e_i+e_n,\,p+1\leq i\leq n-1\ear &
\phi,\,\phi-\al & 1 & \bar{c}1\\n-p-1\ear &
\bar{c}\frac{1}{n+1}\\\frac{1}{2(n+1)}\ear &
\frac{n-p+1}{4(n+1)}\\ \ear$$\ebst

\li

 \bbst\label{cpcn7}
$(\mfsp_n,\mfsp_p\oplus\mfsp_{n-p},\mfu_p\oplus\mfu_{n-p})$,  $1\leq
p\leq n-1$. (Type II)

$$\bar{lll}\R_{\mfp_1}=\{\pm 2e_i:1 \leq i\leq p;\,\pm (e_i+ e_j):1\leq i<j\leq
p\}\\

\R_{\mfp_2}=\{\pm 2e_i:p+1\leq i\leq n;\,\pm (e_i+ e_j):p+1\leq
i<j\leq n\}\\

\mfn\textrm{ irreducible }Ad\,L\textrm{-module}\ear$$

$$\bar{cc} b_1^{\phi} & b_2^{\phi} \\ \hline
 \frac{p+1}{4(n+1)} & \frac{n-p+1}{4(n+1)}
\ear$$\ebst

\li

 \bbst\label{cpcn8}
$(\mfsp_n,\mfsp_p\oplus\mfsp_{n-p},\mfsp_l\oplus\mfsp_{p-l}\oplus\mfu_{n-p})$,
$1\leq p\leq n-1$ and $1\leq l\leq p-1$. (Type II)

$$\bar{lll}\R_{\mfp_1}=\{\pm e_i\pm e_j:1\leq i\leq l,\,l+1\leq j\leq
p\}\\

\R_{\mfp_2}=\{\pm 2e_i:p+1\leq i\leq n;\,\pm (e_i+ e_j):p+1\leq
i<j\leq n\}\\

\mfn=\mfn^1\oplus\mfn^2\ear$$

$$\bar{lcc}\R_{\mfn^i} & b_1^{\phi} & b_2^{\phi} \\ \hline
\{\pm e_i\pm e_j:1\leq i\leq l,\,p+1\leq j\leq
n\} & \frac{p-l}{4(n+1)} & \frac{n-p+1}{4(n+1)}\\

\{\pm e_i\pm e_j:l+1\leq i\leq p,\,p+1\leq j\leq
n\}& \frac{l}{4(n+1)} & \frac{n-p+1}{4(n+1)}
\ear$$\ebst

}}

\newpage

\section{$\mff_4$} In this section we analyze the bisymmetric
triples of the form $(\mff_4,\mfso_9,\mfl)$ and
$(\mff_4,\mfsp_3\oplus\mfsu_2,\mfl)$.

The root system for the simple Lie algebra $\mff_4$ is

\beq\label{rf4}\R=\{\pm e_i,\,\pm e_i\pm e_j, 1\leq i<j\leq 4;\,
\frac{1}{2}\sum_1^4(-1)^{\nu_i}e_i\},\eeq

where $e_1,\ldots,e_4$ is the canonical basis for $\reals^4$ and the
signs are chosen independently. In $\mff_4$, there are two root
lengths. Roots of the form $\pm e_i$ and
$\frac{1}{2}\sum_1^4(-1)^{\nu_i}e_i$ are short, whereas those of the
form $\pm e_i\pm e_j$ are long, and we have

\beqar |\al|^2=\left\{\bar{ll}\frac{1}{18}, \,\al\textrm{ short }\\
\frac{1}{9},\, \al\textrm{ long }\ear\right. .\eeqar

\li

{\footnotesize{

\bsymp $(\mff_4,\mfso_9)$.

$$\bar{ccc}\mfk & \R_{\mfk} & \ga\\ \hline\hline

\mfso_9 & \{\pm e_i,\,\pm e_i\pm e_j, 1\leq i<j\leq 4 \} &
\frac{7}{9}\ear$$

$$\R_{\mfn}=\{\frac{1}{2}\sum_1^4(-1)^{\nu_i}e_i \}$$
\esymp

\li

\bsymp $(\mff_4,\mfsp_3\oplus\mfsu_2)$.

$$\bar{ccc}\mfk_i & \R_{\mfk_i} & \ga_i\\ \hline\hline

\mfsp_3 & \bar{l}<e_4, e_3-e_4,\frac{1}{2}(e_1-e_2-e_3-e_4)>=\\\{\pm
e_3,\pm e_4,\pm e_3\pm
e_4,\,\pm(e_1-e_2),\,\pm\frac{1}{2}(e_1-e_2\pm e_3\pm e_4)\}\ear &
\frac{4}{9}\\\hline \mfsu_2 & \{\pm(e_1+e_2)\} & \frac{2}{9}\ear$$

$$\R_{\mfn}=\{\pm e_i,\,\pm e_i\pm e_j:i=1,2,j=3,4;
\pm\frac{1}{2}(e_1+e_2\pm e_3\pm e_4)\}, \textrm{ with signs chosen
independently.}$$

\esymp

\li

\bbst\label{cpf41} $(\mff_4,\mfso_9, \mfso_p\oplus\mfso_{9-p})$,
$p=2l+1$, $l=0,1,2,3$. (Type I)

$$\bar{l}\R_{\mfp}=\{\pm 2e_i:p+1\leq i\leq n;\,\pm (e_i+ e_j):p+1\leq i<j\leq
n\}\\

C_{\mfp}\textrm{ scalar on }\mfn\ear$$

$$\bar{cccccccc} \phi\in\R_{\mfn} & \al\in\R_{\mfp}^+  & \phi+n\al & d_{\al\phi} & \textrm{No of } \al's & |\al|^2 & b^{\phi}\\ \hline\hline

\xstrut \frac{1}{2}\sum_1^4e_i & \bar{c}e_i,\,l+1\leq i\leq
4\\e_i+e_j,\,1\leq i\leq l,l+1\leq j\leq 4\ear & \phi,\,\phi-\al &
1 & \bar{c}4-l\\l(4-l)\ear &
\bar{c}\frac{1}{18}\\\frac{1}{9}\ear & \frac{p(9-p)}{72}\\
\ear$$\ebst

\li

\bbst\label{cpf42}
 $(\mff_4,\mfsp_3\oplus\mfsu_2,
\mfsp_3\oplus\reals)$. (Type I)

$$\bar{l}\R_{\mfp}=\{\pm (e_1+e_2) \}\\

C_{\mfp}\textrm{ scalar on }\mfn\ear$$

$$\bar{ccccccc} \phi\in\R_{\mfn} & \al\in\R_{\mfp}^+  & \phi+n\al & d_{\al\phi} & \textrm{No of } \al's & |\al|^2 & b^{\phi}\\ \hline\hline

\xstrut  e_1 & e_1+e_2 & \phi,\,\phi-\al & 1 & 1 & \frac{1}{9}&
\frac{1}{18} \ear$$\ebst

\li

\bbst\label{cpf43} $(\mff_4,\mfsp_3\oplus\mfsu_2,
\mfu_3\oplus\mfsu_2)$. (Type I)

$$\bar{l}\R_{\mfp}=\{\pm e_3,\,\pm e_3\pm
e_4,\,\pm(e_1-e_2),\,\pm\frac{1}{2}(e_1-e_2+e_3\pm e_4)\}\\

\mfn=\mfn^1\oplus\mfn^2\\
 \R_{\mfn^1}=  \{\pm(e_1+e_3),\,\pm(e_2-e_3)\}\\
\R_{\mfn^2}=  \{\pm e_1,\,\pm e_2,\,\pm(e_1-e_3),\,\pm(e_1+e_4),
\pm(e_2\pm e_3),\,\pm(e_2\pm e_4),\,\pm\frac{1}{2}(e_1+e_2\pm e_3\pm
e_4)\}\ear$$

$$\bar{cccccccc}\mfn^i & \phi\in\R_{\mfn^i} & \al\in\R_{\mfp}^+  & \phi+n\al & d_{\al\phi} & \textrm{No of } \al's & |\al|^2 & b^{\phi}\\ \hline\hline

\xstrut \mfn^1 & e_1+e_3 &
\bar{c}e_3,\,\frac{1}{2}\pm(e_1-e_2+e_3\pm e_4)\\e_3\pm
e_4,\,e_1-e_2\ear &
\phi,\,\phi-\al & 1 & \bar{c}3\\3\ear & \bar{c}\frac{1}{18}\\\frac{1}{9}\ear & \frac{1}{4}\\
\hline

\xstrut \mfn^2 & e_1 & \bar{c}e_3\\e_1-e_2\\
\frac{1}{2}\pm(e_1-e_2+e_3\pm e_4)\ear &
\bar{c}\phi,\,\phi\pm\al\\\phi,\,\phi-\al\\\phi,\,\phi-\al\ear &
\bar{c}4\\1\\1\ear & \bar{c}1\\1\\2\ear &
\bar{c}\frac{1}{18}\\\frac{1}{9}\\\frac{1}{18}\ear &
\frac{2}{9}\ear$$\ebst

\li

 \bbst\label{cpf44} $(\mff_4,\mfsp_3\oplus\mfsu_2,
\mfsp_2\oplus\mfsu_2\oplus\mfsu_2)$. (Type I)

$$\bar{l}\R_{\mfp}=\{\pm e_3,\,\pm e_4,\,\pm\frac{1}{2}(e_1-e_2\pm (e_3+e_4)\}\\

\mfn=\mfn^1\oplus\mfn^2\\
 \R_{\mfn^1}= \{\pm e_i\pm e_j,\,i=1,2,\,j=3,4;\,\pm\frac{1}{2}(e_1+e_2\pm
(e_3+ e_4))\}\\
\R_{\mfn^2}=  \{\pm e_1,\,\pm e_2,\,\pm\frac{1}{2}(e_1+e_2\pm (e_3-
e_4))\}\ear$$

$$\bar{cccccccc}\mfn^i & \phi\in\R_{\mfn^i} & \al\in\R_{\mfp}^+  & \phi+n\al & d_{\al\phi} & \textrm{No of } \al's & |\al|^2 & b^{\phi}\\ \hline\hline

\xstrut \mfn^1 & \frac{1}{2}(e_1+e_2+e_3+e_4) &
\bar{c}\frac{1}{2}(e_1-e_2-e_3-e_4)\\\textrm{all others}\ear &
\bar{c}\phi,\,\phi+\al\\\phi,\,\phi-\al\ear & 1 & 4 & \frac{1}{18} & \frac{1}{9}\\
\hline

\xstrut \mfn^2 & e_1 & \bar{c}e_3,\,e_4\\\frac{1}{2}(e_1-e_2\pm(e_3+
e_4))\ear & \bar{c}\phi,\,\phi\pm\al\\\phi,\,\phi-\al\ear &
\bar{c}4\\1\ear & \bar{c}2\\2\ear & \frac{1}{18} &
\frac{5}{18}\ear$$\ebst

\li

\bbst\label{cpf45} $(\mff_4,\mfsp_3\oplus\mfsu_2,
\mfu_3\oplus\reals)$. (Type II)

$$\bar{lll}\R_{\mfp_1}=\{\pm e_3,\,\pm e_3\pm
e_4,\,\pm(e_1-e_2),\,\pm\frac{1}{2}(e_1-e_2+e_3\pm e_4)\}\\
 \R_{\mfp_2}=\{\pm (e_1+e_2) \}\\

\mfn=\mfn^1\oplus\mfn^2\ear$$

$$\bar{lcc}\R_{\mfn^i} & b_1^{\phi} & b_2^{\phi} \\ \hline
\{\pm(e_1+e_3),\,\pm(e_2-e_3)\} & \frac{1}{4} & \frac{1}{18}\\

 \{\pm e_1,\,\pm e_2,\,\pm(e_1-e_3),\,\pm(e_1+e_4),
\pm(e_2\pm e_3),\,\pm(e_2\pm e_4),\,\pm\frac{1}{2}(e_1+e_2\pm e_3\pm
e_4)\}& \frac{2}{9} & \frac{1}{18} \ear$$\ebst

\li

\bbst\label{cpf46} $(\mff_4,\mfsp_3\oplus\mfsu_2,
\mfsp_2\oplus\mfsu_2\oplus\reals)$. (Type II)

$$\bar{lll}\R_{\mfp_1}=\{\pm e_3,\,\pm e_4,\,\pm\frac{1}{2}(e_1-e_2\pm
(e_3+e_4)\}
\R_{\mfp_2}=\{\pm (e_1+e_2) \}\\

\mfn=\mfn^1\oplus\mfn^2\ear$$

$$\bar{lcc}\R_{\mfn^i} & b_1^{\phi} & b_2^{\phi} \\ \hline
\{\pm e_i\pm e_j,\,i=1,2,\,j=3,4;\,\pm\frac{1}{2}(e_1+e_2\pm
(e_3+ e_4))\} & \frac{1}{9} & \frac{1}{18}\\

\{\pm e_1,\,\pm e_2,\,\pm\frac{1}{2}(e_1+e_2\pm (e_3- e_4))\} &
\frac{5}{18} & \frac{1}{18} \ear$$\ebst

}}

\newpage

\section{$\mfg_2$}

In this section we analyze all bisymmetric triples
$(\mfg_2,\mfsu_2\oplus\mfsu_2 ,\mfl)$. The symmetric pair
$(\mfg_2,\mfsu_2\oplus\mfsu_2)$ corresponds to a flag manifold of
$G_2$ and thus is obtained by a painted Dynkin diagram. We observe
that each factor $\mfsu_2$ corresponds to a different root length.
We set that $\mfk_1=\mfsu_2$ corresponds to a long root and
$\mfk_2=\mfsu_2$ corresponds to a short root. If we consider the
root system of $\mfg_2$,

$$\bar{rl}\R= & \{\pm(e_1-e_2),\,\pm(e_1-e_3),\,\pm(e_2-e_3),\\
& \pm(2e_1-e_2-e_3),\,\pm(2e_2-e_1-e_3),\,\pm(2e_3-e_1-e_2)\},\ear$$

we can choose $\R_{\mfk_1}=\{\pm(2e_1-e_2-e_3)\}$ and
$\R_{\mfk_2}=\{\pm(e_2-e_3)\}$, the orthogonal of $\R_{\mfk_1}$.

In $\mfg_2$ the length of a root is given by

\beqar |\al|^2=\left\{\bar{ll}\frac{1}{4}, \,\al\textrm{ long }\\
\frac{1}{12},\, \al\textrm{ short }\ear\right. ,\eeqar

where roots of the form $e_a-e_b$ are short and those of the form
$2e_a-e_b-e_c$ are long.

\li

{\footnotesize{

\bsymp  $(\mfg_2,\mfsu_2\oplus\mfsu_2)$.

$$\bar{ccc}\mfk_i & \R_{\mfk_i} & \ga_i\\ \hline\hline

\mfsu_2 & \{\pm(2e_1-e_2-e_3)\} & \frac{1}{2}\\

 \mfsu_2 & \{\pm(e_2-e_3)\} & \frac{1}{6}\ear$$

$$\R_{\mfn}=\{\pm(e_1-e_2),\,\pm(e_1-e_3),\,\pm(2e_2-e_1-e_3),\,\pm(2e_3-e_1-e_2)\}$$

\esymp

\li

\bbst\label{cpg21} $(\mfg_2,\mfsu_2\oplus\mfsu_2
,\reals\oplus\mfsu_2)$. (Type I)

$$\bar{l}\R_{\mfp}=\{\pm(2e_1-e_2-e_3)\}\\

C_{\mfp}\textrm{ scalar on }\mfn\ear$$

$$\bar{ccccccc} \phi\in\R_{\mfn} & \al\in\R_{\mfp}^+  & \phi+n\al & d_{\al\phi} & \textrm{No of } \al's & |\al|^2 & b^{\phi}\\ \hline\hline

\xstrut  e_1-e_2 & 2e_1-e_2-e_3 & \phi,\,\phi-\al & 1 & 1 &
\frac{1}{4}& \frac{1}{8} \ear$$\ebst

\li

\bbst\label{cpg22} $(\mfg_2,\mfsu_2\oplus\mfsu_2
,\mfsu_2\oplus\reals)$. (Type I)

$$\bar{l}\R_{\mfp}=\{\pm(e_2-e_3)\}\\

C_{\mfp}\textrm{ scalar on }\mfn\ear$$

$$\bar{ccccccc} \phi\in\R_{\mfn} & \al\in\R_{\mfp}^+  & \phi+n\al & d_{\al\phi} & \textrm{No of } \al's & |\al|^2 & b^{\phi}\\ \hline\hline

\xstrut  e_1-e_2 & e_2-e_3 & \phi,\,\phi\pm\al & 4 & 1 &
\frac{1}{12}& \frac{1}{6} \ear$$\ebst

\li

\bbst \label{cpg23}$(\mfg_2,\mfsu_2\oplus\mfsu_2 ,\reals\oplus\reals)$. (Type II)

$$\bar{l}\R_{\mfp_1}=\{\pm(2e_1-e_2-e_3)\}\\

\R_{\mfp_2}=\{\pm(e_2-e_3)\}\\

C_{\mfp_i}\textrm{ scalar on }\mfn,\,i=1,2\ear$$

\li

$$\bar{cc}b_{1}^{\phi} & b_2^{\phi}\\ \hline

\frac{1}{8} & \frac{1}{6}\ear$$\ebst

}}

\newpage

\section{$\mfe_8$}In this section we consider the bisymmetric
triples of the form $(\mfe_8,\mfso_{16},\mfl)$ and
$(\mfe_8,\mfe_7\oplus\mfsu_2,\mfl)$. The root system for $\mfe_8$ is

\beq\label{re8}\R=\{\pm e_i\pm e_j, 1\leq i<j\leq 8;\,
\pm\frac{1}{2}\sum_1^8(-1)^{\nu_i}e_i,\sum_1^8\nu_i\textrm{ even }
\},\eeq

where $e_1,\ldots,e_8$ is the canonical basis for $\reals^8$. In
$\mfe_8$ there is only one root length which is
\beq\label{rle8}|\al|^2=\frac{1}{30}.\eeq

\blem \label{sumse8}Let
$\phi=\frac{1}{2}\sum_1^8(-1)^{\nu_i}e_i\in\R$.

(i) Let $\al=\frac{1}{2}\sum_1^8(-1)^{\mu_i}e_i \in\R$. The string
$\phi+n\al$ is either singular, $\phi,\,\phi+\al$ or
$\phi,\,\phi-\al$. So either $d_{\al\phi}=0$ or $1$, respectively.

We have that $\phi+\al$ is a root if and only if $\nu_i=\mu_i$, for
for  two indices $i_1,i_2\in\{1,\ldots,8\}$  and in this
case, $\phi+\al=(-1)^{\nu_{i_1}}e_{i_1}+(-1)^{\nu_{i_2}}e_{i_2}$.

$\phi-\al$ is a root if and only if $\nu_i\neq\mu_i$, for
two indices $i_1,i_2\in\{1,\ldots,8\}$ and in this case,
$\phi-\al=(-1)^{\nu_{i_1}}e_{i_1}+(-1)^{\nu_{i_2}}e_{i_2}$.

(ii) Let $\al'=(-1)^{\mu_j} e_j+(-1)^{\mu_k} e_k\in\R$, $1\leq
j<k\leq 8$. The string $\phi+n\al'$ is either singular or
$\phi,\,\phi-\al'$. $\phi-\al'$ is a root if and only if
$\al'=(-1)^{\nu_j}e_j+(-1)^{\nu_k}e_k$, for $1\leq j<k\leq 8$. In
this case, $\phi-\al'=\frac{1}{2}\big(\sum_{\sbar{l}i=1\\i\neq
j,k\sear}^8(-1)^{\nu_i}e_i+(-1)^{\nu_j+1}e_j+(-1)^{\nu_k+1}e_k\big)$.\footnote{We
may choose $-\al$' in which case we obtain $\phi+\al'$
instead.}\elem

\bproof Consider the roots $\phi=\frac{1}{2}\sum_1^8(-1)^{\nu_i}e_i$
and $\al=\frac{1}{2}\sum_1^8(-1)^{\mu_i}e_i$. Suppose that the
string $\phi+n\al$ is not singular. Then, since $\phi+n\al$ is an
uninterrupted string either $\phi+\al$ or $\phi-\al$ is a root. We
have that
$$\phi+\al=\frac{1}{2}\sum_{1}^8((-1)^{\nu_i}+(-1)^{\mu_i})e_i\textrm{ and }\phi-\al=\frac{1}{2}\sum_{1}^5((-1)^{\nu_i}-(-1)^{\mu_i})e_i.$$

By observing the form of the roots in $\R$ given in (\ref{re7}) we
conclude that $\phi+\al$ is a root if and only if
$(-1)^{\nu_i}+(-1)^{\mu_i}=0$, i.e., $\nu_i=\mu_i$, for  two
indices $i_1,i_2\in\{1,\ldots,8\}$. We observe that this case is
possible since $\sum_1^8\nu_i$ and $\sum_1^8\mu_i$ are even with $8$
even. We thus obtain
$\phi+\al=(-1)^{\nu_{i_1}}e_{i_1}+(-1)^{\nu_{i_2}}e_{i_2}$. Clearly
$\phi+2\al$ is never a root. Similarly, $\phi-\al$ is a root if and
only if $(-1)^{\nu_i}-(-1)^{\mu_i}\neq 0$, i.e.,$\nu_i\neq\mu_i$,
for  two indices $i_1,i_2\in\{1,\ldots,6\}$ and in this case,
$\phi-\al=(-1)^{\nu_{i_1}}e_{i_1}+(-1)^{\nu_{i_2}}e_{i_2}$. The
element $\phi-2\al$ is never a root. Once again we observe that this
case is possible .

\li

Let $\al'=(-1)^{\mu_j} e_j+(-1)^{\mu_k} e_k\in\R$. We have

$$\phi-\al'=\frac{1}{2}\big(\sum_{\sbar{l}i=1\\i\neq
j,k\sear}^8(-1)^{\nu_i}e_i+((-1)^{\nu_j}-2(-1)^{\mu_j})e_j+((-1)^{\nu_k}-2(-1)^{\mu_k})e_k\big).$$

This element is a root if and only if $\mu_j=\nu_j$ and
$\mu_k=\nu_k$ and, in this case,
$$\phi-\al'=\frac{1}{2}\big(\sum_{\sbar{l}i=1\\i\neq
j,k\sear}^8(-1)^{\nu_i}e_i+(-1)^{\nu_j+1}e_j+(-1)^{\nu_k+1}e_k\big).$$

We observe that $\phi-2\al'$ is never a root.

$\Box$\eproof

\li

{\footnotesize{

\bsymp\label{sympe81}   $(\mfe_8,\mfso_{16})$.

$$\bar{ccc}\mfk & \R_{\mfk} & \ga\\ \hline\hline

\mfso_{16} & \{\pm e_i\pm e_j, 1\leq i<j\leq 8 \} &
\frac{7}{15}\ear$$

$$\R_{\mfn}=\R_{\mfn}=\{\pm\frac{1}{2}\sum_1^8(-1)^{\nu_i}e_i,\sum_1^8\nu_i\textrm{
even } \}$$

\esymp

\li

\bsymp \label{sympe82}  $(\mfe_8,\mfe_7\oplus\mfsu_2)$.

$$\bar{ccc}\mfk_i & \R_{\mfk_i} & \ga_i\\ \hline\hline

\mfe_7 & \{\pm e_i\pm e_j, 1\leq i<j\leq 6;\,\pm (e_7-e_8);\,
\pm\frac{1}{2}(e_7-e_8+\sum_1^6(-1)^{\nu_i}e_i),\sum_1^6\nu_i\textrm{
odd } \} & \frac{3}{5}\\\hline

 \mfsu_2 & \{\pm(e_7+e_8)\} & \frac{1}{15}\ear$$

$$\R_{\mfn}=\{\pm e_i\pm e_j, 1\leq i\leq
6,\,j=7,8;\,
\pm\frac{1}{2}(e_7+e_8+\sum_1^6(-1)^{\nu_i}e_i),\sum_1^6\nu_i\textrm{
even }\}$$

\esymp

\li

\bbst\label{cpe81} $(\mfe_8,\mfso_{16}, \mfso_{2p}\oplus\mfso_{2(8-p)})$, $1\leq p\leq 4$. (Type I)

$$\bar{l}\R_{\mfp}=\{\pm e_i\pm e_j, 1\leq i\leq p, p+1\leq j\leq 8 \}\\

C_{\mfp}\textrm{ scalar on }\mfn\ear$$

$$\bar{ccccccc} \phi\in\R_{\mfn} & \al\in\R_{\mfp}^+  & \phi+n\al & d_{\al\phi} & \textrm{No of } \al's & |\al|^2 & b^{\phi}\\ \hline\hline

\xstrut  \frac{1}{2}\sum_{1}^8e_i & e_i+e_j,\,1\leq i\leq p,\,p+1\leq j\leq 8 & \phi,\,\phi-\al & 1 & p(8-p) &
\frac{1}{30} & \frac{p(8-p)}{60} \ear$$\ebst

\li

\bland

 \bbst\label{cpe82} $(\mfe_8,\mfso_{16},
\mfu_8)$. (Type I)

$$\bar{l}\R_{\mfp}=\{\pm (e_i+ e_j), 1\leq i\leq 8 \}\\

\mfn=\oplus_{0,1,2}\mfn^i\\
\R_{\mfn^i}=\{\pm\frac{1}{2}\sum_{1}^8(-1)^{\nu_j}e_j: 2i\textrm{
odd }\nu_j's\},\,i=0,1,2\ear$$

$$\bar{ccccccc} \phi\in\R_{\mfn^i} & \al\in\R_{\mfp}^+  & \phi+n\al & d_{\al\phi} & \textrm{No of } \al's & |\al|^2 & b^{\phi}\\ \hline\hline

\xstrut  \frac{1}{2}\sum_{1}^8(-1)^{\nu_j}e_j,\, 2i\textrm{ odd
}\nu_j's & e_a+e_b,\,\textrm{ s.t. }\nu_a=\nu_b &
\phi,\,\phi+(-1)^{\nu_a+1}\al & 1 & 4i^2-16i+28^{(*)} &
\frac{1}{30} & \frac{i^2-4i+7}{15} \ear$$\ebst

(*) the number of possible pairs $(a,b)$, where $a<b$ and
$\nu_a=\nu_b$, since $2i$ $\nu_j$'s are odd and $8-2i$ $\nu_j$'s
are even, is

$$\bar{rl}& \left(\bar{c}2i\\2\ear\right)+\left(\bar{c}8-2i\\2\ear\right) \\ \\

= & 4i^2-16i+28. \ear$$

\eland

\li

\bbst\label{cpe83} $(\mfe_8,\mfe_7\oplus\mfsu_2,\mfe_7\oplus\reals)$. (Type I)

$$\bar{l}\R_{\mfp}=\{\pm(e_7+e_8)\}\\

C_{\mfp}\textrm{ scalar on }\mfn\ear$$

$$\bar{ccccccc} \phi\in\R_{\mfn} & \al\in\R_{\mfp}^+  & \phi+n\al & d_{\al\phi} & \textrm{No of } \al's & |\al|^2 & b^{\phi}\\ \hline\hline

\xstrut  \bar{l}\frac{1}{2}(e_7+e_8+\sum_1^6(-1)^{\nu_i}e_i),\\\sum_1^6\nu_i\textrm{
even }\ear & e_7+e_8 & \phi,\,\phi-\al & 1 & 1 &
\frac{1}{30} & \frac{1}{60} \ear$$\ebst

\li

 \bbst\label{cpe84} $(\mfe_8,\mfe_7\oplus\mfsu_2,\mfe_6\oplus\reals\oplus\mfsu_2)$. (Type I)

$$\bar{l}\R_{\mfp}=\{\pm e_i\pm e_6, 1\leq i\leq
5;\,\pm(e_7-e_8);\,
\pm\frac{1}{2}(e_8-e_7+e_6+\sum_1^5(-1)^{\nu_i}e_i),\sum_1^5\nu_i\textrm{
odd }\}\\

\mfn=\mfn^1\oplus\mfn^2\oplus\mfn^3\\
 \R_{\mfn^1}=\{\pm e_i\pm e_j:1\leq i\leq 5,\,j=7,8\},\\

\R_{\mfn^2}=\{\pm (e_6-e_7),\,\pm
(e_6+e_8),\,\pm\frac{1}{2}(e_8+e_7-e_6+\sum_1^5(-1)^{\nu_i}e_i),\sum_1^6\nu_i\textrm{
odd }\},\\

\R_{\mfn^3}=\{\pm (e_6+e_7),\,\pm
(e_6-e_8),\,\pm\frac{1}{2}(e_8+e_7+e_6+\sum_1^5(-1)^{\nu_i}e_i),\sum_1^6\nu_i\textrm{
even }\}\ear$$

$$\bar{cccccccc}\mfn^i & \phi\in\R_{\mfn^i} & \al\in\R_{\mfp}^+  & \phi+n\al & d_{\al\phi} & \textrm{No of } \al's & |\al|^2 & b^{\phi}\\ \hline\hline

\xstrut \mfn^1 & e_1+e_8 &
\bar{c}e_1\pm e_6,\,e_7-e_8\\\bar{l}\frac{1}{2}(e_8-e_7+e_6+e_1+\sum_2^5(-1)^{\nu_i}e_i),\\\sum_2^5\nu_i\textrm{ even}\ear\ear &
\phi,\,\phi-\al & 1 & \bar{c}3\\2^3\ear & \frac{1}{30} & \frac{11}{60}\\
\hline

\xstrut \mfn^2 & e_6-e_8 & \bar{c}e_i+ e_6,\,i\leq i\leq
5\\e_i-e_6,\,i\leq i\leq 5\\e_7-e_8\ear &
\bar{c}\phi,\,\phi-\al\\\phi,\,\phi+\al\\\phi,\,\phi+\al\ear & 1 &
\bar{c}5\\5\\1\ear & \frac{1}{30} & \frac{11}{60}\\\hline

\xstrut \mfn^2 & e_6+e_8 & \bar{c}e_i+ e_6,\,i\leq i\leq
5\\e_i-e_6,\,i\leq i\leq
5\\e_7-e_8\\\bar{l}\frac{1}{2}(e_8-e_7+e_6+\sum_1^5(-1)^{\nu_j}e_j,\\\sum_1^5\nu_j\textrm{
odd}\ear\ear &
\bar{c}\phi,\,\phi-\al\\\phi,\,\phi+\al\\\phi,\,\phi+\al\\\\\phi,\,\phi-\al\ear
& 1 & \bar{c}5\\5\\1\\\\2^4\ear & \frac{1}{30} & \frac{9}{20}
\ear$$\ebst

\li

\bbst\label{cpe85}
$(\mfe_8,\mfe_7\oplus\mfsu_2,\mfe_6\oplus\reals\oplus\reals)$. (Type II)

$$\bar{lll}\R_{\mfp_1}=\{\pm(e_7+e_8)\}\\
\R_{\mfp_2}=\{\pm e_i\pm e_6, 1\leq i\leq 5;\,\pm(e_7-e_8);\,
\pm\frac{1}{2}(e_8-e_7+e_6+\sum_1^5(-1)^{\nu_i}e_i),\sum_1^5\nu_i\textrm{
odd }\}\\

\mfn=\mfn^1\oplus\mfn^2\oplus\mfn^3\ear$$

$$\bar{lcc}\R_{\mfn^i} & b_1^{\phi} & b_2^{\phi} \\ \hline
\{\pm e_i\pm e_j:1\leq i\leq 5,\,j=7,8\} & \frac{11}{60} &\frac{1}{60}\\

\{\pm (e_6-e_7),\,\pm
(e_6+e_8),\,\pm\frac{1}{2}(e_8+e_7-e_6+\sum_1^5(-1)^{\nu_i}e_i),\sum_1^6\nu_i\textrm{
odd }\} &\frac{11}{60} & \frac{1}{60}\\

\{\pm (e_6+e_7),\,\pm
(e_6-e_8),\,\pm\frac{1}{2}(e_8+e_7+e_6+\sum_1^5(-1)^{\nu_i}e_i),\sum_1^6\nu_i\textrm{
even }\} & \frac{9}{20} & \frac{1}{60}
\ear$$\ebst

\li

 \bbst\label{cpe86} $(\mfe_8,\mfe_7\oplus\mfsu_2,\mfso_{12}\oplus\mfsu_2\oplus\mfsu_2)$. (Type I)

$$\bar{l}\R_{\mfp}=\{\pm\frac{1}{2}(e_7-e_8+\sum_1^6(-1)^{\nu_i}e_i),\sum_1^6\nu_i\textrm{
odd }\}\\

\mfn=\mfn^1\oplus\mfn^2\\
\R_{\mfn^1}=\{\pm e_i\pm e_j:1\leq i\leq ,\,j=7,8\},\\

\R_{\mfn^2}=\{\pm\frac{1}{2}(e_8+e_7+\sum_1^6(-1)^{\nu_i}e_i),\sum_1^6\nu_i\textrm{
even }\}
\ear$$

$$\bar{cccccccc}\mfn^i & \phi\in\R_{\mfn^i} & \al\in\R_{\mfp}^+  & \phi+n\al & d_{\al\phi} & \textrm{No of } \al's & |\al|^2 & b^{\phi}\\ \hline\hline

\xstrut \mfn^1 & e_1+e_7 &
\bar{l}\frac{1}{2}(e_7-e_8+e_1+\sum_2^6(-1)^{\nu_i}e_i),\\\sum_2^6\nu_i\textrm{ odd}\ear & \phi,\,\phi-\al & 1 & 2^4 & \frac{1}{30} & \frac{4}{15}\\
\hline

\xstrut \mfn^2 & \frac{1}{2}\sum_1^8e_j &
\bar{l}\frac{1}{2}(e_7-e_8+\sum_1^6(-1)^{\nu_i}e_i),\\\textrm{1
even
}\nu_i\\\frac{1}{2}(e_7-e_8+\sum_1^6(-1)^{\nu_i}e_i),\\\textrm{1
odd }\nu_i\ear & \bar{c}\phi,\,\phi+\al\\\\\phi,\,\phi-\al\ear & 1
& \bar{c}6\\\\6\ear & \frac{1}{30} & \frac{1}{5} \ear$$\ebst

\li

\bbst\label{cpe87}
$(\mfe_8,\mfe_7\oplus\mfsu_2,\mfso_{12}\oplus\mfsu_2\oplus\reals)$. (Type II)

$$\bar{lll}\R_{\mfp_1}=\{\pm\frac{1}{2}(e_7-e_8+\sum_1^6(-1)^{\nu_i}e_i),\sum_1^6\nu_i\textrm{
odd }\}\\
\R_{\mfp_2}=\R_{\mfk_2}= \{\pm(e_7+e_8)\}\\

\mfn=\mfn^1\oplus\mfn^2\ear$$

$$\bar{lcc}\R_{\mfn^i} & b_1^{\phi} & b_2^{\phi} \\ \hline
\{\pm e_i\pm e_j:1\leq i\leq ,\,j=7,8\} & \frac{4}{15} & \frac{1}{60}\\

\{\pm\frac{1}{2}(e_8+e_7+\sum_1^6(-1)^{\nu_i}e_i),\sum_1^6\nu_i\textrm{
even }\} & \frac{1}{5} & \frac{1}{60}

\ear$$\ebst

\li

 \bbst\label{cpe88}
$(\mfe_8,\mfe_7\oplus\mfsu_2,\mfsu_8\oplus\mfsu_2)$. (Type I)

$$\bar{l}\R_{\mfp}=\{\pm(e_i+e_j):1\leq i<j\leq 6;\, \pm\frac{1}{2}(e_7-e_8+\sum_{1}^6(-1)^{\nu_i}e_i),\,
\,3\,odd\,\nu_i's\}\\

\mfn \textrm{ irreducible }Ad\,L\textrm{-module}\ear$$

$$\bar{cccccccc} \phi\in\R_{\mfn} & \al\in\R_{\mfp}^+  & \phi+n\al & d_{\al\phi} & \textrm{No of } \al's & |\al|^2 & b^{\phi}\\ \hline\hline

\xstrut  e_1+e_7 & \bar{c}\bar{l}\frac{1}{2}(e_7-e_8+e_1+\sum_2^6(-1)^{\nu_i}e_i),\\\textrm{3 odd }\nu_i's\ear \\ e_1+e_j,\,j\in\{2,\ldots,6\}\ear & \phi,\,\phi-\al & 1 & \bar{c}\left(\bar{c}5\\3\ear\right)\\ 5\ear & \frac{1}{30} & \frac{1}{4}\ear$$\ebst

\li

\bbst\label{cpe89}
$(\mfe_8,\mfe_7\oplus\mfsu_2,\mfsu_8\oplus\reals)$. (Type II)

$$\bar{lll}\R_{\mfp_1}=\{\pm(e_i+e_j):1\leq i<j\leq 6;\, \pm\frac{1}{2}(e_7-e_8+\sum_{1}^6(-1)^{\nu_i}e_i),\,
\,3\,odd\,\nu_i's\}\\
\R_{\mfp_2}=\{\pm(e_7+e_8)\} \\

\mfn \textrm{ irreducible }Ad\,L\textrm{-module}\ear$$

$$\bar{cc} b_1^{\phi} & b_2^{\phi} \\ \hline
 \frac{1}{4} & \frac{1}{60}
\ear$$\ebst

}}

\newpage

\section{$\mfe_7$}In this Section we study the bisymmetric
triples of the form $(\mfe_7, \mfsu_8,\mfl)$,
$(\mfe_7,\mfso_{12}\oplus\mfsu_2, \mfl)$ and $(\mfe_7,
\mfe_6\oplus\reals,\mfl)$.

The root system for $\mfg=\mfe_7$ is

\beq\label{re7}\R=\{\pm e_i\pm e_j, 1\leq i<j\leq 6;\,\pm
(e_7-e_8);\,
\pm\frac{1}{2}(e_7-e_8+\sum_1^6(-1)^{\nu_i}e_i),\sum_1^6\nu_i\textrm{
odd } \},\eeq

where $e_1,\ldots,e_8$ is the canonical basis for $\reals^8$. Throughout all the relations for the $\nu_i$'s  are $mod\,2$. In $\mfe_7$, all the roots have the same length which is

\beq\label{rle7}|\al|^2=\frac{1}{18}.\eeq

\li

\blem \label{sumse7}Let
$\phi=\frac{1}{2}(e_7-e_8+\sum_1^6(-1)^{\nu_i}e_i)\in\R$.

(i) Let $\al=\frac{1}{2}(e_7-e_8+\sum_1^6(-1)^{\mu_i}e_i) \in\R$.
The string $\phi+n\al$ is either singular, $\phi,\,\phi+\al$ or
$\phi,\,\phi-\al$. So either $d_{\al\phi}=0$ or $1$, respectively.
We have that $\phi+\al$ is a root if and only if $\nu_i\neq\mu_i$,
for every $i=1,\ldots,6$ and in this case, $\phi+\al=e_8-e_7$;
$\phi-\al$ is a root if and only if $\nu_i\neq\mu_i$, for two
indices $i_1,i_2\in\{1,\ldots,6\}$ and in this case,
$\phi-\al=(-1)^{\nu_{i_1}}e_{i_1}+(-1)^{\nu_{i_2}}e_{i_2}$.

(ii) Let $\al'=(-1)^{\mu_j} e_j+(-1)^{\mu_k} e_k\in\R$, $1\leq
j<k\leq 6$. The string $\phi+n\al'$ is either singular or
$\phi,\,\phi-\al'$. $\phi-\al'$ is a root if and only if
$\al'=(-1)^{\nu_j}e_j+(-1)^{\nu_k}e_k$, for $1\leq j<k\leq 6$. In
this case,
$\phi-\al'=\frac{1}{2}\big(e_7-e_8+\sum_{\sbar{l}i=1\\i\neq
j,k\sear}^6(-1)^{\nu_i}e_i+(-1)^{\nu_j+1}e_j+(-1)^{\nu_k+1}e_k\big)$.\footnote{We
may choose $-\al$' in which case we obtain $\phi+\al'$
instead.}\elem

(iii) For $\al''=e_7-e_8\in\R$, the string $\phi+n\al''$ is
$\phi,\,\phi-\al''$, with
$\phi-\al''=-\frac{1}{2}(e_7-e_8+\sum_1^6(-1)^{\nu_i+1}e_i)$.

\bproof $\R$ is a subsystem of roots of the root system for
$\mfe_8$. Hence, we use Lemma \ref{sumse8}.

For (i) let us consider the roots
$\phi=\frac{1}{2}(e_7-e_8+\sum_1^6(-1)^{\nu_i}e_i)$ and
$\al=\frac{1}{2}(e_7-e_8+\sum_1^6(-1)^{\mu_i}e_i)$. Suppose that the
string $\phi+n\al$ is not singular. Since $\nu_7=\mu_7$ and
$\nu_8=\mu_8$ at least two indices satisfy $\nu_i=\mu_i$. Hence, if
for every $i\in\{1,\ldots,6\}$, $\nu_i\neq\mu_i$, then
$\phi+\al=e_8-e_7$ is a root; if there is $i\in\{1,\ldots,6\}$,
$\nu_i=\mu_i$, then we obtain a root $\phi-\al$ if and only if
$\nu_i\neq\mu_i$ for precisely two indices
$i_1,i_2\in\{1,\ldots,6\}$. In this case,
$\phi-\al=(-1)^{\nu_{i_1}}e_{i_1}+(-1)^{\nu_{i_2}}e_{i_2}$. We
observe that both conditions on indices in $\{1,\ldots,6\}$ are
possible since $\sum_i^6\nu_i$ and $\sum_i^6\mu_i$ are odd and $6$
is even.

(ii) and (iii) follow directly from (ii) in Lemma \ref{sumse8}.

$\Box$\eproof

\li

{\footnotesize{

\bsymp \label{sympe71}   $(\mfe_7,\mfso_{12}\oplus\mfsu_2)$.

$$\bar{ccc}\mfk_i & \R_{\mfk_i} & \ga_i\\ \hline\hline

\mfso_{12} & \{\pm e_i\pm e_j, 1\leq i<j\leq 6\} & \frac{5}{9}\\

 \mfsu_2 & \{\pm( e_7- e_8)\} & \frac{1}{9}\ear$$

$$\R_{\mfn}=\{\pm\frac{1}{2}(e_7-e_8+\sum_1^6(-1)^{\nu_i}e_i,\sum_1^6\nu_i\textrm{
odd } \}$$

\esymp

\li

\bsymp \label{sympe72}   $(\mfe_7,\mfe_6\oplus\reals)$.

$$\bar{ccc}\mfk & \R_{\mfk} & \ga\\ \hline\hline

\mfe_6 & \{\pm e_i\pm e_j:1\leq i<j\leq 5;\, \pm\frac{1}{2}(e_8-e_7-e_6+\sum_{1}^5(-1)^{\nu_i}e_i):\sum_{1}^5\nu_i\textrm{ is even}\} & \frac{2}{3}\ear$$

$$\R_{\mfn}=\{\pm e_i\pm e_6:1\leq i\leq
5;\,\pm(e_7-e_8)\,\pm\frac{1}{2}(e_7-e_8+e_6\sum_1^5(-1)^{\nu_i}e_i),\sum_1^5\nu_i\textrm{
odd } \}$$
\esymp

\li

\bsymp \label{sympe73}   $(\mfe_7,\mfsu_8,\mfsu_p\oplus\mfsu_{8-p}\oplus\reals)$,
$p=1,\ldots,4$.

$$\bar{ccc}\mfk & \R_{\mfk} & \ga\\ \hline\hline

\mfsu_8 & \{\pm(e_i-e_j):\,1\leq i<j\leq 6;\,
\pm(e_7-e_8);\pm
\frac{1}{2}(e_7-e_8+\sum_{1}^6(-1)^{\nu_i}e_i):\,1\,or\,5\,odd\,\nu_i's\} & \frac{4}{9}\ear$$

$$\R_{\mfn}= \{\pm(e_i+e_j):\,1\leq i<j\leq 6;\,\pm
\frac{1}{2}(e_7-e_8+\sum_{1}^6(-1)^{\nu_i}e_i):3\,odd\,\nu_i's\}$$
\esymp

\li

\bbst\label{cpe71} $(\mfe_7,\mfso_{12}\oplus\mfsu_2,\mfso_{12}\oplus\reals)$. (Type I)

$$\bar{l}\R_{\mfp}=\{\pm( e_7- e_8)\}\\

C_{\mfp}\textrm{ scalar on }\mfn\ear$$

$$\bar{ccccccc} \phi\in\R_{\mfn} & \al\in\R_{\mfp}^+  & \phi+n\al & d_{\al\phi} & \textrm{No of } \al's & |\al|^2 & b^{\phi}\\ \hline\hline

\xstrut  \frac{1}{2}(e_7-e_8-e_1\sum_2^6e_i) & e_7-e_8 & \phi,\,\phi-\al & 1 & 1 &
\frac{1}{18} & \frac{1}{36} \ear$$\ebst

\li

 \bbst\label{cpe72} $(\mfe_7,\mfso_{12}\oplus\mfsu_2,\mfu_6\oplus\mfsu_2)$. (Type I)

$$\bar{l}\R_{\mfp}=\{\pm(e_i+e_j):1\leq i<j\leq 6\}\\

\mfn=\oplus_{0,1,2}\mfn^i\\
\R_{\mfn^i}=\{\pm\frac{1}{2}(e_7-e_8+\sum_1^6(-1)^{\nu_j}e_j:\textrm{  }2i+1\,\nu_j's \textrm{ are odd}\},\,i=0,1,2\ear$$

$$\bar{ccccccc} \phi\in\R_{\mfn^i} & \al\in\R_{\mfp}^+  & \phi+n\al & d_{\al\phi} & \textrm{No of } \al's & |\al|^2 & b^{\phi}\\ \hline\hline

\xstrut  \bar{l}\frac{1}{2}(e_7-e_8+\sum_1^6(-1)^{\nu_j}e_j,\\
2i+1\,\nu_j's\textrm{ odd }\ear & \bar{l}e_a+e_b,\\\nu_a=\nu_b\ear & \phi,\,\phi-(-1)^{\nu_a}\al & 1 &
4i^2-8i+10^{(*)} & \frac{1}{18} & \frac{2i^2-4i+5}{18} \ear$$\ebst

(*) the number of possible pairs $(a,b)$, where $a\neq b$ and $\nu_a=\nu_b$, since $2i+1$ $\nu_j$'s are odd, is

$$\left(\bar{c}2i+1\\2\ear\right)+\left(\bar{c}6-(2i+1)\\2\ear\right) = 4i^2-8i+10. $$

\li

\bbst\label{cpe73}
$(\mfe_7,\mfso_{12}\oplus\mfsu_2,\mfu_6\oplus\reals)$. (Type II)

$$\bar{lll}\R_{\mfp_1}=\{\pm(e_i+e_j):1\leq i<j\leq 6\}\\
\R_{\mfp_2}=\R_{\mfk_2}=\{\pm (e_7-e_8)\}\\

\mfn=\oplus_{0,1,2}\mfn^i\ear$$

$$\bar{lcc}\R_{\mfn^i} & b_1^{\phi} & b_2^{\phi} \\ \hline
\{\pm\frac{1}{2}(e_7-e_8+\sum_1^6(-1)^{\nu_j}e_j:\textrm{  }2i+1\,\nu_j's \textrm{ are odd}\}^{(*)} & \frac{2i^2-4i+5}{18} & \frac{1}{36}\ear$$\ebst

(*) $i=0,1,2$ as in \ref{cpe72}.

\li

\bbst\label{cpe74} $(\mfe_7,\mfso_{12}\oplus\mfsu_2,\mfso_p\oplus\mfso_{12-p}\oplus\mfsu_2)$,
$p=2,4,6$. (Type I)

$$\bar{l}\R_{\mfp}=\{\pm e_i\pm e_j:1\leq i\leq p/2,\,p/2+1\leq j\leq 6\}\\

C_{\mfp}\textrm{ scalar on }\mfn\ear$$

$$\bar{ccccccc} \phi\in\R_{\mfn} & \al\in\R_{\mfp}^+  & \phi+n\al & d_{\al\phi} & \textrm{No of } \al's & |\al|^2 & b^{\phi}\\ \hline\hline

\xstrut  \frac{1}{2}(e_7-e_8\sum_1^6e_i) & e_i-e_j,\,1\leq j\leq 6 & \phi,\,\phi-\al & 1 & \frac{p(12-p)}{4} &
\frac{1}{18} & \frac{p(12-p)}{144} \ear$$\ebst

\li

\bbst \label{cpe75}$(\mfe_7,\mfso_{12}\oplus\mfsu_2,\mfso_p\oplus\mfso_{12-p}\oplus\reals)$,
$p=2,4,6$. (Type II)

$$\bar{l}\R_{\mfp_1}=\{\pm e_i\pm e_j:1\leq i\leq p/2,\,p/2+1\leq j\leq 6\}\\

\R_{\mfp_2}=\{\pm(e_7-e_8)\}\\

C_{\mfp_i}\textrm{ scalar on }\mfn,\,i=1,2\ear$$

\li

$$\bar{cc}b_{1}^{\phi} & b_2^{\phi}\\ \hline

\frac{1}{36} & \frac{p(12-p)}{144} \ear$$\ebst

\li

\bland \bbst\label{cpe76} $(\mfe_7,\mfe_6\oplus\reals,\mfso_{10}\oplus\reals\oplus\reals)$. (Type I)

$$\bar{l}\R_{\mfp}=\{\pm\frac{1}{2}(e_7-e_8-e_6+\sum_1^5(-1)^{\nu_i}e_i,\sum_1^5\nu_i\textrm{
even }\}\\

\mfn=\mfn^1\oplus\mfn^2\oplus\mfn^3\\
\R_{\mfn^1}=\R_{\mfn^1}=\{\pm e_i\pm e_6:1\leq i\leq 5\},\\

\R_{\mfn^2}=\{\pm\frac{1}{2}(e_7-e_8+e_6+\sum_1^5(-1)^{\nu_i}e_i),\sum_1^5\nu_i\,odd\}\\

\R_{\mfn^3}=\{\pm (e_7- e_8)\}\ear$$

$$\bar{cccccccc}\mfn^i & \phi\in\R_{\mfn^i} & \al\in\R_{\mfp}^+  & \phi+n\al & d_{\al\phi} & \textrm{No of } \al's & |\al|^2 & b^{\phi}\\ \hline\hline

\xstrut \mfn^1 & e_1-e_6
 & \frac{1}{2}(e_7-e_8-e_6 +e_1+\sum_2^5(-1)^{\nu_i}e_i),\,\sum_2^5\nu_i\textrm{ even}
 & \phi,\,\phi-\al & 1 & 2^3 & \frac{1}{18} & \frac{2}{9}\\\hline

\xstrut \mfn^2 & \frac{1}{2}(e_7-e_8+e_6-e_5+\sum_1^4e_j) &
\bar{c}\frac{1}{2}(e_7-e_8-e_6+e_5-\sum_1^4e_i)\\
\frac{1}{2}(e_7-e_8-e_6+\sum_1^5e_i)\\
\frac{1}{2}(e_7-e_8-e_6-e_5+\sum_1^4(-1)^{\nu_i}e_i),\textrm{one
odd }\nu_i  \ear & \bar{c}\phi,\,\phi+\al\\\phi,\,\phi-\al\\\\
\phi,\,\phi-\al\ear & 1 & \bar{c}1\\1\\\\4\ear & \frac{1}{18} &
\frac{1}{6}\\ \hline

\xstrut \mfn^3 & e_7-e_8 & all & \phi,\,\phi-\al &
1 & 2^4 & \frac{1}{18} &
\frac{4}{9}\

\ear$$\ebst

\li

\bbst\label{cpe77} $(\mfe_7,\mfe_6\oplus\reals,\mfsu_6\oplus\mfsu_2\oplus\reals)$. (Type I)

$$\bar{l}\R_{\mfp}=\{\pm(e_i+e_j):1\leq i<j\leq
5;\,\pm\frac{1}{2}(e_8-e_7-e_6+\sum_1^5(-1)^{\nu_i}e_i),\textrm{
 2 }\nu_i's\textrm{ odd}\}\\

\mfn=\mfn^1\oplus\mfn^2\\
\R_{\mfn^1}= \{\pm (e_i-e_6):1\leq i\leq 5;\,\pm(e_7-e_8);
\pm\frac{1}{2}(e_8-e_7-e_6+\sum_1^5(-1)^{\nu_i}e_i),\, \,
1\,or\,5\,\nu_i's\,odd \}\\

\R_{\mfn^2}= \{\pm (e_i+e_6):1\leq i\leq
5;\,\pm\frac{1}{2}(e_8-e_7-e_6+\sum_1^5(-1)^{\nu_i}e_i),3\,\nu_i's\,odd
\}\ear$$

$$\bar{cccccccc}\mfn^i & \phi\in\R_{\mfn^i} & \al\in\R_{\mfp}^+  & \phi+n\al & d_{\al\phi} & \textrm{No of } \al's & |\al|^2 & b^{\phi}\\ \hline\hline

\xstrut \mfn^1 &  e_7-e_8 &  \frac{1}{2}(e_8-e_7-e_6+\sum_1^5(-1)^{\nu_i}e_i),\,2\textrm{ odd }\nu_i's & \phi,\,\phi-\al &
1 & \left(\bar{c}5\\2\ear\right) & \frac{1}{18} &
\frac{5}{18}\\\hline

\xstrut \mfn^2 & e_1+e_6
 & \bar{c}e_1-e_i,\,i\leq i\leq 5\\\frac{1}{2}(e_8-e_7-e_6 -e_1+\sum_2^5(-1)^{\nu_i}e_i),\,1\textrm{ odd }\nu_i\ear
 & \bar{c}\phi,\,\phi-\al\\\phi,\,\phi+\al\ear & 1 & \bar{c}4\\4\ear & \frac{1}{18} & \frac{2}{9}
\ear$$\ebst

\li

\bbst\label{cpe78} $(\mfe_7,\mfsu_8,\mfsu_p\oplus\mfsu_{8-p}\oplus\reals)$,
$p=1,\ldots,4$. (Type I)

$$\bar{l}\bar{rl}\R_{\mfl}= \{\pm(e_i-e_j):1\leq i<j\leq p\,or\,p+1\leq i<j\leq 6;\,\pm(e_7-e_8);& \pm\frac{1}{2}(e_7-e_8+\sum_{1}^pe_i+\sum_{p+1}^6(-1)^{\nu_i}e_i):1\,odd\,\nu_i;\\
 & \pm\frac{1}{2}(e_7-e_8-\sum_{1}^pe_i+\sum_{p+1}^6(-1)^{\nu_i}e_i):(5-p)\,odd\,\nu_i's\}\ear \\

\R_{\mfp}= \{\pm(e_i-e_j):1\leq i\leq p,\,p+1\leq j\leq 6;\, \pm\frac{1}{2}(e_7-e_8+\sum_{1}^p(-1)^{\nu_i}e_i+\sum_{p+1}^6e_i):1\,odd\,\nu_i;\,\pm\frac{1}{2}(e_7-e_8+\sum_{1}^p(-1)^{\nu_i}e_i-\sum_{p+1}^6e_i):(p-1)\,odd\,\nu_i's\}\ear$$

\li

\textbf{p=1}

$$\bar{l}\R_{\mfp}=\{\pm(e_1-e_j):2\leq j\leq 6; \pm\frac{1}{2}(e_7-e_8+e_1-\sum_{2}^6e_i);\,\pm\frac{1}{2}(e_7-e_8-e_1+\sum_{2}^6e_i)\}\\

\mfn \textrm{ irreducible }Ad\,L\textrm{-module}\ear$$

\li

$$\bar{ccccccc} \phi\in\R_{\mfn^i} & \al\in\R_{\mfp}^+  & \phi+n\al & d_{\al\phi} & \textrm{No of } \al's & |\al|^2 & b^{\phi}\\ \hline\hline

\xstrut e_1+e_2  & e_1-e_j,\,3\leq j\leq 6  & \phi,\,\phi-\al & 1
& 4 & \frac{1}{18} & \frac{1}{9} \ear$$

\textbf{p=2}

$$\bar{l}\R_{\mfp}= \{\pm(e_i-e_j):1\leq i\leq 2,\,3\leq j\leq 6;
\pm\frac{1}{2}(e_7-e_8+\sum_{1}^2(-1)^{\nu_i}\pm\sum_{3}^6e_i):1\,odd\,\nu_i\}\\

\mfn=\mfn^1\oplus\mfn^2\\

\R_{\mfn^1}= \{\pm(e_1+e_2);\,\pm(e_i+e_j):\,3\leq i<j\leq 6; \pm\frac{1}{2}(e_7-e_8-e_1-e_2+\sum_{3}^6(-1)^{\nu_i}e_i):1\,odd\,\nu_i;\pm\frac{1}{2}(e_7-e_8+e_1+e_2+\sum_{3}^6(-1)^{\nu_i}e_i):3\,odd\,\nu_i\}\\

\R_{\mfn^2}=\{\pm(e_i+e_j):\,4\leq i<j\leq
6;\,\pm\frac{1}{2}(e_7-e_8\pm\sum_1^3e_i\mp\sum_4^6e_i)\}\ear$$

\li

$$\bar{cccccccc}\mfn^i & \phi\in\R_{\mfn^i} & \al\in\R_{\mfp}^+  & \phi+n\al & d_{\al\phi} & \textrm{No of } \al's & |\al|^2 & b^{\phi}\\ \hline\hline

\xstrut \mfn^1 & e_1+e_2  &  e_i-e_j,\,i=1,2,\,3\leq j\leq 6 &
\phi,\,\phi-\al & 1 & 8 & \frac{1}{18} & \frac{2}{9}\\\hline

\xstrut \mfn^2 & e_1+e_6 & \bar{c} e_1-e_j,\,3\leq j\leq 6\\
e_2-e_6\\
\frac{1}{2}(e_7-e_8+e_1-e_2+\sum_3^6e_i)\\\frac{1}{2}(e_7-e_8-e_1+e_2-\sum_3^6e_i)\ear
 & \bar{c}\phi,\,\phi-\al\\ \phi,\,\phi+\al\\\phi,\,\phi-\al\\\phi,\,\phi+\al\ear & 1 & \bar{c}4\\1\\1\\1\ear & \frac{1}{18} & \frac{1}{6}
\ear$$

\newpage
\textbf{p=3}

$$\bar{l}\R_{\mfp}=\{\pm(e_i-e_j):1\leq i\leq 3,\,4\leq j\leq 6; \pm\frac{1}{2}(e_7-e_8+\sum_{1}^3(-1)^{\nu_i}e_i+\sum_{4}^6e_i):1\,odd\,\nu_i;

\pm\frac{1}{2}(e_7-e_8+\sum_{1}^3(-1)^{\nu_i}e_i-\sum_{4}^6e_i):2\,odd\,\nu_i's\}\\

\mfn=\mfn^1\oplus\mfn^2\\

\bar{ll}\R_{\mfn^1}=\{\pm(e_i+e_j):\,1\leq i<j\leq 3\,or\, 1\leq i\leq 3,\,4\leq j\leq 6; &  \pm\frac{1}{2}(e_7-e_8+\sum_{1}^3(-1)^{\mu_i}e_i+\sum_{3}^6(-1)^{\nu_i}e_i):2\,odd\,\mu_i's\,and\,1\,odd\,\nu_i;\\
 & \pm\frac{1}{2}(e_7-e_8+\sum_{1}^3(-1)^{\mu_i}e_i+\sum_{3}^6(-1)^{\nu_i}e_i):1\,odd\,\mu_i's\,and\,2\,odd\,\nu_i\}\ear\\

\R_{\mfn^2}=\{\pm(e_i+e_j):\,4\leq i<j\leq 6;
\pm\frac{1}{2}(e_7-e_8\pm\sum_1^3e_i\mp\sum_4^6e_i) \}\ear$$

\li

$$\bar{cccccccc}\mfn^i & \phi\in\R_{\mfn^i} & \al\in\R_{\mfp}^+  & \phi+n\al & d_{\al\phi} & \textrm{No of } \al's & |\al|^2 & b^{\phi}\\ \hline\hline

\xstrut \mfn^1 & e_1+e_2  & \bar{c}e_i-e_j,\,i=1,2,\,j=4,5,6 \\ \frac{1}{2}(e_7-e_8+e_1+e_2-e_3+\sum_4^6e_i)\\\frac{1}{2}(e_7-e_8-e_1-e_2+e_3-\sum_4^6e_i) \ear  & \bar{c}\phi,\,\phi-\al\\\phi,\,\phi-\al\\\phi,\,\phi+\al\ear &
1 & \bar{c}6\\1\\1\ear & \frac{1}{18} &
\frac{2}{9}\\\hline

\xstrut \mfn^2 & e_4+e_6 &  \bar{c}e_i-e_j,\,i=1,2,3,\,j=4,6 \\ \frac{1}{2}(e_7-e_8+\sum_1^3(-1)^{\nu_i}e_i+\sum_4^6e_i),1\,odd\,\nu_i\\\frac{1}{2}(e_7-e_8+\sum_1^3(-1)^{\nu_i}e_i-\sum_4^6e_i),2\,odd\,\nu_i's \ear
 & \bar{c}\phi,\,\phi-\al\\\phi,\,\phi-\al\\ \phi,\,\phi+\al\ear & 1 & \bar{c}6\\3\\3\ear & \frac{1}{18} & \frac{1}{3}
\ear$$

\textbf{p=4}

$$\bar{l}\R_{\mfp}=  \{\pm(e_i-e_j):1\leq i\leq 4,\,5\leq j\leq 6; \pm\frac{1}{2}(e_7-e_8+\sum_{1}^4(-1)^{\nu_i}e_i+\sum_5^6e_i):1\,odd\,\nu_i;
\pm\frac{1}{2}(e_7-e_8+\sum_{1}^4(-1)^{\nu_i}e_i-\sum_5^6e_i):3\,odd\,\nu_i's\}\\ \\

\mfn=\mfn^1\oplus\mfn^2\oplus\mfn^3\\ \\

\R_{\mfn^1}=\{\pm(e_i+e_j):\,1\leq i<j\leq 4; \pm\frac{1}{2}(e_7-e_8+\sum_{1}^4(-1)^{\nu_i}e_i\pm e_5\mp
e_6):2\,odd\,\nu_i's\}\\

\R_{\mfn^2}= \{\pm(e_5+e_6)\}\\

\R_{\mfn^2}= \{\pm(e_i+e_j):\,1\leq i\leq 4,\,5\leq j\leq 6; \pm\frac{1}{2}(e_7-e_8-e_6-e_5+\sum_{1}^4(-1)^{\nu_i}e_i):1\,odd\,\nu_i; \pm\frac{1}{2}(e_7-e_8+e_6+e_5+\sum_{1}^4(-1)^{\nu_i}e_i):3\,odd\,\nu_i's\}
\ear$$

\li

$$\bar{cccccccc}\mfn^i & \phi\in\R_{\mfn^i} & \al\in\R_{\mfp}^+  & \phi+n\al & d_{\al\phi} & \textrm{No of } \al's & |\al|^2 & b^{\phi}\\ \hline\hline

\xstrut \mfn^1 &  e_1+e_2 & \bar{c}e_i-e_j,\,i=1,2,\,j=5,6 \\ \frac{1}{2}(e_7-e_8+e_6+e_5\mp e_4\pm e_3+e_2+e_1)\\\frac{1}{2}(e_7-e_8-e_6-e_5\mp e_4\pm e_3-e_2-e_1) \ear  & \bar{c}\phi,\,\phi-\al\\\phi,\,\phi-\al\\ \phi,\,\phi+\al\ear &
1 & \bar{c}4\\2\\2\ear & \frac{1}{18} &
\frac{2}{9}\\\hline

\xstrut \mfn^2 & e_5+e_6
 & \bar{c}e_i-e_j,\,i=1,2,3,4,\,j=5,6 \\ \frac{1}{2}(e_7-e_8+e_6+e_5+\sum_1^4(-1)^{\nu_i}e_i),1\,odd\,\nu_i\\\frac{1}{2}(e_7-e_8-e_6-e_5+\sum_1^4(-1)^{\nu_i}e_i),3\,odd\,\nu_i's \ear
 & \bar{c}\phi,\,\phi+\al\\\phi,\,\phi-\al\\ \phi,\,\phi+\al\ear & 1 & \bar{c}8\\4\\4\ear & \frac{1}{18} & \frac{4}{9}\\\hline

 \xstrut \mfn^3 & e_1+e_6
 & \bar{c}e_i-e_6,\,i=1,2,3,4\\e_1-e_5 \\ \frac{1}{2}(e_7-e_8+e_6+e_5+e_1+\sum_2^4(-1)^{\nu_i}e_i),1\,odd\,\nu_i\\\frac{1}{2}(e_7-e_8-e_6-e_5-e_1+\sum_2^4(-1)^{\nu_i}e_i),2\,odd\,\nu_i's \ear
 & \bar{c}\phi,\,\phi+\al\\\phi,\,\phi-\al\\ \phi,\,\phi-\al\\ \phi,\,\phi+\al\ear & 1 & \bar{c}1\\4\\3\\3\ear & \frac{1}{18} & \frac{11}{36}
\ear$$\ebst
\eland

}}

\newpage

\section{$\mfe_6$} In this Section we consider the bisymmetric
triples of the form $(\mfe_6, \mfso_{10}\oplus\reals,\mfl)$ and
$(\mfe_6,\mfsu_6\oplus\mfsu_2,\mfl)$.

If $e_1,\ldots,e_8$ is the canonical basis for $\reals^8$, we can
write the root system for $\mfe_6$ as follows:

$$\R=\{\pm e_i\pm e_j:1\leq i<j\leq 5;\, \pm\frac{1}{2}(e_8-e_7-e_6+\sum_{1}^5(-1)^{\nu_i}e_i):\sum_{1}^5\nu_i\textrm{ is even}\}.$$

Throughout, all the relations for the $\nu_i$'s are $mod\,2$.

We recall that on $\mfe_6$ there is only one root length which is

\beq\label{rle6}|\al|^2=\dfrac{1}{12}\eeq

\blem \label{sumse6}Let
$\phi=\frac{1}{2}(e_8-e_7-e_6\sum_1^5(-1)^{\nu_i}e_i)\in\R$.

(i)Let $\al=\frac{1}{2}(e_8-e_7-e_6\sum_1^5(-1)^{\mu_i}e_i) \in\R$.
The string $\phi+n\al$ is either singular or $\phi,\,\phi-\al$. So
either $d_{\al\phi}=0$ or $1$, respectively. We have that $\phi-\al$
is a root if and only if $\nu_i\neq\mu_i$, for  two indices
$i_1,i_2\in\{1,\ldots,5\}$ and in this case,
$\phi-\al=(-1)^{\nu_{i_1}}e_{i_1}+(-1)^{\nu_{i_2}}e_{i_2}$.

(ii) Let $\al'=(-1)^{\mu_j} e_j+(-1)^{\mu_k} e_k\in\R$, $1\leq
j<k\leq 5$. The string $\phi+n\al'$ is either singular or
$\phi,\,\phi-\al'$. $\phi-\al'$ is a root if and only if
$\al'=(-1)^{\nu_j}e_j+(-1)^{\nu_k}e_k$, for $1\leq j<k\leq 5$. In
this case,
$\phi-\al'=\frac{1}{2}\big(e_8-e_7-e_6\sum_{\sbar{l}i=1\\i\neq
j,k\sear}^5(-1)^{\nu_i}e_i+(-1)^{\nu_j+1}e_j+(-1)^{\nu_k+1}e_k\big)$.\footnote{We
may choose $-\al$' in which case we obtain $\phi+\al'$
instead.}\elem

\bproof $\R$ is a subsystem of roots of the root system for $\mfe_8$
and thus we use Lemma \ref{sumse8}. For (i), since $\nu_i=\mu_i$,
for $i=6,7,8$, the case that $\nu_i=\mu_i$ for precisely two indices
in $\{1,\ldots,8\}$ never happens. Hence, $\phi+\al$ is never a
root. If $\nu_i\neq\mu_i$ for  two indices
$i_1,i_2\in\{1,\ldots,5\}$, then $\phi-\al$ is a root and
$\phi-\al=(-1)^{\nu_{i_1}}e_{i_1}+(-1)^{\nu_{i_2}}e_{i_2}$. This
case may happen since $\sum_{1}^5\nu_i$ and $\sum_{1}^5\mu_i$
are even with $5$ odd.

(ii) follows directly from (ii) in Lemma \ref{sumse8}.

$\Box$\eproof

{\footnotesize{

\bsymp\label{sympe61}   $(\mfe_6, \mfso_{10}\oplus\reals)$.

$$\bar{ccc}\mfk & \R_{\mfk} & \ga\\ \hline\hline

\mfso_{10} & \{\pm e_i\pm e_j:1\leq i<j\leq 5\} &
\frac{2}{3}\ear$$

$$\R_{\mfn}=\{\pm\frac{1}{2}(e_8-e_7-e_6+\sum_{1}^5(-1)^{\nu_i}e_i):\sum_{1}^5\nu_i\,is\,even^{(*)}\}$$

(*) there are  either
$0$, $2$ or $4$ negative signs.

\esymp

\li
\bland

\bsymp \label{sympe62}  $(\mfe_6,\mfsu_6\oplus\mfsu_2)$.

$$\bar{ccc}\mfk_i & \R_{\mfk_i} & \ga_i\\ \hline\hline

\mfsu_6 & \{\pm (e_i- e_j):1\leq i<j\leq 5;\, \pm\frac{1}{2}(e_8-e_7-e_6+\sum_{1}^5(-1)^{\nu_i}e_i): 4\, odd\,\nu_i's \} & \frac{1}{2}\\\hline

 \mfsu_2 & \{\pm\frac{1}{2}(e_8-e_7-e_6+\sum_{1}^5e_i)\} & \frac{1}{6}\ear$$

$$ \R_{\mfn}=\{\pm(e_i+e_j):1\leq i<j\leq
5;\,\pm\frac{1}{2}(e_8-e_7-e_6+\sum_{1}^5(-1)^{\nu_i}e_i):2\,odd\,\nu_i's\}$$

\esymp

\li

\bbst\label{cpe61} $(\mfe_6,
\mfso_{10}\oplus\reals ,\mfu_5 \oplus\reals)$. (Type I)

$$\bar{l}\R_{\mfp}=\{\pm (e_i+e_j):1\leq i<j\leq 5\}\\

\mfn=\mfn^0\oplus\mfn^2\oplus\mfn^4\\
\R_{\mfn^0}=\{\pm\frac{1}{2}(e_8-e_7-e_6+\sum_{1}^5e_i)\},\\

\R_{\mfn^2}=\{\pm\frac{1}{2}(e_8-e_7-e_6+\sum_{1}^5(-1)^{\nu_i}e_i):
 \,2\,odd\,\nu_i's\},\\

\R_{\mfn^4}=\{\pm\frac{1}{2}(e_8-e_7-e_6+\sum_{1}^5(-1)^{\nu_i}e_i):
\,4\,odd\,\nu_i's\}
\ear$$

$$\bar{cccccccc}\mfn^i & \phi\in\R_{\mfn^i} & \al\in\R_{\mfp}^+  & \phi+n\al & d_{\al\phi} & \textrm{No of } \al's & |\al|^2 & b^{\phi}\\ \hline\hline

\xstrut \mfn^0 & \frac{1}{2}(e_8-e_7-e_6+\sum_{1}^5e_i) &
all & \phi,\,\phi-\al & 1 & \left(\bar{c}5\\2\ear\right) & \frac{1}{12} & \frac{5}{12}\\
\hline

\xstrut \mfn^2 & \frac{1}{2}(e_8-e_7-e_6-e_1-e_2+\sum_{3}^5e_i) &
\bar{c}e_1+e_2\\e_i+e_j,\,3\leq i<j\leq 5\ear &
\bar{c}\phi,\,\phi+\al\\\phi,\,\phi-\al\ear & 1 &
\bar{c}1\\3\ear & \frac{1}{12} & \frac{1}{6}\\\hline

\xstrut \mfn^4 & \frac{1}{2}(e_8-e_7-e_6+e_5-\sum_1^4e_i)  & e_i+e_j,\,1\leq i<j\leq 4 & \phi,\,\phi+\al &
1 & \left(\bar{c}4\\2\ear\right) & \frac{1}{12} &
\frac{1}{4}

\ear$$\ebst

\eland

\li

\bland

\bbst\label{cpe62} $(\mfe_6,
\mfso_{10}\oplus\reals ,\mfso_p \oplus\mfso_{10-p}\oplus\reals)$,
$p=2,4$. (Type I)

$$\bar{l}\R_{\mfp}=\{\pm e_i\pm e_j:1\leq i\leq p/2,\,p/2+1<j<5\}\\

C_{\mfp}\textrm{ scalar on }\mfn\ear$$

$$\bar{ccccccc} \phi\in\R_{\mfn} & \al\in\R_{\mfp}^+  & \phi+n\al & d_{\al\phi} & \textrm{No of } \al's & |\al|^2 & b^{\phi}\\ \hline\hline

\xstrut  \frac{1}{2}(e_8-e_7-e_6+\sum_1^5(-1)^{\nu_i}e_i) & \bar{c}(-1)^{\nu_i}e_i-(-1)^{\nu_j}e_j,\\1\leq i\leq p/2,\,p/2+1\leq j\leq 5\ear & \phi,\,\phi-\al & 1 & \frac{p(10-p)}{4} &
\frac{1}{12} & \frac{p(10-p)}{96} \ear$$\ebst

\li

\bbst\label{cpe63} $(\mfe_6,\mfsu_6\oplus\mfsu_2,\mfsu_6\oplus\reals)$. (Type I)

$$\bar{l}\R_{\mfp}=\{\pm\frac{1}{2}(e_8-e_7-e_6+\sum_{1}^5e_i)\}\\

C_{\mfp}\textrm{ scalar on }\mfn\ear$$

$$\bar{ccccccc} \phi\in\R_{\mfn} & \al\in\R_{\mfp}^+  & \phi+n\al & d_{\al\phi} & \textrm{No of } \al's & |\al|^2 & b^{\phi}\\ \hline\hline

\xstrut  \frac{1}{2}(e_8-e_7-e_6-e_5-e_4+\sum_1^3e_i) & \frac{1}{2}(e_8-e_7-e_6+\sum_1^5e_i) & \phi,\,\phi-\al & 1 & 1 &
\frac{1}{12} & \frac{1}{24} \ear$$\ebst

\eland
\li

\bland

\bbst\label{cpe64} $(\mfe_6,\mfsu_6\oplus\mfsu_2,\mfsu_p\oplus\mfsu_{6-p}\oplus\reals\oplus\mfsu_2)$,
$p=1,2,3$. (Type I)

$$\bar{l}\R_{\mfp}= \{\pm(e_i-e_j):1\leq i\leq 6-p,\,7-p\leq j\leq 5;  \pm\frac{1}{2}(e_8-e_7-e_6-\sum_{7-p}^{5}e_i+\sum_{1}^{6-p}(-1)^{\nu_i}e_i):(5-p)\,odd\,\nu_i's\}\\
\mfn=\mfn^1\oplus\mfn^2,\textrm{ for }p=2,3\textrm{ and } \mfn \textrm{ irreducible }Ad\,L\textrm{-module for }p=1\ear$$

$$\bar{cc} & b^{\phi} \\\hline

\mfn^1 & \frac{p+2}{24}\\

\mfn^2 & \frac{p}{8}\ear$$

\li

\textbf{p=1}

$$\bar{l}\R_{\mfl_1}=\{\pm(e_i-e_j):1\leq i<j\leq 5\}\\

\R_{\mfl_2}=\R_{\mfk_2}=\{\pm\frac{1}{2}(e_8-e_7-e_6+\sum_{1}^5e_i)\}\\

\R_{\mfp}=\{\pm\frac{1}{2}(e_8-e_7-e_6+\sum_{1}^5(-1)^{\nu_i}e_i):4\,odd\,\nu_i's\}\\

\mfn  \textrm{ irreducible }Ad\,L\textrm{-module}\ear$$

\li

$$\bar{ccccccc} \phi\in\R_{\mfn} & \al\in\R_{\mfp}^+  & \phi+n\al & d_{\al\phi} & \textrm{No of } \al's & |\al|^2 & b^{\phi}\\ \hline\hline

\xstrut  e_1+e_2 & \frac{1}{2}(e_8-e_7-e_6-e_1-e_2+\sum_1^3(-1)^{\nu_i}e_i),\,2\,odd\,\nu_i's   & \phi,\,\phi-\al & 1 & 3 & \frac{1}{12} & \frac{1}{8}

\ear$$

\li

\textbf{p=2}

$$\bar{l}\R_{\mfl_1}=\{\pm(e_i-e_j): 1\leq i<j\leq 4;\,\pm\frac{1}{2}(e_8-e_7-e_6+e_5-\sum_{1}^4e_i)\}\\

\R_{\mfl_2}=\{\pm\frac{1}{2}(e_8-e_7-e_6+\sum_{1}^5e_i)\}\\

\R_{\mfp}=\{\pm(e_i-e_5):1\leq i\leq
4;\,\pm\frac{1}{2}(e_8-e_7-e_6-e_5+\sum_{1}^4(-1)^{\nu_i}e_i):3\,odd\,\nu_i's\}\\

\R_{\mfn^1}= \{\pm(e_i+e_j):1\leq i<j\leq
4;\pm\frac{1}{2}(e_8-e_7-e_6+e_5+\sum_{1}^4(-1)^{\nu_i}e_i):2\,odd\,\nu_i's\}\\

\R_{\mfn^2}= \{\pm(e_i+e_5):1\leq i\leq 4;\pm\frac{1}{2}(e_8-e_7-e_6-e_5+\sum_{1}^4(-1)^{\nu_i}e_i):1\,odd\,\nu_i\}\ear$$

\li

$$\bar{cccccccc}\mfn^i & \phi\in\R_{\mfn^i} & \al\in\R_{\mfp}^+  & \phi+n\al & d_{\al\phi} & \textrm{No of } \al's & |\al|^2 & b^{\phi}\\ \hline\hline

\xstrut \mfn^1 & e_1+e_2 & \bar{c}e_i-e_5,\,i=1,2\\\frac{1}{2}(e_8-e_7-e_6-e_5\pm e_3\mp e_4-e_1-e_2)   \ear & \bar{c}\phi,\,\phi-\al\\\phi,\,\phi+\al\ear & 1 & \bar{c}2\\2\ear & \frac{1}{12} & \frac{1}{6}\\
\hline

\xstrut \mfn^2 & e_1+e_5 &
\bar{c}e_i-e_5,\,i=2,3,4\\\frac{1}{2}(e_8-e_7-e_6-e_5+\sum_2^4(-1)^{\nu_i}e_i-e_1),\,2\,odd\,\nu_i's\ear
& \phi,\,\phi+\al & 1 & \bar{c}3\\3\ear & \frac{1}{12} &
\frac{1}{4}

\ear$$

\li

\textbf{p=3}

$$\bar{l}\R_{\mfl_1}=  \{\pm(e_i-e_j): 1\leq i<j\leq 3;\,\pm(e_4-e_5); \pm\frac{1}{2}(e_8-e_7-e_6\pm e_5\mp e_4-\sum_{1}^3e_i)\}\\

\R_{\mfl_2}= \{\pm\frac{1}{2}(e_8-e_7-e_6+\sum_{1}^5e_i)\}\\

\R_{\mfp}= \{\pm(e_i-e_j):i=1,2,3,\,j=4,5;\,\pm\frac{1}{2}(e_8-e_7-e_6-e_5-e_4+\sum_{1}^3(-1)^{\nu_i}e_i):\,exactly\,2\,odd\,\nu_i's\}\\

\bar{ll}\R_{\mfn^1}= \{\pm(e_i+e_j):1\leq i<j\leq 3\,or\,\,1\leq
i\leq 3,\,j=4,5; &
\,\pm\frac{1}{2}(e_8-e_7-e_6+e_5+e_4+\sum_{1}^3(-1)^{\nu_i}e_i):2\,odd\,\nu_i's;\\
& \pm\frac{1}{2}(e_8-e_7-e_6\pm e_5\mp
e_4+\sum_{1}^3(-1)^{\nu_i}e_i):1\,odd\,\nu_i\} \ear\\

 \R_{\mfn^2}=\{\pm(e_4+e_5);\,\pm\frac{1}{2}(e_8-e_7-e_6-e_5-e_4+\sum_{1}^3e_i)\}\ear$$

\li

$$\bar{cccccccc}\mfn^i & \phi\in\R_{\mfn^i} & \al\in\R_{\mfp}^+  & \phi+n\al & d_{\al\phi} & \textrm{No of } \al's & |\al|^2 & b^{\phi}\\ \hline\hline

\xstrut \mfn^1 & e_1+e_2 & \bar{l}e_i-e_j,\,i=1,2,j=4,5\\\frac{1}{2}(e_8-e_7-e_6-e_5-e_4+ e_3-e_2-e_1)   \ear & \bar{c}\phi,\,\phi-\al\\\phi,\,\phi+\al\ear & 1 & \bar{c}4\\1\ear & \frac{1}{12} & \frac{5}{24}\\
\hline

\xstrut \mfn^2 & e_4+e_5 &
\bar{c}e_i-e_j,\,i=1,2,3,j=4,5\\\frac{1}{2}(e_8-e_7-e_6-e_5-e_4+\sum_1^3(-1)^{\nu_i}e_i),\,2\,odd\,\nu_i's\ear
& \phi,\,\phi+\al & 1 & \bar{c}6\\3\ear & \frac{1}{12} &
\frac{9}{24}

\ear$$

\ebst

\eland

\li

\bbst \label{cpe65}$(\mfe_6,\mfsu_6\oplus\mfsu_2,\mfsu_p\oplus\mfsu_{6-p}\oplus\reals\oplus\reals)$,
$p=1,2,3$. (Type II)

$$\bar{l}\bar{rl}\R_{\mfp}= &\{\pm(e_i-e_j):1\leq i\leq 6-p,\,7-p\leq j\leq 5; \\ & \pm\frac{1}{2}(e_8-e_7-e_6-\sum_{7-p}^{5}e_i+\sum_{1}^{6-p}(-1)^{\nu_i}e_i):(5-p)\,odd\,\nu_i's\}\ear\\

\R_{\mfp_2}=\{\pm\frac{1}{2}(e_8-e_7-e_6+\sum_{1}^5e_i)\} \\

\textrm{decomposition of }\mfn\textrm{ as in }\ref{cpe64}\ear$$

\li

$$\bar{ccc}\mfn^i & b_1^{\phi} & b_2^{\phi}\\\hline

\mfn^1 & \frac{p+2}{24} & \frac{1}{24}\\

\mfn^2 & \frac{p}{8} & \frac{1}{24}\ear$$

\ebst

}}

\newpage

\section*{Notation List}
\bc{\textit{Common notation:}}\ec

\li

$A_g$ Conjugation by $g\in G$ in a Lie group $G$

\li

$Ad_g$ Adjoint action of $g\in G$ in the Lie algebra $\mfg$ of a Lie
group $G$

\li

$ad_X$ Adjoint linear map of $X\in\mfg$

\li

$Kill$ the Killing form of a Lie algebra $\mfg$

\li

$B=-Kill$

\li

$C_V$ the Casimir operator of a subspace $V$ of a Lie algebra $\mfg$
with respect to the Killing form of $\mfg$

\li

$c_{V,U}$ the eigenvalue of a Casimir operator $C_V$ on a subspace $U$

\li

$R$ the curvature of a Riemannian metric

\li

$Ric$ the Ricci curvature of a Riemannian metric

\li

$K$ the sectional curvature of a Riemannian metric

\li

$h^*(\mfg)$ the dual Coxeter number of a Lie algebra $\mfg$

\li

$\R$ the system of roots of a Lie algebra $\mfg$

\li

$\R_K$ the system of restricted roots to a subgroup of maximal rank
$K$

\li

$S^+$ the subset of positive roots of a set of roots $S$

\li

$\mfg^{\complex}$ the complexification of a Lie algebra $\mfg$

\li

$V^{\complex}$ the complexification of a vector space $V$

\li

$id_V$ the identity map of a vector space $V$

\li

$0_V$ the null map of a vector space $V$

\li

\bc{\textit{Notation for homogeneous fibrations:}}\ec

\li

$G$ compact connected semisimple Lie group

\li

$K$, $L$ compact closed non-trivial subgroups of $G$ such that $L\varsubsetneq K \varsubsetneq G$

\li

$\mfg$, $\mfk$, $\mfl$ the Lie algebras of $G$, $K$ and $L$

\li

$M=G/L$

\li

$N=G/K$

\li

$F=K/L$

\li

$g_M$ an adapted metric metric on $M$ with respect to a fibration $F\rightarrow M\rightarrow N$

\li

$g_N$ the projection of an adapted metric $g_M$ onto the base
space $N$ as above

\li

$g_F$ the restriction of an adapted metric $g_M$ to the fiber
space $F$ as above

\li

$\mfn$ $Ad\,K$-invariant complement of $\mfk$ on $\mfg$

\li

$\mfp$ $Ad\,L$-invariant complement of $\mfl$ on $\mfk$

\li

$\mfm=\mfp\oplus\mfn$ $Ad\,L$-invariant complement of $\mfl$ on $\mfg$

\li

$\mfn=\mfn_1\oplus\ldots\mfn_n$ decomposition of $\mfn$ into pairwise inequivalent irreducible $Ad\,K$-submodules

\li

$\mfn=\mfn^1\oplus\ldots\mfn^{n'}$ decomposition of $\mfn$ into pairwise inequivalent irreducible $Ad\,L$-submodules

\li

$\mfp=\mfp_1\oplus\ldots\mfp_s$ decomposition of $\mfp$ into pairwise inequivalent irreducible $Ad\,L$-submodules

\li

$C_{\mfg}$, $C_{\mfk}$, $C_{\mfl}$ the Casimir operator of $\mfg$, $\mfk$ and $\mfl$, respectively

\li

$C_{\mfp_a}$ the Casimir operator of $\mfp_a$, $a=1,\ldots,s$

\li

$C_{\mfn_i}$ the Casimir operator of $\mfn_i$, $i=1,\ldots,n$

\li

$c_{\mfl,a}$ the eigenvalue of $C_{\mfl}$ on $\mfp_a$, $a=1,\ldots,s$

\li

$c_{\mfl,\mfp}$ the eigenvalue of $C_{\mfl}$ on $\mfp$, when $\mfp$ is $Ad\,L$-irreducible

\li

$c_{\mfk,i}$ the eigenvalue of $C_{\mfk}$ on $\mfn_i$, $i=1,\ldots,n$

\li

$c_{\mfk,\mfn}$ the eigenvalue of $C_{\mfk}$ on $\mfn$, when $\mfn$ is $Ad\,K$-irreducible

\li

$c_{\mfn_i,a}$ the constant defined by $Kill(C_{\mfn_i}\cdot,\cdot)\mid_{\mfp_a\times\mfp_a}=c_{\mfn_i,a}Kill\mid_{\mfp_a\times\mfp_a}$, $a=1,\ldots,s$

\li

$\ga_a$ the constant defined by $Kill_{\mfk}\mid_{\mfp_a\times\mfp_a}=\ga_aKill\mid_{\mfp_a\times\mfp_a}$

\li

$\ga$ the constant defined by $Kill_{\mfk}\mid_{\mfp\times\mfp}=\ga Kill\mid_{\mfp\times\mfp}$, when $\mfp$ is $Ad\,L$-irreducible

\li

$b_a^i$ the eigenvalue of $C_{\mfp_a}$ on $\mfn_i$,  $i=1,\ldots,n$, when this eigenvalue exists ( the indices are dropped in the case of irreducibility as above)

\li

$b_a^{\phi}$ the constant defined by $b_a^{\phi}=B(C_{\mfp_a}X_{\phi},X_{-\phi})=Kill(C_{\mfp_a^{\complex}}E_{\phi},E_{-\phi})$ for a root $\phi$; for $\phi\in\R_{\mfn}$, it represents an eigenvalue of $C_{\mfp_a}$ on $\mfn$

\li

\end{document}